\documentclass{svjour3}

\usepackage{setspace}
\usepackage{amsmath}
\usepackage{geometry}
\usepackage{graphicx}
\usepackage{parskip}
\usepackage{amssymb}
\usepackage{longtable}
\usepackage{caption}
\usepackage{subfig}
\providecommand\phantomcaption{\caption@refstepcounter\@captype}
\usepackage{mathtools}
\usepackage{tabularx}
\usepackage{arydshln}
\usepackage{graphicx}

\usepackage{pgfplots}
\pgfplotsset{compat=1.7}

\usepackage{tikz}
\usetikzlibrary{patterns}
\usetikzlibrary{calc}

\usepackage{placeins}
\usepackage{filecontents}
\usepackage{xparse}
\usepackage{float}
\usepackage{lmodern}
\usepackage{multirow}
\usepackage{makecell}
\usepackage{changepage,calc}
\usepackage{cite}
\usepackage{lscape}

\makeatletter

\newcommand{\Matrix}[1]
{\begin{bmatrix}
		\Matrix@r #1;\@bye;\Matrix@r
\end{bmatrix}}

\def\Matrix@r #1;{\@bye #1\Matrix@z\@bye\Matrix@s #1,\@bye, }%
\def\Matrix@s #1,{#1\Matrix@t }%
\def\Matrix@t #1,{\@bye #1\Matrix@y\@bye\@firstofone {&#1}\Matrix@t}%
\def\Matrix@y #1\Matrix@t{\\ \Matrix@r }%
\def\Matrix@z #1\Matrix@r {}
\def\@bye  #1\@bye   {}

\makeatother

\newcommand{\logLogSlopeTriangle}[5]
{
	\pgfplotsextra
	{
		\pgfkeysgetvalue{/pgfplots/xmin}{\xmin}
		\pgfkeysgetvalue{/pgfplots/xmax}{\xmax}
		\pgfkeysgetvalue{/pgfplots/ymin}{\ymin}
		\pgfkeysgetvalue{/pgfplots/ymax}{\ymax}
		
		\pgfmathsetmacro{\xArel}{#1}
		\pgfmathsetmacro{\yArel}{#3}
		\pgfmathsetmacro{\xBrel}{#1-#2}
		\pgfmathsetmacro{\yBrel}{\yArel}
		\pgfmathsetmacro{\xCrel}{\xBrel}
		
		\pgfmathsetmacro{\lnxB}{\xmin*(1-(#1-#2))+\xmax*(#1-#2)} 
		\pgfmathsetmacro{\lnxA}{\xmin*(1-#1)+\xmax*#1} 
		\pgfmathsetmacro{\lnyA}{\ymin*(1-#3)+\ymax*#3} 
		\pgfmathsetmacro{\lnyC}{\lnyA+#4*(\lnxA-\lnxB)}
		\pgfmathsetmacro{\yCrel}{\lnyC-\ymin)/(\ymax-\ymin)} 
		
		\coordinate (A) at (rel axis cs:\xArel,\yArel);
		\coordinate (B) at (rel axis cs:\xBrel,\yBrel);
		\coordinate (C) at (rel axis cs:\xCrel,\yCrel);
		
		\draw[#5]   (B)-- node[pos=0.5,anchor=north] {1}
		(A)-- 
		(C)-- node[pos=0.5,anchor=east] {#4}
		cycle;
	}
}
\newcommand{\logLogSlopeTriangleUp}[5]
{
	\pgfplotsextra
	{
		\pgfkeysgetvalue{/pgfplots/xmin}{\xmin}
		\pgfkeysgetvalue{/pgfplots/xmax}{\xmax}
		\pgfkeysgetvalue{/pgfplots/ymin}{\ymin}
		\pgfkeysgetvalue{/pgfplots/ymax}{\ymax}
		
		\pgfmathsetmacro{\xArel}{#1}
		\pgfmathsetmacro{\yArel}{#3}
		\pgfmathsetmacro{\xBrel}{#1-#2}
		\pgfmathsetmacro{\yBrel}{\yArel}
		\pgfmathsetmacro{\xCrel}{\xBrel}
		
		\pgfmathsetmacro{\lnxB}{\xmin*(1-(#1-#2))+\xmax*(#1-#2)} 
		\pgfmathsetmacro{\lnxA}{\xmin*(1-#1)+\xmax*#1} 
		\pgfmathsetmacro{\lnyA}{\ymin*(1-#3)+\ymax*#3} 
		\pgfmathsetmacro{\lnyC}{\lnyA-#4*(\lnxA-\lnxB)}
		\pgfmathsetmacro{\yCrel}{\lnyC-\ymin)/(\ymax-\ymin)} 
		
		\coordinate (A) at (rel axis cs:\xArel,\yArel);
		\coordinate (B) at (rel axis cs:\xBrel,\yBrel);
		\coordinate (C) at (rel axis cs:\xCrel,\yCrel);
		
		\draw[#5]   (B)-- node[pos=0.5,anchor=south] {1}
		(A)-- 
		(C)-- node[pos=0.5,anchor=east] {#4}
		cycle;
	}
}

\makeatletter
\def\@seccntformat#1{\@ifundefined{#1@cntformat}%
	{\csname the#1\endcsname\quad}  
	{\csname #1@cntformat\endcsname}
}
\let\oldappendix\appendix 
\renewcommand\appendix{%
	\oldappendix
	\newcommand{\section@cntformat}{\appendixname~\thesection\quad}
}

\usepackage{enumerate}

\title{Taylor-series expansion based numerical methods: a primer, performance benchmarking and new approaches for problems with non-smooth solutions}
\titlerunning{Collocation methods for non-smooth elasticity problems.}

\author{Thibault Jacquemin \and Satyendra Tomar \and Konstantinos Agathos \and Shoya Mohseni-Mofidi \and St\'ephane P.A. Bordas}

\institute{T. Jacquemin \and S. Tomar \and K. Agathos \at
	Institute of Computational Engineering, University of Luxembourg,\\
	Maison du Nombre, 6 Avenue de la Fonte,\\
	L-4364 Esch-sur-Alzette, Luxembourg.\\
	\and S. Mohseni-Mofidi \at
	Fraunhofer-Institute for Mechanics of Materials IWM, Freiburg, Germany.\\
	\and S. Bordas \at
	Institute of Research and Development, Duy Tan University, K7/25 Quang Trung, Danang, Vietnam.\\
	Institute of Computational Engineering, University of Luxembourg, Maison du Nombre, 6 Avenue de la Fonte,\\
	L-4364 Esch-sur-Alzette, Luxembourg.\\
	\email{stephane.bordas@alum.northwestern.edu}\\
}

\date{}

\begin{document}
	
	\maketitle
	\pagenumbering{gobble}
	\pagenumbering{arabic}
	
	\begin{abstract}
		
		We provide a primer to numerical methods based on Taylor series expansions such as generalized finite difference methods and collocation methods. We provide a detailed benchmarking strategy for these methods as well as all data files including input files, boundary conditions, point distribution and solution fields, so as to facilitate future benchmarking of new methods.
		We  review traditional methods and recent ones which appeared in the last decade. We aim to help newcomers to the field understand the main characteristics of these methods and to provide sufficient information to both simplify implementation and benchmarking of new methods. Some of the examples are chosen within a subset of problems where collocation is traditionally known to perform sub-par, namely when the solution sought is non-smooth, i.e. contains discontinuities, singularities or sharp gradients. For such problems and other simpler ones with smooth solutions, we study in depth the influence of the weight function, correction function, and the number of nodes in a given support. We also propose new stabilization approaches  to improve the accuracy of the numerical methods. In particular, we experiment with the use of a Voronoi diagram for weight computation, collocation method stabilization approaches, and support node selection for problems with singular solutions. With an appropriate selection of the above-mentioned parameters, the resulting collocation methods are compared to the moving least-squares method (and variations thereof), the radial basis function finite difference method and the finite element method. Extensive tests involving two and three dimensional problems indicate that the methods perform well in terms of efficiency (accuracy versus computational time), even for non-smooth solutions.
		
	\end{abstract}	
	
	\section{Introduction}
	
	We focus in this paper on Taylor-series expansion based collocation approaches for Partial Differential Equations (PDEs). In these methods, instead of writing the problem in an average sense, as in Galerkin methods, the strong form is written explicitly at a set of computational points, distributed over the domain. Derivative operators are computed through the use of stencils of points, which can be built in different ways. A lot of work has been done in this field since the early 1900's, and collocation methods are regaining interest, due to the advent of massively parallel computing, which lends itself very naturally to these methods. Our goal in this paper is to facilitate the understanding of newcomers to the field, help choose optimal parameters, to benchmark the methods and, finally, to propose novel approaches to deal with non-smooth solutions. Specifically, we aim to : 
	\begin{itemize}
		\item  briefly review approaches to alleviate the mesh burden in computational mechanics;
		\item  provide a gradual, clear and detailed introduction to Taylor-series expansion based collocation approaches;
		\item  investigate the sensitivity of the above approaches to the parameters involved;
		\item  provide recommendations on the methods optimal parameters;
		\item  propose and experiment on a computational approach to handle sharp corners and singularities;
		\item  provide a comprehensive investigation of the relative performance of some of the most popular such approaches;
		\begin{itemize}
			\item for smooth problems;
			\item for rough and singular problems with low solution regularity;
		\end{itemize}
		\item facilitate future benchmarking by providing all data files, including geometries, point distributions, loading and boundary conditions, solution fields, to help benchmarking existing and new methods.
	\end{itemize}
	
	Numerical methods have been under development for approximately 80 years. The first methods which were developed were finite difference methods, which focus on the approximation of the differential operator. The first known reference is the inception of finite difference methods for partial differential equation, in the work of C. Runge in 1908 \cite{Runge1908}. The idea was to use stencils of points in order to approximate differential operators using finite differences. In their initial form, finite difference methods were largely limited to Cartesian domains in space or to time approximations. 
	
	This limitation of finite difference methods to the union of Cartesian domains may have been the motivation for the development of alternative methods including the Ritz method \cite{Ritz1908} and the Galerkin finite element method \cite{Galerkin1915}. Contrary to finite difference schemes, finite element methods were able to handle arbitrarily complex geometries, at the cost of the generation of a mesh, i.e. a cover of the volume with simple shapes including tetrahedral, hexahedral and prismatic elements. 
	
	Shortly after the introduction of the concept of mesh, in 1977 with the creation of the smoothed particle hydrodynamics method \cite{Monaghan1992}, the notion of methods which would later become known as  mesh-free methods came about. SPH enabled the solution of problems which caused difficulties to finite elements, in particular those involving fluid flow, fragmentation and very large deformations. 
	
	The finite element concept of mesh, closely related to that of interpolation and approximation comes with at least five associated challenges:
	\begin{itemize}
		\item the mesh should conform to the potentially complex geometry of the domain and hence be regenerated, at least partially, for each change in the geometry of the component under consideration;
		\item for moving boundaries, the mesh must be regenerated at each geometrical change in the boundary;
		\item the aspect ratios of the elements should be controlled to ensure accuracy, in large deformations, this includes ensuring that the elements do not become too deformed or inverted during deformation;
		\item locking problems have to be accounted for when small parameters appear within the PDE, e.g. for thin plates and shells or incompressible materials, warranting the development of new locking-free formulations;
		\item stability of approximation schemes for coupled multi-field problems must be ensured, leading to the requirement of hybrid methods.
	\end{itemize}

	Some of these challenges may well have motivated the inception of alternative methods known at the time as meshless or meshfree methods \cite{Belytschko1994,Liu1995,Duarte1996,Atluri1998,De2000,Chen2000}. The original idea behind such methods was to decrease the burden posed by the generation and regeneration of a mesh. In particular, the Bubnov-Galerkin or Petrov-Galerkin methods relax some of the constraints associated with locating the points used to construct the approximation and thus simplify local refinement. Nonetheless, these methods rely on non-polynomial approximations which are usually non-interpolating, thus posing additional difficulties associated with enforcing boundary conditions and numerical integration. For many of these methods, numerical integration requires a background mesh or local integration rules on complex domains such as lenses. In their initial formulation, meshfree methods are computationally expensive, which somewhat limits their application to industrial problems. The 2008 review on implementation and recent advances in meshfree methods is a possible reference \cite{Nguyen2008}.
	
	Contemporarily to the birth of meshfree methods, partition of unity approaches see the light of day \cite{Babuska1995,Babuska1997}. In their original form, they enable the introduction of known features about the solution within the finite element approximation. Either this known feature is computed numerically (as in the generalised finite element method \cite{Strouboulis2001}) or they are extracted from analytical knowledge about the solution, as in the extended finite element method (XFEM). These methods, born in parallel to meshfree methods create an intermediate world between finite element methods and meshfree methods, and have similarities with both. For instance, methods such as XFEM enable the simulation of propagating discontinuities in the field variable or its derivative with minimal or no remeshing, whilst some versions of partition of unity methods require special treatment of boundary conditions.
	Some of the most exciting applications of partition of unity methods include fracture mechanics either as enriched finite elements \cite{Mos1999,Sukumar2000,Dolbow2000,Dolbow2001,Sukumar2001,Mos2002,Ji2004,Duflot2008} or as enriched meshfree methods \cite{Rabczuk2007,Rabczuk2007Sec,Rabczuk2007Thi,Bordas2007,Bordas2008,Talebi2011,Natarajan2011}. Note that such partition of unity methods were also used to permit the implicit treatment of (evolving) discontinuities using level set methods, including an implicit description of the boundary of the computational domain \cite{Belytschko2002,Moumnassi2014}. Several recent reviews can be consulted for an overview on partition of unity methods \cite{Rabczuk2010}.
	
	A decade after the appearance of Galerkin meshfree methods, isogeometric analysis (IGA) approaches saw the light \cite{Hughes2005}. Their primary goal is to facilitate the connection between computer aided design (CAD) and computer aided engineering (CAE) with numerical analysis by using the same functions used to describe the geometry of the object to also approximate the unknown field variables. In this way, the method is able to represent complex geometries exactly. During early stage design iterations, any change in the geometry is automatically inherited by the approximation scheme for the field variables, thereby simplifying the iterative design process. Isogeometric analysis boundary element methods (IGABEM) \cite{Simpson2012,Simpson2013,Scott2013,Lian2013,Peng2014,Atroshchenko2015,Lian2016,Peng2017,Lian2017,Atroshchenko2017} transcend the intrinsic limitations of IGA within a finite element context, in particular the requirement of 3D volume parameterisation, akin to hexahedral meshing \cite{Xu2011,Xu2013,Xu2018}. IGA shares many common points with meshfree methods, in particular its natural ability to deal with high order approximations, which makes it suitable to handle Kirchhoff-Love plates and shells and high-order PDEs. Various approaches combining enrichment with IGA were introduced \cite{Nguyen2015}. The reader can refer to the recent overview and computer implementation aspects of IGA presented in \cite{NguyenImplementation2015}.
	
	To overcome the most negative aspects of IGA, i.e. the need for structured Cartesian parameterisation associated with the tensor product nature of the method as well as the consequential difficulties associated with local mesh refinement, the geometry-independent field approximation method (GIFT) was proposed by Atroshchenko and colleagues in a series of papers \cite{Atroshchenko2018}, which relaxes the strict requirement of using the same basis functions to represent the geometry and the field variables, and, hence enables the local refinement of the field approximation independently of the non-uniform rational B-splines (NURBS) representation of the boundary. This therefore maintains the tight coupling between the CAD and the analysis of a given component, without requiring the use of NURBS for field approximations, which has been shown to be suboptimal in certain situations, for example, for problems with corner singularities, or weakly regular solutions.
	
	Contemporarily with IGA, methods based on implicit treatments of boundaries have continued to develop, thanks to the combined efforts of engineers and applied mathematicians \cite{Burman2010,Burman2012,Burman2014,Burman2014_2,Hansbo2014,Burman2015,Claus2015,Claus2017,Claus2018}, and \cite{BordasUnfitted2017}.
	
	In light of the above summary, the finite element methods and the meshfree methods seem to have been abandoned by the computational mechanics community. Collocation methods, however, have been continuously studied from the mid-1950s to date. Collocation methods have been reintroduced into the literature thanks to the recrudescence of advanced computing hardware such as graphical processing units and Xeon Phis, among others. Such computing architectures have memory architectures which are well suited to handling similar data shapes such as the row of a stiffness matrix provided by collocation approaches. 
	
	Now that we have painted an impressionist picture of the path towards mesh-burden reduction, subsequent to the birth of finite difference methods and finite element methods, we proceed to introducing collocation methods, which we classify broadly into two groups. The first group includes all methods which use an approximation of the differential operator to solve the Partial Differential Equation (PDE). In this paper, two methods of the first group, which use a Taylor's series expansion to approximate the field derivatives, are considered. These methods are the Generalized Finite Difference (GFD) method and the Discretization-Corrected Particle Strength Exchange (DC PSE) method, see the discussion below. The second group includes methods which are based on an approximation of the unknown field. The most prominent method in this second group is the Moving Least Squares (MLS) method \cite{Lancaster1981,Shepard1968} that is used in the Element Free Galerkin (EFG) method \cite{Belytschko1994}.
	
	The idea of generalizing the Finite Difference Method (FDM) began in 1953 with MacNeal \cite{Macneal953}, and in 1960 with Forsyth and Wasow \cite{Forsythe1960}. They proposed a method to transform an irregular node distribution over the domain into a regular sub-domain on which the FDM can be applied. In 1962, Jensen \cite{Jensen1972} introduced the basis of the Generalized Finite Difference Method. The method, described for two-dimensional problems, uses a six-node star and a second order Taylor's series expansion to approximate the spatial derivatives up to the second order. In 1980, Liszka and Orkisz \cite{Liszka1980} presented a method based on an eight-node star which allows obtaining a more stable approximation of the derivatives. The method is based on some selected weights and a mean least square approximation of the derivatives. In 1998, Orkisz \cite{Orkisz1998} presented a more complete version of the GFD method covering various subjects, such as the application of the method to the Galerkin framework, and the use of a posteriori error estimators for model adaptivity.
	
	The Particle Strength Exchange (PSE) method was introduced by Degond and Mas-Gallic in 1989  \cite{Degond1989}. Initially developed to approximate the diffusion operator of the convection-diffusion equations, the method has been generalized by Eldredge et al. in 2002 \cite{Eldredge2002} in order to approximate any derivative order. The Discretization-Corrected Particle Strength Exchange (DC PSE) method has been introduced by Schrader et al. \cite{Schrader2010} in 2010 in order to account for the discretization of the domain in the operator calculation. This allows removing the discretization error, which led to the name being ``Discretization-Corrected".
	
	The GFD and the DC PSE methods show many similarities, which are analyzed in this paper. Both methods are based on a set of parameters. In this paper, we study the sensitivity of these methods to these parameters. Some methods, aiming at improving the accuracy of the solution, are presented and analyzed in the paper. The two considered methods are compared to other well known collocation methods. These methods can be classified into two categories. The methods based on an approximation of the differential operator such as the GFD and the DC PSE methods form the first group, and the methods based on an approximation of the field form the second group. For each category, the following methods are considered:
	
	\begin{itemize}
		\item Differential Operator Approximation;
		\begin{itemize}
			\item Generalized Finite Difference Method (GFD);
			\item Discretization-Corrected Particle Strength Exchange Method (DC PSE);
			\item Radial Basis Function Finite Difference Method (RBF-FD).
		\end{itemize}
		\item Field Approximation;
		\begin{itemize}
			\item Moving Least Square Method (MLS);
			\item Interpolating Moving Least Square Method (I-MLS).
		\end{itemize}
	\end{itemize}
		
	A brief outline of the remainder of the paper is as follows. In Section \ref{MethodDescription}, we briefly describe each of the methods considered in this paper. In Section \ref{ProblemAndNorms}, three linear elastic problems, for which an analytical solution is known (i.e. a cylinder under internal pressure, a sphere under internal pressure, and an L-shape domain in mode I loading), are presented. Moreover, the error norms are also introduced in Section \ref{ProblemAndNorms}. The methods are compared for the $L_2$ norm and the $L_{\infty}$ norm of the error for the calculated stress components. In Section \ref{ParametricAnalysis}, we present a parametric sensitivity study of the methods. This includes a study of the weight function, of the correction function (for DC PSE), and of the number of support nodes. In Section \ref{ImprovementMethods_Section}, we present some improvement methods for the GFD and the DC PSE methods, such as Voronoi diagram, stabilization, and criteria for support node selection for singular problems. In Section \ref{Benchmarking}, we present some benchmarking results from the comparison of the various methods listed above. We also present some results on convergence rates and computational expenses of these methods. In Section \ref{3DResults}, we present the results of the GFD method for 3D problems. Moreover, we compare our results with finite elements results obtained using the commercial package ABAQUS \cite{Abaqus2017}. Some conclusions are drawn in Section \ref{Conclusions}. Finally, a detailed comparison of the GFD and DC PSE methods for 1D problems is provided in the Appendix.
	
	As a novelty of our work, we list two main components, namely (1) a detailed comparison of the GFD and DC PSE methods for 2D and 3D linear elastic problems (such as the pressurized cylinder and the L-shape domain in mode I loading), and (2) the assessments of the improvement methods as well as the identification and comparison of variations on DC PSE methods. To the best of the authors' knowledge, these studies are not found in the literature.
	
	\section{Collocation Methods} \label{MethodDescription}
	
	\subsection{Introduction} \label{General_Method}
	
	Solving a problem by collocation methods consists in solving the set of PDEs only at collocation centers. A number of nodes spread over the domain are used to estimate the derivatives at the collocation centers. In most collocation methods, the equations are solved at the nodes. The problem being solved locally, the strong form of the PDEs is considered.
	In this paper, we primarily consider the GFD and DC PSE methods. These methods are compared to the MLS approximation method and to the RBF-FD method, which are among the most popular methods of approximation in the framework of collocation methods. In the remainder of the present section, we present the principles of each of these methods.
	In order to facilitate the comprehension of the methods, the case of a two dimensional problem in a Cartesian coordinate system is considered. The GFD and DC PSE methods are also presented and compared for the case of a 1D problem in Appendix \ref{app:a}. 
	
	In the sections below, the spatial coordinates are denoted by $x$ and $y$. The coordinates of a node $\mathbf{X}$ are then $\mathbf{X}=[x,y]^T$. The subscripts $c$ and $p$ are used to identify, respectively, the collocation node and a particle ``$p$".
	The first and second derivatives in the two spatial directions are denoted by: $\frac{\partial}{\partial x}$, $\frac{\partial}{\partial y}$, $\frac{\partial^{2}}{\partial x^{2}}$, $\frac{\partial^{2}}{\partial x \partial y}$, $\frac{\partial^{2}}{\partial y^{2}}$. In the general case, these derivatives are written as $D^{n_x,n_y}f(\mathbf{X_c})$, where $n_x$ and $n_y$ are, respectively, the derivation orders in the directions $x$ and $y$.
	
	The derivatives at a collocation center are typically approximated based on a defined support. The support is the set of nodes located in the vicinity of the collocation node. Figure \ref{Node_Neighbors} below shows the nodes of the domain $\Omega$ included in the support $\Omega_c$ of a collocation node $\mathbf{X_c}$. In 2D, the support is limited by a circle of radius $R_{\text{sup}}$.
	
	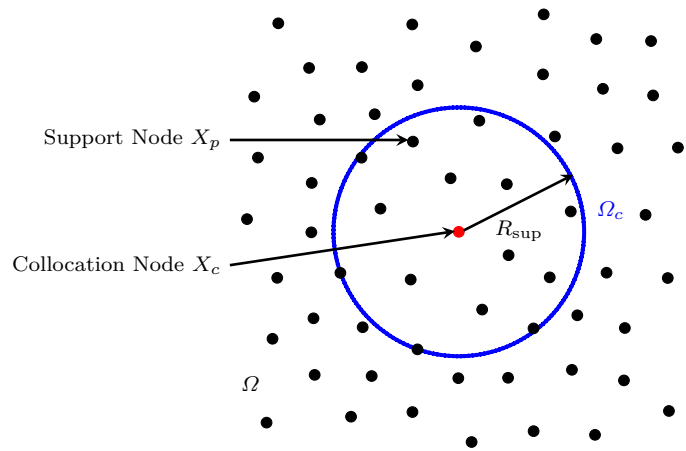
\begin{figure}[H] 
		\centering
		\begin{tikzpicture}
		\definecolor{GreyColor}{rgb}{0.3,0.3,0.3}
		\begin{axis}[height=7.5cm,width=7.5cm, xmin=1.80,xmax=2.26,scaled x ticks = false, ymin=1.80,ymax=2.26,scaled y ticks = false,axis line style={draw=none},tick style={draw=none},xticklabels={,,}, yticklabels={,,}, cells={anchor=west},  font=\footnotesize, rounded corners=1pt]
		\addplot[blue,line width=1.5pt]  table [x=X-Circle, y=Y-Circle, col sep=comma] {SupportDrawing.csv};
		\addplot[only marks,red,mark=*,mark options={fill=red}]  table [x=X-Ref, y=Y-Ref, col sep=comma] {SupportDrawing.csv};
		\addplot[only marks,black,mark=*,mark options={fill=black}]  table [x=X-Nodes, y=Y-Nodes, col sep=comma] {SupportDrawing.csv};
		\end{axis}
		\draw [black,-stealth, line width=1.0pt] (0,2.57) node [left] {Collocation Node $X_c$} -- (2.96,3.02);
		\draw [black,-stealth, line width=1.0pt] (0,4.23) node [left] {Support Node $X_p$} -- (2.35,4.23);
		\draw [black,-stealth, line width=1.0pt] (3.075,3.03) node [anchor=west, align=left] {\quad $R_{\text{sup}}$} --  (4.52,3.77);
		\draw [black] (0.5,1) node [left] {$\Omega$};
		\draw [blue] (5.3,3.3) node [left] {$\Omega_c$};
		\end{tikzpicture}
		\caption{Collocation Node Support}
		\label{Node_Neighbors}
	\end{figure}
	
	\subsection{Generalized Finite Difference Method} \label{GFD_Method}
	
	\subsubsection{Principle}
	
	The FDM is the most simple and one of the oldest methods for derivative approximation. The major drawback of this method is that it requires the use of a regular mesh. In 1972, Jensen \cite{Jensen1972} presented a method to approximate two dimensional derivatives using the Taylor's series approximation on an irregular grid. This method is known as the Generalized Finite Difference (GFD) method. 
	
	For the GFD method, the derivatives are calculated at collocation nodes $\mathbf{X_c}=[x_c,y_c]^T$ using a Taylor's series expansion of the unknown field. The field derivatives at $\mathbf{X_c}$ are computed in order to reproduce the known field values $f(\mathbf{X_{pi}})$ for a number of points $\mathbf{X_{pi}}=[x_{pi},y_{pi}]^T$. The number of selected points depends on the approximated derivative order.
	
	\subsubsection{Differential Operator Approximation}
	
\begin{figure}[H] 
		\centering
		\def\svgwidth{6cm}
		\includegraphics{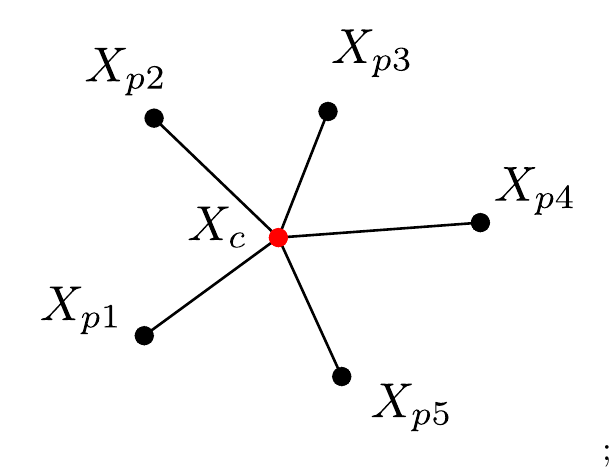}
		\caption{Five Nodes Support of a Collocation Node $X_c$}
		\label{GFD_FiveSupNodes}
	\end{figure}

	Considering a function $f: \rm I\!R^2 \rightarrow \rm  I\!R$, the Taylor's series expansion of this function at $\mathbf{X_{pi}}$ in the vicinity of a collocation node $\mathbf{X_c}$ is written:
	\begin{equation} \label{Taylor2D_AllTerms}
	f(\mathbf{X_{pi}})=\sum_{i=0}^{+\infty} \sum_{j=0}^{+\infty} \frac{\partial^{i+j}f (\mathbf{X_c})}{\partial x^i\partial y^j}  \frac{(x_{pi} - x_c)^i}{i!}  \frac{(y_{pi} - y_c)^j}{j!}.
	\end{equation}
	For ease of notations, we write the second order approximation of the function $f$ at the point $\mathbf{X_{pi}}$ near $\mathbf{X_c}$ as $f_h(\mathbf{X_{pi}})$. For $f_h(\mathbf{X_{pi}})$, Equation (\ref{Taylor2D_AllTerms}) becomes:
	
	\begin{equation} \label{Taylor2D_SecondOrderApprox}
	\begin{aligned}
	f_h(\mathbf{X_{pi}})= &f(\mathbf{X_c}) +(x_{pi} - x_c)\frac{\partial f (\mathbf{X_c})}{\partial x} + (y_{pi} - y_c)\frac{\partial f (\mathbf{X_c})}{\partial y} \\
	&+ \frac{(x_{pi} - x_c)^{2}}{2!}\frac{\partial^2 f (\mathbf{X_c})}{\partial x^2} + (x_{pi} - x_c)(y_{pi} - y_c)\frac{\partial^2 f (\mathbf{X_c})}{\partial x \partial y} + \frac{(y_{pi} - y_c)^{2}}{2!}\frac{\partial^2 f (\mathbf{X_c})}{\partial y^2}.
	\end{aligned}
	\end{equation}
	Equation (\ref{Taylor2D_SecondOrderApprox}) can be cast in a matrix form:
	\begin{equation} \label{Taylor2D_GFD}
	\begin{aligned}
	&\Matrix { x_{pi} - x_c, y_{pi} - y_c, \frac{(x_{pi} - x_c)^{2}}{2!}, (x_{pi} - x_c)(y_{pi} - y_c), \frac{(y_{pi} - y_c)^{2}}{2!}} \Matrix {\frac{\partial f (\mathbf{X_c})}{\partial x} ;\frac{\partial f (\mathbf{X_c})}{\partial y} ; \frac{\partial^2 f (\mathbf{X_c})}{\partial x^2}; \frac{\partial^2 f (\mathbf{X_c})}{\partial x \partial y}; \frac{\partial^2 f (\mathbf{X_c})}{\partial y^2}}
	= f_h(\mathbf{X_{pi}}) - f(\mathbf{X_c}).
	\end{aligned}
	\end{equation}
	In order to determine an approximation of the field derivatives $\mathbf{Df(X)}= \Big[ \frac{\partial f (\mathbf{X})}{\partial x}, \frac{\partial f (\mathbf{X})}{\partial y}, \frac{\partial^2 f (\mathbf{X})}{\partial x^2}, \frac{\partial^2 f (\mathbf{X})}{\partial x \partial y}, \frac{\partial^2 f (\mathbf{X})}{\partial y^2} \Big]^T $ (five unknowns), Equation (\ref{Taylor2D_SecondOrderApprox}) will be written for five nodes $\mathbf{X_{pi}}$ in the vicinity of $\mathbf{X_c}$ (see Figure \ref{GFD_FiveSupNodes}). Thereby, a linear system is obtained.
	\begin{equation} \label{System_GFD}
	\begin{aligned}
	&\Matrix { x_{p1} - x_c, y_{p1} - y_c, \frac{(x_{p1} - x_c)^{2}}{2!}, (x_{p1} - x_c)(y_{p1} - y_c), \frac{(y_{p1} - y_c)^{2}}{2!} ; x_{p2} - x_c, y_{p2} - y_c, \frac{(x_{p2} - x_c)^{2}}{2!}, (x_{p2} - x_c)(y_{p2} - y_c), \frac{(y_{p2} - y_c)^{2}}{2!} ; x_{p3} - x_c, y_{p3} - y_c, \frac{(x_{p3} - x_c)^{2}}{2!}, (x_{p3} - x_c)(y_{p3} - y_c), \frac{(y_{p3} - y_c)^{2}}{2!} ; x_{p4} - x_c, y_{p4} - y_c, \frac{(x_{p4} - x_c)^{2}}{2!}, (x_{p4} - x_c)(y_{p4} - y_c), \frac{(y_{p4} - y_c)^{2}}{2!} ; x_{p5} - x_c, y_{p5} - y_c, \frac{(x_{p5} - x_c)^{2}}{2!}, (x_{p5} - x_c)(y_{p5} - y_c), \frac{(y_{p5} - y_c)^{2}}{2!}} \Matrix {\frac{\partial f (\mathbf{X_c})}{\partial x} ;\frac{\partial f (\mathbf{X_c})}{\partial y} ; \frac{\partial^2 f (\mathbf{X_c})}{\partial x^2}; \frac{\partial^2 f (\mathbf{X_c})}{\partial x \partial y}; \frac{\partial^2 f (\mathbf{X_c})}{\partial y^2}}
	= \Matrix {f_h(\mathbf{X_{p1}}) - f(\mathbf{X_c});f_h(\mathbf{X_{p2}}) - f(\mathbf{X_c});f_h(\mathbf{X_{p3}}) - f(\mathbf{X_c});f_h(\mathbf{X_{p4}}) - f(\mathbf{X_c});f_h(\mathbf{X_{p5}}) - f(\mathbf{X_c})}.
	\end{aligned}
	\end{equation}
	Assuming that $f_h$ is close to $f$ in the vicinity of $\mathbf{X_c}$, the derivatives at the collocation node $\mathbf{X_c}$ can be approximated as a function of $f_h(\mathbf{X_c})$ and $f_h(\mathbf{X_{pi}})$ by solving the above system. If more than five points $\mathbf{X_{pi}}$ are chosen for solving Equation (\ref{Taylor2D_GFD}), the system is overdetermined. In that case, the derivatives at $\mathbf{X_c}$ leading to the minimum error can be determined using the least square method.
	
	\subsubsection{Overdetermined Approximation}
	
	If an arbitrary number of nodes $m$ is selected, the derivatives are determined using the mean least square method. A mean least square functional $B$ is presented below for the two dimensional case for both, the general form (\ref{FunctionalB_AllTerms_GFD}) and the second order approximation (\ref{FunctionalB_SecondOrder_GFD}). A weight function $w$ is typically used to balance the contribution of each node in the approximation. While a wide range of functions can be used as weight, $3^{\text{rd}}$ and $4^{\text{th}}$ order splines are usually preferred.
	\begin{align}
	\begin{split}\label{FunctionalB_AllTerms_GFD}
	B(\mathbf{X_c})= \sum_{i=1}^m { w(\mathbf{X_{pi}} - \mathbf{X_c}) \Big[\sum_{j=0}^{+\infty} \sum_{k=0}^{+\infty} \frac{\partial^{j+k} f (\mathbf{X_c})}{\partial x^j\partial y^k}  \frac{(x_{pi} - x_c)^j}{j!}  \frac{(y_{pi} - y_c)^k}{k!} - f(\mathbf{X_{pi}}) \Big]^2}. \\
	\end{split}\\
	\begin{split} \label{FunctionalB_SecondOrder_GFD}
	B_h(\mathbf{X_c})=\sum_{i=1}^m { w(\mathbf{X_{pi}} - \mathbf{X_c}) \Big[ } & {f(\mathbf{X_c}) - f(\mathbf{X_{pi}}) + (x_{pi} - x_c)\frac{\partial f (\mathbf{X_c})}{\partial x} + (y_{pi} - y_c)\frac{\partial f (\mathbf{X_c})}{\partial y}}\\
	& + {\frac{(x_{pi} - x_c)^{2}}{2!}\frac{\partial^2 f (\mathbf{X_c})}{\partial x^2} + (x_{pi} - x_c)(y_{pi} - y_c)\frac{\partial^2 f (\mathbf{X_c})}{\partial x \partial y}}\\
	& + {\frac{(y_{pi} - y_c)^{2}}{2!}\frac{\partial^2 f (\mathbf{X_c})}{\partial y^2} \Big]^2}. \\
	\end{split}
	\end{align}
	
	The derivatives $\mathbf{Df(X_c)}$, that best approximate the known field values using the Taylor's series expansion, minimize $B_h(\mathbf{X})$ when:
	\begin{equation} \label{DerivativeFunctionalB_SecondOrder_GFD}
	\frac{\partial B_h(\mathbf{X})}{\partial \mathbf{Df(X)}}\biggr\rvert_{\mathbf{X}=\mathbf{X_c}}=0.
	\end{equation}
	Equation (\ref{DerivativeFunctionalB_SecondOrder_GFD}) can be written as a linear system of the form:
	\begin{equation} \label{GFDLinearSystem}
	\mathbf{A(X_c) Df(X_c) = E(X_c) F(X_c)}.
	\end{equation}
	For the two dimensional second order case, the matrices $\mathbf{A(X_c)}$, $\mathbf{E(X_c)}$ and $\mathbf{F(X_c)}$ are:
	\begin{align}
	\begin{split}\label{GFD_MatAx}
	\mathbf{A(X_c)}=& \begin{bmatrix}
	m_{11} & m_{12} & \dots    &  m_{15} \\
	m_{21} & m_{22} & \dots    &  m_{25} \\
	\vdots &        &          & \vdots \\
	m_{51} & m_{52} & \dots    &  m_{55} \\
	\end{bmatrix} \in \rm I\!R^{5 \times 5},
	\end{split}\\[15pt]
	\begin{split}\label{GFD_MatBx}
	\mathbf{E(X_c)}=&\begin{bmatrix}
	-m_{01} & m_{01,1} & \dots    &  m_{01,m} \\
	-m_{02} & m_{02,1} & \dots    &  m_{02,m} \\
	\vdots &        &          & \vdots \\
	-m_{05} & m_{05,1} & \dots    &  m_{05,m} \\
	\end{bmatrix} \in \rm I\!R^{5 \times (m+1)},
	\end{split}\\[15pt]
	\begin{split}\label{GFD_MatFx}
	\mathbf{F(x_c)}=&\begin{bmatrix}
	f(\mathbf{X_c}) & f(\mathbf{X_{p1}}) & f(\mathbf{X_{p2}}) & \dots &	f(\mathbf{X_{pm}}) \\
	\end{bmatrix}^T,
	\end{split}
	\end{align}	
	where the moments $m_{ij,k}$ and $m_{ij}$ correspond to:
	\begin{equation}\label{Moments_GFD}
	\begin{aligned}
	m_{ij,k} &= w(\mathbf{X_{pk}} - \mathbf{X_c}) P_{(i+1),k}(\mathbf{X_c}) P_{(j+1),k}(\mathbf{X_c}), \\
	m_{ij} &= \sum_{k=1}^m {m_{ij,k}}.
	\end{aligned}
	\end{equation}
	The matrix $\mathbf{P(X_c)} \in \rm I\!R^{5 \times m}$ is written as follows:
	\begin{equation}\label{GFD_MatP}
	\mathbf{P(X_c)}=\begin{bmatrix}
	1 &1 & \dots    & 1 \\
	(x_{p1} - x_c) & (x_{p2} - x_c) & \dots    &  (x_{pm} - x_c)  \\
	(y_{p1} - y_c) & (y_{p2} - y_c) & \dots    &  (y_{pm} - y_c)  \\
	\frac{(x_{p1} - x_c)^2}{2!} & \frac{(x_{p2} - x_c)^2}{2!} & \dots    &  \frac{(x_{pm} - x_c)^2}{2!}  \\
	(x_{p1} - x_c)(y_{p1} - y_c) & (x_{p2} - x)(y_{p2} - y_c) & \dots    &  (x_{pm} - x_c)(y_{pm} - y_c)  \\
	\frac{(y_{p1} - y_c)^2}{2!} & \frac{(y_{p2} - y_c)^2}{2!} & \dots    &  \frac{(y_{pm} - y_c)^2}{2!}  \\
	\end{bmatrix}.
	\end{equation}
	The derivative vector $\mathbf{Df(X_c)}$ can then be determined as a function of $\mathbf{F(X_c)}$:
	\begin{equation}\label{FinalSystem_GFD}
	\mathbf{Df(X_c)=A(X_c)^{-1} E(X_c) F(X_c)}.
	\end{equation}
	The approximated derivatives are determined by solving the linear system (\ref{FinalSystem_GFD}). These derivatives are, by definition, consistent with each other as they participate in reproducing the unknown field values based on a Taylor's series expansion.
	
	\subsection{Discretization-Corrected Particle Strength Exchange Method (DC PSE)}\label{DC-PSE_Method}
	
	\subsubsection{General DC PSE Operator}
	
	The DC PSE method is based on a Taylor's series expansion of the unknown field. A convolution function is used to select the approximated derivative term. All the other unknown terms of the expansion are canceled out by the convolution function. The Taylor's series expansions presented in Equation (\ref{Taylor2D_AllTerms}) and Equation (\ref{Taylor2D_SecondOrderApprox}) are convoluted by a function $\eta$ over a domain $\Omega_c$:
	\begin{equation} \label{Convolution_AllTerms}
	\begin{aligned}
	\int_{\Omega_c} {f(\mathbf{X_p})} \eta(\mathbf{X_p}-\mathbf{X_c}) d\mathbf{X_p} = &\sum_{i=0}^{+\infty} \sum_{j=0}^{+\infty} \int_{\Omega_c}{\frac{\partial^{i+j} f (\mathbf{X_c})}{\partial x^{i} \partial x^{j}}  \frac{(x_p - x_c)^i}{i!}  \frac{(y_p - y_c)^j}{j!} \eta(\mathbf{X_p}-\mathbf{X_c}) d\mathbf{X_p}}.
	\end{aligned}
	\end{equation}
	The second order approximation of Equation (\ref{Convolution_AllTerms}) is written as follows:
	\begin{equation} \label{Convolution_SecondOrderApprox}
	\begin{aligned}
	\int_{\Omega_c} {f_h(\mathbf{X_p})} \eta(\mathbf{X_p}-\mathbf{X_c}) d\mathbf{X_p} = &\int_{\Omega_c} {f(\mathbf{X_c})} \eta(\mathbf{X_p}-\mathbf{X_c}) d\mathbf{X_p} \\
	&+ \int_{\Omega_c} {\frac{\partial f(\mathbf{X_c})}{\partial x}} (x_p - x_c) \eta(\mathbf{X_p}-\mathbf{X_c}) d\mathbf{X_p} \\
	&+ \int_{\Omega_c} {\frac{\partial f(\mathbf{X_c})}{\partial y}} (y_p - y_c) \eta(\mathbf{X_p}-\mathbf{X_c}) d\mathbf{X_p} \\
	&+ \int_{\Omega_c} {\frac{\partial^2 f(\mathbf{X_c})}{\partial x^2}} \frac{(x_p - x_c)^{2}}{2!} \eta(\mathbf{X_p}-\mathbf{X_c}) d\mathbf{X_p} \\
	&+ \int_{\Omega_c} {\frac{\partial^2 f(\mathbf{X_c})}{\partial x \partial xy}} (x_p - x_c) (y_p - y_c) \eta(\mathbf{X_p}-\mathbf{X_c}) d\mathbf{X_p} \\
	&+ \int_{\Omega_c} {\frac{\partial^2 f(\mathbf{X_c})}{\partial y^2}} \frac{(y_p - y_c)^{2}}{2!} \eta(\mathbf{X_p}-\mathbf{X_c}) d\mathbf{X_p}.
	\end{aligned}
	\end{equation}
	Equations (\ref{Convolution_AllTerms})  and (\ref{Convolution_SecondOrderApprox}) can be simplified by introducing the moments $M_{i,j}(\mathbf{X_c})$ which are defined as follows:
	\begin{equation} \label{MomentEquation}
	M_{i,j}(\mathbf{X_c})= \int_{\Omega_c} \frac{(x_p - x_c)^i}{i!} \frac{(y_p - y_c)^j}{j!}  \eta(\mathbf{X_p}-\mathbf{X_c}) d\mathbf{X_p}.
	\end{equation}
	Considering that the field is relatively smooth in $\Omega_c$, the integration can be transformed into a discrete summation over the nodes of the domain. Constant values $V_p$ are associated to each of the nodes of the domain. The moments then become:
	\begin{equation} \label{DiscreteMomentEquation}
	M_{i,j}(\mathbf{X_c})= \sum_{p \in \Omega_c} V_p \frac{(x_p - x_c)^i}{i!} \frac{(y_p - y_c)^j}{j!}  \eta(\mathbf{X_p}-\mathbf{X_c}).
	\end{equation}
	The values $V_p$ associated to the particles $p$ are hard to determine in the general case. Assuming a uniform distribution of the particles over the domain, these values are typically set to unity. Equation (\ref{DiscreteMomentEquation}) then becomes:
	\begin{equation} \label{DiscreteMomentEquation_UnitVol}
	M_{i,j}(\mathbf{X_c})= \sum_{p \in \Omega_c} \frac{(x_p - x_c)^i}{i!} \frac{(y_p - y_c)^j}{j!} \eta(\mathbf{X_p}-\mathbf{X_c}).
	\end{equation}
	Using these moments, Equation (\ref{Convolution_AllTerms}) and Equation (\ref{Convolution_SecondOrderApprox}), respectively, become:
	\begin{equation} \label{ConvWithMoment_AllTerms}
	\begin{aligned}
	\sum_{p \in \Omega_c} {f_h(\mathbf{X_p})} \eta(\mathbf{X_p}-\mathbf{X_c}) = &\sum_{i=0}^{+\infty} \sum_{j=0}^{+\infty} {\frac{\partial^{i+j} f(\mathbf{X_c})}{\partial x^i \partial x^j} M_{i,j}(\mathbf{X_c})},
	\end{aligned}
	\end{equation}
	\begin{equation} \label{ConvWithMoment_SecondOrderApprox}
	\begin{aligned}
	\sum_{p \in \Omega_c} {f_h(\mathbf{X_p})} \eta(\mathbf{X_p}-\mathbf{X_c}) = &f(\mathbf{X_c}) M_{0,0}(\mathbf{X_c}) + {\frac{\partial f(\mathbf{X_c})}{\partial x}} M_{1,0}(\mathbf{X_c}) + {\frac{\partial f(\mathbf{X_c})}{\partial y}} M_{0,1}(\mathbf{X_c}) \\
	&+ {\frac{\partial f^2(\mathbf{X_c})}{\partial x^2}} M_{2,0}(\mathbf{X_c}) + {\frac{\partial f^2(\mathbf{X_c})}{\partial x \partial y}} M_{1,1}(\mathbf{X_c})  + {\frac{\partial f^2(\mathbf{X_c})}{\partial y^2}} M_{0,2}(\mathbf{X_c}).
	\end{aligned}
	\end{equation}
	The selection of an appropriate function $\eta$ allows approximating the desired derivative $D^{k,l}f(\mathbf{X_c})=\frac{\partial f^{k+l} (\mathbf{X_c})}{\partial^k x \partial^l y}$ by setting all the moments to zero except the one multiplying $D^{k,l}f(\mathbf{X_c})$, which is set to unity.
	Equation (\ref{ConvWithMoment_AllTerms}) can then be written:
	\begin{equation}\label{DC PSE Operator}
	\left \{
	\begin{aligned}
	&D^{k,l}f(\mathbf{X_c}) = \sum_{p \in \Omega_c} {f_h(\mathbf{X_{p}})} \eta(\mathbf{X_{p}}-\mathbf{X_c})\\
	&\begin{array}{ll}
	\text{with} &M_{k,l}(\mathbf{X_c})=1\\
	&M_{i,j}(\mathbf{X_c})=0 \text{ \quad } \text{if} \ (i,j) \ne (k,l) .\\
	\end{array}\\
	\end{aligned}
	\right.	
	\end{equation}
	
	\subsubsection{The Convolution Function}
	
	In order to satisfy at each node of the domain the moment condition (\ref{DC PSE Operator}), the convolution function $\eta$ needs to be chosen carefully. Schrader et al. \cite{Schrader2012} performed a study of a wide range of functions. In general, the convolution function is composed of the product of two functions: the correction function $K$ and the weight function $w$:
	\begin{equation} \label{KernelForm}
	\eta(\mathbf{X})=K(\mathbf{X}) w(\mathbf{X}).
	\end{equation}
	The correction function is typically derived from a polynomial or an exponential basis. For the case of a two dimensional problem, the polynomial basis $\mathbf{P}=[1, x, y, x^2, xy, y^2]^T$ can be selected. The weight function is a function that returns a scalar based on the distance to a defined origin. It has typically a compact support: the weights are null outside of a defined perimeter. For isotropic weight functions (functions with similar behavior in every direction), the support of a collocation node $X_c$ is limited by a radius $r_c$. The normalized distance to the collocation node is written as $s$. For the node $X_p$ within the support of $X_c$, $s$ is written as $s_p$ and equals to:
	\begin{equation} \label{NormalizedDistance}
	s_p=\frac{\lVert \mathbf{X_p}-\mathbf{X_c} \rVert_2}{r_c}.
	\end{equation}
	The shape of the weight function has a significant impact on the solution as it balances the contribution of each node of the support in the field derivative approximation. Three types of weight functions can be considered in particular. These are:
	\begin{align}
	\intertext{The exponential weight functions:}	
	\begin{split}\label{ExpWeight}
	w(s)=
	\begin{cases}
	e(-s^{\alpha}\epsilon^{-2}) & \text{ \quad if } s \leq 1 \\
	0  & \text{ \quad if } s > 1, \\
	\end{cases}\\
	\end{split}\\
	\intertext{where $\alpha$ is an exponent and $\epsilon$ is a shape parameter,}
	\intertext{3$^\text{rd}$ order spline weight functions:}
	\begin{split}\label{Spline3Weight}
	w(s) =
	\begin{cases}
	\frac{2}{3} - 4 s^2 +	4 s^3 & \text{ \quad if } s \leq 0.5 \\
	\frac{4}{3} - 4 s +	4 s^2 - \frac{4}{3} s^3 & \text{ \quad if } 0.5 < s \leq 1 \\
	0 & \text{ \quad if } s > 1,
	\end{cases}\\
	\end{split}\\
	\intertext{4$^\text{th}$ order spline weight functions:}	
	\begin{split}\label{Spline4Weight}
	w(s)=
	\begin{cases}
	1 - 6 s^2 +	8 s^3 - 3 s^4 & \text{ \quad if } s \leq 1 \\
	0  & \text{ \quad if } s > 1. \\
	\end{cases}\\
	\end{split}
	\end{align}
	A typical convolution function, composed of a polynomial correction function $\mathbf{P}$ and a vector of coefficients $\mathbf{a}$, is written as:
	\begin{equation} \label{SimpleKernel}
	\eta(\mathbf{X_p}-\mathbf{X_c}) ={\mathbf{P(X_p-X_c)}}^T \mathbf{a} \ w(s_p).
	\end{equation}
	In order for this convolution function to satisfy the moment condition, the polynomial order shall be of at least the derivation order.
	
	\subsubsection{Correction Function Calculation}
	
	The coefficient vector $\mathbf{a}$ is the solution of a linear system $\mathbf{A_M(X_c)a }=\mathbf{B_M}$, where the left side of the equation corresponds to the moments calculation for the unknown convolution function $\eta$. The vector $\mathbf{B_M}$ corresponds to the moment condition which needs to be satisfied to obtain the desired derivative approximation.
	For instance, in order to approximate the derivative $D^{2,0}f(\mathbf{X_c})$, the system is:
	\begin{equation} \label{DCPSE_LinerSyst}
	\left \{
	\begin{array}{ll}
	\begin{aligned}
	&M_{0,0}(\mathbf{X_c})=0 &\Leftrightarrow \quad &\sum_{p \in \Omega_c} {\mathbf{P(X_p-X_c)}}^T \mathbf{a} w(s_p) = 0 \\
	&M_{1,0}(\mathbf{X_c})=0 &\Leftrightarrow \quad &\sum_{p \in \Omega_c} (x_p - x_c) {\mathbf{P(X_p-X_c)}}^T \mathbf{a} w(s_p) = 0 \\
	&M_{0,1}(\mathbf{X_c})=0 &\Leftrightarrow \quad &\sum_{p \in \Omega_c} (y_p - y_c) {\mathbf{P(X_p-X_c)}}^T \mathbf{a} w(s_p) = 0 \\
	&M_{2,0}(\mathbf{X_c})=1 &\Leftrightarrow \quad &\sum_{p \in \Omega_c} \frac{(x_p - x_c)^2}{2!} {\mathbf{P(X_p-X_c)}}^T \mathbf{a} w(s_p) = 1 \\
	&M_{1,1}(\mathbf{X_c})=0 &\Leftrightarrow \quad &\sum_{p \in \Omega_c} (x_p - x_c)(y_p - y_c) {\mathbf{P(X_p-X_c)}}^T \mathbf{a} w(s_p) = 0 \\
	&M_{0,2}(\mathbf{X_c})=0 &\Leftrightarrow \quad &\sum_{p \in \Omega_c} \frac{(y_p - y_c)^2}{2!} {\mathbf{P(X_p-X_c)}}^T \mathbf{a} w(s_p) = 0. \\
	\end{aligned}
	\end{array}
	\right.
	\end{equation}
	Considering the vector $\mathbf{Q(X_c,X_p)}=[1,(x_p - x_c),(y_p - y_c),\frac{(x_p - x_c)^2}{2!},(x_p - x_c)(y_p - y_c),\frac{(y_p - y_c)^2}{2!}]^T$, the correction function basis $\mathbf{P}$ and the weight function $w$, the coefficients of the matrix $\mathbf{A_M} \in \rm I\!R^{6 \times 6}$ can be written:
	\begin{equation}\label{DCPSE_ACoefficients}
	A_{M(i,j)}(\mathbf{X_c})=\sum_{p \in \Omega_c}Q_i(\mathbf{X_c},\mathbf{X_p})P_j(\mathbf{X_p}-\mathbf{X_c})w(s_p).
	\end{equation}
	Having solved the system of equations ($\mathbf{a}=\mathbf{A_M^{-1}(X_c) B_M}$), the derivative $D^{2,0}f(\mathbf{X_c})$ can be approximated with Equation (\ref{DC PSE Operator}) as a function of $f_h(\mathbf{X_p}), \: p \in  \Omega_c$. From a computational point of view, it shall be noted that the inversion of the matrix $\mathbf{A_M(X_c)}$ only needs to be performed once per collocation node $X_c$. If the approximation of another derivative is required for the solution of the partial differential equation, only the moment condition set by the vector $\mathbf{B_M}$ is updated. 
	
	\subsubsection{Identified Variations of DC PSE Method}\label{DCPEVariationsSec}
	
	It can be observed from the DC PSE method presented above that setting the moment $M_{0,0}(\mathbf{X_c})$ to zero is not necessary. The value $f(\mathbf{X_c})$ is determined in the global problem, and does not need to be canceled by a null moment $M_{0,0}(\mathbf{X_c})$. Based on this remark, three approaches can be considered. These approaches have respectively been labeled DCPSE0, DCPSE1 and DCPSE2.
	\begin{itemize}
		\item \underline{DCPSE0:} \quad $M_{0,0}(\mathbf{X_c})$ is set to zero (case presented above);
		\item \underline{DCPSE1:} \quad $M_{0,0}(\mathbf{X_c})$ is set to a constant value (e.g. 1), thereby reducing the sparsity of the global matrix;
		\item \underline{DCPSE2:} \quad The $M_{0,0}(\mathbf{X_c})$ moment is not introduced in the polynomial coefficient calculation. The dimension of the polynomial basis is then reduced by one. The matrix $\mathbf{A_M}$ belongs to $\rm I\!R^{5 \times 5}$.
	\end{itemize}
	All of these methods are compared in Section \ref{DCPSEVariationComp}.
	
	\subsection{Radial Basis Function Finite Difference Method}
	
	\subsubsection{Principle}
	
	In \cite{Kansa1990a, Kansa1990b}, Kansa introduced the idea of using Radial Basis Functions for solving differential equations over a domain. Contributions to the RBF-FD method were made later by Driscoll and Fornberg \cite{Driscoll2002}, Shu \cite{Shu2003}, Fornberg \cite{Fornberg2011,Fornberg2013} and Davydov \cite{Davydov2011,Davydov2011a}. The principle of the RBF-FD method is to determine an approximation of the differential operator at a collocation node $\mathbf{X_c}$ based on a linear combination of the field values at the nodes $\mathbf{X_p}$ nearby. Considering $m$ nodes in the support of $\mathbf{X_c}$, the aim is to determine a set of weights $\lambda$, written in a vector form as $\mathbf{W(X_c)}=[\lambda_{p1}\ \dots \ \lambda_{pm}]$, so that:
	\begin{equation} \label{DiffApprox_RBF-FD}
	D^{n_x,n_y}f(\mathbf{X_c}) = \mathbf{W(X_c)} \begin{bmatrix}
	f(\mathbf{X_{p1}}) \\
	\vdots \\
	f(\mathbf{X_{pm}})\\
	\end{bmatrix}.
	\end{equation}
	The RBF-FD method assumes that Equation (\ref{DiffApprox_RBF-FD}) is exact for all radial basis functions $\varphi$ centered at each node of the support.
	\subsubsection{Radial Basis Functions}
	
	Various classes of RBFs can be used for the purpose of the RBF-FD method. Depending on the type of RBF, one or two shape parameters need to be selected. The selection of the shape parameter(s) is critical to the accuracy of the solution as it balances the contribution of the neighboring nodes in the derivative approximation. The main classes of RBFs are presented in Table \ref{RBF_Types} below.
	\begin{table}[h]
		\centering
		\caption{Radial Basis Functions Types \cite{Kee2008}}
		\label{RBF_Types}
		\renewcommand{\arraystretch}{1.5}
		\begin{tabular}{|l|l|c|}
			\hline
			\multicolumn{1}{|c|}{\textbf{Type}} & \multicolumn{1}{c|}{\textbf{Expression}} & \textbf{Shape Parameters} \\
			\hline
			Multi-quadratics (MQ) & $\varphi(s_p)=({s_p}^2+c^2)^q$ & c, q \\
			Gaussian (EXP) & $\varphi(s_p)=e^{-c {s_p}^2}$ & c \\
			Thin plate spline (TPS) & $\varphi(s_p)={s_p}^\eta$ & $\eta$ \\
			Logarithmic & $\varphi(s_p)={s_p}^\eta log(s_p)$ & $\eta$ \\
			3$^\text{rd}$ order spline & \multicolumn{1}{c|}{Equation (\ref{Spline3Weight})} & - \\
			4$^\text{th}$ order spline &\multicolumn{1}{c|}{Equation (\ref{Spline4Weight})} & - \\
			\hline
		\end{tabular}
	\end{table}
	\FloatBarrier
	In order for the overall system matrix to be sparse, the radial basis functions are chosen with a compact support.
	
	\subsubsection{Obtaining the Differentiation Matrix}
	
	Replacing the unknown field $f$ by the selected radial basis function $\varphi_{X_c}$ centered in $X_c$, Equation (\ref{DiffApprox_RBF-FD}) becomes:
	\begin{equation} \label{DiffRBFApprox_RBF-FD}
	D^{n_x,n_y}\varphi_{X_c}(\mathbf{X_c}) = \mathbf{W(X_c)} \begin{bmatrix}
	\varphi_{X_c}(\mathbf{X_{p1}}) \\
	\vdots \\
	\varphi_{X_c}(\mathbf{X_{pm}})\\
	\end{bmatrix}.
	\end{equation}
	The weights are determined so that the approximation is exact for every radial basis function centered at each support node of the collocation node support. A linear system of equations is thus obtained:
	\begin{equation} \label{0OrderApprox_RBF-FD}
	\begin{bmatrix}
	\varphi_{X_{p1}}(\mathbf{X_{p1}}) & \varphi_{X_{p1}}(\mathbf{X_{p2}}) & \dots  & \varphi_{X_{p1}}(\mathbf{X_{pm}}) \\
	\varphi_{X_{p2}}(\mathbf{X_{p1}}) & \varphi_{X_{p2}}(\mathbf{X_{p2}}) & \dots  & \varphi_{X_{p2}}(\mathbf{X_{pm}}) \\
	\vdots & \vdots & & \vdots \\
	\varphi_{X_{pm}}(\mathbf{X_{p1}}) & \varphi_{X_{pm}}(\mathbf{X_{p2}}) & \dots  & \varphi_{X_{pm}}(\mathbf{X_{pm}}) \\
	\end{bmatrix} \begin{bmatrix}
	\lambda_{p1} \\ \lambda_{p2} \\ \vdots \\ \lambda_{pm}\\
	\end{bmatrix} = \begin{bmatrix}
	D\varphi_{X_{p1}}(\mathbf{X_c}) \\ D\varphi_{X_{p2}}(\mathbf{X_c}) \\ \vdots \\ D\varphi_{X_{pm}}(\mathbf{X_c})\\
	\end{bmatrix}.
	\end{equation}
	Additional constraints can be added to the system in order for the radial basis functions to reproduce exactly polynomials of at least the derivative order. This also ensures a certain regularity of the solution. For instance, for the case of a 2D problem and for a first order derivation in the $y$ direction, a first order polynomial basis can be added to the set of RBFs. Thereby, the system presented in Equation (\ref{0OrderApprox_RBF-FD}) becomes:
	\begin{equation} \label{1OrderApprox_RBF-FD}
	\begin{bmatrix}\begin{array}{ccc|ccc}
	\varphi_{X_{p1}}(\mathbf{X_{p1}}) & \dots  & \varphi_{X_{p1}}(\mathbf{X_{pm}}) & 1 & x_{p1} & y_{p1}  \\
	\varphi_{X_{p2}}(\mathbf{X_{p1}}) & \dots  & \varphi_{X_{p2}}(\mathbf{X_{pm}}) & 1 & x_{p2} & y_{p2}  \\
	\vdots  &  & \vdots & \vdots & \vdots & \vdots  \\
	\varphi_{X_{pm}}(\mathbf{X_{p1}}) & \dots  & \varphi_{X_{pm}}(\mathbf{X_{pm}}) & 1 & x_{pm} & y_{pm}  \\ \hline
	1 & \dots  & 1 & 0 & 0 & 0 \\
	x_{p1} & \dots  & x_{pm} & 0 & 0 & 0  \\
	y_{p1} & \dots  & y_{pm} & 0 & 0 & 0  \\
	\end{array}
	\end{bmatrix}
	\begin{bmatrix}\begin{array}{c}
	\lambda_{p1} \\ \lambda_{p2} \\ \vdots \\ \lambda_{pm}\\ \hline
	\lambda_{m+1} \\ \lambda_{m+2}\\ \lambda_{m+3}\\
	\end{array}\end{bmatrix} =
	\begin{bmatrix}
	\begin{array}{c}
	D^{0,1}\varphi_{X_{p1}}(\mathbf{X_c}) \\ D^{0,1}\varphi_{X_{p2}}(\mathbf{X_c}) \\ \vdots \\ D^{0,1}\varphi_{X_{pm}}(\mathbf{X_c})\\ \hline
	D^{0,1}1=0\\ D^{0,1}x=0\\ D^{0,1}y=1\\
	\end{array}\end{bmatrix}.
	\end{equation}
	Once the weights $W(\mathbf{X_c})$ determined by the solution of Equation (\ref{0OrderApprox_RBF-FD}) or Equation (\ref{1OrderApprox_RBF-FD}), the derivative $D^{n_x,n_y}f(\mathbf{X_c})$ can be approximated.
	
	\subsection{Moving Least Square Approximation}
	
	\subsubsection{Field Approximation}
	
	The MLS method has been introduced by Lancaster and Salkauskas in 1981 \cite{Lancaster1981}. Interpolation for the lowest order has been introduced by Shepard in 1968 \cite{Shepard1968}. The method has been widely used in the context of the Element-Free Galerkin (EFG) method \cite{Belytschko1994} and in the framework of collocation for the Finite Point Method \cite{Onate1996}. The MLS method consists in approximating the unknown field using a function basis. Differentiation of the approximated field and solution of a partial differential equation then becomes possible. The unknown field can be approximated with various types of functions depending on the considered application. Polynomial functions are typically used for linear elasticity problems.
	
	Considering a polynomial basis $P(\mathbf{X})$ and a coefficient vector $\mathbf{a(X_c)}$, an approximation of the field $f$ around a collocation node $\mathbf{X_c}$ can be written as follows:
	\begin{equation} \label{FieldApprox_MLS}
	f_h(\mathbf{X},\mathbf{X_c}) = \mathbf{P(X)^T a(X_c)}.
	\end{equation}
	The coefficients $\mathbf{a(X_c)}$ are determined by minimizing the error of the approximated field over a set of $m$ nodes around the collocation node. The error is weighted by a function $w$ centered in $X_c$. The minimization problem can be expressed by a functional $B(\mathbf{X_c})$:
	\begin{equation} \label{FunctionalB_MLS}
	B(\mathbf{X_c})= \sum_{i=1}^m {w(\mathbf{X_c} - \mathbf{X_{pi}}) \Big[\mathbf{P(X_{pi})}^T \mathbf{a(X_c)} - f(\mathbf{X_{pi}})} \Big]^2.
	\end{equation}
	The resulting error represented by the functional $B(\mathbf{X_c})$ is minimal when:
	\begin{equation} \label{DerivativeFunctionalB_MLS}
	\frac{\partial B(\mathbf{X_c})}{\partial \mathbf{a(X_c)}}=0.
	\end{equation}
	This problem is a linear system of the form:
	\begin{equation} \label{DerivativeSystem_MLS}
	\mathbf{A(X_c) a(X_c) = E(X_c) F(X_c)}.
	\end{equation}
	For a polynomial vector of size $n$, the matrices $\mathbf{A(X_c)}$, $\mathbf{E(X_c)}$ and $\mathbf{F(X_c)}$ correspond to:
	\begin{align}
	\begin{split}\label{MLS_MatAx}
	\mathbf{A(X_c)}=\begin{bmatrix}
	m_{11} & m_{12} & \dots    &  m_{1n} \\
	m_{21} & m_{22} & \dots    &  m_{2n} \\
	\vdots &        &          & \vdots \\
	m_{n1} & m_{n2} & \dots    &  m_{nn} \\
	\end{bmatrix} \ \in \rm I\!R^{n \times n},
	\end{split}\\[15pt]
	\begin{split}\label{MLS_MatBx}
	\mathbf{E(X_c)}=\begin{bmatrix}
	m_{01,1} & m_{01,2} & \dots &  m_{01,m} \\
	m_{02,1} & m_{02,2} & \dots &  m_{02,m} \\
	\vdots &        &          & \vdots \\
	m_{0n,1} & m_{0n,2} & \dots &  m_{0n,m} \\
	\end{bmatrix} \ \in \rm I\!R^{n \times m},
	\end{split}\\[15pt]
	\begin{split}\label{MLS_MatFx}
	\mathbf{F(x_c)}=\begin{bmatrix}
	f\mathbf{(X_{p1}}) & f(\mathbf{X_{p2}}) &	\dots &	f(\mathbf{X_{pm}}) \\
	\end{bmatrix}^T,
	\end{split}\\
	\intertext{where}
	\begin{split}\label{MLS_Moments0}
	m_{ij,k}= w(\mathbf{X_c} - \mathbf{X_{pk}}) P_i^{X_c}(\mathbf{X_{pk}})P_j^{X_c}(\mathbf{X_{pk}}),\\
	\end{split} \\[15pt]
	\begin{split}\label{MLS_Moments}
	m_{ij}= \sum_{k=1}^m {m_{ij,k}}. \\
	\end{split} \\[15pt]
	\intertext{For a two dimensional second order case with a polynomial basis, $\mathbf{P}$ is chosen as:}
	\begin{split}\label{MLS_P-Values}
	\mathbf{P^{X_c}(X_{pk})}=\begin{bmatrix}
	1 \\
	(x_{pk} - x_c) \\	(y_{pk} - y_c) \\	(x_{pk} - x_c)^2 \\
	(x_{pk} - x_c)(y_{p1} - y_c) \\	(y_{p1} - y_c)^2 \\
	\end{bmatrix}. \\
	\end{split}
	\end{align}
	As for the RBF-FD method, the dimension of the function basis $\mathbf{P}$ can be augmented in order to improve the regularity of the solution.
	
	\subsubsection{Boundary Condition Enforcement}
	
	The MLS method does not interpolate the field values. This impacts the Dirichlet boundary condition enforcement. Unlike the methods presented in the previous sections, the Dirichlet boundary condition is applied to the approximated field rather than directly to the degree of freedom solved in the linear system. The shape functions of the approximated field need to be calculated at the Dirichlet boundary nodes and are used to set the boundary condition. This leads to a slightly denser linear system as the rows of the matrix corresponding to the Dirichlet degree of freedoms are filled with coefficients allowing the approximation of the field at the boundary condition location.
	
	\subsubsection{Interpolating Moving Least Square Method}
	
	The Interpolating Moving Least Square (IMLS) method is a variation of the Moving Least Square method that allows the approximated field to interpolate the solution. The method has been presented in \cite{Lancaster1986}, and analyzed in a number of papers \cite{Ishida1999}, \cite{Maisuradze2003}. Interpolation of the approximated field can be achieved by various means. One of which consists in choosing a near singular weight function. The weight function assigns to the reference node a very large weight compared to the other nodes of the support. This makes the system (\ref{DerivativeSystem_MLS}) nearly singular but allows interpolation of the field. In this paper, the following weight function is considered:
	\begin{equation} \label{IMLSWeight}
	w(s)=
	\begin{cases}
	e^{-s^{2}}(s^n-\epsilon)^{-1} & \text{ \quad if } s \leq 1 \\
	0  & \text{ \quad if } s > 1. \\
	\end{cases}
	\end{equation}
	
	Here the parameters $n$ and $\epsilon$ control the singularity of the function. Based on the analysis of Maisuradze et al. \cite{Maisuradze2003}, the following parameters have been selected: $n=4$ and $\epsilon=10^{-15}$.
	
	\section{Problems and Error Norms}\label{ProblemAndNorms}
	
	\subsection{Problems Considered}
	\label{RefProblems}
	
	For our comparisons, we have selected two two-dimensional and one three-dimensional linear elastic problems. These are:
	\begin{itemize}
		\item A cylinder under internal pressure (2D - plane stress model);
		\item An L-shape domain in Mode I loading (2D - plane stress model);
		\item A sphere under internal pressure (3D).
	\end{itemize}
	An analytical solution is known for each of these problems. The 2D and 3D problems are, respectively, presented in Figure \ref{2DProblemsConsidered} and Figure \ref{SphereDwg} below. A Cartesian coordinate system has been used. Due to the symmetries of the problems in the Cartesian coordinate system, we have only considered $\text{1/4}^\text{th}$ of the cylinder and $\text{1/8}^\text{th}$ of the sphere.
	
	For all the problems presented in this section, a regular node discretizations has been selected. 
	
	\pagebreak
	
	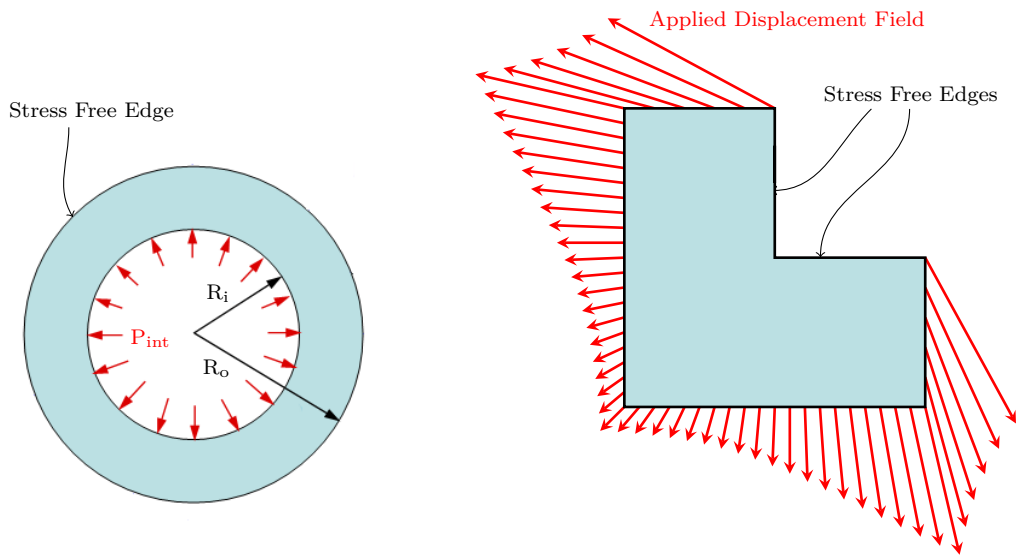
\begin{figure}[H] 
		\centering
		\begin{tikzpicture}
		\node at (0,0) {\includegraphics[scale=0.45]{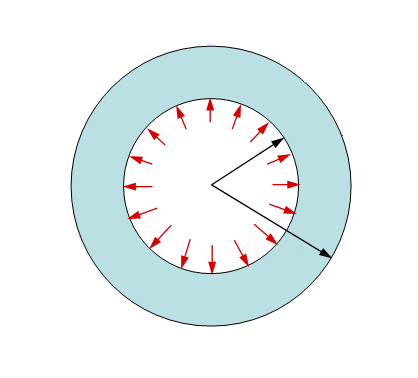}};
		\node[color=red] at (-0.9,0) [right] {$\text{P}_{\text{int}}$};
		\node[color=black] at (0.10,0.60) [right] {$\text{R}_{\text{i}}$};
		\node[color=black] at (0.05,-0.4) [right] {$\text{R}_{\text{o}}$};
		\node[color=black] at (-2.5,3) [right] {Stress Free Edge};
		\node[anchor=east] at (-1.5,2.9) (Start2) {};
		\node[anchor=west] at (-1.6,1.5) (HorizontalLine) {};
		\draw (Start2) edge[out=-90,in=120,->] (HorizontalLine);
		\end{tikzpicture}
		\begin{tikzpicture}
	
	\node[color=black] at (4.5,6.1) [right] {Stress Free Edges};
	\node[color=red] at (2.2,7.1) [right] {Applied Displacement Field};
	
	\definecolor{LightBlue}{rgb}{0.733,0.878,0.890}
	\begin{axis}[height=9.5cm,width=9.5cm, xmin=-4,xmax=4,scaled x ticks = false, ymin=-4,ymax=4,scaled y ticks = false,anchor=south west,axis line style={draw=none},tick style={draw=none},xticklabels={,,}, yticklabels={,,}, cells={anchor=west},  font=\footnotesize, rounded corners=2pt]
	
	\node[anchor=east] at (550,610) (Start1) {};
	\node[anchor=west] at (375,490) (VerticalLine) {};
	\draw (Start1) edge[out=-135,in=0,->] (VerticalLine);
	\node[anchor=east] at (590,610) (Start2) {};
	\node[anchor=west] at (450,390) (HorizontalLine) {};
	\draw (Start2) edge[out=-90,in=90,->] (HorizontalLine);
	
	\addplot[-stealth,line width=1pt,red]  table [x=X-1, y=Y-1, col sep=comma] {LShapedModel.csv};
	\addplot[-stealth,line width=1pt,red]  table [x=X-3, y=Y-3, col sep=comma] {LShapedModel.csv};
	\addplot[-stealth,line width=1pt,red]  table [x=X-5, y=Y-5, col sep=comma] {LShapedModel.csv};
	\addplot[-stealth,line width=1pt,red]  table [x=X-7, y=Y-7, col sep=comma] {LShapedModel.csv};
	\addplot[-stealth,line width=1pt,red]  table [x=X-9, y=Y-9, col sep=comma] {LShapedModel.csv};
	\addplot[-stealth,line width=1pt,red]  table [x=X-11, y=Y-11, col sep=comma] {LShapedModel.csv};
	\addplot[-stealth,line width=1pt,red]  table [x=X-12, y=Y-12, col sep=comma] {LShapedModel.csv};
	\addplot[-stealth,line width=1pt,red]  table [x=X-13, y=Y-13, col sep=comma] {LShapedModel.csv};
	\addplot[-stealth,line width=1pt,red]  table [x=X-14, y=Y-14, col sep=comma] {LShapedModel.csv};
	\addplot[-stealth,line width=1pt,red]  table [x=X-15, y=Y-15, col sep=comma] {LShapedModel.csv};
	\addplot[-stealth,line width=1pt,red]  table [x=X-16, y=Y-16, col sep=comma] {LShapedModel.csv};
	\addplot[-stealth,line width=1pt,red]  table [x=X-17, y=Y-17, col sep=comma] {LShapedModel.csv};
	\addplot[-stealth,line width=1pt,red]  table [x=X-18, y=Y-18, col sep=comma] {LShapedModel.csv};
	\addplot[-stealth,line width=1pt,red]  table [x=X-19, y=Y-19, col sep=comma] {LShapedModel.csv};
	\addplot[-stealth,line width=1pt,red]  table [x=X-20, y=Y-20, col sep=comma] {LShapedModel.csv};
	\addplot[-stealth,line width=1pt,red]  table [x=X-21, y=Y-21, col sep=comma] {LShapedModel.csv};
	\addplot[-stealth,line width=1pt,red]  table [x=X-22, y=Y-22, col sep=comma] {LShapedModel.csv};
	\addplot[-stealth,line width=1pt,red]  table [x=X-23, y=Y-23, col sep=comma] {LShapedModel.csv};
	\addplot[-stealth,line width=1pt,red]  table [x=X-24, y=Y-24, col sep=comma] {LShapedModel.csv};
	\addplot[-stealth,line width=1pt,red]  table [x=X-25, y=Y-25, col sep=comma] {LShapedModel.csv};
	\addplot[-stealth,line width=1pt,red]  table [x=X-26, y=Y-26, col sep=comma] {LShapedModel.csv};
	\addplot[-stealth,line width=1pt,red]  table [x=X-27, y=Y-27, col sep=comma] {LShapedModel.csv};
	\addplot[-stealth,line width=1pt,red]  table [x=X-28, y=Y-28, col sep=comma] {LShapedModel.csv};
	\addplot[-stealth,line width=1pt,red]  table [x=X-29, y=Y-29, col sep=comma] {LShapedModel.csv};
	\addplot[-stealth,line width=1pt,red]  table [x=X-30, y=Y-30, col sep=comma] {LShapedModel.csv};
	\addplot[-stealth,line width=1pt,red]  table [x=X-31, y=Y-31, col sep=comma] {LShapedModel.csv};
	\addplot[-stealth,line width=1pt,red]  table [x=X-32, y=Y-32, col sep=comma] {LShapedModel.csv};
	\addplot[-stealth,line width=1pt,red]  table [x=X-33, y=Y-33, col sep=comma] {LShapedModel.csv};
	\addplot[-stealth,line width=1pt,red]  table [x=X-34, y=Y-34, col sep=comma] {LShapedModel.csv};
	\addplot[-stealth,line width=1pt,red]  table [x=X-35, y=Y-35, col sep=comma] {LShapedModel.csv};
	\addplot[-stealth,line width=1pt,red]  table [x=X-36, y=Y-36, col sep=comma] {LShapedModel.csv};
	\addplot[-stealth,line width=1pt,red]  table [x=X-37, y=Y-37, col sep=comma] {LShapedModel.csv};
	\addplot[-stealth,line width=1pt,red]  table [x=X-38, y=Y-38, col sep=comma] {LShapedModel.csv};
	\addplot[-stealth,line width=1pt,red]  table [x=X-39, y=Y-39, col sep=comma] {LShapedModel.csv};
	\addplot[-stealth,line width=1pt,red]  table [x=X-40, y=Y-40, col sep=comma] {LShapedModel.csv};
	\addplot[-stealth,line width=1pt,red]  table [x=X-41, y=Y-41, col sep=comma] {LShapedModel.csv};
	\addplot[-stealth,line width=1pt,red]  table [x=X-42, y=Y-42, col sep=comma] {LShapedModel.csv};
	\addplot[-stealth,line width=1pt,red]  table [x=X-43, y=Y-43, col sep=comma] {LShapedModel.csv};
	\addplot[-stealth,line width=1pt,red]  table [x=X-44, y=Y-44, col sep=comma] {LShapedModel.csv};
	\addplot[-stealth,line width=1pt,red]  table [x=X-45, y=Y-45, col sep=comma] {LShapedModel.csv};
	\addplot[-stealth,line width=1pt,red]  table [x=X-46, y=Y-46, col sep=comma] {LShapedModel.csv};
	\addplot[-stealth,line width=1pt,red]  table [x=X-47, y=Y-47, col sep=comma] {LShapedModel.csv};
	\addplot[-stealth,line width=1pt,red]  table [x=X-48, y=Y-48, col sep=comma] {LShapedModel.csv};
	\addplot[-stealth,line width=1pt,red]  table [x=X-49, y=Y-49, col sep=comma] {LShapedModel.csv};
	\addplot[-stealth,line width=1pt,red]  table [x=X-50, y=Y-50, col sep=comma] {LShapedModel.csv};
	\addplot[-stealth,line width=1pt,red]  table [x=X-51, y=Y-51, col sep=comma] {LShapedModel.csv};
	\addplot[-stealth,line width=1pt,red]  table [x=X-53, y=Y-53, col sep=comma] {LShapedModel.csv};
	\addplot[-stealth,line width=1pt,red]  table [x=X-55, y=Y-55, col sep=comma] {LShapedModel.csv};
	\addplot[-stealth,line width=1pt,red]  table [x=X-57, y=Y-57, col sep=comma] {LShapedModel.csv};
	\addplot[-stealth,line width=1pt,red]  table [x=X-59, y=Y-59, col sep=comma] {LShapedModel.csv};
	\addplot[-stealth,line width=1pt,red]  table [x=X-61, y=Y-61, col sep=comma] {LShapedModel.csv};
	
	\addplot[patch,patch type=rectangle,color=LightBlue,faceted color=none] coordinates { (0,0) (-2,0) (-2,-2) (0,-2) };
	\addplot[patch,patch type=rectangle,color=LightBlue,faceted color=none] coordinates { (0,0) (2,0) (2,-2) (0,-2) };
	\addplot[patch,patch type=rectangle,color=LightBlue,faceted color=none] coordinates { (0,0) (-2,0) (-2,2) (0,2) };
	\addplot[black, sharp corners,line width=1.0pt]  table [x=X-Out, y=Y-Out, col sep=comma] {LShapedModel.csv};
	\end{axis}
\end{tikzpicture}
		\caption{Pressurized Cylinder (Left), L-Shape Domain in Mode I Loading (Right)}
		\label{2DProblemsConsidered}
	\end{figure}
	
	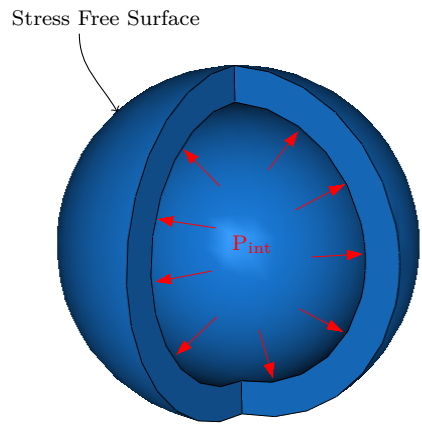
\begin{figure}[H] 
		\centering
		\begin{tikzpicture}
		\node at (0,0){\includegraphics[scale=0.5]{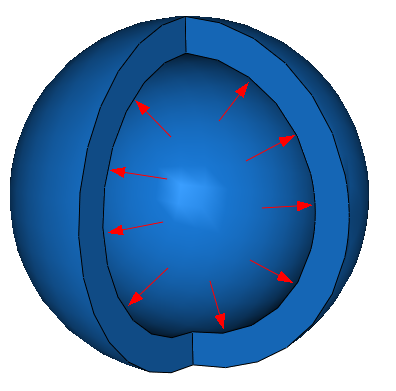}};
		\node[color=red] at (-0.3,0) [right] {$\text{P}_{\text{int}}$};
		\node[color=black] at (-3.2,3) [right] {Stress Free Surface};
		\node[anchor=east] at (-2.1,2.9) (Start2) {};
		\node[anchor=west] at (-1.75,1.65) (HorizontalLine) {};
		\draw (Start2) edge[out=-90,in=120,->] (HorizontalLine);
		\end{tikzpicture}
		\caption{Pressurized Sphere Model}
		\label{SphereDwg}
	\end{figure}
	
	The stress solution in terms of $\sigma_{11}$, $\sigma_{12}$, $\sigma_{22}$ for each problem is presented in Figure \ref{SolutionPipe}, Figure \ref{SolutionL-Shaped} and Figure \ref{SolutionSphere}, respectively, for the pressurized cylinder, the L-shape and the pressurized sphere. The equivalent von Mises stress (noted $\sigma_{VM}$) is also presented.
	\begin{figure}[H] 
		\centering
		\begin{tabular}{c}
			\text{Stress $\sigma_{11}$} \\
			\includegraphics[scale=0.4]{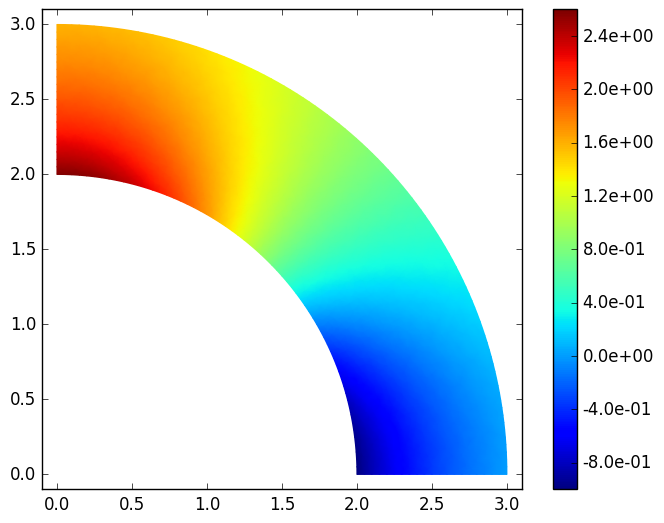}\\
		\end{tabular}
		\begin{tabular}{c}
			\text{Stress $\sigma_{12}$} \\
			\includegraphics[scale=0.4]{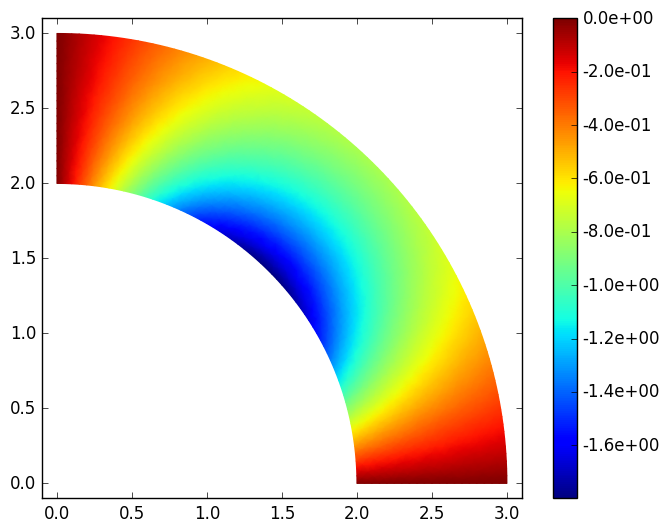}\\
		\end{tabular}
		\phantomcaption
	\end{figure}
	\begin{figure}[H] 
		\centering	
		\ContinuedFloat
		\begin{tabular}{c}
			\text{Stress $\sigma_{22}$}\\
			\includegraphics[scale=0.4]{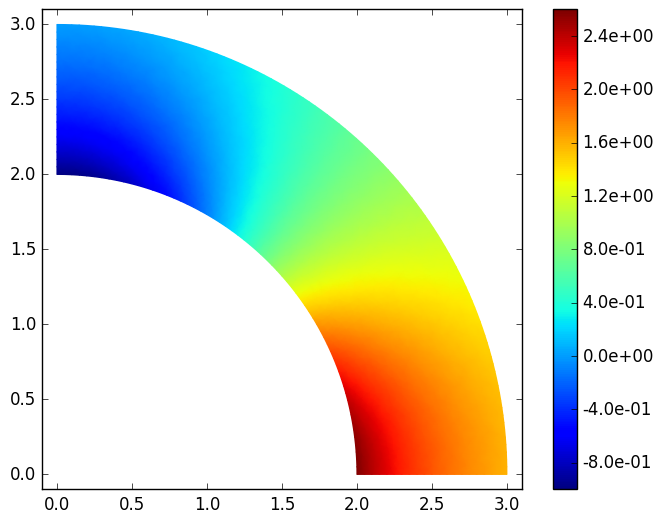} \\
		\end{tabular}
		\begin{tabular}{c}
			\text{Stress $\sigma_{VM}$ } \\
			\includegraphics[scale=0.4]{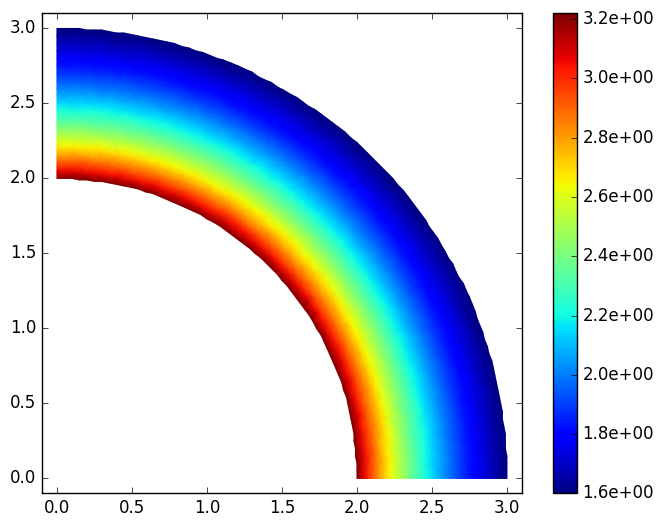}\\
		\end{tabular}
		\caption{Pressurized Cylinder - Stress Solution}
		\label{SolutionPipe}
	\end{figure}
	
	\setcounter{subfigure}{0}
	\begin{figure}[H] 
		\centering
		\subfloat[][]{%
			\begin{tabular}{c}
				\text{Stress $\sigma_{11}$} \\
				\includegraphics[scale=0.4]{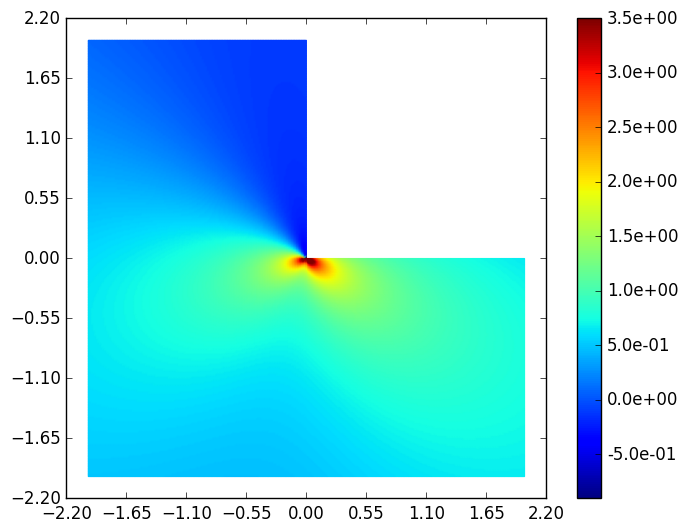}
			\end{tabular}
		}
		\subfloat[][]{%
			\begin{tabular}{c}
				\text{Stress $\sigma_{12}$} \\
				\includegraphics[scale=0.4]{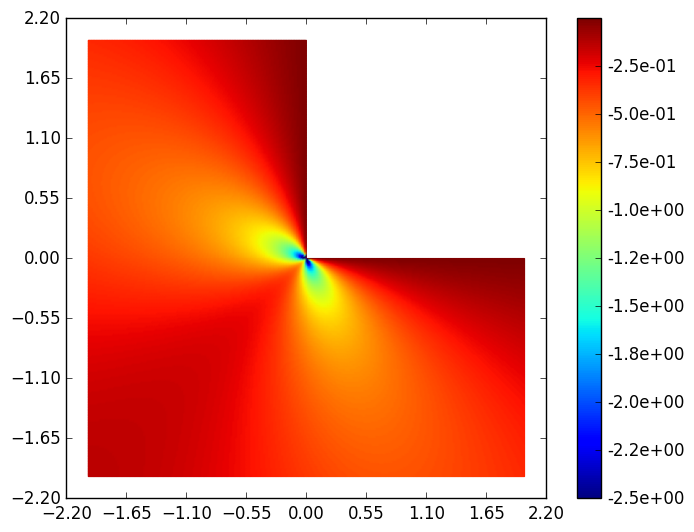}\\
			\end{tabular}
		}
		\phantomcaption
	\end{figure}
	
	\begin{figure}[H] 
		\centering
		\ContinuedFloat
		\subfloat[][]{%
			\begin{tabular}{c}
				\text{Stress $\sigma_{22}$} \\
				\includegraphics[scale=0.4]{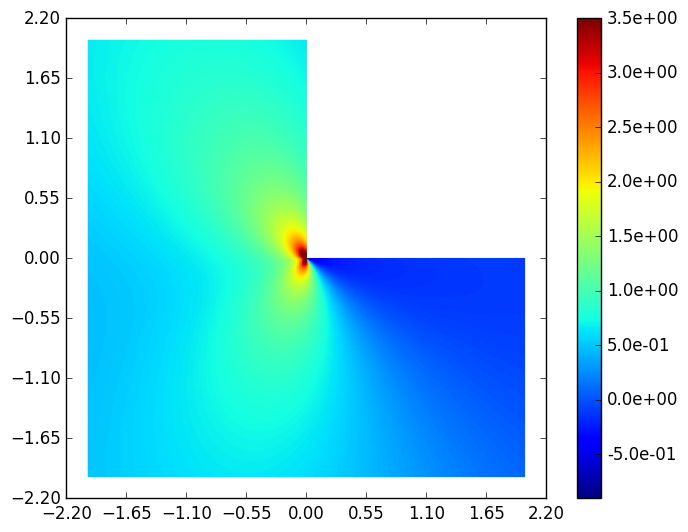} 
			\end{tabular} 
		} \qquad
		\subfloat[][]{%
			\begin{tabular}{c}
				\text{Stress $\sigma_{VM}$} \\
				\includegraphics[scale=0.4]{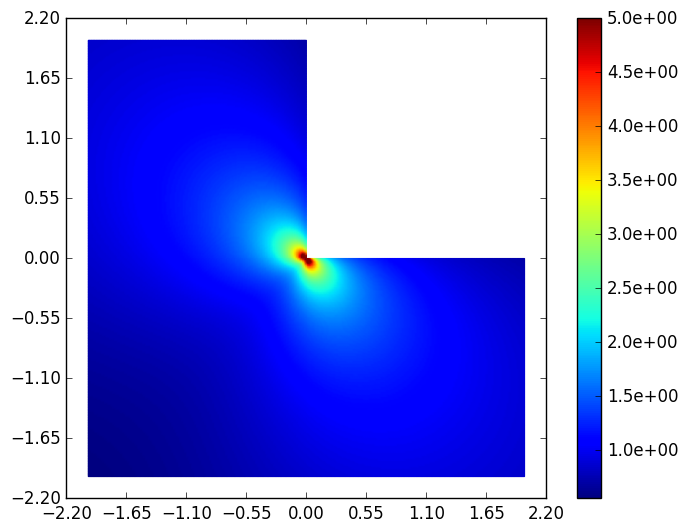}\\
			\end{tabular}
		}
		\caption{L-Shape Plate in Mode I Loading - Stress Solution}
		\label{SolutionL-Shaped}
	\end{figure}
	It should be noted that the stress solution tends to infinity at the interior corner of the L-shape domain. In order to represent the solution in Figure \ref{SolutionL-Shaped}, we have truncated the stress results around the singularity and represented truncated values in the same color as the selected threshold values.
	\setcounter{subfigure}{0}
	\begin{figure}[H] 
		\centering
		\subfloat[][]{%
			\begin{tabular}{c}
				\text{Stress $\sigma_{11}$} \\
				\includegraphics[scale=0.4]{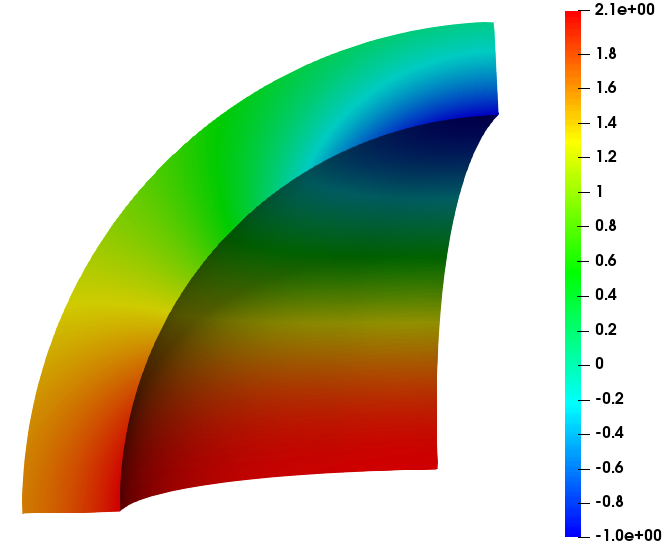} \\
			\end{tabular} 
		} 
		\subfloat[][]{%
			\begin{tabular}{c}
				\text{Stress $\sigma_{12}$} \\
				\includegraphics[scale=0.4]{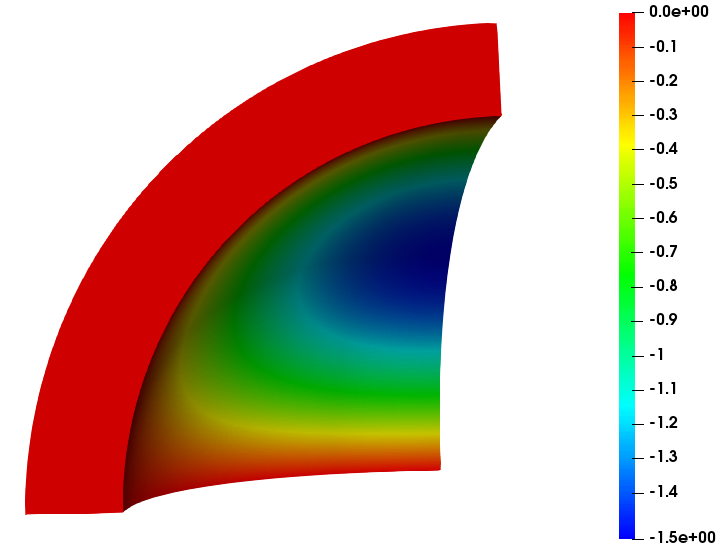} \\
			\end{tabular}
		}
		\phantomcaption
	\end{figure}
	
	\begin{figure}[H] 
		\centering
		\ContinuedFloat
		\subfloat[][]{%
			\begin{tabular}{c}
				\text{Stress $\sigma_{22}$} \\
				\includegraphics[scale=0.4]{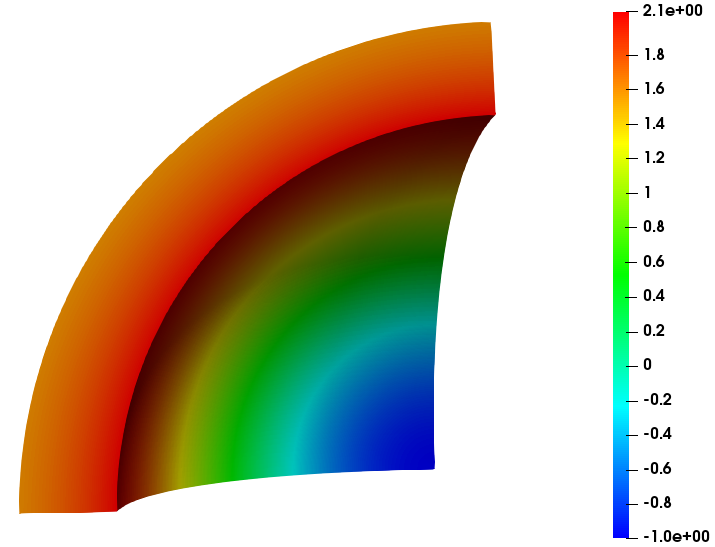} \\
			\end{tabular} 
		} 
		\subfloat[][]{%
			\begin{tabular}{c}
				\text{Stress $\sigma_{VM}$} \\
				\includegraphics[scale=0.4]{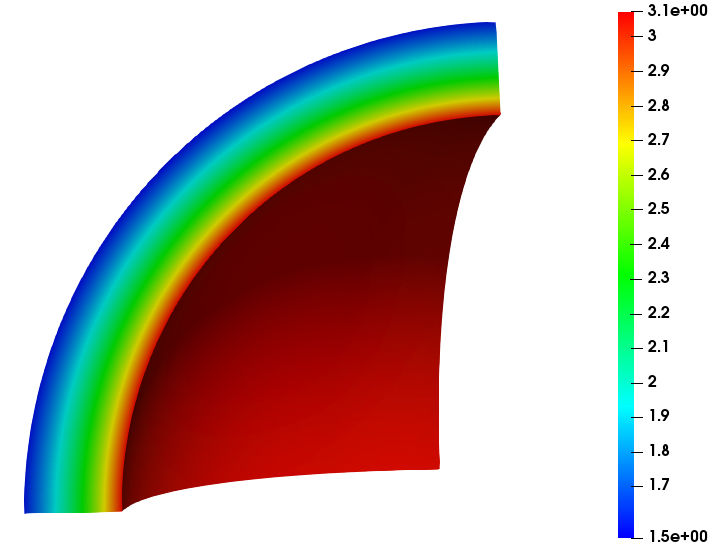} \\
			\end{tabular}
		}
		\caption{Pressurized Sphere - Stress Solution}
		\label{SolutionSphere}
	\end{figure}
	
	\subsection{Error Estimation}\label{ErrorNorms}
	
	To assess the influence of the considered parameters on the solution, and to compare the accuracy of the methods, a set of errors on the stress results has been calculated. For this purpose, a properly scaled $L_2$ error norm and a $L_{\infty}$ error norm have been selected. The $L_2$ error norm averages out the error over all the collocation nodes.
	
	Considering a domain $\Omega$ discretized with $n$ collocation nodes, where $\sigma_{ij}^{e}(\mathbf{X_k})$ and $\sigma_{ij}^{h}(\mathbf{X_k})$, respectively, represent the exact and the approximated stress values at a node $\mathbf{X_k}$, the $L_2$ error norm is calculated as follows:
	\begin{equation} \label{L_2Norm}
	L_{2}(\sigma_{ij})= \frac{\sqrt{\sum_{k=1}^{n}{(\sigma_{ij}^{e}(\mathbf{X_k}) - \sigma_{ij}^{h}(\mathbf{X_k}))^2}}}{n}.
	\end{equation}
	For the L-shape problem, the singular point has not been included in the $L_2$ error norm as the analytical solution diverges at this point.	
	The $L_{\infty}$ error norm corresponds to the maximum absolute error observed over the considered domain. 
	\begin{equation} \label{L_InfinityNorm}
	L_{\infty}(\sigma_{ij})= \max_{k \in \Omega} \big(|\sigma_{ij}^{e}(\mathbf{X_k}) - \sigma_{ij}^{h}(\mathbf{X_k}) | \big).
	\end{equation}

	\section{Parametric Study}\label{ParametricAnalysis}
	
	\subsection{General}
	
	In order to better understand the GFD and the DC PSE methods, and to assess the impact of the various parameters on the solution, we performed a parametric study. The purpose of this study is to select a single set of parameters that can be applied to any problem without knowing a priori the solution type. The following items have been considered:
	\begin{itemize}
		\item The selected weight function;
		\item The selected correction function basis for the DC PSE method;
		\item The number of support nodes.
	\end{itemize}
	The 2D cylinder, the 2D L-shape and the 3D sphere models have been considered for this study. A regular distribution of 5,372 nodes and of 13,735 nodes has been selected for the 2D problems, respectively. A regular distribution of 83,174 nodes has been selected for the 3D problem. The DCPSE1 variation of the DC PSE method has been selected for the purpose of this sensitivity study. The sensitivity analyses for the weight function and the DC PSE correction function have been performed for the 2D model only. The impact of the number of support nodes on the error has been assessed for both, the 2D and 3D problems. Due to the symmetries of the models, results are only presented in terms of error on the $\sigma_{11}$ and $\sigma_{12}$ components of the stress tensor.
	
	\subsection{Weight Function Sensitivity}
	
	\paragraph{GFD}\
	
\begin{figure}[H] 
		\centering
		\begin{tabular}{c:c}
			\begin{tikzpicture}
			\begin{axis}[height=7cm,width=7.5cm, ymin=0,ymax=1.4,ytick={0,0.2,0.4,0.6,0.8,1,1.2},xmin=0,xmax=1.3, legend entries={($\text{3}^{\text{rd}}$ Order Spline)$^{0.4}$,($\text{3}^{\text{rd}}$ Order Spline)$^{0.6}$,($\text{3}^{\text{rd}}$ Order Spline)$^{0.8}$,($\text{3}^{\text{rd}}$ Order Spline)$^{1.0}$,($\text{3}^{\text{rd}}$ Order Spline)$^{1.2}$,Linear Function},legend style={legend columns=1, cells={anchor=west},  font=\footnotesize, rounded corners=2pt,}, legend pos=north east,xlabel=Normalized Distance to Ref. Node,ylabel=Weight Value]
			\addplot+[green,mark=square*,mark options={fill=green}]   table [x=DIST-NORM, y=S3-0.4
			, col sep=comma] {WeightFunctionsGFD.csv};
			\addplot+[red,mark=diamond*,mark options={fill=red}]   table [x=DIST-NORM, y=S3-0.6
			, col sep=comma] {WeightFunctionsGFD.csv};
			\addplot+[blue,mark=triangle*,mark options={fill=blue}]   table [x=DIST-NORM, y=S3-0.8
			, col sep=comma] {WeightFunctionsGFD.csv};
			\addplot+[yellow,mark=*,mark options={fill=yellow}]   table [x=DIST-NORM, y=S3-1.0
			, col sep=comma] {WeightFunctionsGFD.csv};
			\addplot+[olive,mark=star,mark options={fill=olive}]   table [x=DIST-NORM, y=S3-1.2
			, col sep=comma] {WeightFunctionsGFD.csv};
			\addplot+[black,mark=pentagon*,mark options={fill=black}] table [x=DIST-NORM, y=LINEAR-FUNC
			, col sep=comma] {WeightFunctionsGFD.csv};
			\end{axis}
			\end{tikzpicture} &
			
			\begin{tikzpicture}
			\begin{axis}[height=7cm,width=7.5cm, ymin=0,ymax=1.4,ytick={0,0.2,0.4,0.6,0.8,1,1.2},xmin=0,xmax=1.3, legend entries={($\text{4}^{\text{rd}}$ Order Spline)$^{0.4}$,($\text{4}^{\text{rd}}$ Order Spline)$^{0.6}$,($\text{4}^{\text{rd}}$ Order Spline)$^{0.8}$,($\text{4}^{\text{rd}}$ Order Spline)$^{1.0}$,($\text{4}^{\text{rd}}$ Order Spline)$^{1.2}$},legend style={legend columns=1, cells={anchor=west},  font=\footnotesize, rounded corners=2pt,}, legend pos=north east,xlabel=Normalized Distance to Ref. Node,ylabel=Weight Value]
			\addplot+[green,mark=square*,mark options={fill=green}]   table [x=DIST-NORM, y=S4-0.4
			, col sep=comma] {WeightFunctionsGFD.csv};
			\addplot+[red,mark=diamond*,mark options={fill=red}]   table [x=DIST-NORM, y=S4-0.6
			, col sep=comma] {WeightFunctionsGFD.csv};
			\addplot+[blue,mark=triangle*,mark options={fill=blue}]   table [x=DIST-NORM, y=S4-0.8
			, col sep=comma] {WeightFunctionsGFD.csv};
			\addplot+[yellow,mark=*,mark options={fill=yellow}]   table [x=DIST-NORM, y=S4-1.0
			, col sep=comma] {WeightFunctionsGFD.csv};
			\addplot+[olive,mark=star,mark options={fill=olive}]   table [x=DIST-NORM, y=S4-1.2
			, col sep=comma] {WeightFunctionsGFD.csv};
			\end{axis}
			\end{tikzpicture}
		\end{tabular}
		\caption{GFD Weight Functions: 3$^\text{rd}$ Order Spline (Left) and 4$^\text{th}$ Order Spline (Right) for power parameters ranging from 0.4 to 1.2. The linear weight function is also given for reference purposes.}
		\label{WeightFunctionsGFD}
	\end{figure}
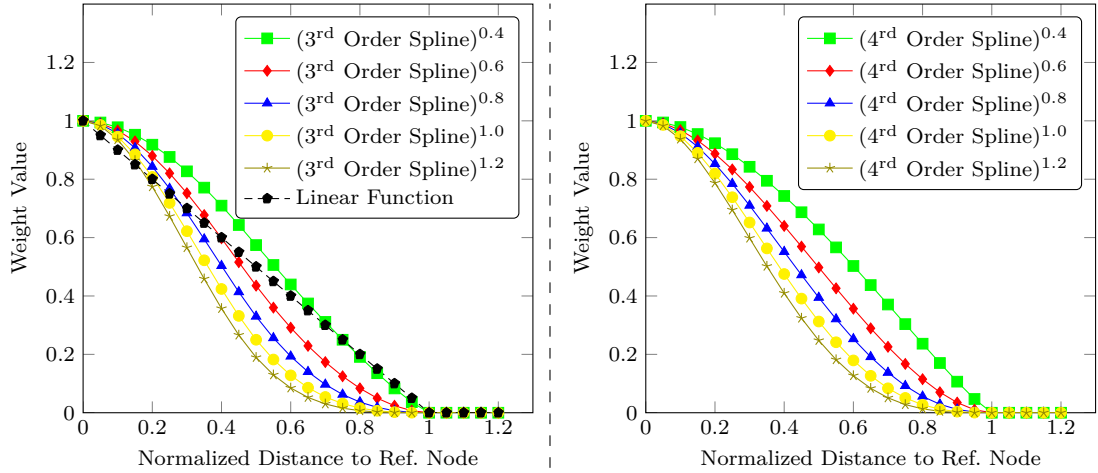

	Various functions can be considered for the weight function introduced in Equation (\ref{FunctionalB_AllTerms_GFD}) and Equation (\ref{Moments_GFD}). The $3^\text{rd}$ and $4^\text{th}$ order splines are the preferred types (see Equation (\ref{Spline3Weight}) and  Equation (\ref{Spline4Weight})). In order to vary the shape of the weight functions, we have composed the splines with the following power function
	\begin{equation} \label{GFDWeightForm}
	W(s)=(w(s))^{\gamma},
	\end{equation}
	where $w$ is the spline function and $W$ is the modified weight function. We have compared the results obtained for power parameters $\gamma$ between 0.4 and 1.2.
	We have also considered a linear weight function for comparison purposes. The equation of this function is given by
	\begin{equation} \label{GFDLinearWeight}
	w(s)=
	\begin{cases}
	1-s & \text{ \quad if } s \leq 1 \\
	0  & \text{ \quad if } s > 1. \\
	\end{cases}\\
	\end{equation}
	The shapes of the considered functions are presented in Figure \ref{WeightFunctionsGFD}. 
	
	The $L_2$ and $L_{\infty}$ errors obtained with the considered weight functions are presented in Figure \ref{WeightFunctionsResultsGFD} and Figure \ref{WeightFunctionsResultsGFD_LShaped} for the 2D cylinder and for the 2D L-shape, respectively.
	The error for the linear weight function is also presented in these figures for comparison purposes, even though a power parameter is not used in this function.
	\begin{figure}[H] 
		\centering
		
		\begin{tabular}{c:c}
			\begin{tikzpicture}
			\begin{axis}[height=6cm,width=7.5cm,ymode=log, ymin=0.000001,ymax=0.01,xmin=0.15,xmax=1.25, legend entries={$\text{3}^{\text{rd}}$ Order Spline,$\text{4}^{\text{th}}$ Order Spline,Linear Function},legend style={legend columns=1, cells={anchor=west},  font=\footnotesize, rounded corners=2pt,}, legend pos=north east,xlabel=Power Parameter,ylabel=$L_2$ Error - $\sigma_{11}$]
			\addplot+[green,mark=square*,mark options={fill=green}]   table [x=Exponent, y=S3-L2-REL-S11
			, col sep=comma] {WeightFunctionsGFD.csv};
			\addplot+[red,mark=diamond*,mark options={fill=red}]   table [x=Exponent, y=S4-L2-REL-S11
			, col sep=comma] {WeightFunctionsGFD.csv};
			\addplot+[blue,mark=triangle*,mark options={fill=blue}]   table [x=Exponent-Linear, y=LIN-L2-REL-S11
			, col sep=comma] {WeightFunctionsGFD.csv};
			\end{axis}
			\end{tikzpicture} & 
			
			\begin{tikzpicture}
			\begin{axis}[height=6cm,width=7.5cm,ymode=log, ymin=0.0001,ymax=10,xmin=0.15,xmax=1.25, legend entries={$\text{3}^{\text{rd}}$ Order Spline,$\text{4}^{\text{th}}$ Order Spline,Linear Function},legend style={legend columns=1, cells={anchor=west},  font=\footnotesize, rounded corners=2pt,}, legend pos=north east,xlabel=Power Parameter,ylabel=$L_{\infty}$ Error - $\sigma_{11}$]
			\addplot+[green,mark=square*,mark options={fill=green}]   table [x=Exponent, y=S3-L-INF-S11
			, col sep=comma] {WeightFunctionsGFD.csv};
			\addplot+[red,mark=diamond*,mark options={fill=red}]   table [x=Exponent, y=S4-L-INF-S11
			, col sep=comma] {WeightFunctionsGFD.csv};
			\addplot+[blue,mark=triangle*,mark options={fill=blue}]   table [x=Exponent-Linear, y=LIN-L-INF-S11
			, col sep=comma] {WeightFunctionsGFD.csv};
			\end{axis}
			\end{tikzpicture} \\

			\begin{tikzpicture}
			\begin{axis}[height=6cm,width=7.5cm,ymode=log, ymin=0.000001,ymax=0.01,xmin=0.15,xmax=1.25, legend entries={$\text{3}^{\text{rd}}$ Order Spline,$\text{4}^{\text{th}}$ Order Spline,Linear Function},legend style={legend columns=1, cells={anchor=west},  font=\footnotesize, rounded corners=2pt,}, legend pos=north east,xlabel=Power Parameter,ylabel=$L_2$ Error - $\sigma_{12}$]
			\addplot+[green,mark=square*,mark options={fill=green}]   table [x=Exponent, y=S3-L2-REL-S12
			, col sep=comma] {WeightFunctionsGFD.csv};
			\addplot+[red,mark=diamond*,mark options={fill=red}]   table [x=Exponent, y=S4-L2-REL-S12
			, col sep=comma] {WeightFunctionsGFD.csv};
			\addplot+[blue,mark=triangle*,mark options={fill=blue}]   table [x=Exponent-Linear, y=LIN-L2-REL-S12
			, col sep=comma] {WeightFunctionsGFD.csv};
			\end{axis}
			\end{tikzpicture} &
			
			\begin{tikzpicture}
			\begin{axis}[height=6cm,width=7.5cm,ymode=log, ymin=0.0001,ymax=10,xmin=0.15,xmax=1.25, legend entries={$\text{3}^{\text{rd}}$ Order Spline,$\text{4}^{\text{th}}$ Order Spline,Linear Function},legend style={legend columns=1, cells={anchor=west},  font=\footnotesize, rounded corners=2pt,}, legend pos=north east,xlabel=Power Parameter,ylabel=$L_{\infty}$ Error - $\sigma_{12}$]
			\addplot+[green,mark=square*,mark options={fill=green}]   table [x=Exponent, y=S3-L-INF-S12
			, col sep=comma] {WeightFunctionsGFD.csv};
			\addplot+[red,mark=diamond*,mark options={fill=red}]   table [x=Exponent, y=S4-L-INF-S12
			, col sep=comma] {WeightFunctionsGFD.csv};
			\addplot+[blue,mark=triangle*,mark options={fill=blue}]   table [x=Exponent-Linear, y=LIN-L-INF-S12
			, col sep=comma] {WeightFunctionsGFD.csv};
			\end{axis}
			\end{tikzpicture}\\
			
		\end{tabular}
		
		\caption{GFD Weight Sensitivity - 2D Cylinder. $L_2$ (Left) and $L_{\infty}$ (Right) errors. Comparison for $\text{3}^{\text{rd}}$ and $\text{4}^{\text{th}}$ order splines weight functions composed with a power function of various exponents. The linear weight function is also given for reference purposes. The 4$^\text{th}$ order spline consistently leads to a low error.}
		\label{WeightFunctionsResultsGFD}
	\end{figure}
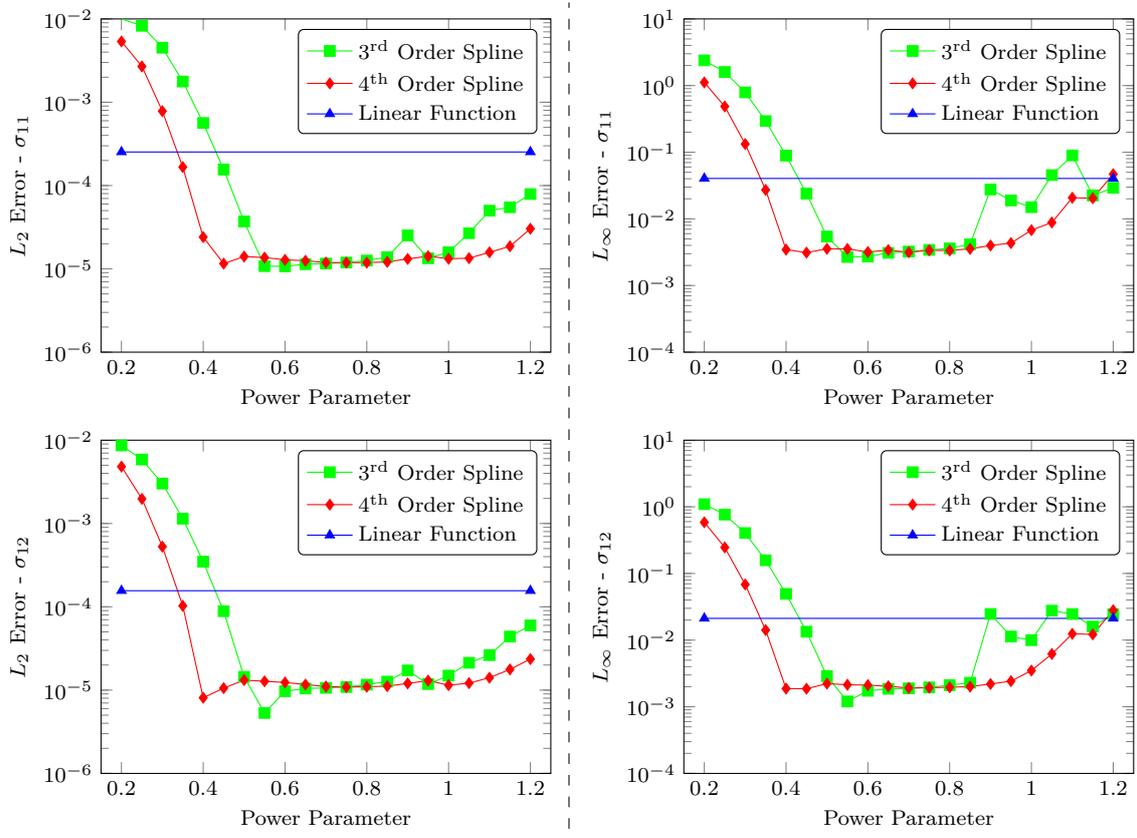
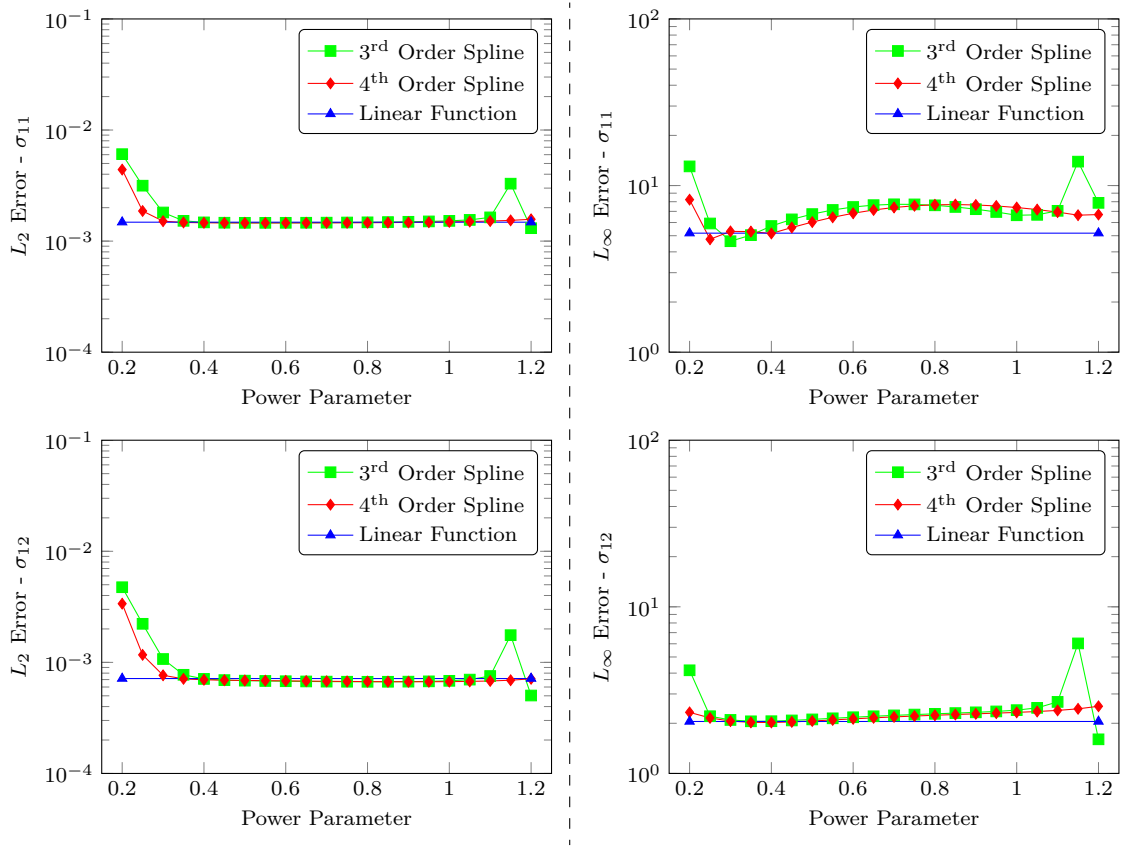
\begin{figure}[H] 
	\centering
	
	\begin{tabular}{c:c}
		\begin{tikzpicture}
		\begin{axis}[height=6cm,width=7.5cm,ymode=log, ymin=0.0001,ymax=0.1,xmin=0.15,xmax=1.25, legend entries={$\text{3}^{\text{rd}}$ Order Spline,$\text{4}^{\text{th}}$ Order Spline,Linear Function},legend style={legend columns=1, cells={anchor=west},  font=\footnotesize, rounded corners=2pt,}, legend pos=north east,xlabel=Power Parameter,ylabel=$L_2$ Error - $\sigma_{11}$]
		\addplot+[green,mark=square*,mark options={fill=green}]   table [x=Exponent, y=S3-L2-REL-S11
		, col sep=comma] {WeightFunctionsGFD-LShaped.csv};
		\addplot+[red,mark=diamond*,mark options={fill=red}]   table [x=Exponent, y=S4-L2-REL-S11
		, col sep=comma] {WeightFunctionsGFD-LShaped.csv};
		\addplot+[blue,mark=triangle*,mark options={fill=blue}]   table [x=Exponent-Linear, y=LIN-L2-REL-S11
		, col sep=comma] {WeightFunctionsGFD-LShaped.csv};
		\end{axis}
		\end{tikzpicture} & 
		
		\begin{tikzpicture}
		\begin{axis}[height=6cm,width=7.5cm,ymode=log, ymin=1,ymax=100,xmin=0.15,xmax=1.25, legend entries={$\text{3}^{\text{rd}}$ Order Spline,$\text{4}^{\text{th}}$ Order Spline,Linear Function},legend style={legend columns=1, cells={anchor=west},  font=\footnotesize, rounded corners=2pt,}, legend pos=north east,xlabel=Power Parameter,ylabel=$L_{\infty}$ Error - $\sigma_{11}$]
		\addplot+[green,mark=square*,mark options={fill=green}]   table [x=Exponent, y=S3-L-INF-S11
		, col sep=comma] {WeightFunctionsGFD-LShaped.csv};
		\addplot+[red,mark=diamond*,mark options={fill=red}]   table [x=Exponent, y=S4-L-INF-S11
		, col sep=comma] {WeightFunctionsGFD-LShaped.csv};
		\addplot+[blue,mark=triangle*,mark options={fill=blue}]   table [x=Exponent-Linear, y=LIN-L-INF-S11
		, col sep=comma] {WeightFunctionsGFD-LShaped.csv};
		\end{axis}
		\end{tikzpicture} \\

		\begin{tikzpicture}
		\begin{axis}[height=6cm,width=7.5cm,ymode=log, ymin=0.0001,ymax=0.1,xmin=0.15,xmax=1.25, legend entries={$\text{3}^{\text{rd}}$ Order Spline,$\text{4}^{\text{th}}$ Order Spline,Linear Function},legend style={legend columns=1, cells={anchor=west},  font=\footnotesize, rounded corners=2pt,}, legend pos=north east,xlabel=Power Parameter,ylabel=$L_2$ Error - $\sigma_{12}$]
		\addplot+[green,mark=square*,mark options={fill=green}]   table [x=Exponent, y=S3-L2-REL-S12
		, col sep=comma] {WeightFunctionsGFD-LShaped.csv};
		\addplot+[red,mark=diamond*,mark options={fill=red}]   table [x=Exponent, y=S4-L2-REL-S12
		, col sep=comma] {WeightFunctionsGFD-LShaped.csv};
		\addplot+[blue,mark=triangle*,mark options={fill=blue}]   table [x=Exponent-Linear, y=LIN-L2-REL-S12
		, col sep=comma] {WeightFunctionsGFD-LShaped.csv};
		\end{axis}
		\end{tikzpicture} &
		
		\begin{tikzpicture}
		\begin{axis}[height=6cm,width=7.5cm,ymode=log, ymin=1,ymax=100,xmin=0.15,xmax=1.25, legend entries={$\text{3}^{\text{rd}}$ Order Spline,$\text{4}^{\text{th}}$ Order Spline,Linear Function},legend style={legend columns=1, cells={anchor=west},  font=\footnotesize, rounded corners=2pt,}, legend pos=north east,xlabel=Power Parameter,ylabel=$L_{\infty}$ Error - $\sigma_{12}$]
		\addplot+[green,mark=square*,mark options={fill=green}]   table [x=Exponent, y=S3-L-INF-S12
		, col sep=comma] {WeightFunctionsGFD-LShaped.csv};
		\addplot+[red,mark=diamond*,mark options={fill=red}]   table [x=Exponent, y=S4-L-INF-S12
		, col sep=comma] {WeightFunctionsGFD-LShaped.csv};
		\addplot+[blue,mark=triangle*,mark options={fill=blue}]   table [x=Exponent-Linear, y=LIN-L-INF-S12
		, col sep=comma] {WeightFunctionsGFD-LShaped.csv};
		\end{axis}
		\end{tikzpicture}\\
		
	\end{tabular}
	
	\caption{GFD Weight Sensitivity - 2D L-Shape. $L_2$ (Left) and $L_{\infty}$ (Right) errors. Comparison for $\text{3}^{\text{rd}}$ and $\text{4}^{\text{th}}$ order splines weight functions composed with a power function of various exponents. The linear weight function is also given for reference purposes. The $\text{4}^{\text{th}}$ order spline leads to a low error than the  $\text{3}^{\text{rd}}$ order spline for both stress components and both error norms. The linear function leads to the lowest error in terms of $L_{\infty}$ error norm.}
	\label{WeightFunctionsResultsGFD_LShaped}
\end{figure}
	We can see from Figure \ref{WeightFunctionsResultsGFD} that within the range [0.6; 0.9] the power parameter has little impact on the error. In this range, the type of spline used does not significantly impact the error either. From Figure \ref{WeightFunctionsResultsGFD_LShaped} we can see that the power parameter has little impact on the error in terms of $L_2$ norm. A linear weight function leads to similar results as the spline weight functions is the range of power parameter [0.3;1.1]. In terms of $L_{\infty}$ error norm, the linear function leads to a lower error than the spline functions in the range of power parameter [0.4;1.2].

	A 4$^\text{th}$ order spline with a power parameter of 0.75 appears to be a reasonable choice as it leads to a minimum error for both considered problem. It leads to a low error for the 2D cylinder without being too close to the rapid error increase that is observed when the power parameter decreases below 0.6. It also leads to a reasonably low error for the singular problem both in terms of $L_2$ and $L_{\infty}$ error norms. A unique set of parameters has been selected for both problems in order to be applied to a wide variety of problems in the domain of linear elasticity.
	
	\paragraph{DC PSE} \
	
\begin{figure}[H] 
		\centering
		\begin{tikzpicture}
		\begin{axis}[height=7cm,width=10cm, ymin=0,ymax=1,xmin=0,xmax=1, legend entries={Exp.=1.0 - Shape=0.30,Exp.=2.0 - Shape=0.33,Exp.=3.0 - Shape=0.40,$\text{3}^{\text{rd}}$ Order Spline,$\text{4}^{\text{th}}$ Order Spline},legend style={legend columns=1, cells={anchor=west},  font=\footnotesize, rounded corners=2pt,}, legend pos=north east,xlabel=Normalized Distance to Ref. Node,ylabel=Weight Value]
		\addplot+[green,mark=square*,mark options={fill=green}]   table [x=DIST-NORM, y=EXP1-0.30-W
		, col sep=comma] {WeightFunctionsDCPSE.csv};
		\addplot+[red,mark=diamond*,mark options={fill=red}]   table [x=DIST-NORM, y=EXP2-0.33-W
		, col sep=comma] {WeightFunctionsDCPSE.csv};
		\addplot+[blue,mark=triangle*,mark options={fill=blue}]   table [x=DIST-NORM, y=EXP3-0.4-W
		, col sep=comma] {WeightFunctionsDCPSE.csv};
		\addplot+[yellow,mark=*,mark options={fill=yellow}]   table [x=DIST-NORM, y=SPLINE3-W
		, col sep=comma] {WeightFunctionsDCPSE.csv};
		\addplot+[olive,mark=star,mark options={fill=olive}]   table [x=DIST-NORM, y=SPLINE4-W
		, col sep=comma] {WeightFunctionsDCPSE.csv};
		\end{axis}
		\end{tikzpicture}
		\caption{DC PSE Weight Functions: Comparison of the profile of typical exponential and spline functions used as weight in the DC PSE approximation.}
		\label{WindowFunctionShapes}
	\end{figure}
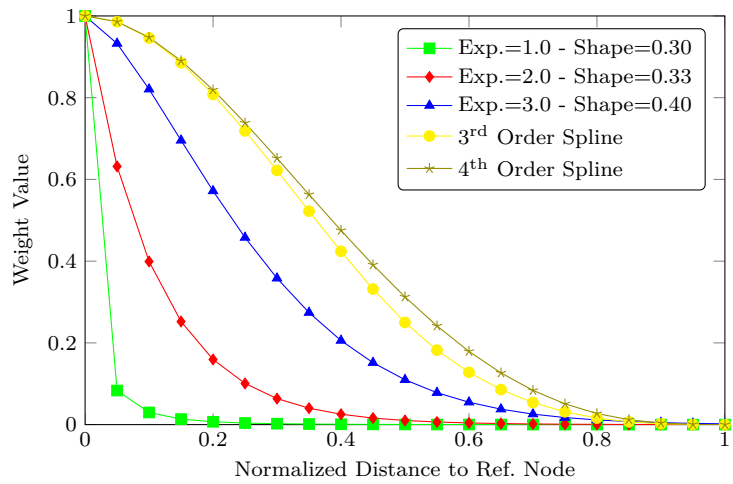

	In this Section, we assess the influence of the selected weight function on the error for the DC PSE method. The exponential weight function presented in Equation (\ref{ExpWeight}) is compared to the 3$^\text{rd}$ and 4$^\text{th}$ order splines, of equations (\ref{Spline3Weight}) and (\ref{Spline4Weight}) respectively. We have considered various combinations of shape parameters ($\epsilon$) and exponents ($\alpha$). In Figure \ref{WindowFunctionShapes}, the profile of three exponential weight functions, along with the two splines functions, is presented.

	In Figure \ref{WeightFunctionsResultsDCPSE} and Figure \ref{WeightFunctionsResultsDCPSE_LShaped}, the error in terms of $L_2$ and $L_{\infty}$ norms for the $\sigma_{11}$ and $\sigma_{12}$ stress components is presented for various combinations of exponents and shape parameters. The error obtained with the $3^\text{rd}$ and $4^\text{th}$ order splines is also presented on the graphs for comparison purposes, even though the shape parameter does not apply for these functions.

	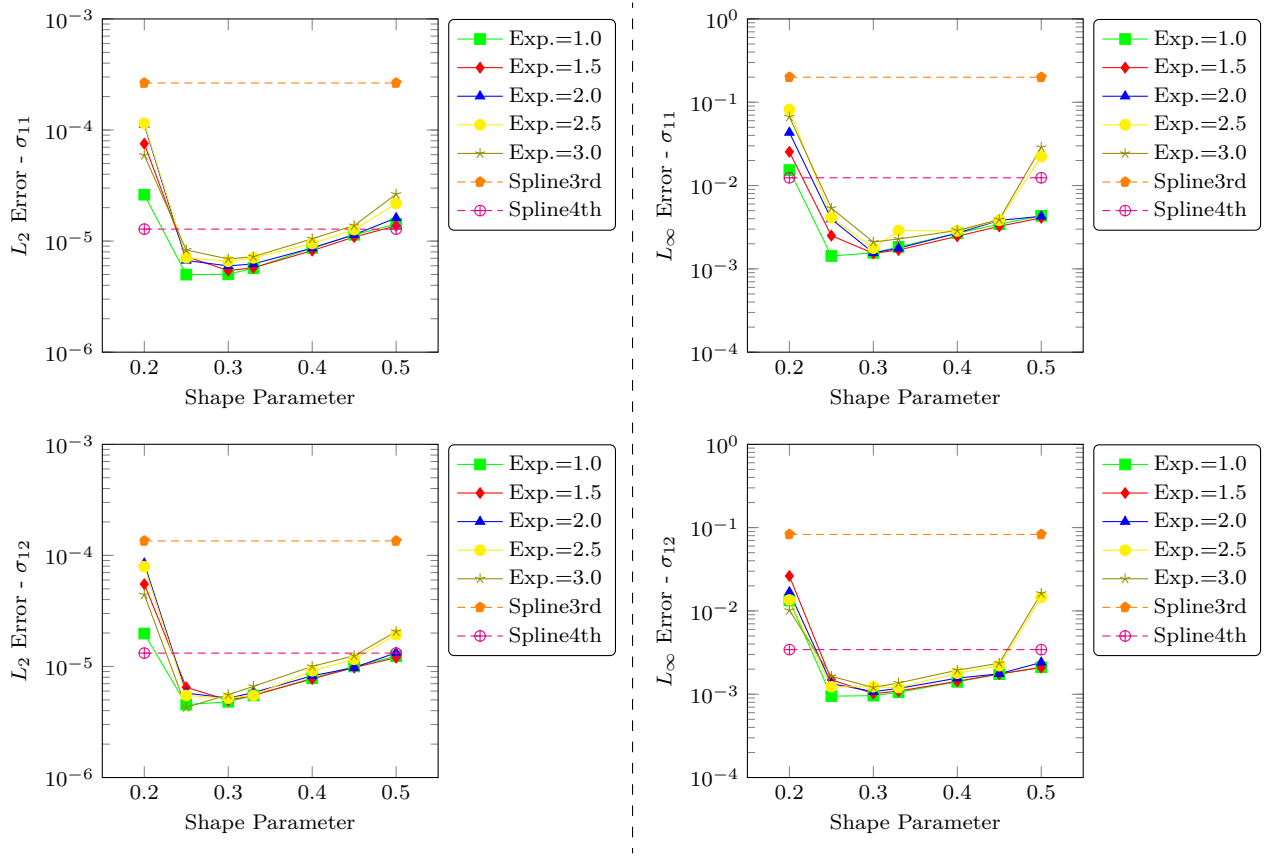
\begin{figure}[H] 
		\centering
		\begin{tabular}{c:c}
			\begin{tikzpicture}
			\begin{axis}[height=6cm,width=6cm,ymode=log, ymin=0.000001,ymax=0.001,xmin=0.15,xmax=0.55, legend entries={Exp.=1.0,Exp.=1.5,Exp.=2.0,Exp.=2.5,Exp.=3.0,Spline3rd ,Spline4th},legend style={legend columns=1, cells={anchor=west},  font=\footnotesize, rounded corners=2pt,}, legend pos=outer north east,xlabel=Shape Parameter,ylabel=$L_2$ Error - $\sigma_{11}$]
			\addplot+[green,mark=square*,mark options={fill=green}]   table [x=Exponent, y=EXP-1.0-L2-REL-S11
			, col sep=comma] {WeightFunctionsDCPSE.csv};
			\addplot+[red,mark=diamond*,mark options={fill=red}]   table [x=Exponent, y=EXP-1.5-L2-REL-S11
			, col sep=comma] {WeightFunctionsDCPSE.csv};
			\addplot+[blue,mark=triangle*,mark options={fill=blue}]   table [x=Exponent, y=EXP-2.0-L2-REL-S11
			, col sep=comma] {WeightFunctionsDCPSE.csv};
			\addplot+[yellow,mark=*,mark options={fill=yellow}]   table [x=Exponent, y=EXP-2.5-L2-REL-S11
			, col sep=comma] {WeightFunctionsDCPSE.csv};
			\addplot+[olive,mark=star,mark options={fill=olive}]   table [x=Exponent, y=EXP-3.0-L2-REL-S11
			, col sep=comma] {WeightFunctionsDCPSE.csv};
			\addplot+[orange,mark=pentagon*,mark options={fill=orange}]   table [x=Exponent-Spline, y=SPLINE3-L2-REL-S11
			, col sep=comma] {WeightFunctionsDCPSE.csv};
			\addplot+[magenta,mark=oplus]   table [x=Exponent-Spline, y=SPLINE4-L2-REL-S11
			, col sep=comma] {WeightFunctionsDCPSE.csv};
			\end{axis}
			\end{tikzpicture} & 
			
			\begin{tikzpicture}
			\begin{axis}[height=6cm,width=6cm,ymode=log, ymin=0.0001,ymax=1,xmin=0.15,xmax=0.55, legend entries={Exp.=1.0,Exp.=1.5,Exp.=2.0,Exp.=2.5,Exp.=3.0,Spline3rd ,Spline4th},legend style={legend columns=1, cells={anchor=west},  font=\footnotesize, rounded corners=2pt,}, legend pos=outer north east,xlabel=Shape Parameter,ylabel=$L_{\infty}$ Error - $\sigma_{11}$]
			\addplot+[green,mark=square*,mark options={fill=green}]   table [x=Exponent, y=EXP-1.0-L-INF-S11
			, col sep=comma] {WeightFunctionsDCPSE.csv};
			\addplot+[red,mark=diamond*,mark options={fill=red}]   table [x=Exponent, y=EXP-1.5-L-INF-S11
			, col sep=comma] {WeightFunctionsDCPSE.csv};
			\addplot+[blue,mark=triangle*,mark options={fill=blue}]   table [x=Exponent, y=EXP-2.0-L-INF-S11
			, col sep=comma] {WeightFunctionsDCPSE.csv};
			\addplot+[yellow,mark=*,mark options={fill=yellow}]   table [x=Exponent, y=EXP-2.5-L-INF-S11
			, col sep=comma] {WeightFunctionsDCPSE.csv};
			\addplot+[olive,mark=star,mark options={fill=olive}]   table [x=Exponent, y=EXP-3.0-L-INF-S11
			, col sep=comma] {WeightFunctionsDCPSE.csv};
			\addplot+[orange,mark=pentagon*,mark options={fill=orange}]   table [x=Exponent-Spline, y=SPLINE3-L-INF-S11
			, col sep=comma] {WeightFunctionsDCPSE.csv};
			\addplot+[magenta,mark=oplus]   table [x=Exponent-Spline, y=SPLINE4-L-INF-S11
			, col sep=comma] {WeightFunctionsDCPSE.csv};
			\end{axis}
			\end{tikzpicture} \\
			
			\begin{tikzpicture}
			\begin{axis}[height=6cm,width=6cm,ymode=log, ymin=0.000001,ymax=0.001,xmin=0.15,xmax=0.55, legend entries={Exp.=1.0,Exp.=1.5,Exp.=2.0,Exp.=2.5,Exp.=3.0,Spline3rd ,Spline4th},legend style={legend columns=1, cells={anchor=west},  font=\footnotesize, rounded corners=2pt,}, legend pos=outer north east,xlabel=Shape Parameter,ylabel=$L_2$ Error - $\sigma_{12}$]
			\addplot+[green,mark=square*,mark options={fill=green}]   table [x=Exponent, y=EXP-1.0-L2-REL-S12
			, col sep=comma] {WeightFunctionsDCPSE.csv};
			\addplot+[red,mark=diamond*,mark options={fill=red}]   table [x=Exponent, y=EXP-1.5-L2-REL-S12
			, col sep=comma] {WeightFunctionsDCPSE.csv};
			\addplot+[blue,mark=triangle*,mark options={fill=blue}]   table [x=Exponent, y=EXP-2.0-L2-REL-S12
			, col sep=comma] {WeightFunctionsDCPSE.csv};
			\addplot+[yellow,mark=*,mark options={fill=yellow}]   table [x=Exponent, y=EXP-2.5-L2-REL-S12
			, col sep=comma] {WeightFunctionsDCPSE.csv};
			\addplot+[olive,mark=star,mark options={fill=olive}]   table [x=Exponent, y=EXP-3.0-L2-REL-S12
			, col sep=comma] {WeightFunctionsDCPSE.csv};
			\addplot+[orange,mark=pentagon*,mark options={fill=orange}]   table [x=Exponent-Spline, y=SPLINE3-L2-REL-S12
			, col sep=comma] {WeightFunctionsDCPSE.csv};
			\addplot+[magenta,mark=oplus]   table [x=Exponent-Spline, y=SPLINE4-L2-REL-S12
			, col sep=comma] {WeightFunctionsDCPSE.csv};
			\end{axis}
			\end{tikzpicture} & 
			
			\begin{tikzpicture}
			\begin{axis}[height=6cm,width=6cm,ymode=log, ymin=0.0001,ymax=1,xmin=0.15,xmax=0.55, legend entries={Exp.=1.0,Exp.=1.5,Exp.=2.0,Exp.=2.5,Exp.=3.0,Spline3rd ,Spline4th},legend style={legend columns=1, cells={anchor=west},  font=\footnotesize, rounded corners=2pt,}, legend pos=outer north east,xlabel=Shape Parameter,ylabel=$L_{\infty}$ Error - $\sigma_{12}$]
			\addplot+[green,mark=square*,mark options={fill=green}]   table [x=Exponent, y=EXP-1.0-L-INF-S12
			, col sep=comma] {WeightFunctionsDCPSE.csv};
			\addplot+[red,mark=diamond*,mark options={fill=red}]   table [x=Exponent, y=EXP-1.5-L-INF-S12
			, col sep=comma] {WeightFunctionsDCPSE.csv};
			\addplot+[blue,mark=triangle*,mark options={fill=blue}]   table [x=Exponent, y=EXP-2.0-L-INF-S12
			, col sep=comma] {WeightFunctionsDCPSE.csv};
			\addplot+[yellow,mark=*,mark options={fill=yellow}]   table [x=Exponent, y=EXP-2.5-L-INF-S12
			, col sep=comma] {WeightFunctionsDCPSE.csv};
			\addplot+[olive,mark=star,mark options={fill=olive}]   table [x=Exponent, y=EXP-3.0-L-INF-S12
			, col sep=comma] {WeightFunctionsDCPSE.csv};
			\addplot+[orange,mark=pentagon*,mark options={fill=orange}]   table [x=Exponent-Spline, y=SPLINE3-L-INF-S12
			, col sep=comma] {WeightFunctionsDCPSE.csv};
			\addplot+[magenta,mark=oplus]   table [x=Exponent-Spline, y=SPLINE4-L-INF-S12
			, col sep=comma] {WeightFunctionsDCPSE.csv};
			\end{axis}
			\end{tikzpicture} \\
			
		\end{tabular}
		\caption{DC PSE Weight Function Sensitivity - 2D Cylinder. $L_2$ (Left) and $L_{\infty}$ (Right) errors as a function of the shape parameters for exponential functions of various exponents. Comparison to results obtained with $3^\text{rd}$ and $4^\text{th}$ order splines. A shape parameter of 0.3 leads to a low error for all the exponents considered.}
		\label{WeightFunctionsResultsDCPSE}
	\end{figure}
	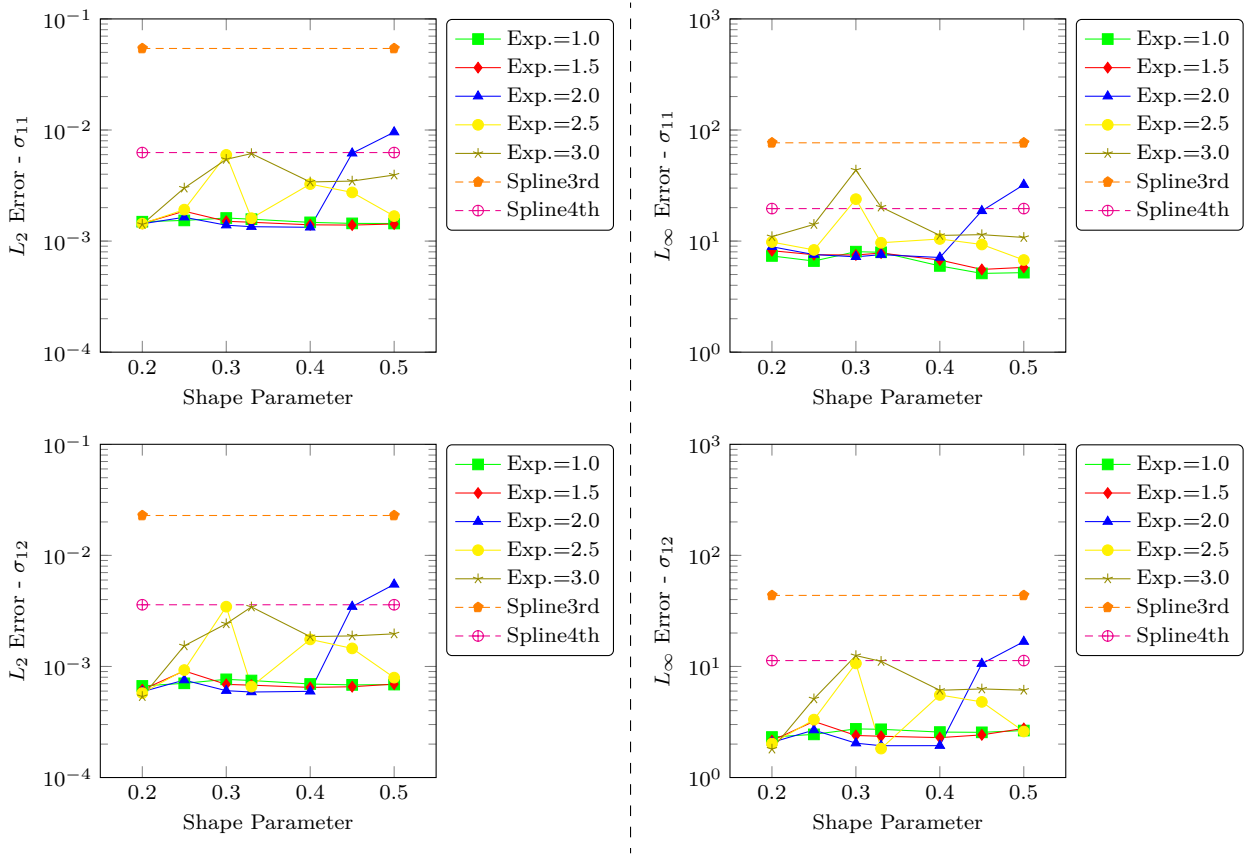
\begin{figure}[H] 
		\centering
		\begin{tabular}{c:c}
			\begin{tikzpicture}
			\begin{axis}[height=6cm,width=6cm,ymode=log, ymin=0.0001,ymax=0.1,xmin=0.15,xmax=0.55, legend entries={Exp.=1.0,Exp.=1.5,Exp.=2.0,Exp.=2.5,Exp.=3.0,Spline3rd ,Spline4th},legend style={legend columns=1, cells={anchor=west},  font=\footnotesize, rounded corners=2pt,}, legend pos=outer north east,xlabel=Shape Parameter,ylabel=$L_2$ Error - $\sigma_{11}$]
			\addplot+[green,mark=square*,mark options={fill=green}]   table [x=Exponent, y=EXP-1.0-L2-REL-S11
			, col sep=comma] {WeightFunctionsDCPSE-LShaped.csv};
			\addplot+[red,mark=diamond*,mark options={fill=red}]   table [x=Exponent, y=EXP-1.5-L2-REL-S11
			, col sep=comma] {WeightFunctionsDCPSE-LShaped.csv};
			\addplot+[blue,mark=triangle*,mark options={fill=blue}]   table [x=Exponent, y=EXP-2.0-L2-REL-S11
			, col sep=comma] {WeightFunctionsDCPSE-LShaped.csv};
			\addplot+[yellow,mark=*,mark options={fill=yellow}]   table [x=Exponent, y=EXP-2.5-L2-REL-S11
			, col sep=comma] {WeightFunctionsDCPSE-LShaped.csv};
			\addplot+[olive,mark=star,mark options={fill=olive}]   table [x=Exponent, y=EXP-3.0-L2-REL-S11
			, col sep=comma] {WeightFunctionsDCPSE-LShaped.csv};
			\addplot+[orange,mark=pentagon*,mark options={fill=orange}]   table [x=Exponent-Spline, y=SPLINE3-L2-REL-S11
			, col sep=comma] {WeightFunctionsDCPSE-LShaped.csv};
			\addplot+[magenta,mark=oplus]   table [x=Exponent-Spline, y=SPLINE4-L2-REL-S11
			, col sep=comma] {WeightFunctionsDCPSE-LShaped.csv};
			\end{axis}
			\end{tikzpicture} & 
			
			\begin{tikzpicture}
			\begin{axis}[height=6cm,width=6cm,ymode=log, ymin=1,ymax=1000,xmin=0.15,xmax=0.55, legend entries={Exp.=1.0,Exp.=1.5,Exp.=2.0,Exp.=2.5,Exp.=3.0,Spline3rd ,Spline4th},legend style={legend columns=1, cells={anchor=west},  font=\footnotesize, rounded corners=2pt,}, legend pos=outer north east,xlabel=Shape Parameter,ylabel=$L_{\infty}$ Error - $\sigma_{11}$]
			\addplot+[green,mark=square*,mark options={fill=green}]   table [x=Exponent, y=EXP-1.0-L-INF-S11
			, col sep=comma] {WeightFunctionsDCPSE-LShaped.csv};
			\addplot+[red,mark=diamond*,mark options={fill=red}]   table [x=Exponent, y=EXP-1.5-L-INF-S11
			, col sep=comma] {WeightFunctionsDCPSE-LShaped.csv};
			\addplot+[blue,mark=triangle*,mark options={fill=blue}]   table [x=Exponent, y=EXP-2.0-L-INF-S11
			, col sep=comma] {WeightFunctionsDCPSE-LShaped.csv};
			\addplot+[yellow,mark=*,mark options={fill=yellow}]   table [x=Exponent, y=EXP-2.5-L-INF-S11
			, col sep=comma] {WeightFunctionsDCPSE-LShaped.csv};
			\addplot+[olive,mark=star,mark options={fill=olive}]   table [x=Exponent, y=EXP-3.0-L-INF-S11
			, col sep=comma] {WeightFunctionsDCPSE-LShaped.csv};
			\addplot+[orange,mark=pentagon*,mark options={fill=orange}]   table [x=Exponent-Spline, y=SPLINE3-L-INF-S11
			, col sep=comma] {WeightFunctionsDCPSE-LShaped.csv};
			\addplot+[magenta,mark=oplus]   table [x=Exponent-Spline, y=SPLINE4-L-INF-S11
			, col sep=comma] {WeightFunctionsDCPSE-LShaped.csv};
			\end{axis}
			\end{tikzpicture} \\
			
			\begin{tikzpicture}
			\begin{axis}[height=6cm,width=6cm,ymode=log, ymin=0.0001,ymax=0.1,xmin=0.15,xmax=0.55, legend entries={Exp.=1.0,Exp.=1.5,Exp.=2.0,Exp.=2.5,Exp.=3.0,Spline3rd ,Spline4th},legend style={legend columns=1, cells={anchor=west},  font=\footnotesize, rounded corners=2pt,}, legend pos=outer north east,xlabel=Shape Parameter,ylabel=$L_2$ Error - $\sigma_{12}$]
			\addplot+[green,mark=square*,mark options={fill=green}]   table [x=Exponent, y=EXP-1.0-L2-REL-S12
			, col sep=comma] {WeightFunctionsDCPSE-LShaped.csv};
			\addplot+[red,mark=diamond*,mark options={fill=red}]   table [x=Exponent, y=EXP-1.5-L2-REL-S12
			, col sep=comma] {WeightFunctionsDCPSE-LShaped.csv};
			\addplot+[blue,mark=triangle*,mark options={fill=blue}]   table [x=Exponent, y=EXP-2.0-L2-REL-S12
			, col sep=comma] {WeightFunctionsDCPSE-LShaped.csv};
			\addplot+[yellow,mark=*,mark options={fill=yellow}]   table [x=Exponent, y=EXP-2.5-L2-REL-S12
			, col sep=comma] {WeightFunctionsDCPSE-LShaped.csv};
			\addplot+[olive,mark=star,mark options={fill=olive}]   table [x=Exponent, y=EXP-3.0-L2-REL-S12
			, col sep=comma] {WeightFunctionsDCPSE-LShaped.csv};
			\addplot+[orange,mark=pentagon*,mark options={fill=orange}]   table [x=Exponent-Spline, y=SPLINE3-L2-REL-S12
			, col sep=comma] {WeightFunctionsDCPSE-LShaped.csv};
			\addplot+[magenta,mark=oplus]   table [x=Exponent-Spline, y=SPLINE4-L2-REL-S12
			, col sep=comma] {WeightFunctionsDCPSE-LShaped.csv};
			\end{axis}
			\end{tikzpicture} & 
			
			\begin{tikzpicture}
			\begin{axis}[height=6cm,width=6cm,ymode=log, ymin=1,ymax=1000,xmin=0.15,xmax=0.55, legend entries={Exp.=1.0,Exp.=1.5,Exp.=2.0,Exp.=2.5,Exp.=3.0,Spline3rd ,Spline4th},legend style={legend columns=1, cells={anchor=west},  font=\footnotesize, rounded corners=2pt,}, legend pos=outer north east,xlabel=Shape Parameter,ylabel=$L_{\infty}$ Error - $\sigma_{12}$]
			\addplot+[green,mark=square*,mark options={fill=green}]   table [x=Exponent, y=EXP-1.0-L-INF-S12
			, col sep=comma] {WeightFunctionsDCPSE-LShaped.csv};
			\addplot+[red,mark=diamond*,mark options={fill=red}]   table [x=Exponent, y=EXP-1.5-L-INF-S12
			, col sep=comma] {WeightFunctionsDCPSE-LShaped.csv};
			\addplot+[blue,mark=triangle*,mark options={fill=blue}]   table [x=Exponent, y=EXP-2.0-L-INF-S12
			, col sep=comma] {WeightFunctionsDCPSE-LShaped.csv};
			\addplot+[yellow,mark=*,mark options={fill=yellow}]   table [x=Exponent, y=EXP-2.5-L-INF-S12
			, col sep=comma] {WeightFunctionsDCPSE-LShaped.csv};
			\addplot+[olive,mark=star,mark options={fill=olive}]   table [x=Exponent, y=EXP-3.0-L-INF-S12
			, col sep=comma] {WeightFunctionsDCPSE-LShaped.csv};
			\addplot+[orange,mark=pentagon*,mark options={fill=orange}]   table [x=Exponent-Spline, y=SPLINE3-L-INF-S12
			, col sep=comma] {WeightFunctionsDCPSE-LShaped.csv};
			\addplot+[magenta,mark=oplus]   table [x=Exponent-Spline, y=SPLINE4-L-INF-S12
			, col sep=comma] {WeightFunctionsDCPSE-LShaped.csv};
			\end{axis}
			\end{tikzpicture} \\
			
		\end{tabular}
		\caption{DC PSE Weight Function Sensitivity - 2D L-Shape. $L_2$ (Left) and $L_{\infty}$ (Right) errors as a function of the shape parameters for exponential functions of various exponents. Comparison to results obtained with $3^\text{rd}$ and $4^\text{th}$ order splines. An exponent of 2.0 and a shape parameter of 0.33 lead to the lowest error in terms of $L_2$ norm.}
		\label{WeightFunctionsResultsDCPSE_LShaped}
	\end{figure}	
	
	The analysis of Figure \ref{WeightFunctionsResultsDCPSE} and Figure \ref{WeightFunctionsResultsDCPSE_LShaped} shows that the exponential weight functions lead to smaller errors than the spline functions for shape parameters between 0.25 and 0.45 and for both problems. The combination of an exponent of 1.0 and a shape parameter of 0.25 leads to the smallest error for the 2D cylinder. A shape parameter of 0.30, associated to an exponent of 1.0, leads to similar results for this problem without being too close to a rapid increase in the observed error.
	
	The analysis of the results for the 2D L-shape shows that an exponent of 2.0 associated with a shape parameter of 0.33 leads relatively constantly to the lowest error.
	
	In order to be applied to most of the problems where the type of solution is a priori unknown, a single set of parameters is selected. This set of parameters is: a shape parameter of 0.30 and an exponent of 1.0. This combination leads to a more significant error reduction than the set of parameters leading to the minimum error for the 2D L-shape problem.
	
	\subsection{DC PSE Correction Function}
	
	The basis functions used to build the correction function in the DC PSE method can be selected from various function types. The most commons bases are the polynomial basis and the exponential basis. For the case of a two-dimensional problem, the polynomial basis is $\mathbf{P}=[1, x, y, x^2, x y, y^2]^T$, and the exponential basis: $\mathbf{P}=[1, e^x, e^y, e^{2x}, e^{x+y}, e^{2y}]^T$. In order for these bases to be independent from the node density, the functions have been scaled according to the support radius. For a node $\mathbf{X_{pi}}$ in the support of the collocation node $\mathbf{X_c}$, the scaling parameters are $Sx_i=\frac{x_c-x_{pi}}{r_c}$ and $Sy_i=\frac{y_c-y_{pi}}{r_c}$, where $r_c$ is the support radius of the collocation node $\mathbf{X_c}$. The polynomial basis becomes:
	\begin{equation} \label{ExpBasisScaled}
	\mathbf{P}=\Big[ 1, Sx_i, Sy_i, {Sx_i}^2, {Sx_i}{Sy_i}, {Sy_i}^2 \Big]^T \\
	\end{equation}
	
	The errors obtained for each correction function basis are presented in Figure \ref{DCPSE_CorrFunction_Comparison} below for the 2D cylinder and in Figure \ref{DCPSE_CorrFunction_Comparison_LShaped} for the 2D L-shape.
	\begin{figure}[H] 
		\centering
		\begin{tabular}{c:c}
			\begin{tikzpicture}
			\begin{axis}[height=6cm,width=7.5cm, ymin=0.0000002,ymax=0.003,ymode=log,xmin=1000,xmax=30000, xmode=log, legend entries={Polynomial Basis,Exponential Basis},legend style={ at={(0.5,-0.2)},anchor=south west,legend columns=1, cells={anchor=west},  font=\footnotesize, rounded corners=2pt,}, legend pos=north east,xlabel=Number of Nodes,ylabel=$L_2$ Error - $\sigma_{11}$]
			\addplot+[green,mark=square*,mark options={fill=green}]  table [x=NodeNum, y=L2-REL-S11-DCPSE1, col sep=comma] {DCPSECorrFunctionCylinder.csv};
			\addplot+[red,mark=diamond*,mark options={fill=red}]   table [x=NodeNum, y=L2-REL-S11-DCPSE1-EXP, col sep=comma] {DCPSECorrFunctionCylinder.csv};
			\logLogSlopeTriangle{0.85}{0.1}{0.13}{1.6}{green};
			\logLogSlopeTriangle{0.85}{0.1}{0.48}{1.0}{red};
			\end{axis}
			\end{tikzpicture} &
			
			\begin{tikzpicture}
			\begin{axis}[height=6cm,width=7.5cm, ymin=0.0001,ymax=1,ymode=log,xmin=1000,xmax=30000, xmode=log, legend entries={Polynomial Basis,Exponential Basis},legend style={ at={(0.5,-0.2)},anchor=south west,legend columns=1, cells={anchor=west},  font=\footnotesize, rounded corners=2pt,}, legend pos=north east,xlabel=Number of Nodes,ylabel=$L_{\infty}$ Error - $\sigma_{11}$]
			\addplot+[green,mark=square*,mark options={fill=green}]  table [x=NodeNum, y=L-INF-S11-DCPSE1, col sep=comma] {DCPSECorrFunctionCylinder.csv};
			\addplot+[red,mark=diamond*,mark options={fill=red}]   table [x=NodeNum, y=L-INF-S11-DCPSE1-EXP, col sep=comma] {DCPSECorrFunctionCylinder.csv};
			\logLogSlopeTriangle{0.85}{0.1}{0.14}{1.1}{green};
			\logLogSlopeTriangle{0.85}{0.1}{0.48}{0.5}{red};
			\end{axis}
			\end{tikzpicture} \\
			
			\begin{tikzpicture}
			\begin{axis}[height=6cm, width=7.5cm, ymin=0.0000002, ymax=0.003, ymode=log, xmin=1000,xmax=30000, xmode=log, legend entries={Polynomial Basis,Exponential Basis},legend style={ at={(0.5,-0.2)}, anchor=south west,legend columns=1, cells={anchor=west},  font=\footnotesize, rounded corners=2pt,}, legend pos=north east,xlabel=Number of Nodes, ylabel=$L_2$ Error - $\sigma_{12}$]
			\addplot+[green,mark=square*,mark options={fill=green}]  table [x=NodeNum, y=L2-REL-S12-DCPSE1, col sep=comma] {DCPSECorrFunctionCylinder.csv};
			\addplot+[red,mark=diamond*,mark options={fill=red}]   table [x=NodeNum, y=L2-REL-S12-DCPSE1-EXP, col sep=comma] {DCPSECorrFunctionCylinder.csv};
			\logLogSlopeTriangle{0.85}{0.1}{0.11}{1.6}{green};
			\logLogSlopeTriangle{0.85}{0.1}{0.42}{1.0}{red};
			\end{axis}
			\end{tikzpicture} &
			
			\begin{tikzpicture}
			\begin{axis}[height=6cm,width=7.5cm, ymin=0.0001,ymax=1,ymode=log,xmin=1000,xmax=30000, xmode=log, legend entries={Polynomial Basis,Exponential Basis},legend style={ at={(0.5,-0.2)},anchor=south west,legend columns=1, cells={anchor=west},  font=\footnotesize, rounded corners=2pt,}, legend pos=north east,xlabel=Number of Nodes,ylabel=$L_{\infty}$ Error - $\sigma_{12}$]
			\addplot+[green,mark=square*,mark options={fill=green}]  table [x=NodeNum, y=L-INF-S12-DCPSE1, col sep=comma] {DCPSECorrFunctionCylinder.csv};
			\addplot+[red,mark=diamond*,mark options={fill=red}]   table [x=NodeNum, y=L-INF-S12-DCPSE1-EXP, col sep=comma] {DCPSECorrFunctionCylinder.csv};
			\logLogSlopeTriangle{0.85}{0.1}{0.09}{0.9}{green};
			\logLogSlopeTriangle{0.85}{0.1}{0.34}{0.5}{red};
			\end{axis}
			\end{tikzpicture}
		\end{tabular}
		\caption{DC PSE Correction Function Basis Comparison - 2D Cylinder. $L_2$ (Left) and $L_{\infty}$ (Right) errors for polynomial and exponential bases functions as a function of the number of nodes in the model. The use of a polynomial basis leads to a lower error and a faster convergence.}
		\label{DCPSE_CorrFunction_Comparison}
	\end{figure}
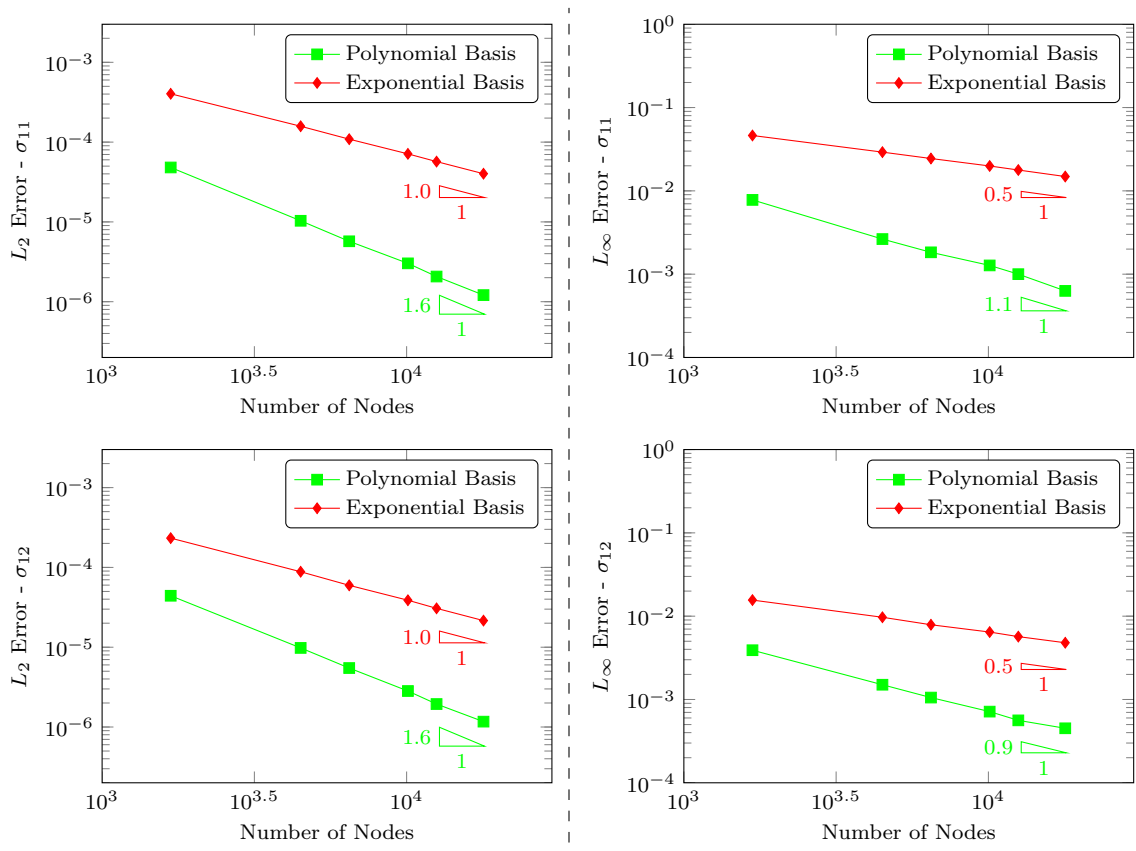
\begin{figure}[H] 
	\centering
	\begin{tabular}{c:c}
		\begin{tikzpicture}
		\begin{axis}[height=6cm,width=7.5cm, ymin=0.00001,ymax=0.01,ymode=log,xmin=3000,xmax=1000000, xmode=log, legend entries={Polynomial Basis,Exponential Basis},legend style={ at={(0.5,-0.2)},anchor=south west,legend columns=1, cells={anchor=west},  font=\footnotesize, rounded corners=2pt,}, legend pos=north east,xlabel=Number of Nodes,ylabel=$L_2$ Error - $\sigma_{11}$]
		\addplot+[green,mark=square*,mark options={fill=green}]  table [x=NodeNum, y=L2-REL-S11-DCPSE1, col sep=comma] {DCPSECorrFunctionLShaped.csv};
		\addplot+[red,mark=diamond*,mark options={fill=red}]   table [x=NodeNum, y=L2-REL-S11-DCPSE1-EXP, col sep=comma] {DCPSECorrFunctionLShaped.csv};
		\logLogSlopeTriangle{0.85}{0.1}{0.6}{0.6}{green};
		\logLogSlopeTriangle{0.85}{0.1}{0.3}{0.6}{red};
		\end{axis}
		\end{tikzpicture} &
		
		\begin{tikzpicture}
		\begin{axis}[height=6cm,width=7.5cm, ymin=1,ymax=100,ymode=log,xmin=3000,xmax=1000000, xmode=log, legend entries={Polynomial Basis,Exponential Basis},legend style={ at={(0.5,-0.2)},anchor=south west,legend columns=1, cells={anchor=west},  font=\footnotesize, rounded corners=2pt,}, legend pos=north east,xlabel=Number of Nodes,ylabel=$L_{\infty}$ Error - $\sigma_{11}$]
		\addplot+[green,mark=square*,mark options={fill=green}]  table [x=NodeNum, y=L-INF-S11-DCPSE1, col sep=comma] {DCPSECorrFunctionLShaped.csv};
		\addplot+[red,mark=diamond*,mark options={fill=red}]   table [x=NodeNum, y=L-INF-S11-DCPSE1-EXP, col sep=comma] {DCPSECorrFunctionLShaped.csv};
		\logLogSlopeTriangleUp{0.85}{0.1}{0.67}{0.2}{red};
		\logLogSlopeTriangleUp{0.85}{0.1}{0.45}{0.2}{green};
		\end{axis}
		\end{tikzpicture} \\
		
		\begin{tikzpicture}
		\begin{axis}[height=6cm, width=7.5cm, ymin=0.00001, ymax=0.01, ymode=log, xmin=3000,xmax=1000000, xmode=log, legend entries={Polynomial Basis,Exponential Basis},legend style={ at={(0.5,-0.2)}, anchor=south west,legend columns=1, cells={anchor=west},  font=\footnotesize, rounded corners=2pt,}, legend pos=north east,xlabel=Number of Nodes, ylabel=$L_2$ Error - $\sigma_{12}$]
		\addplot+[green,mark=square*,mark options={fill=green}]  table [x=NodeNum, y=L2-REL-S12-DCPSE1, col sep=comma] {DCPSECorrFunctionLShaped.csv};
		\addplot+[red,mark=diamond*,mark options={fill=red}]   table [x=NodeNum, y=L2-REL-S12-DCPSE1-EXP, col sep=comma] {DCPSECorrFunctionLShaped.csv};
		\logLogSlopeTriangle{0.85}{0.1}{0.5}{0.6}{green};
		\logLogSlopeTriangle{0.85}{0.1}{0.15}{0.6}{red};
		\end{axis}
		\end{tikzpicture} &
		
		\begin{tikzpicture}
		\begin{axis}[height=6cm,width=7.5cm, ymin=1,ymax=100,ymode=log,xmin=3000,xmax=1000000, xmode=log, legend entries={Polynomial Basis,Exponential Basis},legend style={ at={(0.5,-0.2)},anchor=south west,legend columns=1, cells={anchor=west},  font=\footnotesize, rounded corners=2pt,}, legend pos=north east,xlabel=Number of Nodes,ylabel=$L_{\infty}$ Error - $\sigma_{12}$]
		\addplot+[green,mark=square*,mark options={fill=green}]  table [x=NodeNum, y=L-INF-S12-DCPSE1, col sep=comma] {DCPSECorrFunctionLShaped.csv};
		\addplot+[red,mark=diamond*,mark options={fill=red}]   table [x=NodeNum, y=L-INF-S12-DCPSE1-EXP, col sep=comma] {DCPSECorrFunctionLShaped.csv};
		\logLogSlopeTriangleUp{0.85}{0.1}{0.4}{0.2}{green};
		\logLogSlopeTriangleUp{0.85}{0.1}{0.15}{0.2}{red};
		\end{axis}
		\end{tikzpicture}
	\end{tabular}
	\caption{DC PSE Correction Function Basis Comparison - 2D L-Shape. $L_2$ (Left) and $L_{\infty}$ (Right) errors for polynomial and exponential bases functions as a function of the number of nodes in the model. Both function basis lead to similar convergence rates. The use of a exponential basis leads to a slightly lower error for the $L_2$ error norm. For the $L_{\infty}$ error norm, the polynomial basis lead to the lowest error for the $\sigma_{11}$ stress component while the exponential basis lead to the lowest error for the $\sigma_{12}$ stress component.}
	\label{DCPSE_CorrFunction_Comparison_LShaped}
\end{figure}
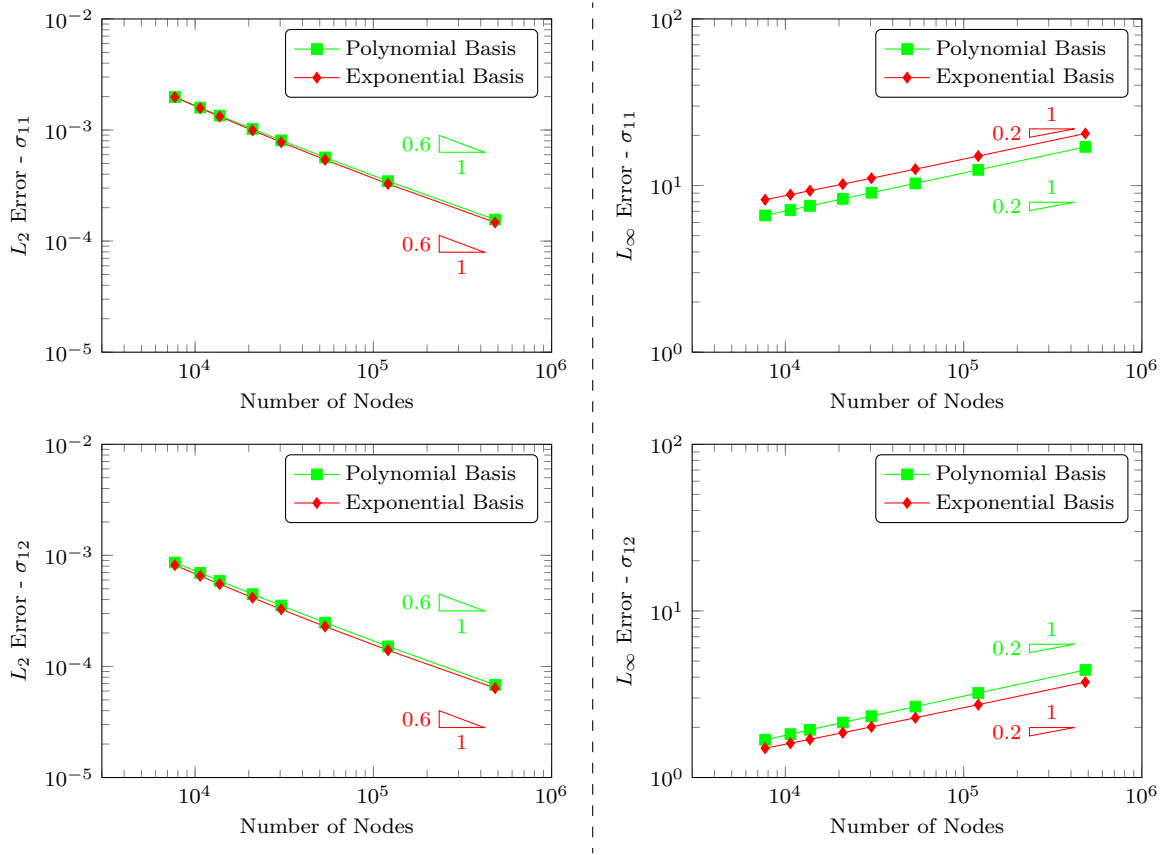
	
	It can be observed from Figure \ref{DCPSE_CorrFunction_Comparison} that the polynomial function basis constantly leads to a much lower error than the exponential basis. Depending on the number of nodes, the error increases by a factor between 5 and 30 when the exponential basis is used.
	
	In Figure \ref{DCPSE_CorrFunction_Comparison_LShaped} we can see that both correction function bases lead to similar results. In terms of $L_2$ error norm the exponential basis leads to a reduction of around 5\% compared to the results where a polynomial basis is used.
	
	Based on the results presented in this section, the polynomial basis is preferred as it leads to a much higher error reduction than the exponential basis for the considered problems. This basis function is expected to give a reasonably low error for most of the problems in the domain of linear elasticity.
	
	\subsection{Number of Support Nodes}
	
	\paragraph{Support Radius Selection} \
	
	In this work, we have selected the support radii of the collocation nodes based on the number of nodes within the support they define. The number of nodes in the support of a collocation node shall be of at least the number of approximated derivatives. In practice, in order to account for the node distribution, a larger number of nodes is used.
	
	\paragraph{GFD} \
	
	In this Section, we study the impact of the number of nodes in the support of a collocation node for the GFD method. Various combinations of inner node and boundary node support sizes have been considered. Results are presented in Figure \ref{BoundaryNodeResultsGFD}, Figure \ref{BoundaryNodeResultsGFD_LShaped} and Figure \ref{BoundaryNode_3DResultsGFD} below, respectively, for the 2D cylinder, the 2D L-shape and the sphere under internal pressure. The $L_2$ and the $L_{\infty}$ norms are presented as a function of the number of support nodes for collocation nodes located on the boundary. The results are presented for three sizes of inner node supports for the 2D problems and four sizes for the 3D problem.
	\begin{figure}[H] 
		\centering
		\begin{tabular}{c:c}
			\begin{tikzpicture}
			\begin{axis}[height=6cm,width=7.5cm, ymode=log, ymin=0.000001,ymax=0.01,xmin=12,xmax=25, legend entries={Inn. Sup=11,Inn. Sup=13,Inn. Sup=15},legend style={ at={(0.5,-0.2)},anchor=south west,legend columns=1, cells={anchor=west},  font=\footnotesize, rounded corners=2pt,}, legend pos=north east,xlabel=Support Size on Boundary,ylabel=$L_2$ Error - $\sigma_{11}$]
			\addplot+[green,mark=square*,mark options={fill=green}]  table [x=X-ALL, y=L2-REL-S11-11.0, col sep=comma] {SupportSizeGFD.csv};
			\addplot+[red,mark=diamond*,mark options={fill=red}]   table [x=X-ALL, y=L2-REL-S11-13.0, col sep=comma] {SupportSizeGFD.csv};
			\addplot+[blue,mark=triangle*,mark options={fill=blue}]   table [x=X-ALL, y=L2-REL-S11-15.0, col sep=comma] {SupportSizeGFD.csv};
			\end{axis}
			\end{tikzpicture} &
			
			\begin{tikzpicture}
			\begin{axis}[height=6cm,width=7.5cm,ymode=log, ymin=0.001,ymax=10,xmin=12,xmax=25, legend entries={Inn. Sup=11,Inn. Sup=13,Inn. Sup=15},legend style={ at={(0.5,-0.2)},anchor=south west,legend columns=1, cells={anchor=west},  font=\footnotesize, rounded corners=2pt,}, legend pos=north east,xlabel=Support Size on Boundary,ylabel=$L_{\infty}$ Error - $\sigma_{11}$]
			\addplot+[green,mark=square*,mark options={fill=green}]  table [x=X-ALL, y=L-INF-S11-11.0, col sep=comma] {SupportSizeGFD.csv};
			\addplot+[red,mark=diamond*,mark options={fill=red}]   table [x=X-ALL, y=L-INF-S11-13.0, col sep=comma] {SupportSizeGFD.csv};
			\addplot+[blue,mark=triangle*,mark options={fill=blue}]   table [x=X-ALL, y=L-INF-S11-15.0, col sep=comma] {SupportSizeGFD.csv};
			\end{axis}
			\end{tikzpicture} \\
			
			\begin{tikzpicture}
			\begin{axis}[height=6cm,width=7.5cm, ymode=log, ymin=0.000001,ymax=0.01,xmin=12,xmax=25, legend entries={Inn. Sup=11,Inn. Sup=13,Inn. Sup=15},legend style={ at={(0.5,-0.2)},anchor=south west,legend columns=1, cells={anchor=west},  font=\footnotesize, rounded corners=2pt,}, legend pos=north east,xlabel=Support Size on Boundary,ylabel=$L_2$ Error - $\sigma_{12}$]
			\addplot+[green,mark=square*,mark options={fill=green}]  table [x=X-ALL, y=L2-REL-S12-11.0, col sep=comma] {SupportSizeGFD.csv};
			\addplot+[red,mark=diamond*,mark options={fill=red}]   table [x=X-ALL, y=L2-REL-S12-13.0, col sep=comma] {SupportSizeGFD.csv};
			\addplot+[blue,mark=triangle*,mark options={fill=blue}]   table [x=X-ALL, y=L2-REL-S12-15.0, col sep=comma] {SupportSizeGFD.csv};
			\end{axis}
			\end{tikzpicture} &
			
			\begin{tikzpicture}
			\begin{axis}[height=6cm,width=7.5cm,ymode=log, ymin=0.001,ymax=10,xmin=12,xmax=25, legend entries={Inn. Sup=11,Inn. Sup=13,Inn. Sup=15},legend style={ at={(0.5,-0.2)},anchor=south west,legend columns=1, cells={anchor=west},  font=\footnotesize, rounded corners=2pt,}, legend pos=north east,xlabel=Support Size on Boundary,ylabel=$L_{\infty}$ Error - $\sigma_{12}$]
			\addplot+[green,mark=square*,mark options={fill=green}]  table [x=X-ALL, y=L-INF-S12-11.0, col sep=comma] {SupportSizeGFD.csv};
			\addplot+[red,mark=diamond*,mark options={fill=red}]   table [x=X-ALL, y=L-INF-S12-13.0, col sep=comma] {SupportSizeGFD.csv};
			\addplot+[blue,mark=triangle*,mark options={fill=blue}]   table [x=X-ALL, y=L-INF-S12-15.0, col sep=comma] {SupportSizeGFD.csv};
			\end{axis}
			\end{tikzpicture} \\	
		\end{tabular}
		\caption{GFD Support Node Number Sensitivity - 2D Cylinder. $L_2$ (Left) and $L_{\infty}$ (Right) errors for various combinations of inner nodes and boundary nodes support sizes. Inner collocation nodes with 11 support nodes lead to the lowest observed error. The error stops decreasing for boundary nodes supports larger than 18 nodes.}
		\label{BoundaryNodeResultsGFD}
	\end{figure}
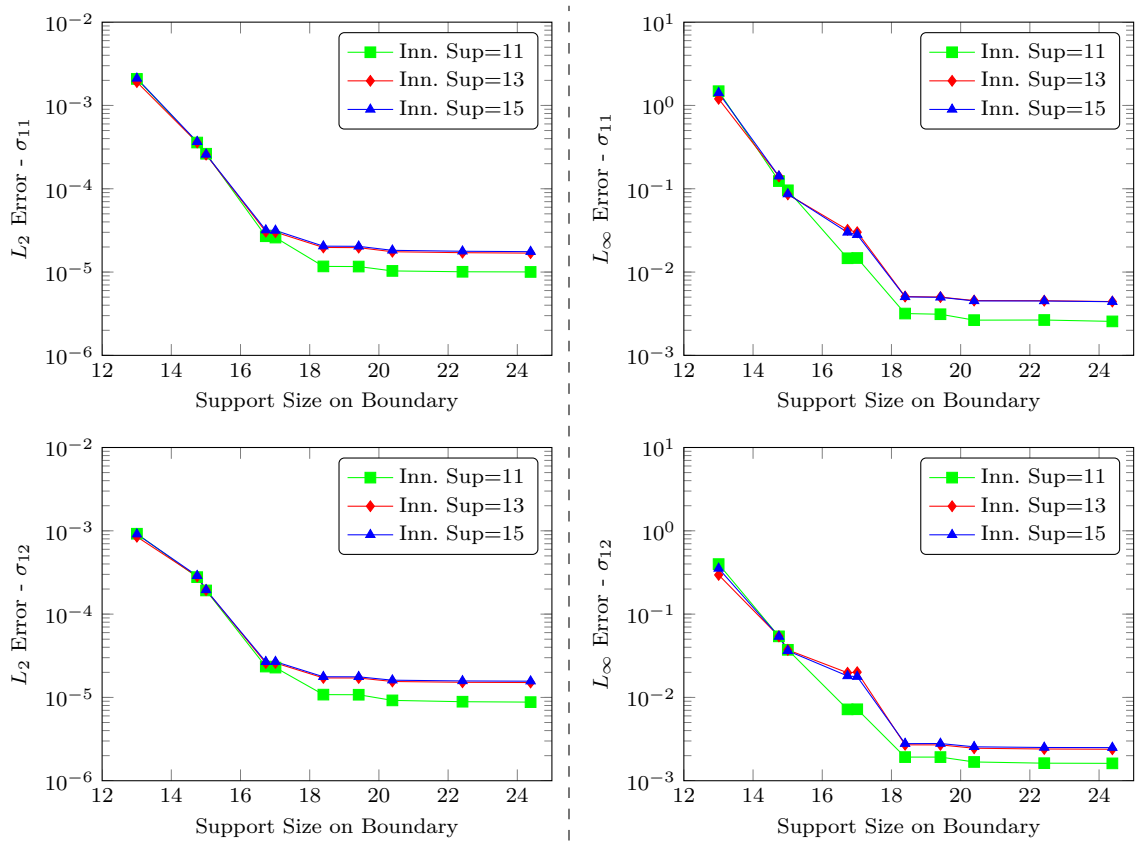
	\begin{figure}[H] 
		\centering
		\begin{tabular}{c:c}
			\begin{tikzpicture}
			\begin{axis}[height=6cm,width=7.5cm, ymode=log, ymin=0.0001,ymax=0.01,xmin=12,xmax=27, legend entries={Inn. Sup=11,Inn. Sup=13,Inn. Sup=15},legend style={ at={(0.5,-0.2)},anchor=south west,legend columns=1, cells={anchor=west},  font=\footnotesize, rounded corners=2pt,}, legend pos=north east,xlabel=Support Size on Boundary,ylabel=$L_2$ Error - $\sigma_{11}$]
			\addplot+[green,mark=square*,mark options={fill=green}]  table [x=X-ALL, y=L2-REL-S11-11.0, col sep=comma] {SupportSizeGFD-LShaped.csv};
			\addplot+[red,mark=diamond*,mark options={fill=red}]   table [x=X-ALL, y=L2-REL-S11-13.0, col sep=comma] {SupportSizeGFD-LShaped.csv};
			\addplot+[blue,mark=triangle*,mark options={fill=blue}]   table [x=X-ALL, y=L2-REL-S11-15.0, col sep=comma] {SupportSizeGFD-LShaped.csv};
			\end{axis}
			\end{tikzpicture} &
			
			\begin{tikzpicture}
			\begin{axis}[height=6cm,width=7.5cm,ymode=log, ymin=1,ymax=100,xmin=12,xmax=27, legend entries={Inn. Sup=11,Inn. Sup=13,Inn. Sup=15},legend style={ at={(0.5,-0.2)},anchor=south west,legend columns=1, cells={anchor=west},  font=\footnotesize, rounded corners=2pt,}, legend pos=north east,xlabel=Support Size on Boundary,ylabel=$L_{\infty}$ Error - $\sigma_{11}$]
			\addplot+[green,mark=square*,mark options={fill=green}]  table [x=X-ALL, y=L-INF-S11-11.0, col sep=comma] {SupportSizeGFD-LShaped.csv};
			\addplot+[red,mark=diamond*,mark options={fill=red}]   table [x=X-ALL, y=L-INF-S11-13.0, col sep=comma] {SupportSizeGFD-LShaped.csv};
			\addplot+[blue,mark=triangle*,mark options={fill=blue}]   table [x=X-ALL, y=L-INF-S11-15.0, col sep=comma] {SupportSizeGFD-LShaped.csv};
			\end{axis}
			\end{tikzpicture} \\
			
			\begin{tikzpicture}
			\begin{axis}[height=6cm,width=7.5cm, ymode=log, ymin=0.0001,ymax=0.01,xmin=12,xmax=27, legend entries={Inn. Sup=11,Inn. Sup=13,Inn. Sup=15},legend style={ at={(0.5,-0.2)},anchor=south west,legend columns=1, cells={anchor=west},  font=\footnotesize, rounded corners=2pt,}, legend pos=north east,xlabel=Support Size on Boundary,ylabel=$L_2$ Error - $\sigma_{12}$]
			\addplot+[green,mark=square*,mark options={fill=green}]  table [x=X-ALL, y=L2-REL-S12-11.0, col sep=comma] {SupportSizeGFD-LShaped.csv};
			\addplot+[red,mark=diamond*,mark options={fill=red}]   table [x=X-ALL, y=L2-REL-S12-13.0, col sep=comma] {SupportSizeGFD-LShaped.csv};
			\addplot+[blue,mark=triangle*,mark options={fill=blue}]   table [x=X-ALL, y=L2-REL-S12-15.0, col sep=comma] {SupportSizeGFD-LShaped.csv};
			\end{axis}
			\end{tikzpicture} &
			
			\begin{tikzpicture}
			\begin{axis}[height=6cm,width=7.5cm,ymode=log, ymin=1,ymax=100,xmin=12,xmax=27, legend entries={Inn. Sup=11,Inn. Sup=13,Inn. Sup=15},legend style={ at={(0.5,-0.2)},anchor=south west,legend columns=1, cells={anchor=west},  font=\footnotesize, rounded corners=2pt,}, legend pos=north east,xlabel=Support Size on Boundary,ylabel=$L_{\infty}$ Error - $\sigma_{12}$]
			\addplot+[green,mark=square*,mark options={fill=green}]  table [x=X-ALL, y=L-INF-S12-11.0, col sep=comma] {SupportSizeGFD-LShaped.csv};
			\addplot+[red,mark=diamond*,mark options={fill=red}]   table [x=X-ALL, y=L-INF-S12-13.0, col sep=comma] {SupportSizeGFD-LShaped.csv};
			\addplot+[blue,mark=triangle*,mark options={fill=blue}]   table [x=X-ALL, y=L-INF-S12-15.0, col sep=comma] {SupportSizeGFD-LShaped.csv};
			\end{axis}
			\end{tikzpicture} \\	
		\end{tabular}
		\caption{GFD Support Node Number Sensitivity - 2D L-Shape. $L_2$ (Left) and $L_{\infty}$ (Right) errors for various combinations of inner nodes and boundary nodes support sizes. All the combinations of inner and boundary nodes support size lead to similar errors. This is due to the system being loaded via Dirichlet boundary conditions.}
		\label{BoundaryNodeResultsGFD_LShaped}
	\end{figure}
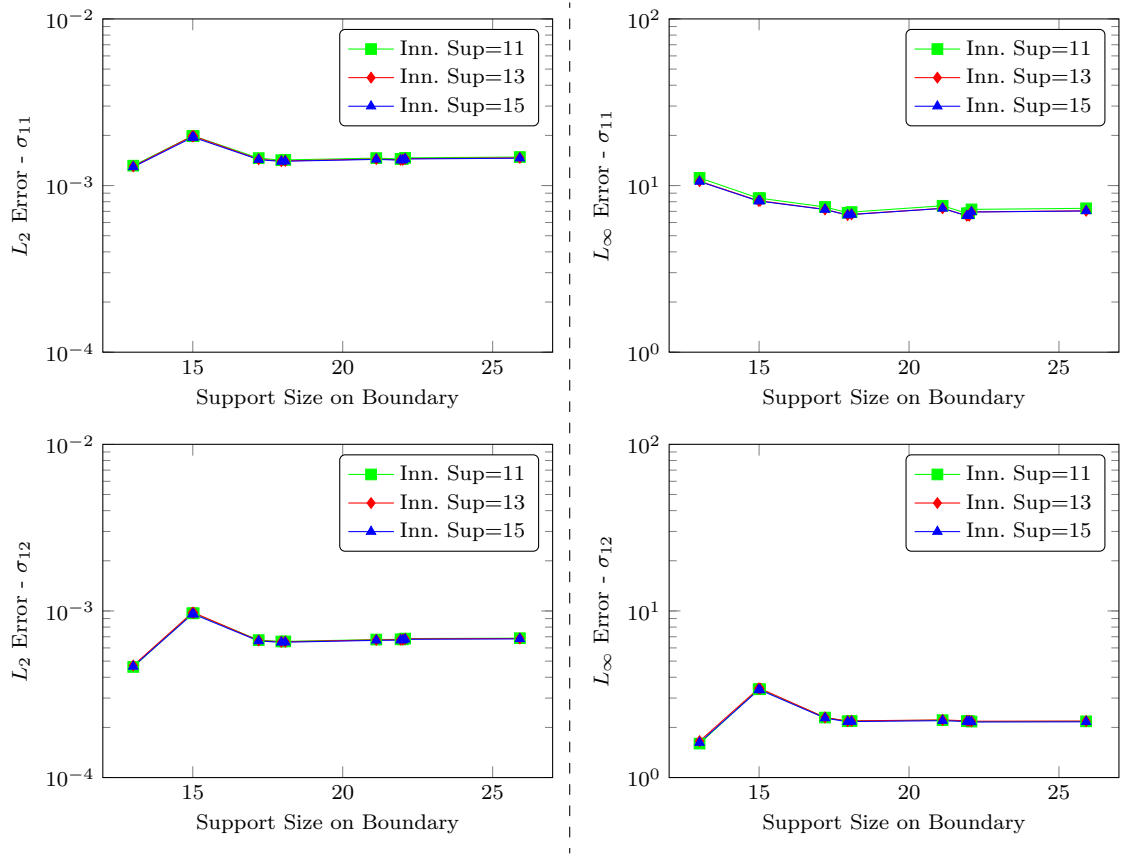

	The results in terms of $L_2$ and $L_{\infty}$ errors present a similar trend for the $\sigma_{11}$ and $\sigma_{12}$ stress components. We can see that increasing the number of nodes in the inner nodes supports does not necessarily lead to an error reduction. The loss in terms of resolution is not compensated by the gain in solution smoothness. Increasing the number of support nodes for the boundary nodes steadily (and rapidly) reduces the error for the 2D cylinder problem. An error reduction of a factor of approximately one hundred is observed when the number of support nodes for boundary collocation nodes is increased from 13 to 18.
	
	The number of support nodes on the boundary does not affect much the observed error for the 2D L-shape problem as the model is loaded via Dirichlet boundary conditions. The size of the inner nodes support has also little impact on the error for this problem.
	
	A combination of 11 support nodes for interior nodes and 19 support nodes for boundary nodes is selected as it leads to a low error for both problems while maintaining the fill of the system matrix reasonably low.
	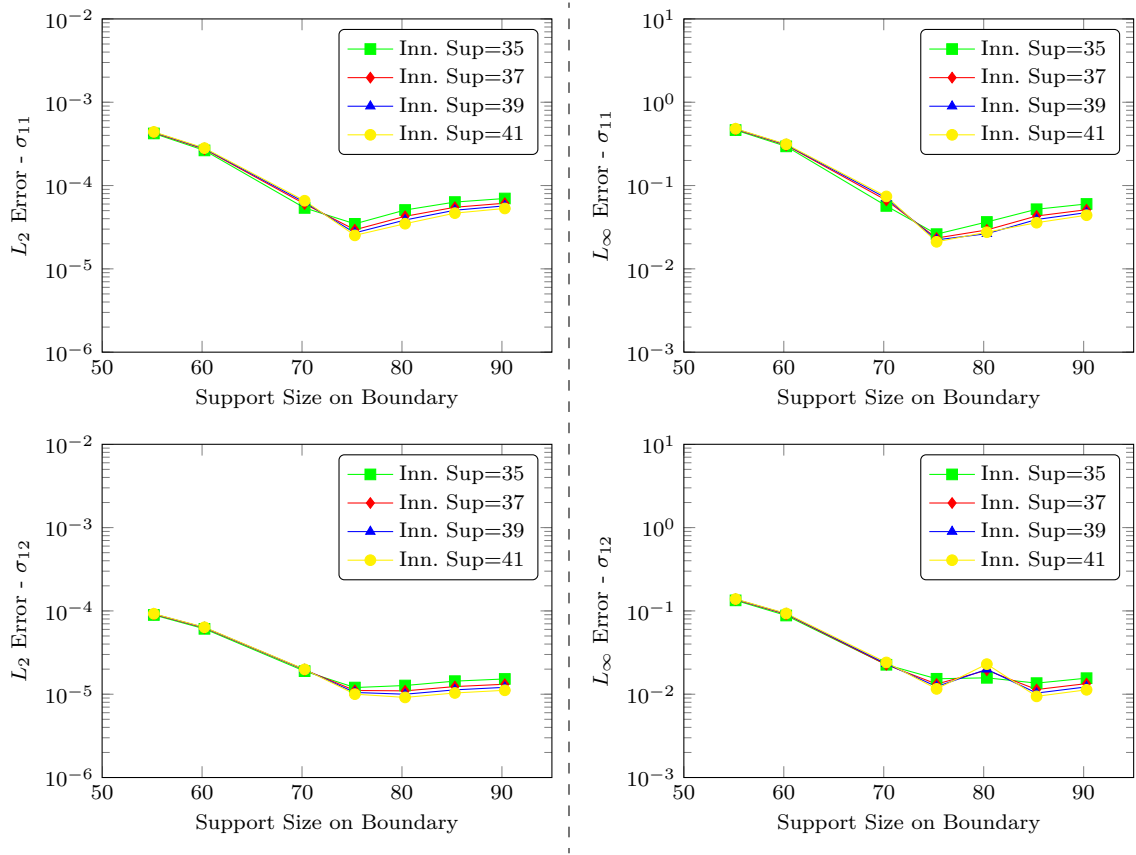
\begin{figure}[H] 
		\centering
		\begin{tabular}{c:c}
			\begin{tikzpicture}
			\begin{axis}[height=6cm,width=7.5cm, ymode=log, ymin=0.000001,ymax=0.01,xmin=50,xmax=95, legend entries={Inn. Sup=35,Inn. Sup=37,Inn. Sup=39,Inn. Sup=41},legend style={ at={(0.5,-0.2)},anchor=south west,legend columns=1, cells={anchor=west},  font=\footnotesize, rounded corners=2pt,}, legend pos=north east,xlabel=Support Size on Boundary,ylabel=$L_2$ Error - $\sigma_{11}$]
			\addplot+[green,mark=square*,mark options={fill=green}]  table [x=X-3D, y=L2-REL-S11-35, col sep=comma] {SupportSizeGFD3D.csv};
			\addplot+[red,mark=diamond*,mark options={fill=red}]   table [x=X-3D, y=L2-REL-S11-37, col sep=comma] {SupportSizeGFD3D.csv};
			\addplot+[blue,mark=triangle*,mark options={fill=blue}]   table [x=X-3D, y=L2-REL-S11-39, col sep=comma] {SupportSizeGFD3D.csv};
			\addplot+[yellow,mark=*,mark options={fill=yellow}]   table [x=X-3D, y=L2-REL-S11-41, col sep=comma] {SupportSizeGFD3D.csv};
			\end{axis}
			\end{tikzpicture} &
			
			\begin{tikzpicture}
			\begin{axis}[height=6cm,width=7.5cm,ymode=log, ymin=0.001,ymax=10,xmin=50,xmax=95, legend entries={Inn. Sup=35,Inn. Sup=37,Inn. Sup=39,Inn. Sup=41},legend style={ at={(0.5,-0.2)},anchor=south west,legend columns=1, cells={anchor=west},  font=\footnotesize, rounded corners=2pt,}, legend pos=north east,xlabel=Support Size on Boundary,ylabel=$L_{\infty}$ Error - $\sigma_{11}$]
			\addplot+[green,mark=square*,mark options={fill=green}]  table [x=X-3D, y=L-INF-S11-35, col sep=comma] {SupportSizeGFD3D.csv};
			\addplot+[red,mark=diamond*,mark options={fill=red}]   table [x=X-3D, y=L-INF-S11-37, col sep=comma] {SupportSizeGFD3D.csv};
			\addplot+[blue,mark=triangle*,mark options={fill=blue}]   table [x=X-3D, y=L-INF-S11-39, col sep=comma] {SupportSizeGFD3D.csv};
			\addplot+[yellow,mark=*,mark options={fill=yellow}]   table [x=X-3D, y=L-INF-S11-41, col sep=comma] {SupportSizeGFD3D.csv};
			\end{axis}
			\end{tikzpicture} \\
			
			\begin{tikzpicture}
			\begin{axis}[height=6cm,width=7.5cm, ymode=log, ymin=0.000001,ymax=0.01,xmin=50,xmax=95, legend entries={Inn. Sup=35,Inn. Sup=37,Inn. Sup=39,Inn. Sup=41},legend style={ at={(0.5,-0.2)},anchor=south west,legend columns=1, cells={anchor=west},  font=\footnotesize, rounded corners=2pt,}, legend pos=north east,xlabel=Support Size on Boundary,ylabel=$L_2$ Error - $\sigma_{12}$]
			\addplot+[green,mark=square*,mark options={fill=green}]  table [x=X-3D, y=L2-REL-S12-35, col sep=comma] {SupportSizeGFD3D.csv};
			\addplot+[red,mark=diamond*,mark options={fill=red}]   table [x=X-3D, y=L2-REL-S12-37, col sep=comma] {SupportSizeGFD3D.csv};
			\addplot+[blue,mark=triangle*,mark options={fill=blue}]   table [x=X-3D, y=L2-REL-S12-39, col sep=comma] {SupportSizeGFD3D.csv};
			\addplot+[yellow,mark=*,mark options={fill=yellow}]   table [x=X-3D, y=L2-REL-S12-41, col sep=comma] {SupportSizeGFD3D.csv};
			\end{axis}
			\end{tikzpicture} &
			
			\begin{tikzpicture}
			\begin{axis}[height=6cm,width=7.5cm,ymode=log, ymin=0.001,ymax=10,xmin=50,xmax=95, legend entries={Inn. Sup=35,Inn. Sup=37,Inn. Sup=39,Inn. Sup=41},legend style={ at={(0.5,-0.2)},anchor=south west,legend columns=1, cells={anchor=west},  font=\footnotesize, rounded corners=2pt,}, legend pos=north east,xlabel=Support Size on Boundary,ylabel=$L_{\infty}$ Error - $\sigma_{12}$]
			\addplot+[green,mark=square*,mark options={fill=green}]  table [x=X-3D, y=L-INF-S12-35, col sep=comma] {SupportSizeGFD3D.csv};
			\addplot+[red,mark=diamond*,mark options={fill=red}]   table [x=X-3D, y=L-INF-S12-37, col sep=comma] {SupportSizeGFD3D.csv};
			\addplot+[blue,mark=triangle*,mark options={fill=blue}]   table [x=X-3D, y=L-INF-S12-39, col sep=comma] {SupportSizeGFD3D.csv};
			\addplot+[yellow,mark=*,mark options={fill=yellow}]   table [x=X-3D, y=L-INF-S12-41, col sep=comma] {SupportSizeGFD3D.csv};
			\end{axis}
			\end{tikzpicture} \\	
		\end{tabular}
		\caption{GFD Support Node Number Sensitivity - 3D Sphere. $L_2$ (Left) and $L_{\infty}$ (Right) errors for various combinations of inner nodes and boundary nodes support sizes. Boundary collocation nodes with 75 support nodes lead to the lowest error. The size of the support of inner collocation nodes has little impact on the error.}
		\label{BoundaryNode_3DResultsGFD}
	\end{figure}
		
	It can be observed from Figure \ref{BoundaryNode_3DResultsGFD} that, for the 3D sphere, a minimum error is obtained for 75 support nodes for boundary collocation nodes. Increasing the size of the support from 55 to 75 for boundary nodes reduces in average by a factor 10 the observed error both in terms of $L_2$ and $L_{\infty}$ norms. The number of support nodes for the collocation nodes located in the domain has a smaller impact on the error. 37 support nodes appears to be a reasonable choice as it leads to a low error while keeping the fill of the matrix reasonable.
	
	The sparsity of the problem matrix is reduced when the number of support nodes increases. However, the impact of an increase in the number of support nodes on the boundary is limited as it only affects a fraction of the nodes of the domain.
	
	\paragraph{DC PSE} \
	
	As for the GFD method, the impact of the support size on the observed error is presented in this section for the DC PSE method. The results are presented in Figure \ref{BoundaryNodeResultsDCPSE}, Figure \ref{BoundaryNodeResultsDCPSE_LShaped} and Figure \ref{BoundaryNode_3DResultsDCPSE} below, respectively, for the 2D cylinder, the 2D L-shape and the sphere under internal pressure for various combinations of inner node and boundary node support sizes.
	\begin{figure}[H] 
		\centering
		\begin{tabular}{c:c}
			\begin{tikzpicture}
			\begin{axis}[height=6cm,width=7.5cm, ymode=log, ymin=0.000001,ymax=0.0001,xmin=12,xmax=25, legend entries={Inn. Sup=9,Inn. Sup=11,Inn. Sup=13,Inn. Sup=15},legend style={ at={(0.5,-0.2)},anchor=south west,legend columns=1, cells={anchor=west},  font=\footnotesize, rounded corners=2pt,}, legend pos=north east,xlabel=Support Size on Boundary,ylabel=$L_2$ Error - $\sigma_{11}$]
			\addplot+[green,mark=square*,mark options={fill=green}]  table [x=X-9.0, y=L2-REL-S11-9.0, col sep=comma] {SupportSizeDCPSE.csv};
			\addplot+[red,mark=diamond*,mark options={fill=red}]   table [x=X-11.0, y=L2-REL-S11-11.0, col sep=comma] {SupportSizeDCPSE.csv};
			\addplot+[blue,mark=triangle*,mark options={fill=blue}]   table [x=X-13.0, y=L2-REL-S11-13.0, col sep=comma] {SupportSizeDCPSE.csv};
			\addplot+[yellow,mark=*,mark options={fill=yellow}]   table [x=X-15.0, y=L2-REL-S11-15.0, col sep=comma] {SupportSizeDCPSE.csv};
			\end{axis}
			\end{tikzpicture} &
			
			\begin{tikzpicture}
			\begin{axis}[height=6cm,width=7.5cm,ymode=log, ymin=0.0001,ymax=0.06,xmin=12,xmax=25, legend entries={Inn. Sup=9,Inn. Sup=11,Inn. Sup=13,Inn. Sup=15},legend style={ at={(0.5,-0.2)},anchor=south west,legend columns=1, cells={anchor=west},  font=\footnotesize, rounded corners=2pt,}, legend pos=north east,xlabel=Support Size on Boundary,ylabel=$L_{\infty}$ Error - $\sigma_{11}$]
			\addplot+[green,mark=square*,mark options={fill=green}]  table [x=X-9.0, y=L-INF-S11-9.0, col sep=comma] {SupportSizeDCPSE.csv};
			\addplot+[red,mark=diamond*,mark options={fill=red}]   table [x=X-11.0, y=L-INF-S11-11.0, col sep=comma] {SupportSizeDCPSE.csv};
			\addplot+[blue,mark=triangle*,mark options={fill=blue}]   table [x=X-13.0, y=L-INF-S11-13.0, col sep=comma] {SupportSizeDCPSE.csv};
			\addplot+[yellow,mark=*,mark options={fill=yellow}]   table [x=X-15.0, y=L-INF-S11-15.0, col sep=comma] {SupportSizeDCPSE.csv};
			\end{axis}
			\end{tikzpicture} \\
			
			\begin{tikzpicture}
			\begin{axis}[height=6cm,width=7.5cm, ymode=log, ymin=0.000001,ymax=0.0001,xmin=12,xmax=25, legend entries={Inn. Sup=9,Inn. Sup=11,Inn. Sup=13,Inn. Sup=15},legend style={ at={(0.5,-0.2)},anchor=south west,legend columns=1, cells={anchor=west},  font=\footnotesize, rounded corners=2pt,}, legend pos=north east,xlabel=Support Size on Boundary,ylabel=$L_2$ Error - $\sigma_{12}$]
			\addplot+[green,mark=square*,mark options={fill=green}]  table [x=X-9.0, y=L2-REL-S12-9.0, col sep=comma] {SupportSizeDCPSE.csv};
			\addplot+[red,mark=diamond*,mark options={fill=red}]   table [x=X-11.0, y=L2-REL-S12-11.0, col sep=comma] {SupportSizeDCPSE.csv};
			\addplot+[blue,mark=triangle*,mark options={fill=blue}]   table [x=X-13.0, y=L2-REL-S12-13.0, col sep=comma] {SupportSizeDCPSE.csv};
			\addplot+[yellow,mark=*,mark options={fill=yellow}]   table [x=X-15.0, y=L2-REL-S12-15.0, col sep=comma] {SupportSizeDCPSE.csv};
			\end{axis}
			\end{tikzpicture} &
			
			\begin{tikzpicture}
			\begin{axis}[height=6cm,width=7.5cm,ymode=log, ymin=0.0001,ymax=0.06,xmin=12,xmax=25, legend entries={Inn. Sup=9,Inn. Sup=11,Inn. Sup=13,Inn. Sup=15},legend style={ at={(0.5,-0.2)},anchor=south west,legend columns=1, cells={anchor=west},  font=\footnotesize, rounded corners=2pt,}, legend pos=north east,xlabel=Support Size on Boundary,ylabel=$L_{\infty}$ Error - $\sigma_{12}$]
			\addplot+[green,mark=square*,mark options={fill=green}]  table [x=X-9.0, y=L-INF-S12-9.0, col sep=comma] {SupportSizeDCPSE.csv};
			\addplot+[red,mark=diamond*,mark options={fill=red}]   table [x=X-11.0, y=L-INF-S12-11.0, col sep=comma] {SupportSizeDCPSE.csv};
			\addplot+[blue,mark=triangle*,mark options={fill=blue}]   table [x=X-13.0, y=L-INF-S12-13.0, col sep=comma] {SupportSizeDCPSE.csv};
			\addplot+[yellow,mark=*,mark options={fill=yellow}]   table [x=X-15.0, y=L-INF-S12-15.0, col sep=comma] {SupportSizeDCPSE.csv};
			\end{axis}
			\end{tikzpicture} \\	
		\end{tabular}
		\caption{DC PSE Support Node Number Sensitivity - 2D Cylinder. $L_2$ (Left) and $L_{\infty}$ (Right) errors for various combinations of inner nodes and boundary nodes support sizes. Inner collocation nodes with 13 support nodes lead relatively constantly to a low error. The error starts to increase for the 13 inner nodes case when the number of boundary nodes is larger than 19.}
		\label{BoundaryNodeResultsDCPSE}
	\end{figure}
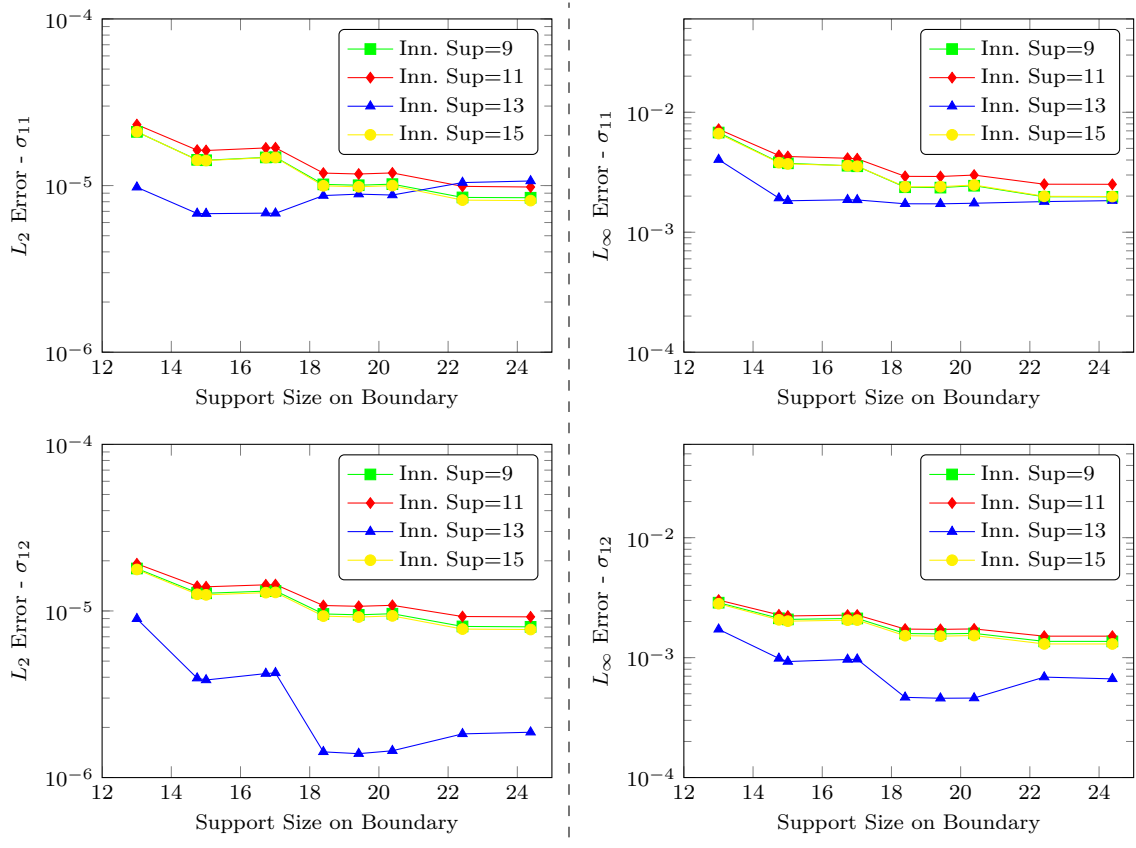
	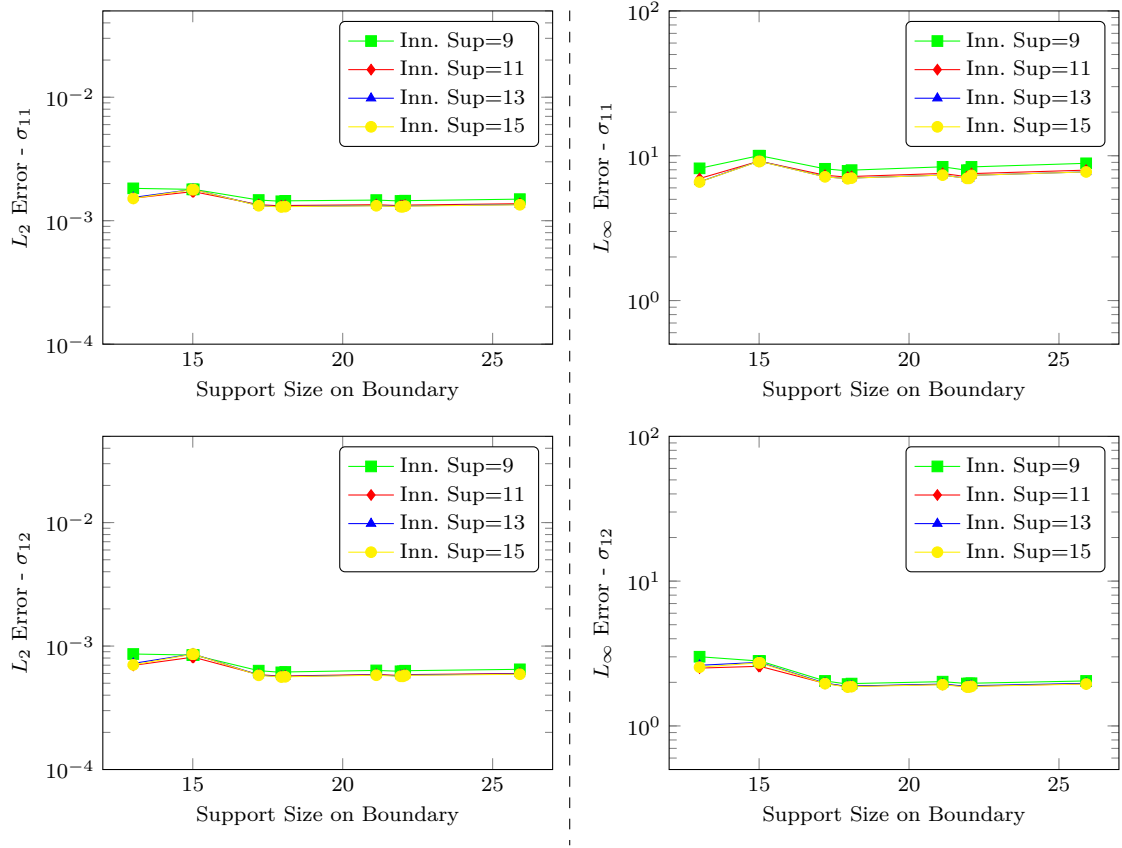
\begin{figure}[H] 
		\centering
		\begin{tabular}{c:c}
			\begin{tikzpicture}
			\begin{axis}[height=6cm,width=7.5cm, ymode=log, ymin=0.0001,ymax=0.05,xmin=12,xmax=27, legend entries={Inn. Sup=9,Inn. Sup=11,Inn. Sup=13,Inn. Sup=15},legend style={ at={(0.5,-0.2)},anchor=south west,legend columns=1, cells={anchor=west},  font=\footnotesize, rounded corners=2pt,}, legend pos=north east,xlabel=Support Size on Boundary,ylabel=$L_2$ Error - $\sigma_{11}$]
			\addplot+[green,mark=square*,mark options={fill=green}]  table [x=X-9.0, y=L2-REL-S11-9.0, col sep=comma] {SupportSizeDCPSE-LShaped.csv};
			\addplot+[red,mark=diamond*,mark options={fill=red}]   table [x=X-11.0, y=L2-REL-S11-11.0, col sep=comma] {SupportSizeDCPSE-LShaped.csv};
			\addplot+[blue,mark=triangle*,mark options={fill=blue}]   table [x=X-13.0, y=L2-REL-S11-13.0, col sep=comma] {SupportSizeDCPSE-LShaped.csv};
			\addplot+[yellow,mark=*,mark options={fill=yellow}]   table [x=X-15.0, y=L2-REL-S11-15.0, col sep=comma] {SupportSizeDCPSE-LShaped.csv};
			\end{axis}
			\end{tikzpicture} &
			
			\begin{tikzpicture}
			\begin{axis}[height=6cm,width=7.5cm,ymode=log, ymin=0.5,ymax=100,xmin=12,xmax=27, legend entries={Inn. Sup=9,Inn. Sup=11,Inn. Sup=13,Inn. Sup=15},legend style={ at={(0.5,-0.2)},anchor=south west,legend columns=1, cells={anchor=west},  font=\footnotesize, rounded corners=2pt,}, legend pos=north east,xlabel=Support Size on Boundary,ylabel=$L_{\infty}$ Error - $\sigma_{11}$]
			\addplot+[green,mark=square*,mark options={fill=green}]  table [x=X-9.0, y=L-INF-S11-9.0, col sep=comma] {SupportSizeDCPSE-LShaped.csv};
			\addplot+[red,mark=diamond*,mark options={fill=red}]   table [x=X-11.0, y=L-INF-S11-11.0, col sep=comma] {SupportSizeDCPSE-LShaped.csv};
			\addplot+[blue,mark=triangle*,mark options={fill=blue}]   table [x=X-13.0, y=L-INF-S11-13.0, col sep=comma] {SupportSizeDCPSE-LShaped.csv};
			\addplot+[yellow,mark=*,mark options={fill=yellow}]   table [x=X-15.0, y=L-INF-S11-15.0, col sep=comma] {SupportSizeDCPSE-LShaped.csv};
			\end{axis}
			\end{tikzpicture} \\
			
			\begin{tikzpicture}
			\begin{axis}[height=6cm,width=7.5cm, ymode=log, ymin=0.0001,ymax=0.05,xmin=12,xmax=27, legend entries={Inn. Sup=9,Inn. Sup=11,Inn. Sup=13,Inn. Sup=15},legend style={ at={(0.5,-0.2)},anchor=south west,legend columns=1, cells={anchor=west},  font=\footnotesize, rounded corners=2pt,}, legend pos=north east,xlabel=Support Size on Boundary,ylabel=$L_2$ Error - $\sigma_{12}$]
			\addplot+[green,mark=square*,mark options={fill=green}]  table [x=X-9.0, y=L2-REL-S12-9.0, col sep=comma] {SupportSizeDCPSE-LShaped.csv};
			\addplot+[red,mark=diamond*,mark options={fill=red}]   table [x=X-11.0, y=L2-REL-S12-11.0, col sep=comma] {SupportSizeDCPSE-LShaped.csv};
			\addplot+[blue,mark=triangle*,mark options={fill=blue}]   table [x=X-13.0, y=L2-REL-S12-13.0, col sep=comma] {SupportSizeDCPSE-LShaped.csv};
			\addplot+[yellow,mark=*,mark options={fill=yellow}]   table [x=X-15.0, y=L2-REL-S12-15.0, col sep=comma] {SupportSizeDCPSE-LShaped.csv};
			\end{axis}
			\end{tikzpicture} &
			
			\begin{tikzpicture}
			\begin{axis}[height=6cm,width=7.5cm,ymode=log, ymin=0.5,ymax=100,xmin=12,xmax=27, legend entries={Inn. Sup=9,Inn. Sup=11,Inn. Sup=13,Inn. Sup=15},legend style={ at={(0.5,-0.2)},anchor=south west,legend columns=1, cells={anchor=west},  font=\footnotesize, rounded corners=2pt,}, legend pos=north east,xlabel=Support Size on Boundary,ylabel=$L_{\infty}$ Error - $\sigma_{12}$]
			\addplot+[green,mark=square*,mark options={fill=green}]  table [x=X-9.0, y=L-INF-S12-9.0, col sep=comma] {SupportSizeDCPSE-LShaped.csv};
			\addplot+[red,mark=diamond*,mark options={fill=red}]   table [x=X-11.0, y=L-INF-S12-11.0, col sep=comma] {SupportSizeDCPSE-LShaped.csv};
			\addplot+[blue,mark=triangle*,mark options={fill=blue}]   table [x=X-13.0, y=L-INF-S12-13.0, col sep=comma] {SupportSizeDCPSE-LShaped.csv};
			\addplot+[yellow,mark=*,mark options={fill=yellow}]   table [x=X-15.0, y=L-INF-S12-15.0, col sep=comma] {SupportSizeDCPSE-LShaped.csv};
			\end{axis}
			\end{tikzpicture} \\	
		\end{tabular}
		\caption{DC PSE Support Node Number Sensitivity - 2D L-Shape. $L_2$ (Left) and $L_{\infty}$ (Right) errors for various combinations of inner nodes and boundary nodes support sizes. All the combinations of inner and boundary nodes support size lead to similar errors. This is due to the system being loaded via Dirichlet boundary conditions.}
		\label{BoundaryNodeResultsDCPSE_LShaped}
	\end{figure}
	It can be observed from Figure \ref{BoundaryNodeResultsDCPSE} that an inner node support composed of 13 nodes leads almost always to the minimum error. It can also be observed that an increasing number of nodes in the support of the boundary nodes reduces relatively steadily the error for the 2D cylinder. An error reduction of a factor two is observed when increasing the number of support nodes from 13 to 19 for most inner nodes support sizes.

	It can be observed from Figure \ref{BoundaryNodeResultsDCPSE_LShaped} that the number of support nodes on the boundary has little effect on the observed error for the L-shape problem. This is because the model is loaded via Dirichlet boundary conditions. The number of support nodes for inner collocation nodes has little impact on the error.
	
	Based the results from Figure \ref{BoundaryNodeResultsDCPSE} and Figure \ref{BoundaryNodeResultsDCPSE_LShaped}, 13 support nodes for inner collocation nodes and 19 support nodes for boundary collocation nodes is a reasonable choice for most problems as it leads to a minimum error.
	\begin{figure}[H] 
		\centering
		\begin{tabular}{c:c}
			\begin{tikzpicture}
			\begin{axis}[height=6cm,width=7.5cm, ymin=0.00003,ymax=0.00005,xmin=55,xmax=95, legend entries={Inn. Sup=35,Inn. Sup=37,Inn. Sup=39,Inn. Sup=41},legend style={ at={(0.5,-0.2)},anchor=south west,legend columns=1, cells={anchor=west},  font=\footnotesize, rounded corners=2pt,}, legend pos=north east,xlabel=Support Size on Boundary,ylabel=$L_2$ Error - $\sigma_{11}$]
			\addplot+[green,mark=square*,mark options={fill=green}]  table [x=X-3D, y=L2-REL-S11-35, col sep=comma] {SupportSizeDCPSE3D.csv};
			\addplot+[red,mark=diamond*,mark options={fill=red}]   table [x=X-3D, y=L2-REL-S11-37, col sep=comma] {SupportSizeDCPSE3D.csv};
			\addplot+[blue,mark=triangle*,mark options={fill=blue}]   table [x=X-3D, y=L2-REL-S11-39, col sep=comma] {SupportSizeDCPSE3D.csv};
			\addplot+[yellow,mark=*,mark options={fill=yellow}]   table [x=X-3D, y=L2-REL-S11-41, col sep=comma] {SupportSizeDCPSE3D.csv};
			\end{axis}
			\end{tikzpicture} &
			
			\begin{tikzpicture}
			\begin{axis}[height=6cm,width=7.5cm,xmin=55,xmax=95,ymin=0.03,ymax=0.05, legend entries={Inn. Sup=35,Inn. Sup=37,Inn. Sup=39,Inn. Sup=41},legend style={ at={(0.5,-0.2)},anchor=south west,legend columns=1, cells={anchor=west},  font=\footnotesize, rounded corners=2pt,}, legend pos=north east,xlabel=Support Size on Boundary,ylabel=$L_{\infty}$ Error - $\sigma_{11}$]
			\addplot+[green,mark=square*,mark options={fill=green}]  table [x=X-3D, y=L-INF-S11-35, col sep=comma] {SupportSizeDCPSE3D.csv};
			\addplot+[red,mark=diamond*,mark options={fill=red}]   table [x=X-3D, y=L-INF-S11-37, col sep=comma] {SupportSizeDCPSE3D.csv};
			\addplot+[blue,mark=triangle*,mark options={fill=blue}]   table [x=X-3D, y=L-INF-S11-39, col sep=comma] {SupportSizeDCPSE3D.csv};
			\addplot+[yellow,mark=*,mark options={fill=yellow}]   table [x=X-3D, y=L-INF-S11-41, col sep=comma] {SupportSizeDCPSE3D.csv};
			\end{axis}
			\end{tikzpicture} \\
			
			\begin{tikzpicture}
			\begin{axis}[height=6cm,width=7.5cm,xmin=55,xmax=95,ymin=0.000006,ymax=0.00001, legend entries={Inn. Sup=35,Inn. Sup=37,Inn. Sup=39,Inn. Sup=41},legend style={ at={(0.5,-0.2)},anchor=south west,legend columns=1, cells={anchor=west},  font=\footnotesize, rounded corners=2pt,}, legend pos=north east,xlabel=Support Size on Boundary,ylabel=$L_2$ Error - $\sigma_{12}$]
			\addplot+[green,mark=square*,mark options={fill=green}]  table [x=X-3D, y=L2-REL-S12-35, col sep=comma] {SupportSizeDCPSE3D.csv};
			\addplot+[red,mark=diamond*,mark options={fill=red}]   table [x=X-3D, y=L2-REL-S12-37, col sep=comma] {SupportSizeDCPSE3D.csv};
			\addplot+[blue,mark=triangle*,mark options={fill=blue}]   table [x=X-3D, y=L2-REL-S12-39, col sep=comma] {SupportSizeDCPSE3D.csv};
			\addplot+[yellow,mark=*,mark options={fill=yellow}]   table [x=X-3D, y=L2-REL-S12-41, col sep=comma] {SupportSizeDCPSE3D.csv};
			\end{axis}
			\end{tikzpicture} &
			
			\begin{tikzpicture}
			\begin{axis}[height=6cm,width=7.5cm,xmin=55,xmax=95,ymin=0.008,ymax=0.013, legend entries={Inn. Sup=35,Inn. Sup=37,Inn. Sup=39,Inn. Sup=41},legend style={ at={(0.5,-0.2)},anchor=south west,legend columns=1, cells={anchor=west},  font=\footnotesize, rounded corners=2pt,}, legend pos=north east,xlabel=Support Size on Boundary,ylabel=$L_{\infty}$ Error - $\sigma_{12}$]
			\addplot+[green,mark=square*,mark options={fill=green}]  table [x=X-3D, y=L-INF-S12-35, col sep=comma] {SupportSizeDCPSE3D.csv};
			\addplot+[red,mark=diamond*,mark options={fill=red}]   table [x=X-3D, y=L-INF-S12-37, col sep=comma] {SupportSizeDCPSE3D.csv};
			\addplot+[blue,mark=triangle*,mark options={fill=blue}]   table [x=X-3D, y=L-INF-S12-39, col sep=comma] {SupportSizeDCPSE3D.csv};
			\addplot+[yellow,mark=*,mark options={fill=yellow}]   table [x=X-3D, y=L-INF-S12-41, col sep=comma] {SupportSizeDCPSE3D.csv};
			\end{axis}
			\end{tikzpicture} \\	
		\end{tabular}
		\caption{DC PSE Support Node Number Sensitivity - 3D Sphere. $L_2$ (Left) and $L_{\infty}$ (Right) errors for various combinations of inner node and boundary node support sizes. Inner collocation nodes with 37 support nodes lead to the lowest observed error.}
		\label{BoundaryNode_3DResultsDCPSE}
	\end{figure}
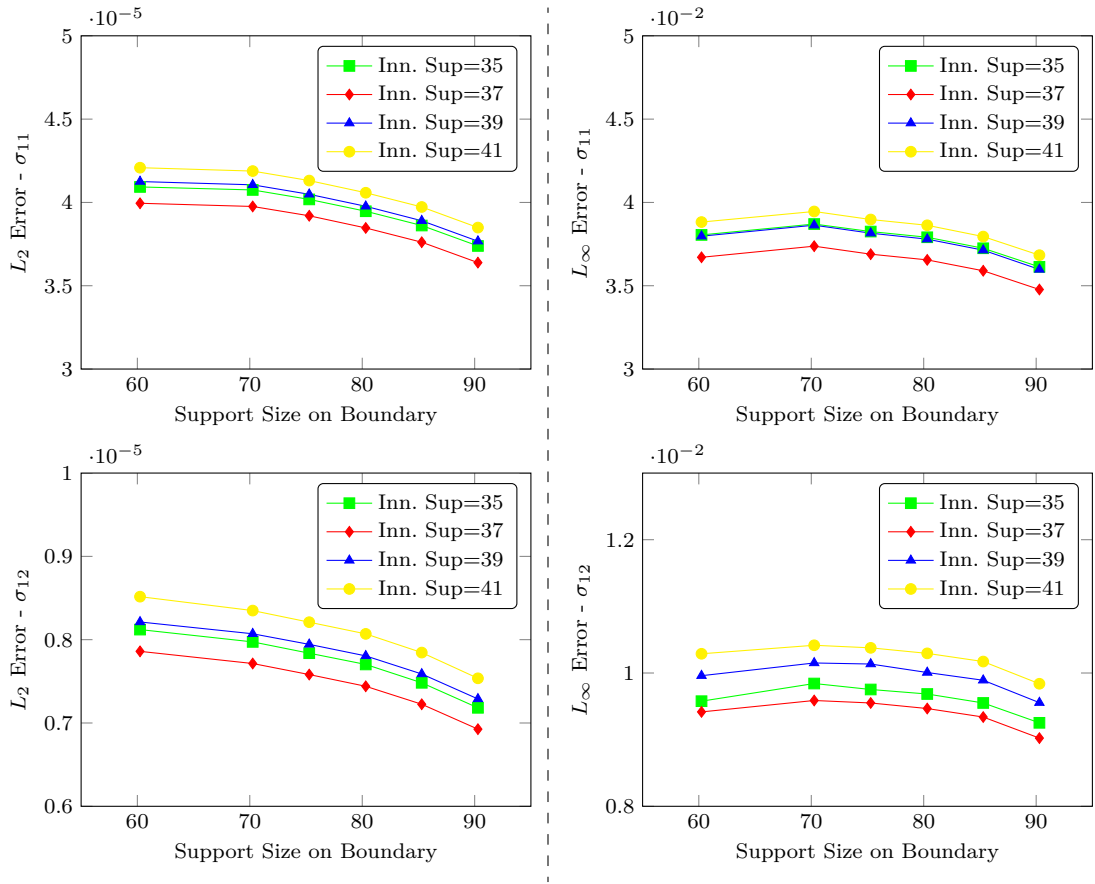
	It can be observed from Figure \ref{BoundaryNode_3DResultsDCPSE} that, for the 3D sphere, a minimum error is observed for inner node supports composed of 37 nodes. The number of boundary support nodes has little impact on the error. Increasing the number of support nodes from 60 to 90 reduces by only 8\% in average the observed error. 75 support nodes for boundary collocation nodes has been selected as for the GFD method in order to keep the fill of the system matrix as low as possible while maintaining the error low.
	
	\subsection{Results Summary}
	
	Based on the results presented in the above sections, the parameters that lead to a minimum error, while maintaining the computational expense reasonably low, are presented in Table \ref{ParametersSummary} below. These parameters are expected to lead to a low error for a wide variety of linear elasticity problems, including singular problems.	They have been used as a base case for the studies presented in the next sections of the paper.
	\begin{table}[h]
		\centering
		\caption{Summary of the results from the parametric study}
		\label{ParametersSummary}
		\renewcommand{\arraystretch}{1.5}
		\begin{tabular}{|l|c|c|}
			\hline
			\multicolumn{1}{|c|}{\textbf{Parameter}} & \multicolumn{1}{c|}{\textbf{GFD}} &
			\multicolumn{1}{c|}{\textbf{DC PSE}} \\
			\hline
			Weight Function Type & 4$^\text{th}$ Order Spline & Exponential \\
			Weight Function Parameter & $\gamma=0.75$ &  $\alpha=1$, $\epsilon$=0.30\\
			Correction Function & N/A & Polynomial \\
			Size of Inner Nodes Support (2D/3D) & 11/37 & 13/37 \\
			Size of Boundary Nodes Support (2D/3D) & 19/75 & 19/75 \\
			\hline
		\end{tabular}
	\end{table}
	
	\section{Improvement Methods} \label{ImprovementMethods_Section}
	
	In this section, we present three methods which are expected to improve the accuracy of the GFD and DC PSE methods.
	
	\subsection{Use of a Voronoi Diagram in Collocation} \label{Voronoi_Section}
	
	\subsubsection{General}
	
\begin{figure}[H] 
		\centering
		\begin{tikzpicture}
		\definecolor{GreyColor}{rgb}{0.3,0.3,0.3}
		\begin{axis}[height=7.5cm,width=7.5cm, xmin=-0.1,xmax=0.1,scaled x ticks = false, ymin=-0.1,ymax=0.1,scaled y ticks = false,anchor=south west,axis line style={draw=none},tick style={draw=none},xticklabels={,,}, yticklabels={,,}, cells={anchor=west},  font=\footnotesize, rounded corners=2pt]
		\addplot[GreyColor,line width=2pt]  table [x=X-1, y=Y-1, col sep=comma] {VoroCellX1.csv};
		\addplot[GreyColor,line width=2pt]  table [x=X-2, y=Y-2, col sep=comma] {VoroCellX1.csv};
		\addplot[GreyColor,line width=2pt]  table [x=X-3, y=Y-3, col sep=comma] {VoroCellX1.csv};
		\addplot[GreyColor,line width=2pt]  table [x=X-4, y=Y-4, col sep=comma] {VoroCellX1.csv};
		\addplot[GreyColor,line width=2pt]  table [x=X-5, y=Y-5, col sep=comma] {VoroCellX1.csv};
		\addplot[GreyColor,line width=2pt]  table [x=X-6, y=Y-6, col sep=comma] {VoroCellX6.csv};
		\addplot[blue,line width=3pt]  table [x=X-Circle, y=Y-Circle, col sep=comma] {VoroCellCircle.txt};
		\addplot[only marks,red,mark=*,mark options={fill=red}]  table [x=X-Scatt, y=Y-Scatt, col sep=comma] {VoroCellScatt.txt};
		\end{axis}
		\draw [black,-stealth, line width=1.0pt] (0,3.5) node [left] {Collocation Node} -- (2.9,3);
		\draw [black,-stealth, line width=1.0pt] (0,5.5) node [left] {Support Node} -- (2.4,4.4);
		\end{tikzpicture}
		\caption{2D Voronoi Diagram on the disc support of a collocation node. The cells associated to each node are delimited by gray lines, while the boundary of the support is drawn in blue.}
		\label{VoronoiCell}
	\end{figure}
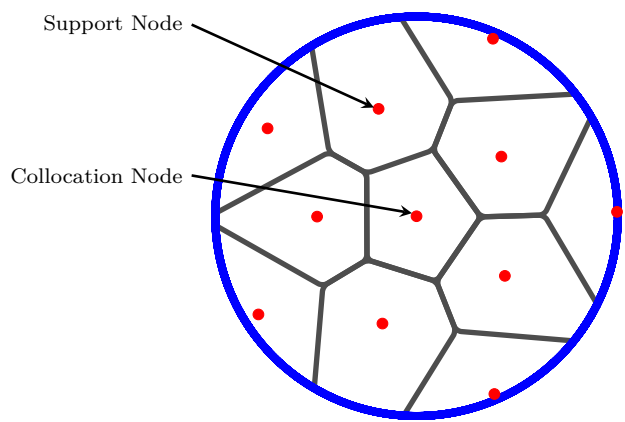

	A Voronoi diagram is a partition of a selected region over which nodes are distributed. A cell is associated to each node. The boundaries of the cell are defined so that all the points contained in it are closer to the cell reference node than to any other node of the domain. Figure \ref{VoronoiCell} shows a typical 2D Voronoi diagram drawn on the support of an inner node of the domain. The boundary of the support is drawn in blue, the nodes in red and the Voronoi cells are limited by grey lines. Sukumar \cite{Sukumar2003} and Zhou et al. \cite{Zhou2007}, respectively, used Voronoi diagrams for node selection, and body integration. The purpose of this section is to assess if using a Voronoi diagram on the support of a collocation node helps reducing the error of the considered methods. 
	
	\subsubsection{Application to the GFD and DC PSE Methods} \label{VoroTheory}
	
	\paragraph{GFD} \
	
	The principle of the GFD method has been presented in Section \ref{GFD_Method}. When more nodes than derivatives are present in the node support, a mean least square approximation is used to determine the field derivatives that best fit the distribution.
	The contribution of each node in the least square approximation is weighted by a function which only depends on the distance between the reference node and the support node. A Voronoi diagram can be used to determine an additional weight based on the spatial arrangement of the nodes. This weight is the area or volume $v$ of the considered Voronoi cell and is multiplied with the distance based weight $w$. Equation (\ref{Moments_GFD}) becomes:
	\begin{equation} \label{Moments_GFD_Voro}
	m_{ij}= \sum_{k=1}^m {w(\mathbf{X_{pk}} - \mathbf{X_c}) v(\mathbf{X_{pk}}^c) P_{ik}(\mathbf{X_c}) P_{jk}(\mathbf{X_c})}.
	\end{equation}
	
	\paragraph{DC PSE} \
	
	One of the key aspects of the DC PSE method presented in Section \ref{DC-PSE_Method} is the convolution of the Taylor's series expansion with a correction function $\eta$. The domain integral is transformed into a discrete summation with a volume $V_p$ associated to each particle $\mathbf{X_p}$ of the support. In a first approximation, all $V_p$ values are set to unity. In order to improve the accuracy of the method, a Voronoi diagram can be used to set $V_p$ equal to the volume of the Voronoi cell associated to each node $\mathbf{X_p}$. 
	
	\subsubsection{Results}
	
\begin{figure}[H] 
		\centering
		\begin{tabular}{c:c}
			\begin{tikzpicture}
			\begin{axis}[height=8cm,width=8cm, xmin=-0.2,xmax=3.2, ymin=-0.2,ymax=3.2,mark size=0.5pt,legend style={ at={(0.5,-0.2)},anchor=south west,legend columns=1, cells={anchor=west},  font=\footnotesize, rounded corners=2pt,}, legend pos=north east,xlabel=X,ylabel=Y]
			\addplot+[black, only marks,mark=*,mark options={fill=black}]  table [x=X-STRUCT, y=Y-STRUCT, col sep=comma] {CylStructFree.csv};
			\end{axis}
			\end{tikzpicture} &
			
			\begin{tikzpicture}
			\begin{axis}[height=8cm,width=8cm, xmin=-0.2,xmax=3.2, ymin=-0.2,ymax=3.2,mark size=0.5pt,legend style={ at={(0.5,-0.2)},anchor=south west,legend columns=1, cells={anchor=west},  font=\footnotesize, rounded corners=2pt,}, legend pos=north east,xlabel=X,ylabel=Y]
			\addplot+[black, only marks,mark=*,mark options={fill=black}]  table [x=X-FREE, y=Y-FREE, col sep=comma] {CylStructFree.csv};
			\end{axis}
			\end{tikzpicture}
		\end{tabular}
		\caption{2D Cylinder Node Distribution - Structured consisting of 1680 Nodes (Left) and Free consisting of 1762 Nodes (Right). The structured node distribution is based on constant angle and radius increments while the free node distribution uses a Delaunay triangulation of the domain.}
		\label{NodeDispositionVoro}
	\end{figure}
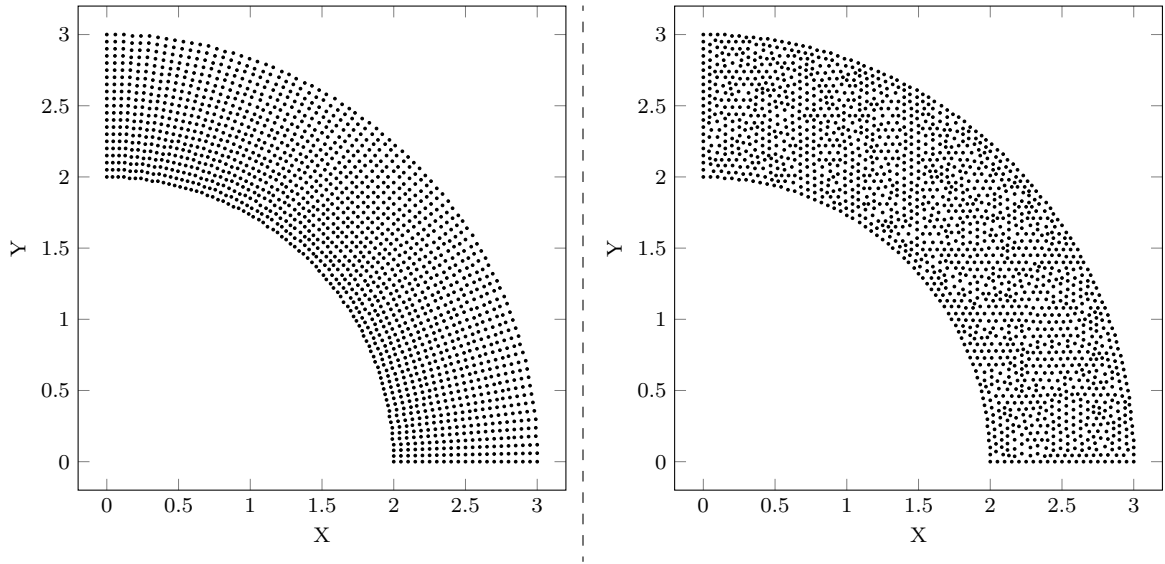

	In this section, we present the results for the 2D cylinder and the 2D L-shape problems. The error are compared between a model where Voroni weights are used, and a model where these weights are not considered. Two types of node distributions are considered: a structured and a free node distribution. The structured node distribution is created using a constant angle and radius increment for the 2D cylinder. For the 2D L-shape, a grid-type arrangement is used. The free node distribution uses a Delaunay triangulation of the domain for both problems. The two types of node arrangements are presented in Figure \ref{NodeDispositionVoro} for the 2D cylinder problem.
	
	\paragraph{GFD} \
	
	The results obtained with the GFD method are presented in Figure \ref{ResultsVoroGFD} and Figure \ref{ResultsVoroGFD_LShape} respectively for the 2D cylinder and the 2D L-shape problems for the  both node distributions. A slight error reduction can be observed for the $L_2$ and the $L_{\infty}$ error norms when Voronoi weights are used but this reduction is not observed for all node densities. For the 2D cylinder and the 2D L-shape, an error reduction of around 2\% is observed when Voronoi based weights are used with a regular discretization of the domain. A more significant error reduction is observed for the 2D cylinder with a free discretization of the domain. The error reduction is of around 17\%. For the 2D L-shape, an error increase of 3\% is observed when Voronoi weights are used with a free discretization of the domain. 
	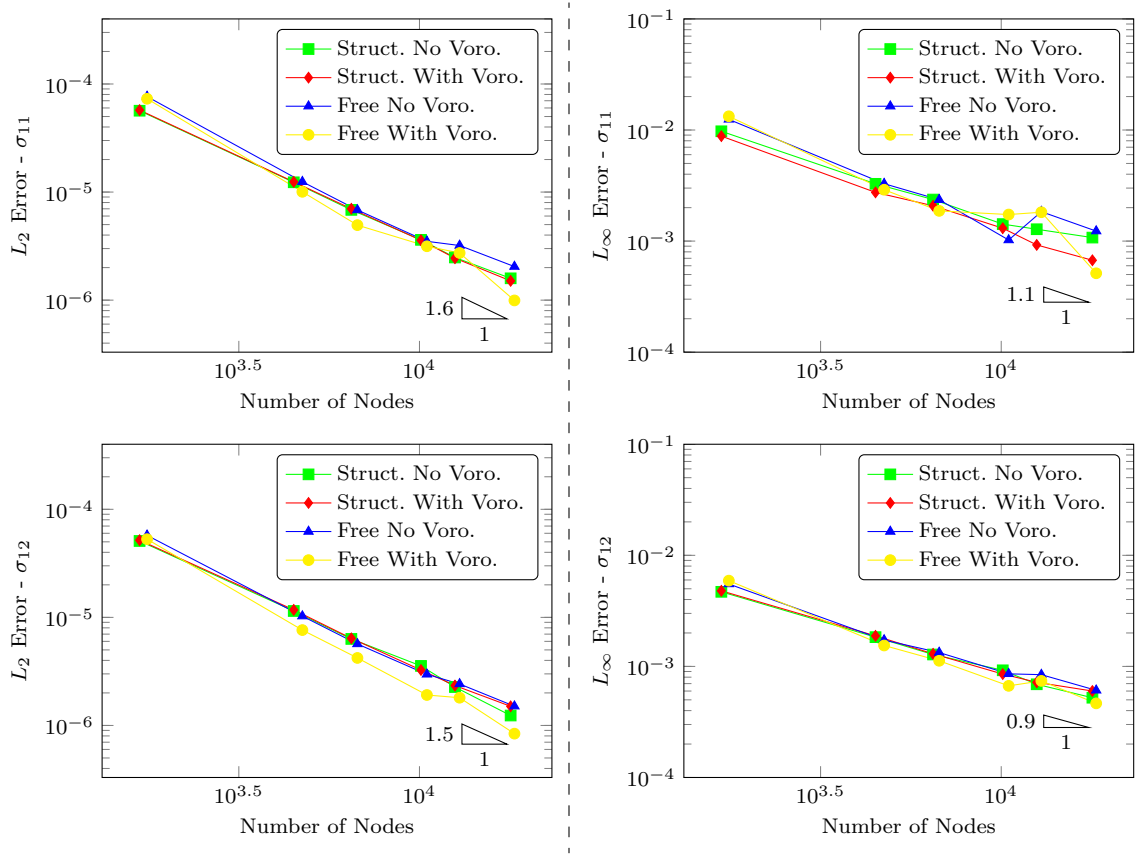
\begin{figure}[H] 
		\centering
		\begin{tabular}{c:c}
			\begin{tikzpicture}
			\begin{axis}[height=6cm,width=7.5cm, xmode=log, ymode=log, ymin=0.00000033,ymax=0.0004, legend entries={Struct. No Voro.,Struct. With Voro.,Free No Voro.,Free With Voro.},legend style={ at={(0.5,-0.2)},anchor=south west,legend columns=1, cells={anchor=west},  font=\footnotesize, rounded corners=2pt,}, legend pos=north east,xlabel=Number of Nodes,ylabel=$L_2$ Error - $\sigma_{11}$]
			\addplot+[green,mark=square*,mark options={fill=green}]  table [x=NODE-NUM-STRUC, y=S11-L2-REL-GFD-NoVoro-Struc, col sep=comma] {VoroResults.csv};
			\addplot+[red,mark=diamond*,mark options={fill=red}]   table [x=NODE-NUM-STRUC, y=S11-L2-REL-GFD-WithVoro-Struc, col sep=comma] {VoroResults.csv};
			\addplot+[blue,mark=triangle*,mark options={fill=blue}]   table [x=NODE-NUM-FREE, y=S11-L2-REL-GFD-NoVoro-Free, col sep=comma] {VoroResults.csv};
			\addplot+[yellow,mark=*,mark options={fill=yellow}]   table [x=NODE-NUM-FREE, y=S11-L2-REL-GFD-WithVoro-Free, col sep=comma] {VoroResults.csv};
			\logLogSlopeTriangle{0.9}{0.1}{0.1}{1.6}{black};
			\end{axis}
			\end{tikzpicture} &
			
			\begin{tikzpicture}
			\begin{axis}[height=6cm,width=7.5cm, xmode=log, ymode=log, ymin=0.0001,ymax=0.1, legend entries={Struct. No Voro.,Struct. With Voro.,Free No Voro.,Free With Voro.},legend style={ at={(0.5,-0.2)},anchor=south west,legend columns=1, cells={anchor=west},  font=\footnotesize, rounded corners=2pt,}, legend pos=north east,xlabel=Number of Nodes,ylabel=$L_{\infty}$ Error - $\sigma_{11}$]
			\addplot+[green,mark=square*,mark options={fill=green}]  table [x=NODE-NUM-STRUC, y=S11-LINF-GFD-NoVoro-Struc, col sep=comma] {VoroResults.csv};
			\addplot+[red,mark=diamond*,mark options={fill=red}]   table [x=NODE-NUM-STRUC, y=S11-LINF-GFD-WithVoro-Struc, col sep=comma] {VoroResults.csv};
			\addplot+[blue,mark=triangle*,mark options={fill=blue}]   table [x=NODE-NUM-FREE, y=S11-LINF-GFD-NoVoro-Free, col sep=comma] {VoroResults.csv};
			\addplot+[yellow,mark=*,mark options={fill=yellow}]   table [x=NODE-NUM-FREE, y=S11-LINF-GFD-WithVoro-Free, col sep=comma] {VoroResults.csv};
			\logLogSlopeTriangle{0.9}{0.1}{0.15}{1.1}{black};
			\end{axis}
			\end{tikzpicture} \\
			
			\begin{tikzpicture}
			\begin{axis}[height=6cm,width=7.5cm, xmode=log, ymode=log, ymin=0.00000033,ymax=0.0004, legend entries={Struct. No Voro.,Struct. With Voro.,Free No Voro.,Free With Voro.},legend style={ at={(0.5,-0.2)},anchor=south west,legend columns=1, cells={anchor=west},  font=\footnotesize, rounded corners=2pt,}, legend pos=north east,xlabel=Number of Nodes,ylabel=$L_2$ Error - $\sigma_{12}$]
			\addplot+[green,mark=square*,mark options={fill=green}]  table [x=NODE-NUM-STRUC, y=S12-L2-REL-GFD-NoVoro-Struc, col sep=comma] {VoroResults.csv};
			\addplot+[red,mark=diamond*,mark options={fill=red}]   table [x=NODE-NUM-STRUC, y=S12-L2-REL-GFD-WithVoro-Struc, col sep=comma] {VoroResults.csv};
			\addplot+[blue,mark=triangle*,mark options={fill=blue}]   table [x=NODE-NUM-FREE, y=S12-L2-REL-GFD-NoVoro-Free, col sep=comma] {VoroResults.csv};
			\addplot+[yellow,mark=*,mark options={fill=yellow}]   table [x=NODE-NUM-FREE, y=S12-L2-REL-GFD-WithVoro-Free, col sep=comma] {VoroResults.csv};
			\logLogSlopeTriangle{0.9}{0.1}{0.1}{1.5}{black};
			\end{axis}
			\end{tikzpicture} &
			
			\begin{tikzpicture}
			\begin{axis}[height=6cm,width=7.5cm, xmode=log, ymode=log, ymin=0.0001,ymax=0.1, legend entries={Struct. No Voro.,Struct. With Voro.,Free No Voro.,Free With Voro.},legend style={ at={(0.5,-0.2)},anchor=south west,legend columns=1, cells={anchor=west},  font=\footnotesize, rounded corners=2pt,}, legend pos=north east,xlabel=Number of Nodes,ylabel=$L_{\infty}$ Error - $\sigma_{12}$]
			\addplot+[green,mark=square*,mark options={fill=green}]  table [x=NODE-NUM-STRUC, y=S12-LINF-GFD-NoVoro-Struc, col sep=comma] {VoroResults.csv};
			\addplot+[red,mark=diamond*,mark options={fill=red}]   table [x=NODE-NUM-STRUC, y=S12-LINF-GFD-WithVoro-Struc, col sep=comma] {VoroResults.csv};
			\addplot+[blue,mark=triangle*,mark options={fill=blue}]   table [x=NODE-NUM-FREE, y=S12-LINF-GFD-NoVoro-Free, col sep=comma] {VoroResults.csv};
			\addplot+[yellow,mark=*,mark options={fill=yellow}]   table [x=NODE-NUM-FREE, y=S12-LINF-GFD-WithVoro-Free, col sep=comma] {VoroResults.csv};
			\logLogSlopeTriangle{0.9}{0.1}{0.15}{0.9}{black};
			\end{axis}
			\end{tikzpicture} 
		\end{tabular}
		\caption{Impact of Voronoi based Weights on the Errors for the GFD Method - 2D Cylinder. $L_2$ (Left) and $L_{\infty}$ (Right) errors for structured and free node distributions. A reduction in the error is only observed for the free node distribution.}
		\label{ResultsVoroGFD}
	\end{figure}
	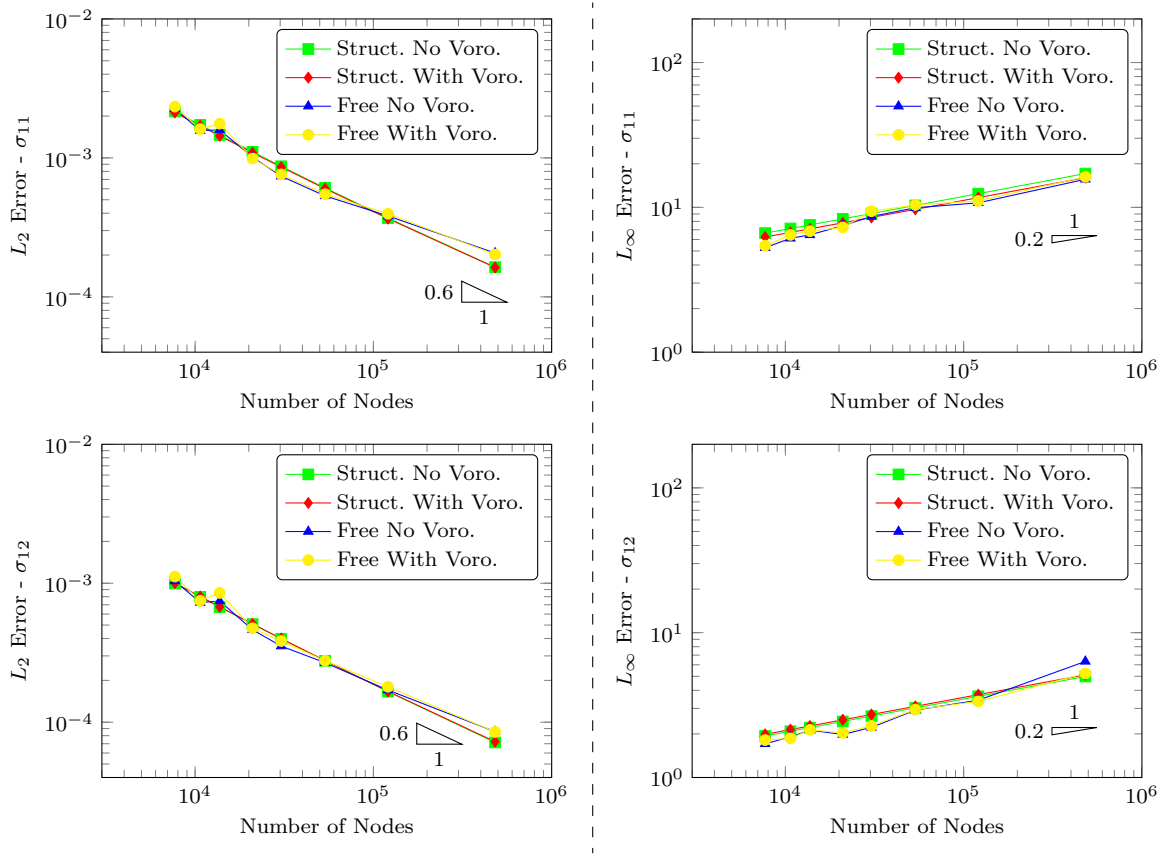
\begin{figure}[H] 
		\centering
		\begin{tabular}{c:c}
			\begin{tikzpicture}
			\begin{axis}[height=6cm,width=7.5cm, xmode=log,xmin=3000,xmax=1000000, ymode=log, ymin=0.00004,ymax=0.01, legend entries={Struct. No Voro.,Struct. With Voro.,Free No Voro.,Free With Voro.},legend style={ at={(0.5,-0.2)},anchor=south west,legend columns=1, cells={anchor=west},  font=\footnotesize, rounded corners=2pt,}, legend pos=north east,xlabel=Number of Nodes,ylabel=$L_2$ Error - $\sigma_{11}$]
			\addplot+[green,mark=square*,mark options={fill=green}]  table [x=NODE-NUM-STRUC, y=S11-L2-REL-GFD-NoVoro-Struc, col sep=comma] {VoroResults-LShape.csv};
			\addplot+[red,mark=diamond*,mark options={fill=red}]   table [x=NODE-NUM-STRUC, y=S11-L2-REL-GFD-WithVoro-Struc, col sep=comma] {VoroResults-LShape.csv};
			\addplot+[blue,mark=triangle*,mark options={fill=blue}]   table [x=NODE-NUM-FREE, y=S11-L2-REL-GFD-NoVoro-Free, col sep=comma] {VoroResults-LShape.csv};
			\addplot+[yellow,mark=*,mark options={fill=yellow}]   table [x=NODE-NUM-FREE, y=S11-L2-REL-GFD-WithVoro-Free, col sep=comma] {VoroResults-LShape.csv};
			\logLogSlopeTriangle{0.9}{0.1}{0.15}{0.6}{black};
			\end{axis}
			\end{tikzpicture} &
			
			\begin{tikzpicture}
			\begin{axis}[height=6cm,width=7.5cm, xmode=log,xmin=3000,xmax=1000000, ymode=log, ymin=1,ymax=200, legend entries={Struct. No Voro.,Struct. With Voro.,Free No Voro.,Free With Voro.},legend style={ at={(0.5,-0.2)},anchor=south west,legend columns=1, cells={anchor=west},  font=\footnotesize, rounded corners=2pt,}, legend pos=north east,xlabel=Number of Nodes,ylabel=$L_{\infty}$ Error - $\sigma_{11}$]
			\addplot+[green,mark=square*,mark options={fill=green}]  table [x=NODE-NUM-STRUC, y=S11-LINF-GFD-NoVoro-Struc, col sep=comma] {VoroResults-LShape.csv};
			\addplot+[red,mark=diamond*,mark options={fill=red}]   table [x=NODE-NUM-STRUC, y=S11-LINF-GFD-WithVoro-Struc, col sep=comma] {VoroResults-LShape.csv};
			\addplot+[blue,mark=triangle*,mark options={fill=blue}]   table [x=NODE-NUM-FREE, y=S11-LINF-GFD-NoVoro-Free, col sep=comma] {VoroResults-LShape.csv};
			\addplot+[yellow,mark=*,mark options={fill=yellow}]   table [x=NODE-NUM-FREE, y=S11-LINF-GFD-WithVoro-Free, col sep=comma] {VoroResults-LShape.csv};
			\logLogSlopeTriangleUp{0.9}{0.1}{0.35}{0.2}{black};
			\end{axis}
			\end{tikzpicture} \\
			
			\begin{tikzpicture}
			\begin{axis}[height=6cm,width=7.5cm, xmode=log,xmin=3000,xmax=1000000, ymode=log, ymin=0.00004,ymax=0.01, legend entries={Struct. No Voro.,Struct. With Voro.,Free No Voro.,Free With Voro.},legend style={ at={(0.5,-0.2)},anchor=south west,legend columns=1, cells={anchor=west},  font=\footnotesize, rounded corners=2pt,}, legend pos=north east,xlabel=Number of Nodes,ylabel=$L_2$ Error - $\sigma_{12}$]
			\addplot+[green,mark=square*,mark options={fill=green}]  table [x=NODE-NUM-STRUC, y=S12-L2-REL-GFD-NoVoro-Struc, col sep=comma] {VoroResults-LShape.csv};
			\addplot+[red,mark=diamond*,mark options={fill=red}]   table [x=NODE-NUM-STRUC, y=S12-L2-REL-GFD-WithVoro-Struc, col sep=comma] {VoroResults-LShape.csv};
			\addplot+[blue,mark=triangle*,mark options={fill=blue}]   table [x=NODE-NUM-FREE, y=S12-L2-REL-GFD-NoVoro-Free, col sep=comma] {VoroResults-LShape.csv};
			\addplot+[yellow,mark=*,mark options={fill=yellow}]   table [x=NODE-NUM-FREE, y=S12-L2-REL-GFD-WithVoro-Free, col sep=comma] {VoroResults-LShape.csv};
			\logLogSlopeTriangle{0.8}{0.1}{0.1}{0.6}{black};
			\end{axis}
			\end{tikzpicture} &
			
			\begin{tikzpicture}
			\begin{axis}[height=6cm,width=7.5cm, xmode=log,xmin=3000,xmax=1000000, ymode=log, ymin=1,ymax=200, legend entries={Struct. No Voro.,Struct. With Voro.,Free No Voro.,Free With Voro.},legend style={ at={(0.5,-0.2)},anchor=south west,legend columns=1, cells={anchor=west},  font=\footnotesize, rounded corners=2pt,}, legend pos=north east,xlabel=Number of Nodes,ylabel=$L_{\infty}$ Error - $\sigma_{12}$]
			\addplot+[green,mark=square*,mark options={fill=green}]  table [x=NODE-NUM-STRUC, y=S12-LINF-GFD-NoVoro-Struc, col sep=comma] {VoroResults-LShape.csv};
			\addplot+[red,mark=diamond*,mark options={fill=red}]   table [x=NODE-NUM-STRUC, y=S12-LINF-GFD-WithVoro-Struc, col sep=comma] {VoroResults-LShape.csv};
			\addplot+[blue,mark=triangle*,mark options={fill=blue}]   table [x=NODE-NUM-FREE, y=S12-LINF-GFD-NoVoro-Free, col sep=comma] {VoroResults-LShape.csv};
			\addplot+[yellow,mark=*,mark options={fill=yellow}]   table [x=NODE-NUM-FREE, y=S12-LINF-GFD-WithVoro-Free, col sep=comma] {VoroResults-LShape.csv};
			\logLogSlopeTriangleUp{0.9}{0.1}{0.15}{0.2}{black};
			\end{axis}
			\end{tikzpicture} 
		\end{tabular}
		\caption{Impact of Voronoi based Weights on the Errors for the GFD Method - 2D L-Shape. $L_2$ (Left) and $L_{\infty}$ (Right) errors for structured and free node distributions. A slight reduction in the error is observed for the structured node distribution.}
		\label{ResultsVoroGFD_LShape}
	\end{figure}

	\paragraph{DC PSE} \
	
	As for the GFD method, we now assess the impact of Voronoi based volumes for the two node distributions for the 2D cylinder and the 2D L-shape problems. The results are presented in Figure \ref{ResultsVoroDCPSE} and Figure \ref{ResultsVoroDCPSE_LShape} below.
	
	It can be observed from Figure \ref{ResultsVoroDCPSE} that, for the 2D cylinder, the use Voronoi based volumes leads to a large error increase for the structured node distribution. For the free node distribution, an average reduction of 10\% is observed.
	
	From Figure \ref{ResultsVoroDCPSE_LShape} we can see that the trend for the 2D L-shape is the opposite than for the 2D cylinder. An error reduction of around 5\% is observed when Voronoi based volumes are used with the structured node distribution. A slight error increase (less than 1\%) is observed for the free node distribution when Voronoi based volumes are used.
	\begin{figure}[H] 
		\centering
		\begin{tabular}{c:c}
			\begin{tikzpicture}
			\begin{axis}[height=6cm,width=7.5cm, xmode=log, ymode=log, ymin=0.00000003,ymax=0.007, legend entries={Struct. No Voro.,Struct. With Voro.,Free No Voro.,Free With Voro.},legend style={ at={(0.5,-0.2)},anchor=south west,legend columns=1, cells={anchor=west},  font=\footnotesize, rounded corners=2pt,}, legend pos=north east,xlabel=Number of Nodes,ylabel=$L_2$ Error - $\sigma_{11}$]
			\addplot+[green,mark=square*,mark options={fill=green}]  table [x=NODE-NUM-STRUC, y=S11-L2-REL-DCPSE-NoVoro-Struc, col sep=comma] {VoroResults.csv};
			\addplot+[red,mark=diamond*,mark options={fill=red}]   table [x=NODE-NUM-STRUC, y=S11-L2-REL-DCPSE-WithVoro-Struc, col sep=comma] {VoroResults.csv};
			\addplot+[blue,mark=triangle*,mark options={fill=blue}]   table [x=NODE-NUM-FREE, y=S11-L2-REL-DCPSE-NoVoro-Free, col sep=comma] {VoroResults.csv};
			\addplot+[yellow,mark=*,mark options={fill=yellow}]   table [x=NODE-NUM-FREE, y=S11-L2-REL-DCPSE-WithVoro-Free, col sep=comma] {VoroResults.csv};
			\logLogSlopeTriangle{0.9}{0.1}{0.20}{1.3}{black};
			\end{axis}
			\end{tikzpicture} &
			
			\begin{tikzpicture}
			\begin{axis}[height=6cm,width=7.5cm, xmode=log, ymode=log, ymin=0.00004,ymax=0.3, legend entries={Struct. No Voro.,Struct. With Voro.,Free No Voro.,Free With Voro.},legend style={ at={(0.5,-0.2)},anchor=south west,legend columns=1, cells={anchor=west},  font=\footnotesize, rounded corners=2pt,}, legend pos=north east,xlabel=Number of Nodes,ylabel=$L_{\infty}$ Error - $\sigma_{11}$]
			\addplot+[green,mark=square*,mark options={fill=green}]  table [x=NODE-NUM-STRUC, y=S11-LINF-DCPSE-NoVoro-Struc, col sep=comma] {VoroResults.csv};
			\addplot+[red,mark=diamond*,mark options={fill=red}]   table [x=NODE-NUM-STRUC, y=S11-LINF-DCPSE-WithVoro-Struc, col sep=comma] {VoroResults.csv};
			\addplot+[blue,mark=triangle*,mark options={fill=blue}]   table [x=NODE-NUM-FREE, y=S11-LINF-DCPSE-NoVoro-Free, col sep=comma] {VoroResults.csv};
			\addplot+[yellow,mark=*,mark options={fill=yellow}]   table [x=NODE-NUM-FREE, y=S11-LINF-DCPSE-WithVoro-Free, col sep=comma] {VoroResults.csv};
			\logLogSlopeTriangle{0.9}{0.1}{0.2}{0.8}{black};
			\end{axis}
			\end{tikzpicture} \\
			
			\begin{tikzpicture}
			\begin{axis}[height=6cm,width=7.5cm, xmode=log, ymode=log, ymin=0.00000003,ymax=0.007, legend entries={Struct. No Voro.,Struct. With Voro.,Free No Voro.,Free With Voro.},legend style={ at={(0.5,-0.2)},anchor=south west,legend columns=1, cells={anchor=west},  font=\footnotesize, rounded corners=2pt,}, legend pos=north east,xlabel=Number of Nodes,ylabel=$L_2$ Error - $\sigma_{12}$]
			\addplot+[green,mark=square*,mark options={fill=green}]  table [x=NODE-NUM-STRUC, y=S12-L2-REL-DCPSE-NoVoro-Struc, col sep=comma] {VoroResults.csv};
			\addplot+[red,mark=diamond*,mark options={fill=red}]   table [x=NODE-NUM-STRUC, y=S12-L2-REL-DCPSE-WithVoro-Struc, col sep=comma] {VoroResults.csv};
			\addplot+[blue,mark=triangle*,mark options={fill=blue}]   table [x=NODE-NUM-FREE, y=S12-L2-REL-DCPSE-NoVoro-Free, col sep=comma] {VoroResults.csv};
			\addplot+[yellow,mark=*,mark options={fill=yellow}]   table [x=NODE-NUM-FREE, y=S12-L2-REL-DCPSE-WithVoro-Free, col sep=comma] {VoroResults.csv};
			\logLogSlopeTriangle{0.9}{0.1}{0.12}{1.3}{black};
			\end{axis}
			\end{tikzpicture} &
			
			\begin{tikzpicture}
			\begin{axis}[height=6cm,width=7.5cm, xmode=log, ymode=log, ymin=0.00004,ymax=0.3, legend entries={Struct. No Voro.,Struct. With Voro.,Free No Voro.,Free With Voro.},legend style={ at={(0.5,-0.2)},anchor=south west,legend columns=1, cells={anchor=west},  font=\footnotesize, rounded corners=2pt,}, legend pos=north east,xlabel=Number of Nodes,ylabel=$L_{\infty}$ Error - $\sigma_{12}$]
			\addplot+[green,mark=square*,mark options={fill=green}]  table [x=NODE-NUM-STRUC, y=S12-LINF-DCPSE-NoVoro-Struc, col sep=comma] {VoroResults.csv};
			\addplot+[red,mark=diamond*,mark options={fill=red}]   table [x=NODE-NUM-STRUC, y=S12-LINF-DCPSE-WithVoro-Struc, col sep=comma] {VoroResults.csv};
			\addplot+[blue,mark=triangle*,mark options={fill=blue}]   table [x=NODE-NUM-FREE, y=S12-LINF-DCPSE-NoVoro-Free, col sep=comma] {VoroResults.csv};
			\addplot+[yellow,mark=*,mark options={fill=yellow}]   table [x=NODE-NUM-FREE, y=S12-LINF-DCPSE-WithVoro-Free, col sep=comma] {VoroResults.csv};
			\logLogSlopeTriangle{0.9}{0.1}{0.12}{0.8}{black};
			\end{axis}
			\end{tikzpicture} 
		\end{tabular}
		\caption{Impact of Voronoi Integration on the Errors for the DC PSE Method - 2D Cylinder. $L_2$ (Left) and $L_{\infty}$ (Right) errors for structured and free node distributions. A slight reduction in the error is only observed for the free node distribution.}
		\label{ResultsVoroDCPSE}
	\end{figure}
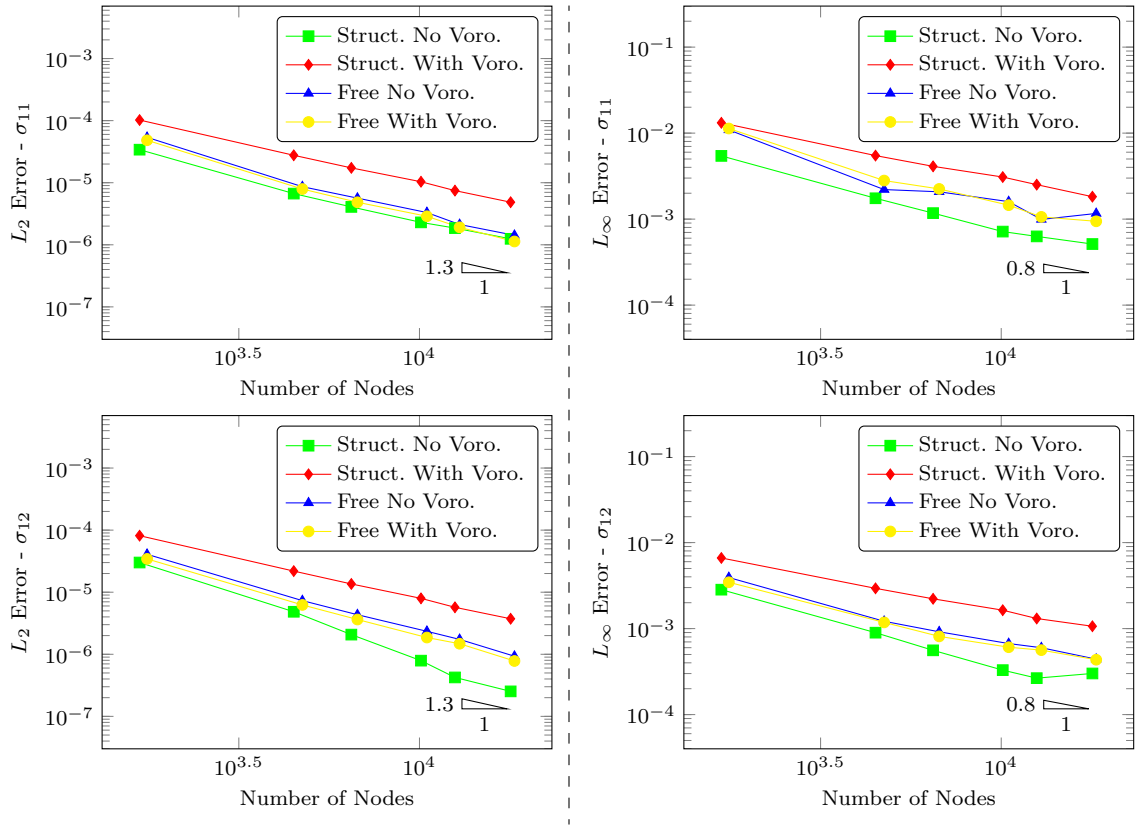
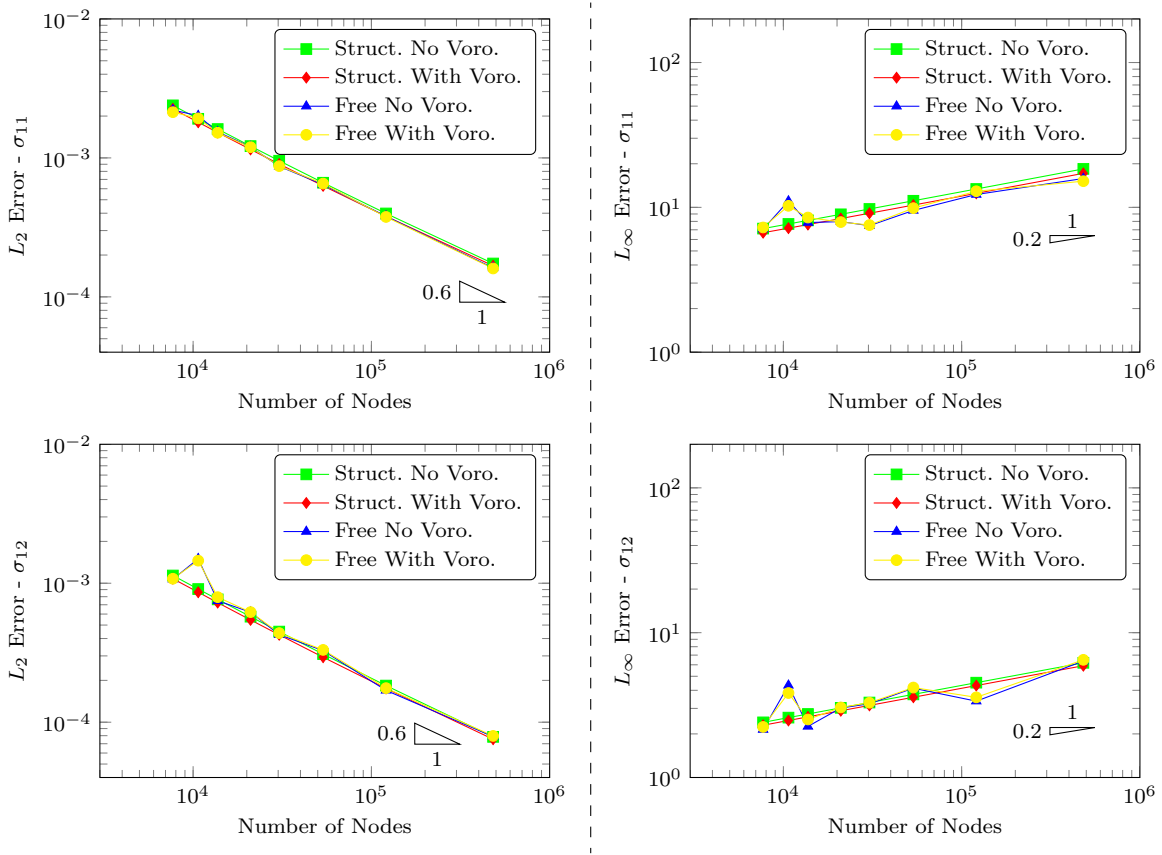
\begin{figure}[H] 
	\centering
	\begin{tabular}{c:c}
		\begin{tikzpicture}
		\begin{axis}[height=6cm,width=7.5cm, xmode=log,xmin=3000,xmax=1000000, ymode=log, ymin=0.00004,ymax=0.01, legend entries={Struct. No Voro.,Struct. With Voro.,Free No Voro.,Free With Voro.},legend style={ at={(0.5,-0.2)},anchor=south west,legend columns=1, cells={anchor=west},  font=\footnotesize, rounded corners=2pt,}, legend pos=north east,xlabel=Number of Nodes,ylabel=$L_2$ Error - $\sigma_{11}$]
		\addplot+[green,mark=square*,mark options={fill=green}]  table [x=NODE-NUM-STRUC, y=S11-L2-REL-DCPSE-NoVoro-Struc, col sep=comma] {VoroResults-LShape.csv};
		\addplot+[red,mark=diamond*,mark options={fill=red}]   table [x=NODE-NUM-STRUC, y=S11-L2-REL-DCPSE-WithVoro-Struc, col sep=comma] {VoroResults-LShape.csv};
		\addplot+[blue,mark=triangle*,mark options={fill=blue}]   table [x=NODE-NUM-FREE, y=S11-L2-REL-DCPSE-NoVoro-Free, col sep=comma] {VoroResults-LShape.csv};
		\addplot+[yellow,mark=*,mark options={fill=yellow}]   table [x=NODE-NUM-FREE, y=S11-L2-REL-DCPSE-WithVoro-Free, col sep=comma] {VoroResults-LShape.csv};
		\logLogSlopeTriangle{0.9}{0.1}{0.15}{0.6}{black};
		\end{axis}
		\end{tikzpicture} &
		
		\begin{tikzpicture}
		\begin{axis}[height=6cm,width=7.5cm, xmode=log,xmin=3000,xmax=1000000, ymode=log, ymin=1,ymax=200, legend entries={Struct. No Voro.,Struct. With Voro.,Free No Voro.,Free With Voro.},legend style={ at={(0.5,-0.2)},anchor=south west,legend columns=1, cells={anchor=west},  font=\footnotesize, rounded corners=2pt,}, legend pos=north east,xlabel=Number of Nodes,ylabel=$L_{\infty}$ Error - $\sigma_{11}$]
		\addplot+[green,mark=square*,mark options={fill=green}]  table [x=NODE-NUM-STRUC, y=S11-LINF-DCPSE-NoVoro-Struc, col sep=comma] {VoroResults-LShape.csv};
		\addplot+[red,mark=diamond*,mark options={fill=red}]   table [x=NODE-NUM-STRUC, y=S11-LINF-DCPSE-WithVoro-Struc, col sep=comma] {VoroResults-LShape.csv};
		\addplot+[blue,mark=triangle*,mark options={fill=blue}]   table [x=NODE-NUM-FREE, y=S11-LINF-DCPSE-NoVoro-Free, col sep=comma] {VoroResults-LShape.csv};
		\addplot+[yellow,mark=*,mark options={fill=yellow}]   table [x=NODE-NUM-FREE, y=S11-LINF-DCPSE-WithVoro-Free, col sep=comma] {VoroResults-LShape.csv};
		\logLogSlopeTriangleUp{0.9}{0.1}{0.35}{0.2}{black};
		\end{axis}
		\end{tikzpicture} \\
		
		\begin{tikzpicture}
		\begin{axis}[height=6cm,width=7.5cm, xmode=log,xmin=3000,xmax=1000000, ymode=log, ymin=0.00004,ymax=0.01, legend entries={Struct. No Voro.,Struct. With Voro.,Free No Voro.,Free With Voro.},legend style={ at={(0.5,-0.2)},anchor=south west,legend columns=1, cells={anchor=west},  font=\footnotesize, rounded corners=2pt,}, legend pos=north east,xlabel=Number of Nodes,ylabel=$L_2$ Error - $\sigma_{12}$]
		\addplot+[green,mark=square*,mark options={fill=green}]  table [x=NODE-NUM-STRUC, y=S12-L2-REL-DCPSE-NoVoro-Struc, col sep=comma] {VoroResults-LShape.csv};
		\addplot+[red,mark=diamond*,mark options={fill=red}]   table [x=NODE-NUM-STRUC, y=S12-L2-REL-DCPSE-WithVoro-Struc, col sep=comma] {VoroResults-LShape.csv};
		\addplot+[blue,mark=triangle*,mark options={fill=blue}]   table [x=NODE-NUM-FREE, y=S12-L2-REL-DCPSE-NoVoro-Free, col sep=comma] {VoroResults-LShape.csv};
		\addplot+[yellow,mark=*,mark options={fill=yellow}]   table [x=NODE-NUM-FREE, y=S12-L2-REL-DCPSE-WithVoro-Free, col sep=comma] {VoroResults-LShape.csv};
		\logLogSlopeTriangle{0.8}{0.1}{0.1}{0.6}{black};
		\end{axis}
		\end{tikzpicture} &
		
		\begin{tikzpicture}
		\begin{axis}[height=6cm,width=7.5cm, xmode=log,xmin=3000,xmax=1000000, ymode=log, ymin=1,ymax=200, legend entries={Struct. No Voro.,Struct. With Voro.,Free No Voro.,Free With Voro.},legend style={ at={(0.5,-0.2)},anchor=south west,legend columns=1, cells={anchor=west},  font=\footnotesize, rounded corners=2pt,}, legend pos=north east,xlabel=Number of Nodes,ylabel=$L_{\infty}$ Error - $\sigma_{12}$]
		\addplot+[green,mark=square*,mark options={fill=green}]  table [x=NODE-NUM-STRUC, y=S12-LINF-DCPSE-NoVoro-Struc, col sep=comma] {VoroResults-LShape.csv};
		\addplot+[red,mark=diamond*,mark options={fill=red}]   table [x=NODE-NUM-STRUC, y=S12-LINF-DCPSE-WithVoro-Struc, col sep=comma] {VoroResults-LShape.csv};
		\addplot+[blue,mark=triangle*,mark options={fill=blue}]   table [x=NODE-NUM-FREE, y=S12-LINF-DCPSE-NoVoro-Free, col sep=comma] {VoroResults-LShape.csv};
		\addplot+[yellow,mark=*,mark options={fill=yellow}]   table [x=NODE-NUM-FREE, y=S12-LINF-DCPSE-WithVoro-Free, col sep=comma] {VoroResults-LShape.csv};
		\logLogSlopeTriangleUp{0.9}{0.1}{0.15}{0.2}{black};
		\end{axis}
		\end{tikzpicture} 
	\end{tabular}
	\caption{Impact of Voronoi based Weights on the Errors for the DC PSE Method - 2D L-Shape. $L_2$ (Left) and $L_{\infty}$ (Right) errors for structured and free node distributions. A slight reduction in the error is observed for the structured node distribution.}
	\label{ResultsVoroDCPSE_LShape}
\end{figure}

	\paragraph{Discussion} \
	
	From Figure \ref{ResultsVoroGFD} to Figure \ref{ResultsVoroDCPSE_LShape} we can see that the error reduction achieved by the use of Voronoi based volumes is not guaranteed. The use of such volumes for the 2D cylinder problem lead to an error reduction of up to 17\% for the 2D cylinder problem. However, for the 2D L-shape, an error increase of 3\% has been observed for the free node distribution. From this study, we can conclude that the use of the Voronoi based weights shall be used with care as is may lead to a significant error.
	
	\subsection{Collocation Method Stabilization}
	
	\subsubsection{General}
	
	Within the framework of collocation, the boundary conditions are applied at the nodes using the strong form of the partial differential equations. This may lead to ill conditioning of the linear system of equations, as both the boundary conditions and the equilibrium equation cannot be enforced simultaneously at a boundary node. To overcome this issue, a stabilization method, known as the \textit{Finite Increment Calculus} (FIC), has been presented by E. O\~nate \cite{Oate1998} for structural problems that are solved with the Finite Point Method. The method is presented in this section.
	Results for the 2D cylinder and for the 2D L-shape problems are presented with and without stabilization. This stabilization approach is used for both the methods, GFD and DC PSE.
	
	\subsubsection{Stabilized Equations}
	
	Considering an unknown field $f$, a partial differential problem is defined by a differential operator $\mathcal{A}$ applied to the interior domain $\Omega$, a field $\overline{f}$ set to the boundary $\Gamma_u$, and a differential operator $\mathcal{B}$ applied to the boundary $\Gamma_t$ (see Figure \ref{FICDrawing}).
	
	\begin{figure}[H] 
		\centering
		\begin{tikzpicture}
		\def\svgwidth{8cm}
		\node at (0,0) {\includegraphics{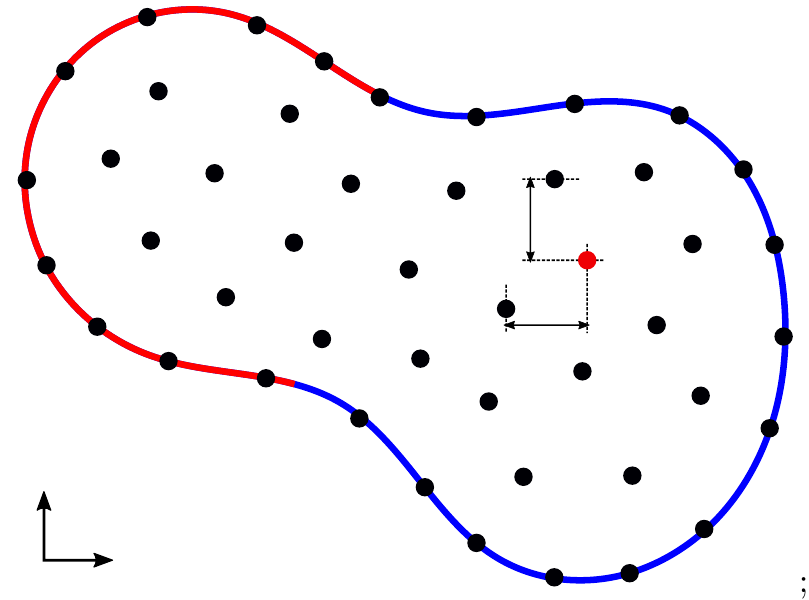}};
		\node[color=black] at (-0.3,0.5) [left] {$\Omega$};
		\node[color=red] at (-4,1) [left] {$\Gamma_u$};
		\node[color=blue] at (4.4,1) [left] {$\Gamma_t$};
		\node[color=red] at (2.6,0.6) [left] {$X_c$};
		\node[color=black] at (1.75,-0.6) [left] {$h_1$};
		\node[color=black] at (1.25,0.7) [left] {$h_2$};
		\node[color=black] at (-2.5,-2.5) [left] {$x_1$};
		\node[color=black] at (-3.1,-1.85) [left] {$x_2$};
		\end{tikzpicture}
		\caption{2D domain $\Omega$ on which Dirichlet boundary conditions are applied to the boundary $\Gamma_u$ and Neumann boundary conditions  to $\Gamma_t$. The characteristic lengths $h_1$ and $h_2$ are presented for the collocation node $\mathbf{X_c}$.}
		\label{FICDrawing}
	\end{figure}
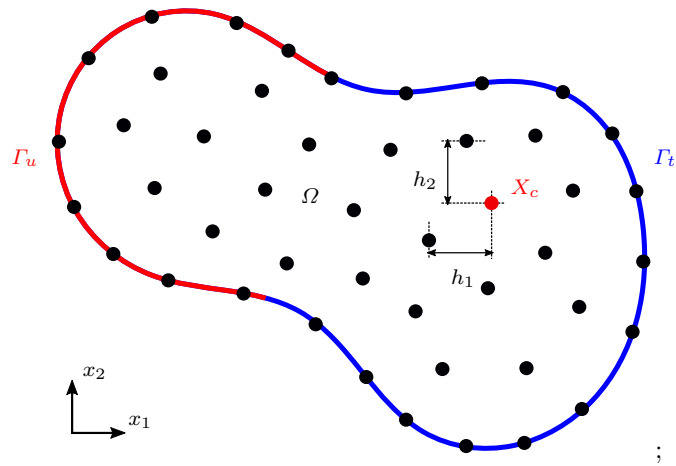
	\begin{equation} \label{ProblemEquations}
	\begin{aligned}
	\mathbf{\mathcal{A}(f)} &=0 \text{ \quad in \quad} \Omega, \\
	\mathbf{f-\overline{f}}&=0 \text{ \quad on \quad} \Gamma_u, \\
	\mathbf{\mathcal{B}(f)} &=0 \text{ \quad on \quad} \Gamma_t. \\
	\end{aligned}
	\end{equation}
	Based on \cite{Oate1998} and \cite{Oate2001}, the stabilized system of equations is:
	\begin{equation} \label{StabProblemEquations}
	\begin{aligned}
	\mathbf{\mathcal{A}(f)}- \frac{1}{2} \sum_{j=1}^m {h_j \frac{\mathbf{\partial \mathcal{A}(f)}}{\partial x_j}}&=0 \text{ \quad in \quad} \Omega, \\
	\mathbf{f-\overline{f}}&=0 \text{ \quad on \quad} \Gamma_u, \\
	\mathbf{\mathcal{B}(f)}- \sum_{j=1}^m {h_j n_j \mathbf{\mathcal{A}(f)}} &=0 \text{ \quad on \quad} \Gamma_t,
	\end{aligned}
	\end{equation}
	where $m$ is the dimension of the domain, $h_j$ is the characteristic length of the domain in the direction $j$, and $n$ is the unit normal. 
	
	For isotropic support weight functions, $h_j$ reduces to $h$ and can be expressed as follows:
	\begin{equation} \label{CharctLenCalc}
	\begin{aligned}
	h=&R_{\text{Sup}}\left(\frac{\pi}{N_{\text{Sup}}}\right)^{\frac{1}{2}} & \text{ \quad for 2D problems, \quad}  \\
	h=&R_{\text{Sup}}\left(\frac{4\pi}{3N_{\text{Sup}}}\right)^{\frac{1}{3}} & \text{ \quad for 3D problems. \quad}, \\
	\end{aligned}
	\end{equation}
	where $R_{\text{Sup}}$ and $N_{\text{Sup}}$, respectively, represent the radius of the node support, and the number of nodes in the support. 
	
	\subsubsection{Results}
	
	Equation (\ref{StabProblemEquations}) has only been applied to the boundary nodes where the maximum error is usually observed. Also, the stabilized equation on the boundary does not require the approximation of an additional derivative order. The results obtained are presented from Figure \ref{ResultsStab-Cyl-GFD} to Figure \ref{ResultsStab-LShape-DCPSE}. $L_{\infty}$ error results have not been presented in this section as this error highly depends on the proximity of the closest node to the singularity.
	\begin{figure}[H] 
		\centering
		\begin{tabular}{c:c}
			\begin{tikzpicture}
			\begin{axis}[height=6cm,width=7.5cm, ymode=log, ymin=0.00000002,ymax=0.01,xmin=1000, xmode=log, legend entries={No Stabilization,With Stabilization},legend style={ at={(0.5,-0.2)},anchor=south west,legend columns=1, cells={anchor=west},  font=\footnotesize, rounded corners=2pt,}, legend pos=north east,xlabel=Number of Nodes,ylabel=$L_2$ Error - $\sigma_{11}$]
			\addplot+[green,mark=square*,mark options={fill=green}]  table [x=NODE-NUM, y=S11-L2-REL-GFD-NoStab, col sep=comma] {StabResults-Cyl.csv};
			\addplot+[red,mark=diamond*,mark options={fill=red}]   table [x=NODE-NUM, y=S11-L2-REL-GFD-WithStab, col sep=comma] {StabResults-Cyl.csv};
			\logLogSlopeTriangle{0.9}{0.1}{0.23}{1.5}{green};
			\logLogSlopeTriangle{0.9}{0.1}{0.5}{1.6}{red};
			\end{axis}
			\end{tikzpicture} &
			
			\begin{tikzpicture}
			\begin{axis}[height=6cm,width=7.5cm, ymode=log, ymin=0.00000002,ymax=0.01,xmin=1000, xmode=log, legend entries={No Stabilization,With Stabilization},legend style={ at={(0.5,-0.2)},anchor=south west,legend columns=1, cells={anchor=west},  font=\footnotesize, rounded corners=2pt,}, legend pos=north east,xlabel=Number of Nodes,ylabel=$L_2$ Error - $\sigma_{12}$]
			\addplot+[green,mark=square*,mark options={fill=green}]  table [x=NODE-NUM, y=S12-L2-REL-GFD-NoStab, col sep=comma] {StabResults-Cyl.csv};
			\addplot+[red,mark=diamond*,mark options={fill=red}]   table [x=NODE-NUM, y=S12-L2-REL-GFD-WithStab, col sep=comma] {StabResults-Cyl.csv};
			\logLogSlopeTriangle{0.9}{0.1}{0.25}{1.6}{green};
			\logLogSlopeTriangle{0.9}{0.1}{0.45}{1.5}{red};
			\end{axis}
			\end{tikzpicture} \\
		\end{tabular}
		\caption{Stabilization results comparison - 2D Cylinder - GFD Method. $L_2$ error for stabilized and non-stabilized PDE for increasing node numbers. A lower error is observed for the non-stabilized PDE.}
		\label{ResultsStab-Cyl-GFD}
	\end{figure}
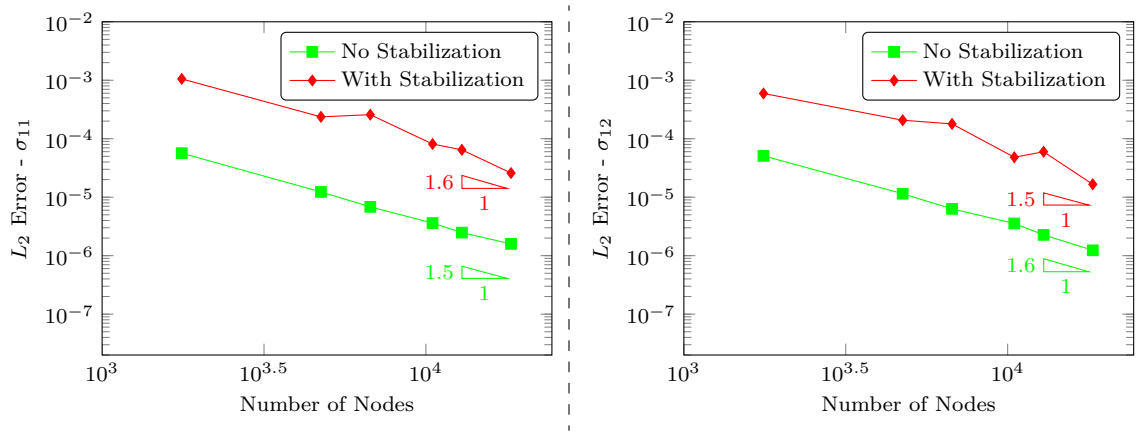
	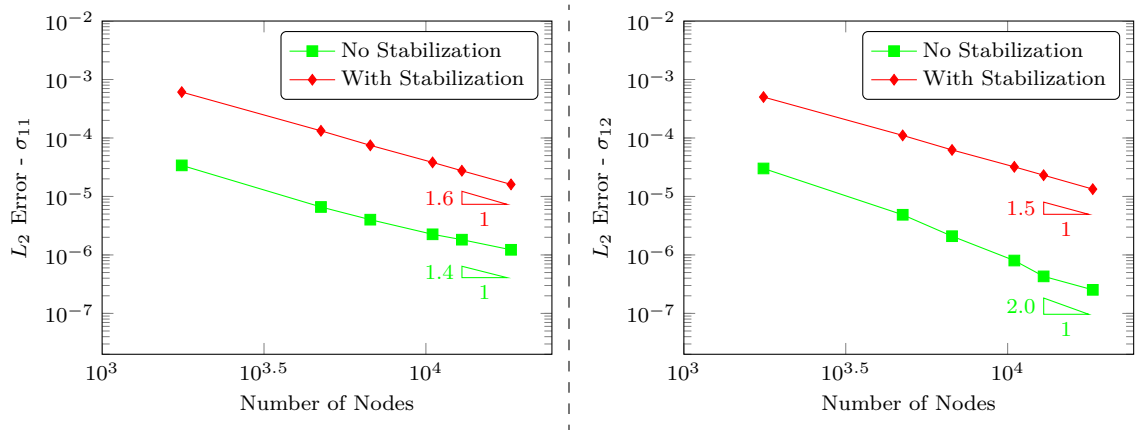
\begin{figure}[H] 
		\centering
		\begin{tabular}{c:c}
			\begin{tikzpicture}
			\begin{axis}[height=6cm,width=7.5cm, ymode=log, ymin=0.00000002,ymax=0.01,xmin=1000, xmode=log, legend entries={No Stabilization,With Stabilization},legend style={ at={(0.5,-0.2)},anchor=south west,legend columns=1, cells={anchor=west},  font=\footnotesize, rounded corners=2pt,}, legend pos=north east,xlabel=Number of Nodes,ylabel=$L_2$ Error - $\sigma_{11}$]
			\addplot+[green,mark=square*,mark options={fill=green}]  table [x=NODE-NUM, y=S11-L2-REL-DCPSE-NoStab, col sep=comma] {StabResults-Cyl.csv};
			\addplot+[red,mark=diamond*,mark options={fill=red}]   table [x=NODE-NUM, y=S11-L2-REL-DCPSE-WithStab, col sep=comma] {StabResults-Cyl.csv};
			\logLogSlopeTriangle{0.9}{0.1}{0.23}{1.4}{green};
			\logLogSlopeTriangle{0.9}{0.1}{0.45}{1.6}{red};
			\end{axis}
			\end{tikzpicture} &
			
			\begin{tikzpicture}
			\begin{axis}[height=6cm,width=7.5cm, ymode=log, ymin=0.00000002,ymax=0.01,xmin=1000, xmode=log, legend entries={No Stabilization,With Stabilization},legend style={ at={(0.5,-0.2)},anchor=south west,legend columns=1, cells={anchor=west},  font=\footnotesize, rounded corners=2pt,}, legend pos=north east,xlabel=Number of Nodes,ylabel=$L_2$ Error - $\sigma_{12}$]
			\addplot+[green,mark=square*,mark options={fill=green}]  table [x=NODE-NUM, y=S12-L2-REL-DCPSE-NoStab, col sep=comma] {StabResults-Cyl.csv};
			\addplot+[red,mark=diamond*,mark options={fill=red}]   table [x=NODE-NUM, y=S12-L2-REL-DCPSE-WithStab, col sep=comma] {StabResults-Cyl.csv};
			\logLogSlopeTriangle{0.9}{0.1}{0.12}{2.0}{green};
			\logLogSlopeTriangle{0.9}{0.1}{0.42}{1.5}{red};
			\end{axis}
			\end{tikzpicture} \\
		\end{tabular}
		\caption{Stabilization results comparison - 2D Cylinder - DC PSE Method. $L_2$ error for stabilized and non-stabilized PDE for increasing node numbers. A lower error is observed for the non-stabilized PDE.}
		\label{ResultsStab-Cyl-DCPSE}
	\end{figure}
	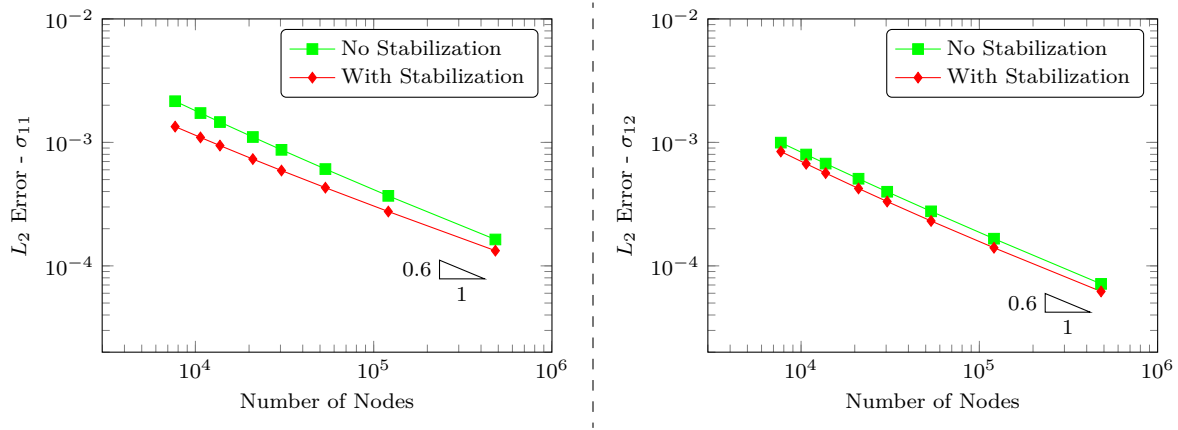
\begin{figure}[H] 
		\centering
		\begin{tabular}{c:c}
			\begin{tikzpicture}
			\begin{axis}[height=6cm,width=7.5cm, ymode=log, ymin=0.00002,ymax=0.01,xmin=3000,xmax=1000000, xmode=log, legend entries={No Stabilization,With Stabilization},legend style={ at={(0.5,-0.2)},anchor=south west,legend columns=1, cells={anchor=west},  font=\footnotesize, rounded corners=2pt,}, legend pos=north east,xlabel=Number of Nodes,ylabel=$L_2$ Error - $\sigma_{11}$]
			\addplot+[green,mark=square*,mark options={fill=green}]  table [x=NODE-NUM, y=S11-L2-REL-GFD-NoStab, col sep=comma] {StabResults-LShape.csv};
			\addplot+[red,mark=diamond*,mark options={fill=red}]   table [x=NODE-NUM, y=S11-L2-REL-GFD-WithStab, col sep=comma] {StabResults-LShape.csv};
			\logLogSlopeTriangle{0.85}{0.1}{0.22}{0.6}{black};
			\end{axis}
			\end{tikzpicture} &
			
			\begin{tikzpicture}
			\begin{axis}[height=6cm,width=7.5cm, ymode=log, ymin=0.00002,ymax=0.01,xmin=3000,xmax=1000000, xmode=log, legend entries={No Stabilization,With Stabilization},legend style={ at={(0.5,-0.2)},anchor=south west,legend columns=1, cells={anchor=west},  font=\footnotesize, rounded corners=2pt,}, legend pos=north east,xlabel=Number of Nodes,ylabel=$L_2$ Error - $\sigma_{12}$]
			\addplot+[green,mark=square*,mark options={fill=green}]  table [x=NODE-NUM, y=S12-L2-REL-GFD-NoStab, col sep=comma] {StabResults-LShape.csv};
			\addplot+[red,mark=diamond*,mark options={fill=red}]   table [x=NODE-NUM, y=S12-L2-REL-GFD-WithStab, col sep=comma] {StabResults-LShape.csv};
			\logLogSlopeTriangle{0.85}{0.1}{0.12}{0.6}{black};
			\end{axis}
			\end{tikzpicture} \\
		\end{tabular}
		\caption{Stabilization Results Comparison - 2D L-Shape - GFD Method. $L_2$ error for stabilized and non-stabilized PDE for increasing node numbers. A lower error is observed for the stabilized PDE.}
		\label{ResultsStab-LShape-GFD}
	\end{figure}
	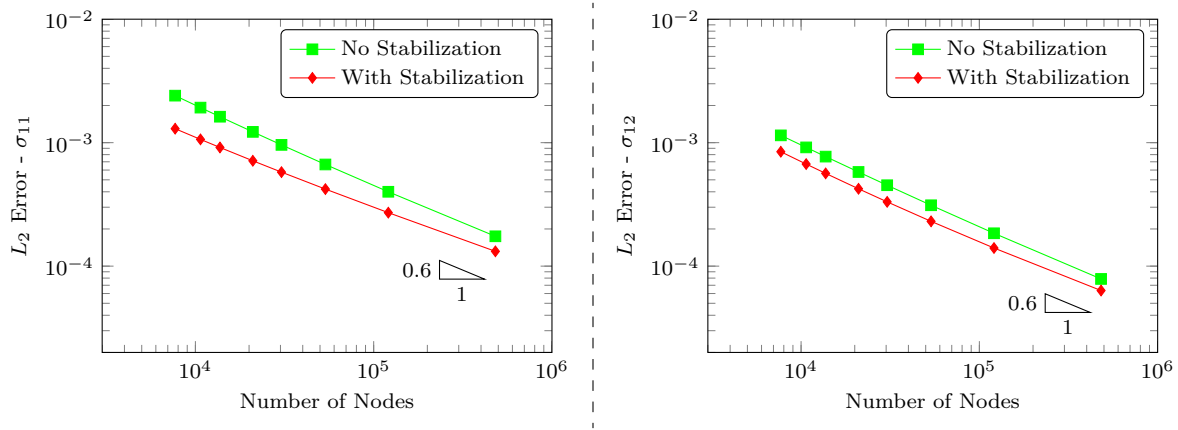
\begin{figure}[H] 
		\centering
		\begin{tabular}{c:c}
			\begin{tikzpicture}
			\begin{axis}[height=6cm,width=7.5cm, ymode=log, ymin=0.00002,ymax=0.01,xmin=3000,xmax=1000000, xmode=log, legend entries={No Stabilization,With Stabilization},legend style={ at={(0.5,-0.2)},anchor=south west,legend columns=1, cells={anchor=west},  font=\footnotesize, rounded corners=2pt,}, legend pos=north east,xlabel=Number of Nodes,ylabel=$L_2$ Error - $\sigma_{11}$]
			\addplot+[green,mark=square*,mark options={fill=green}]  table [x=NODE-NUM, y=S11-L2-REL-DCPSE-NoStab, col sep=comma] {StabResults-LShape.csv};
			\addplot+[red,mark=diamond*,mark options={fill=red}]   table [x=NODE-NUM, y=S11-L2-REL-DCPSE-WithStab, col sep=comma] {StabResults-LShape.csv};
			\logLogSlopeTriangle{0.85}{0.1}{0.22}{0.6}{black};
			\end{axis}
			\end{tikzpicture} &
			
			\begin{tikzpicture}
			\begin{axis}[height=6cm,width=7.5cm, ymode=log, ymin=0.00002,ymax=0.01,xmin=3000,xmax=1000000, xmode=log, legend entries={No Stabilization,With Stabilization},legend style={ at={(0.5,-0.2)},anchor=south west,legend columns=1, cells={anchor=west},  font=\footnotesize, rounded corners=2pt,}, legend pos=north east,xlabel=Number of Nodes,ylabel=$L_2$ Error - $\sigma_{12}$]
			\addplot+[green,mark=square*,mark options={fill=green}]  table [x=NODE-NUM, y=S12-L2-REL-DCPSE-NoStab, col sep=comma] {StabResults-LShape.csv};
			\addplot+[red,mark=diamond*,mark options={fill=red}]   table [x=NODE-NUM, y=S12-L2-REL-DCPSE-WithStab, col sep=comma] {StabResults-LShape.csv};
			\logLogSlopeTriangle{0.85}{0.1}{0.12}{0.6}{black};
			\end{axis}
			\end{tikzpicture} \\
		\end{tabular}
		\caption{Stabilization results comparison - 2D L-Shape - DC PSE Method. $L_2$ error for stabilized and non-stabilized PDE for increasing node numbers. A lower error is observed for the stabilized PDE.}
		\label{ResultsStab-LShape-DCPSE}
	\end{figure}
	
	It can be observed from the results presented from Figure \ref{ResultsStab-Cyl-GFD} and Figure \ref{ResultsStab-Cyl-DCPSE} that the stabilization equations lead to an error increase for the 2D cylinder problem. Using the stabilization method increases by a factor 30 the error for the GFD method and by a factor 20 for the DC PSE method. We can see from Figure \ref{ResultsStab-LShape-GFD} and Figure \ref{ResultsStab-LShape-DCPSE} that the error is reduced by the use of the stabilization method for the L-shape problem. An average error reduction of 25\% is observed for the GFD method and of 35\% for the DC PSE method.
	
	For the 2D cylinder, the loading is applied via Neumann boundary conditions, which represent the pressure loading. The L-shape problem on the other hand is loaded using Dirichlet boundary conditions. It can be concluded from this study that stabilization of the Neumann boundary conditions does not necessarily lead to a reduction of the observed error. Thus, stabilization of Neumann loaded problems does not seem to be an effective solution for the considered methods.
	
	\subsection{Support Node Selection for Singular Problems}\label{SupNodeSelection}
	
	\subsubsection{General} \label{SupNodeSelectionGle}
	
	For singular problems, such as the L-shape problem presented in Section \ref{RefProblems}, the selection of the support nodes in the vicinity of the singularity impacts the solution.
	
\begin{figure}[H] 
		\centering
		\def\svgwidth{7cm}
		\includegraphics{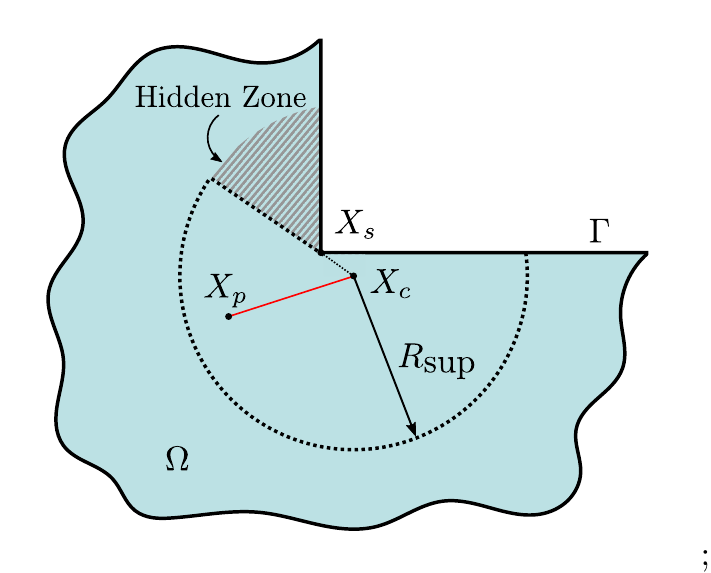}
		\caption{Visibility criterion. Only the nodes that can be connected to the collocation node by a segment which does not intersect the boundary of the domain are included in the support.}
		\label{VisibilityDrawing}
	\end{figure}

	In 1994, with the Element-Free Galerkin (EFG) method, Belytschko et al. \cite{Belytschko1994} introduced the visibility criterion for support nodes selection. This criterion has been widely used in the context of EFG fracture mechanics, see e.g., Duflot \cite{Duflot2004}. The support of a collocation node $\mathbf{X_c}$ in the domain $\Omega$ is selected so that any point $\mathbf{X_p}$ of the support, defined by a radius $R_{Sup}$, can be connected to the collocation node by a segment which does not intersect the domain boundary $\Gamma$ (see Figure \ref{VisibilityDrawing}). The ``hidden" zone is the zone within the support radius for which the segment between the collocation node and the support node intersects the boundary.
	
	The diffraction criterion, introduced by Organ in 1996 \cite{Organ1996}, is based on the same principle as the visibility criterion. The nodes in the zone ``hidden" from the collocation node are only included in the support if the sum of the length between the support node $\mathbf{X_p}$ and the singularity $\mathbf{X_s}$ and the length between the singularity and the collocation node $\mathbf{X_c}$ is smaller than the support radius $R_{Sup}$. The weights associated to the nodes in this zone are based on this increased distance to the collocation node.
	\begin{figure}[H] 
		\centering
		\def\svgwidth{7cm}
		\includegraphics{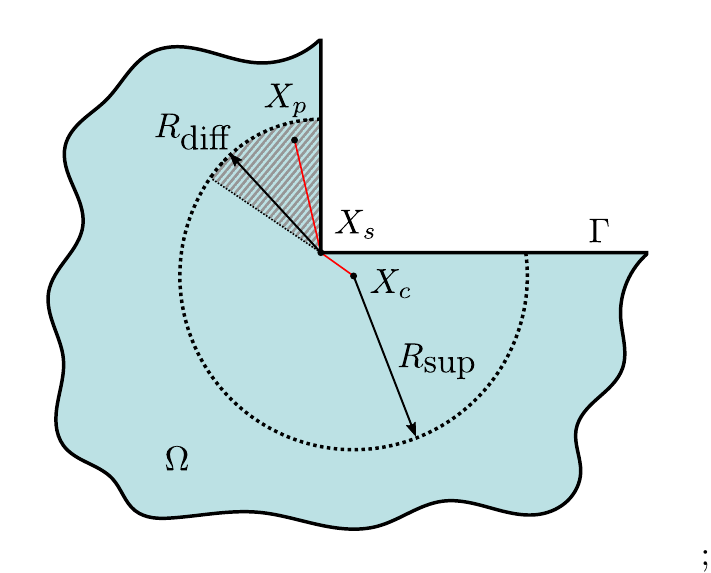}
		\caption{Diffraction Criterion. The nodes in the ``hidden" zone according to the visibility criterion are included in the support only if the sum of the distances [$\mathbf{X_p}$;$\mathbf{X_s}$] and [$\mathbf{X_s}$;$\mathbf{X_c}$] is smaller than the support radius.}
		\label{DiffractionDrawing}
	\end{figure}
	
	\subsubsection{Results}
	
	In this section, we assess the impact of the support node selection criterion on the $L_2$ error for the L-shape problem. Results obtained with the visibility criterion and with the diffraction criterion are compared to results where no criterion is considered (i.e. every node within the support radius of a collocation node is included in the support).
	
	Both the GFD and DC PSE methods are considered in this section.
	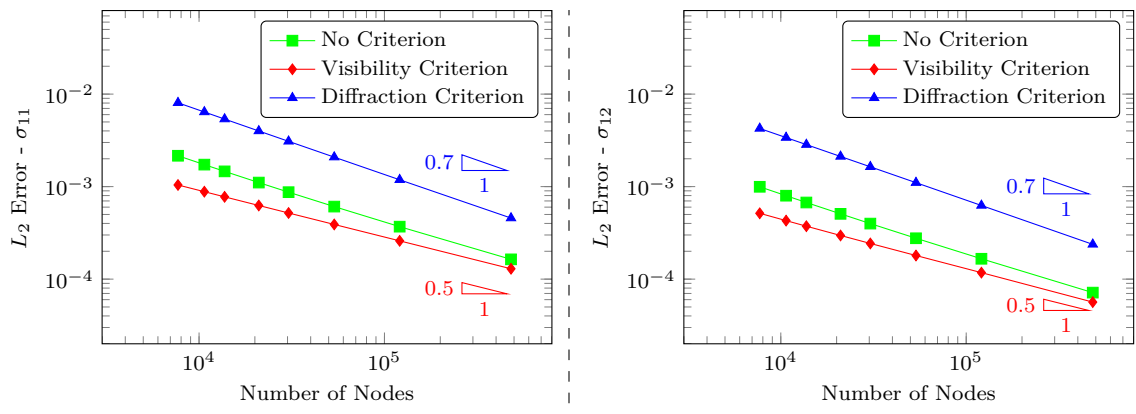
\begin{figure}[H] 
		\centering
		\begin{tabular}{c:c}
			\begin{tikzpicture}
			\begin{axis}[height=6cm,width=7.5cm, ymode=log, ymin=0.00002,ymax=0.08,xmin=3000, xmode=log, legend entries={No Criterion,Visibility Criterion,Diffraction Criterion},legend style={ at={(0.5,-0.2)},anchor=south west,legend columns=1, cells={anchor=west},  font=\footnotesize, rounded corners=2pt,}, legend pos=north east,xlabel=Number of Nodes,ylabel=$L_2$ Error - $\sigma_{11}$]
			\addplot+[green,mark=square*,mark options={fill=green}]  table [x=NODE-NUM, y=S11-L2-REL-GFD-NoVis, col sep=comma] {VisResults-LShape.csv};
			\addplot+[red,mark=diamond*,mark options={fill=red}]   table [x=NODE-NUM, y=S11-L2-REL-GFD-Vis, col sep=comma] {VisResults-LShape.csv};
			\addplot+[blue,mark=triangle*,mark options={fill=blue}]   table [x=NODE-NUM, y=S11-L2-REL-GFD-Diff, col sep=comma] {VisResults-LShape.csv};
			\logLogSlopeTriangle{0.9}{0.1}{0.15}{0.5}{red};
			\logLogSlopeTriangle{0.9}{0.1}{0.52}{0.7}{blue};
			\end{axis}
			\end{tikzpicture} &
			
			\begin{tikzpicture}
			\begin{axis}[height=6cm,width=7.5cm, ymode=log, ymin=0.00002,ymax=0.08,xmin=3000, xmode=log, legend entries={No Criterion,Visibility Criterion,Diffraction Criterion},legend style={ at={(0.5,-0.2)},anchor=south west,legend columns=1, cells={anchor=west},  font=\footnotesize, rounded corners=2pt,}, legend pos=north east,xlabel=Number of Nodes,ylabel=$L_2$ Error - $\sigma_{12}$]
			\addplot+[green,mark=square*,mark options={fill=green}]  table [x=NODE-NUM, y=S12-L2-REL-GFD-NoVis, col sep=comma] {VisResults-LShape.csv};
			\addplot+[red,mark=diamond*,mark options={fill=red}]   table [x=NODE-NUM, y=S12-L2-REL-GFD-Vis, col sep=comma] {VisResults-LShape.csv};
			\addplot+[blue,mark=triangle*,mark options={fill=blue}]   table [x=NODE-NUM, y=S12-L2-REL-GFD-Diff, col sep=comma] {VisResults-LShape.csv};
			\logLogSlopeTriangle{0.9}{0.1}{0.10}{0.5}{red};
			\logLogSlopeTriangle{0.9}{0.1}{0.45}{0.7}{blue};
			\end{axis}
			\end{tikzpicture} \\
		\end{tabular}
	\caption{Support node selection results comparison - 2D L-Shape - GFD Method. $L_2$ error obtained with no node selection criterion, with the visibility criterion, and the diffraction criterion. The lowest error is observed for the visibility criterion.}
		\label{ResultsVis-LShape-GFD}
	\end{figure}
	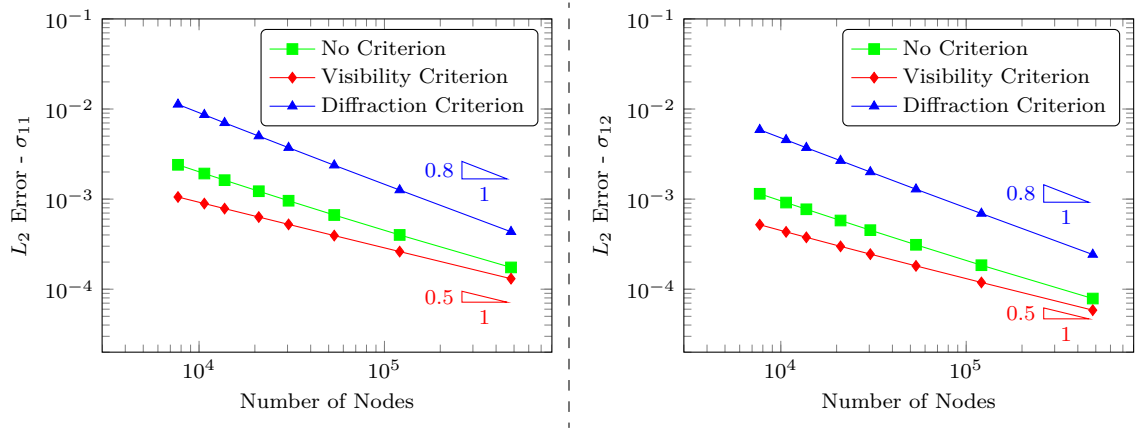
\begin{figure}[H] 
		\centering
		\begin{tabular}{c:c}
			\begin{tikzpicture}
			\begin{axis}[height=6cm,width=7.5cm, ymode=log, ymin=0.00002,ymax=0.1,xmin=3000, xmode=log, legend entries={No Criterion,Visibility Criterion,Diffraction Criterion},legend style={ at={(0.5,-0.2)},anchor=south west,legend columns=1, cells={anchor=west},  font=\footnotesize, rounded corners=2pt,}, legend pos=north east,xlabel=Number of Nodes,ylabel=$L_2$ Error - $\sigma_{11}$]
			\addplot+[green,mark=square*,mark options={fill=green}]  table [x=NODE-NUM, y=S11-L2-REL-DCPSE-NoVis, col sep=comma] {VisResults-LShape.csv};
			\addplot+[red,mark=diamond*,mark options={fill=red}]   table [x=NODE-NUM, y=S11-L2-REL-DCPSE-Vis, col sep=comma] {VisResults-LShape.csv};
			\addplot+[blue,mark=triangle*,mark options={fill=blue}]   table [x=NODE-NUM, y=S11-L2-REL-DCPSE-Diff, col sep=comma] {VisResults-LShape.csv};
			\logLogSlopeTriangle{0.9}{0.1}{0.15}{0.5}{red};
			\logLogSlopeTriangle{0.9}{0.1}{0.52}{0.8}{blue};
			\end{axis}
			\end{tikzpicture} &
			
			\begin{tikzpicture}
			\begin{axis}[height=6cm,width=7.5cm, ymode=log, ymin=0.00002,ymax=0.1,xmin=3000, xmode=log, legend entries={No Criterion,Visibility Criterion,Diffraction Criterion},legend style={ at={(0.5,-0.2)},anchor=south west,legend columns=1, cells={anchor=west},  font=\footnotesize, rounded corners=2pt,}, legend pos=north east,xlabel=Number of Nodes,ylabel=$L_2$ Error - $\sigma_{12}$]
			\addplot+[green,mark=square*,mark options={fill=green}]  table [x=NODE-NUM, y=S12-L2-REL-DCPSE-NoVis, col sep=comma] {VisResults-LShape.csv};
			\addplot+[red,mark=diamond*,mark options={fill=red}]   table [x=NODE-NUM, y=S12-L2-REL-DCPSE-Vis, col sep=comma] {VisResults-LShape.csv};
			\addplot+[blue,mark=triangle*,mark options={fill=blue}]   table [x=NODE-NUM, y=S12-L2-REL-DCPSE-Diff, col sep=comma] {VisResults-LShape.csv};
			\logLogSlopeTriangle{0.9}{0.1}{0.1}{0.5}{red};
			\logLogSlopeTriangle{0.9}{0.1}{0.45}{0.8}{blue};
			\end{axis}
			\end{tikzpicture} \\
			
		\end{tabular}
	\caption{Support Node Selection Results Comparison - 2D L-Shape - DC PSE Method. $L_2$ error obtained with no node selection criterion, with the visibility criterion, and the diffraction criterion. The lowest error is observed for the visibility criterion.}
		\label{ResultsVis-LShape-DCPSE}
	\end{figure}	
	It can be observed that, for both methods, the use of the visibility criterion leads to a significant error reduction compared to results where no criterion is applied. For the GFD method, the error reduction ranges from 50\% to 20\% depending on the node density. For the DC PSE method, the error reduction ranges from 55\% to 25\%. These results are expected as the singularity of the domain is more accurately represented with the visibility criterion. The diffraction criterion leads to an error increase ranging from a factor 2 to factor 10 for the GFD and the DC PSE methods. For both methods, the convergence rate is larger with the diffraction criterion. The visibility criterion appears to be a sensible choice for convex and singular problems.	
	
	\section{Benchmarking}\label{Benchmarking}
	
	In Section \ref{ParametricAnalysis}, we presented the studies of the parameters influencing the GFD and DC PSE methods and have selected a set of optimum parameters presented in Table \ref{ParametersSummary}. Based on these parameters, we compare in this section first the three variations of the DC PSE method presented in Section \ref{DCPEVariationsSec}, and then the methods presented in Section \ref{MethodDescription}. The results are assessed in terms of the $L_2$ and $L_{\infty}$ error norms for the 2D cylinder and the 2D L-shape problems. None of the improvement methods presented in Section \ref{ImprovementMethods_Section} have been used to derive the results presented in this section.
	
	\subsection{DC PSE Variations Comparison}\label{DCPSEVariationComp}
	
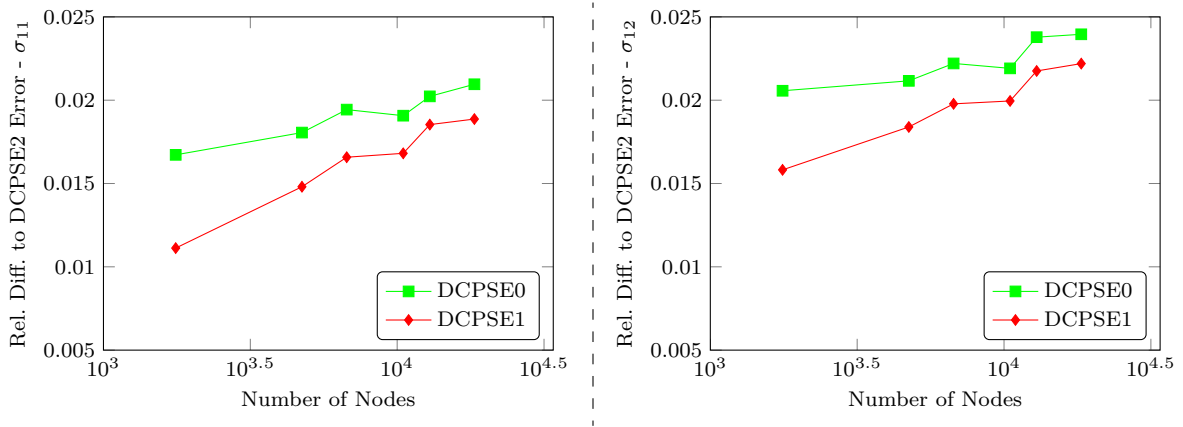
\begin{figure}[H] 
		\centering
		\begin{tabular}{c:c}
			\begin{tikzpicture}
			\begin{axis}[height=6cm,width=7.5cm, ymin=0.005,ymax=0.025,xmin=1000,xmax=34000, xmode=log, legend entries={DCPSE0,DCPSE1},legend style={ at={(0.5,-0.2)},anchor=south west,legend columns=1, cells={anchor=west},  font=\footnotesize, rounded corners=2pt,}, legend pos=south east,xlabel=Number of Nodes,ylabel=Rel. Diff. to DCPSE2 Error - $\sigma_{11}$, ytick={0,0.005,0.010,0.015,0.020,0.025}, scaled y ticks=false,yticklabel style={/pgf/number format/fixed, /pgf/number format/precision=5}]
			\addplot+[green,mark=square*,mark options={fill=green}]  table [x=NodeNum, y=L2-REL-S11-DCPSE0-FRAC, col sep=comma] {DCPSECompCylinder.csv};
			\addplot+[red,mark=diamond*,mark options={fill=red}]   table [x=NodeNum, y=L2-REL-S11-DCPSE1-FRAC, col sep=comma] {DCPSECompCylinder.csv};
			\end{axis}
			\end{tikzpicture} &
			
			\begin{tikzpicture}
			\begin{axis}[height=6cm,width=7.5cm, ymin=0.005,ymax=0.025,xmin=1000,xmax=34000, xmode=log, legend entries={DCPSE0,DCPSE1},legend style={ at={(0.5,-0.2)},anchor=south west,legend columns=1, cells={anchor=west},  font=\footnotesize, rounded corners=2pt,}, legend pos=south east,xlabel=Number of Nodes,ylabel=Rel. Diff. to DCPSE2 Error - $\sigma_{12}$, ytick={0,0.005,0.010,0.015,0.020,0.025}, scaled y ticks=false,yticklabel style={/pgf/number format/fixed, /pgf/number format/precision=5}]
			\addplot+[green,mark=square*,mark options={fill=green}]  table [x=NodeNum, y=L2-REL-S12-DCPSE0-FRAC, col sep=comma] {DCPSECompCylinder.csv};
			\addplot+[red,mark=diamond*,mark options={fill=red}]   table [x=NodeNum, y=L2-REL-S12-DCPSE1-FRAC, col sep=comma] {DCPSECompCylinder.csv};
			\end{axis}
			\end{tikzpicture} 
		\end{tabular}
		\caption{DC PSE method variations comparison - 2D Cylinder. Relative difference to the DCPSE2 $L_2$ error for the DCPSE0 and DCPSE1 methods. The DCPSE2 methods leads to the lowest error followed by the DCPSE1 method.}
		\label{DCPSE_Comparison2}
	\end{figure}

\begin{figure}[H] 
		\centering
		\begin{tabular}{c:c}
			\begin{tikzpicture}
			\begin{axis}[height=6cm,width=7.5cm, ymin=-0.007,ymax=-0.0025,ytick={-0.007,-0.005,-0.003},xmin=3000,xmax=1000000, xmode=log, legend entries={DCPSE0,DCPSE1},legend style={ at={(0.5,-0.2)},anchor=south west,legend columns=1, cells={anchor=west},  font=\footnotesize, rounded corners=2pt,}, legend pos=south east,xlabel=Number of Nodes,ylabel=Rel. Diff. to DCPSE2 Error - $\sigma_{11}$, scaled y ticks=false,yticklabel style={/pgf/number format/fixed, /pgf/number format/precision=5}]
			\addplot+[green,mark=square*,mark options={fill=green}]  table [x=NodeNum, y=L2-REL-S11-DCPSE0-FRAC, col sep=comma] {DCPSECompLShape.csv};
			\addplot+[red,mark=diamond*,mark options={fill=red}]   table [x=NodeNum, y=L2-REL-S11-DCPSE1-FRAC, col sep=comma] {DCPSECompLShape.csv};
			\end{axis}
			\end{tikzpicture} &
			
			\begin{tikzpicture}
			\begin{axis}[height=6cm,width=7.5cm, ymin=-0.015,ymax=0.01,xmin=3000,xmax=1000000, xmode=log, legend entries={DCPSE0,DCPSE1},legend style={ at={(0.5,-0.2)},anchor=south west,legend columns=1, cells={anchor=west},  font=\footnotesize, rounded corners=2pt,}, legend pos=south east,xlabel=Number of Nodes,ylabel=Rel. Diff. to DCPSE2 Error - $\sigma_{12}$, scaled y ticks=false,yticklabel style={/pgf/number format/fixed, /pgf/number format/precision=5}]
			\addplot+[green,mark=square*,mark options={fill=green}]  table [x=NodeNum, y=L2-REL-S12-DCPSE0-FRAC, col sep=comma] {DCPSECompLShape.csv};
			\addplot+[red,mark=diamond*,mark options={fill=red}]   table [x=NodeNum, y=L2-REL-S12-DCPSE1-FRAC, col sep=comma] {DCPSECompLShape.csv};
			\end{axis}
			\end{tikzpicture} 
		\end{tabular}
		\caption{DC PSE method variations comparison - 2D L-Shape. Relative difference to the DCPSE2 $L_2$ error for the DCPSE0 and DCPSE1 methods. The DCPSE1 methods leads to the lowest error followed by the DCPSE0 method.}
		\label{DCPSE_Comparison2_LShaped}
	\end{figure}
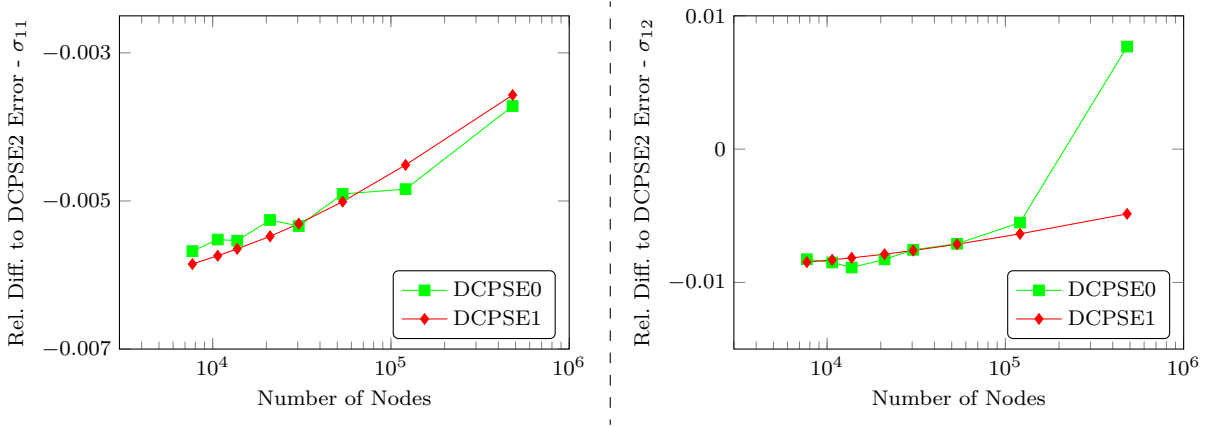

	The three variations of the DC PSE method are studied in this section. The results in terms of $L_2$ error are compared in Figure \ref{DCPSE_Comparison} and in Figure \ref{DCPSE_Comparison_L-Shape}, respectively, for the 2D cylinder and the 2D L-shape problems.
	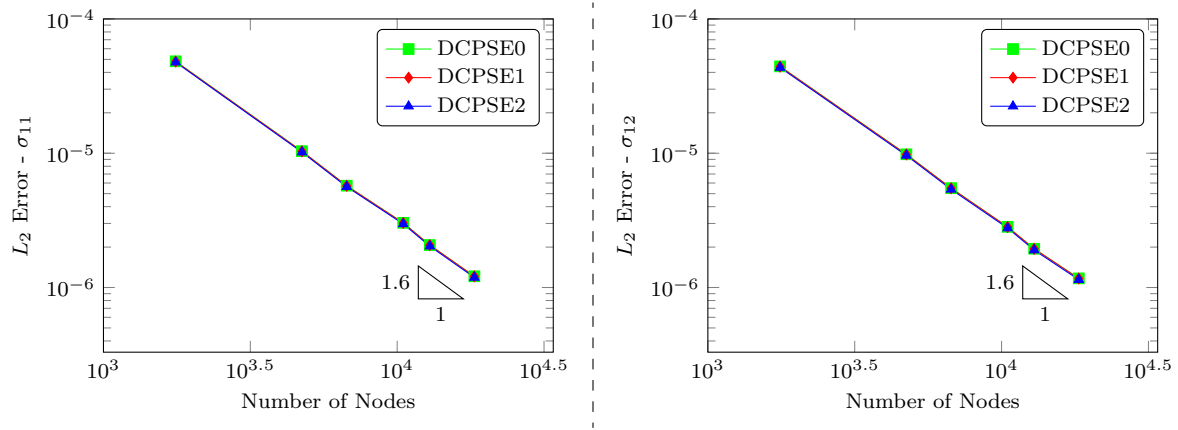
\begin{figure}[H] 
		\centering
		\begin{tabular}{c:c}
			\begin{tikzpicture}
			\begin{axis}[height=6cm,width=7.5cm, ymin=0.00000033,ymax=0.0001,ymode=log,xmin=1000,xmax=34000, xmode=log, legend entries={DCPSE0,DCPSE1,DCPSE2},legend style={ at={(0.5,-0.2)},anchor=south west,legend columns=1, cells={anchor=west},  font=\footnotesize, rounded corners=2pt,}, legend pos=north east,xlabel=Number of Nodes,ylabel=$L_2$ Error - $\sigma_{11}$]
			\addplot+[green,mark=square*,mark options={fill=green}]  table [x=NodeNum, y=L2-REL-S11-DCPSE0, col sep=comma] {DCPSECompCylinder.csv};
			\addplot+[red,mark=diamond*,mark options={fill=red}]   table [x=NodeNum, y=L2-REL-S11-DCPSE1, col sep=comma] {DCPSECompCylinder.csv};
			\addplot+[blue,mark=triangle*,mark options={fill=blue}]   table [x=NodeNum, y=L2-REL-S11-DCPSE2, col sep=comma] {DCPSECompCylinder.csv};
			\logLogSlopeTriangle{0.8}{0.1}{0.16}{1.6}{black};
			\end{axis}
			\end{tikzpicture} &
			
			\begin{tikzpicture}
			\begin{axis}[height=6cm,width=7.5cm, ymin=0.00000033,ymax=0.0001,ymode=log,xmin=1000,xmax=34000, xmode=log, legend entries={DCPSE0,DCPSE1,DCPSE2},legend style={ at={(0.5,-0.2)},anchor=south west,legend columns=1, cells={anchor=west},  font=\footnotesize, rounded corners=2pt,}, legend pos=north east,xlabel=Number of Nodes,ylabel=$L_2$ Error - $\sigma_{12}$]
			\addplot+[green,mark=square*,mark options={fill=green}]  table [x=NodeNum, y=L2-REL-S12-DCPSE0, col sep=comma] {DCPSECompCylinder.csv};
			\addplot+[red,mark=diamond*,mark options={fill=red}]   table [x=NodeNum, y=L2-REL-S12-DCPSE1, col sep=comma] {DCPSECompCylinder.csv};
			\addplot+[blue,mark=triangle*,mark options={fill=blue}]   table [x=NodeNum, y=L2-REL-S12-DCPSE2, col sep=comma] {DCPSECompCylinder.csv};
			\logLogSlopeTriangle{0.8}{0.1}{0.16}{1.6}{black};
			\end{axis}
			\end{tikzpicture} 
		\end{tabular}
		\caption{DC PSE method variations comparison - 2D Cylinder. $L_2$ error as a function of the number of nodes in the model for the DCPSE0, DCPSE1 and DCPSE2 variations of the DC PSE method. No distinction can be observed between the different methods.}
		\label{DCPSE_Comparison}
	\end{figure}
	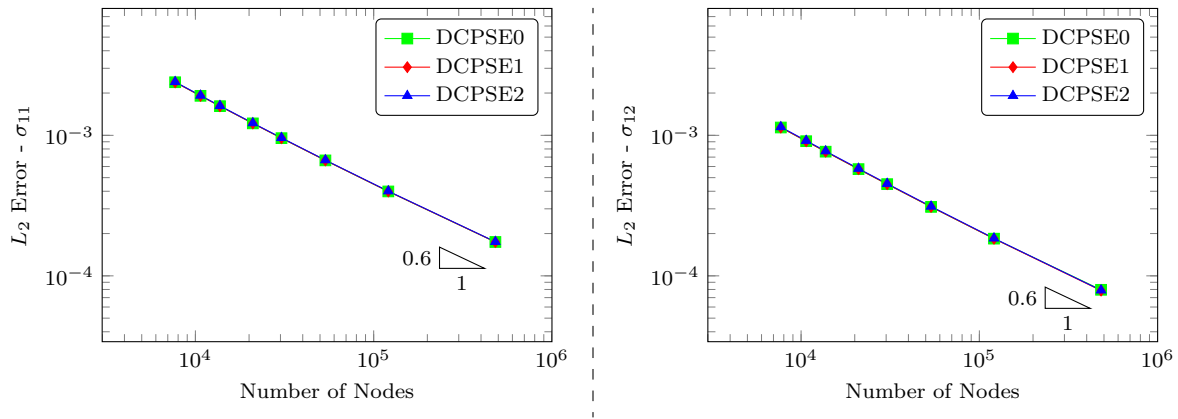
\begin{figure}[H] 
		\centering
		\begin{tabular}{c:c}
			\begin{tikzpicture}
			\begin{axis}[height=6cm,width=7.5cm, ymin=0.000034,ymax=0.008,ymode=log, xmin=3000,xmax=1000000, xmode=log, legend entries={DCPSE0,DCPSE1,DCPSE2},legend style={ at={(0.5,-0.2)},anchor=south west,legend columns=1, cells={anchor=west},  font=\footnotesize, rounded corners=2pt,}, legend pos=north east,xlabel=Number of Nodes,ylabel=$L_2$ Error - $\sigma_{11}$]
			\addplot+[green,mark=square*,mark options={fill=green}]  table [x=NodeNum, y=L2-REL-S11-DCPSE0, col sep=comma] {DCPSECompLShape.csv};
			\addplot+[red,mark=diamond*,mark options={fill=red}]   table [x=NodeNum, y=L2-REL-S11-DCPSE1, col sep=comma] {DCPSECompLShape.csv};
			\addplot+[blue,mark=triangle*,mark options={fill=blue}]   table [x=NodeNum, y=L2-REL-S11-DCPSE2, col sep=comma] {DCPSECompLShape.csv};
			\logLogSlopeTriangle{0.85}{0.1}{0.22}{0.6}{black};
			\end{axis}
			\end{tikzpicture} &
			
			\begin{tikzpicture}
			\begin{axis}[height=6cm,width=7.5cm, ymin=0.000034,ymax=0.008,ymode=log,xmin=3000,xmax=1000000, xmode=log, legend entries={DCPSE0,DCPSE1,DCPSE2},legend style={ at={(0.5,-0.2)},anchor=south west,legend columns=1, cells={anchor=west},  font=\footnotesize, rounded corners=2pt,}, legend pos=north east,xlabel=Number of Nodes,ylabel=$L_2$ Error - $\sigma_{12}$]
			\addplot+[green,mark=square*,mark options={fill=green}]  table [x=NodeNum, y=L2-REL-S12-DCPSE0, col sep=comma] {DCPSECompLShape.csv};
			\addplot+[red,mark=diamond*,mark options={fill=red}]   table [x=NodeNum, y=L2-REL-S12-DCPSE1, col sep=comma] {DCPSECompLShape.csv};
			\addplot+[blue,mark=triangle*,mark options={fill=blue}]   table [x=NodeNum, y=L2-REL-S12-DCPSE2, col sep=comma] {DCPSECompLShape.csv};
			\logLogSlopeTriangle{0.85}{0.1}{0.1}{0.6}{black};
			\end{axis}
			\end{tikzpicture} 
		\end{tabular}
		\caption{DC PSE method variations comparison - 2D L-Shape. $L_2$ error as a function of the number of nodes in the model for the DCPSE0, DCPSE1 and DCPSE2 variations of the DC PSE method. No distinction can be observed between the different methods.}
		\label{DCPSE_Comparison_L-Shape}
	\end{figure}
	It can be observed from Figure \ref{DCPSE_Comparison} and Figure \ref{DCPSE_Comparison_L-Shape} that all three methods lead to very similar results. In order to quantify the difference, the relative difference compared to the DCPSE2 is presented in Figure \ref{DCPSE_Comparison2} and Figure \ref{DCPSE_Comparison2_LShaped} for the DCPSE0 and DCPSE1 methods.
	
	It can be observed from Figure \ref{DCPSE_Comparison2} that the minimum error is obtained with the DCPSE2 method. The DCPSE1 method leads to a slightly lower error than the DCPSE0 method. The errors obtained with the DCPSE0 and DCPSE1 methods are between 1\% and 2.5\% larger than the error obtained with the DCPSE2 method. For the 2D L-shape problem, the results presented in Figure \ref{DCPSE_Comparison2_LShaped} show that the DCPSE2 method leads to a larger error for most node densities. The DCPSE0 and DCPSE1 methods lead to similar errors.
	
	It can be concluded from this study that all variations of the DC PSE method lead to very similar results. The assembly of the linear problem is slightly faster with the DCPSE2 method as the coefficients of the correction function are obtained with the inversion of a linear problem of a lower dimension. The DCPSE1 method leads in general to a lower error than the DCPSE0 method, and thus, has been selected for the comparison to other collocation methods presented in Subsection \ref{MethodComparison}.
	
	\subsection{GFD, DC PSE and Other Methods Comparison}\label{MethodComparison}
	
	\subsubsection{General}
	
	In this section, the GFD and the DC PSE methods are compared to the MLS, IMLS and RBF-FD collocation methods. The same number of support nodes has been chosen for all methods. A 3$^\text{rd}$ order spline weight function has been chosen for the MLS method. The weight function presented in Equation (\ref{IMLSWeight}) has been considered for the IMLS method, and a Gaussian radial basis function has been selected for the RBF-FD methods. The GFD and the DC PSE methods are based on the parameters presented in Table \ref{ParametersSummary}. For reference purpose, results for the finite element method (FEM), obtained using the commercial software ABAQUS \cite{Abaqus2017}, are also included in the comparison.
	
	The results from FEM are extrapolated to the nodes. This allows the results from the FEM to be compared with the same error norms to the results obtained with the collocation methods. The same discretization as for the collocation model has been selected for the FE models. The adjacent nodes of the regular distribution have been grouped into bilinear quadrilateral elements with four integration points.
	
	\subsubsection{Results for the 2D Cylinder Under Internal Pressure}
	
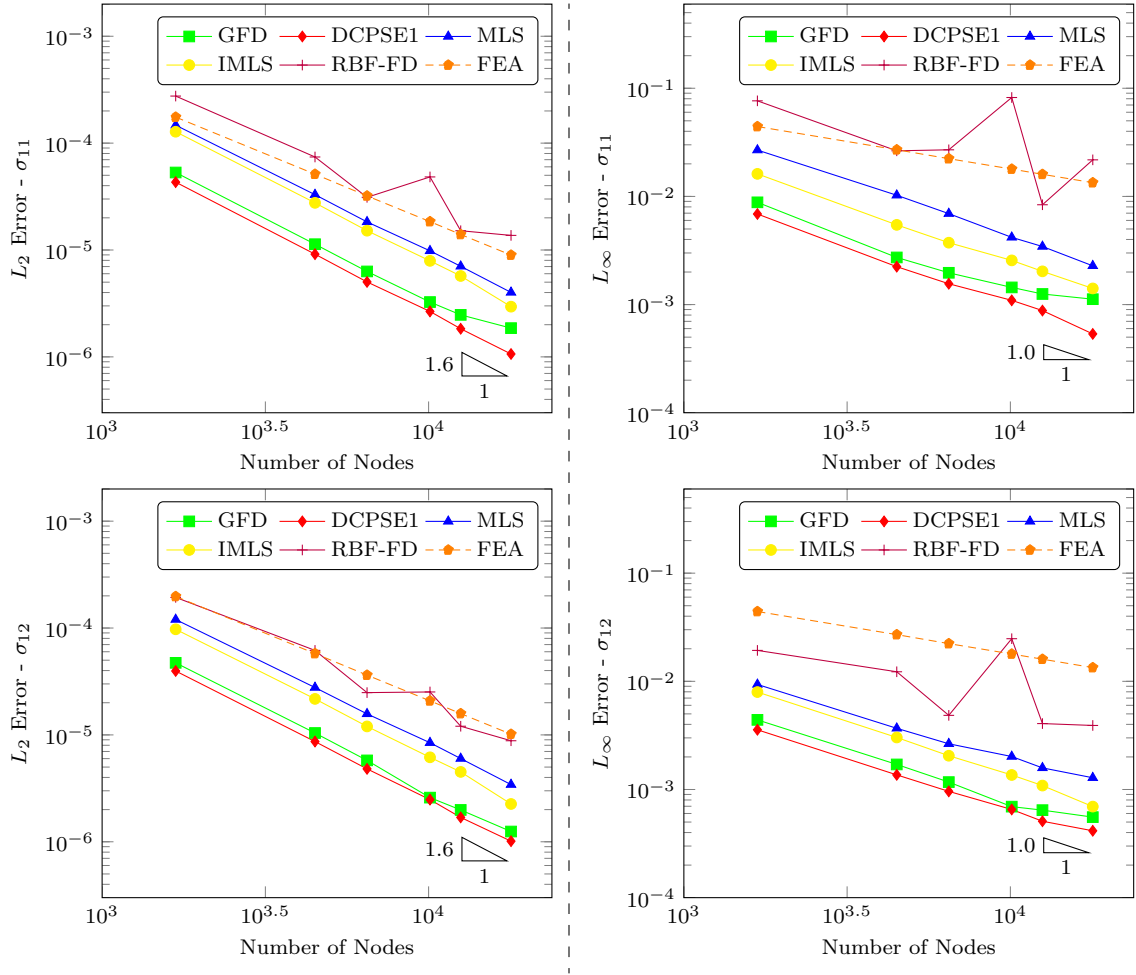
\begin{figure}[H] 
		\centering
		\begin{tabular}{c:c}
			\begin{tikzpicture}
			\begin{axis}[height=7cm,width=7.5cm,ymode=log, ymin=0.0000003,ymax=0.002, xmin=1000,xmode=log, legend entries={GFD,DCPSE1,MLS,IMLS,RBF-FD,FEA},legend style={ at={(0.5,-0.2)},anchor=south west,legend columns=3, cells={anchor=west},  font=\footnotesize, rounded corners=2pt,}, legend pos=north east,xlabel=Number of Nodes,ylabel=$L_2$ Error - $\sigma_{11}$]
			\addplot+[green,mark=square*,mark options={fill=green}]   table [x=NodeNum, y=L2-REL-S11-GFD
			, col sep=comma] {MethodCompCylinder.csv};
			\addplot+[red,mark=diamond*,mark options={fill=red}]   table [x=NodeNum, y=L2-REL-S11-DCPSE1
			, col sep=comma] {MethodCompCylinder.csv};
			\addplot+[blue,mark=triangle*,mark options={fill=blue}]   table [x=NodeNum, y=L2-REL-S11-MLS1
			, col sep=comma] {MethodCompCylinder.csv};
			\addplot+[yellow,mark=*,mark options={fill=yellow}]   table [x=NodeNum, y=L2-REL-S11-IMLS
			, col sep=comma] {MethodCompCylinder.csv};
			\addplot+[purple,mark=+,mark options={fill=purple}]   table [x=NodeNum, y=L2-REL-S11-RBF-FD
			, col sep=comma] {MethodCompCylinder.csv};
			\addplot+[orange,mark=pentagon*,mark options={fill=orange}]   table [x=NodeNum, y=L2-REL-S11-FEA
			, col sep=comma] {MethodCompCylinder.csv};
			\logLogSlopeTriangle{0.9}{0.1}{0.09}{1.6}{black};
			\end{axis}
			\end{tikzpicture} & 
			
			\begin{tikzpicture}
			\begin{axis}[height=7cm,width=7.5cm,ymode=log, ymin=0.0001,ymax=0.6, xmin=1000,xmode=log, legend entries={GFD,DCPSE1,MLS,IMLS,RBF-FD,FEA},legend style={ at={(0.5,-0.2)},anchor=south west,legend columns=3, cells={anchor=west},  font=\footnotesize, rounded corners=2pt,}, legend pos=north east,xlabel=Number of Nodes,ylabel=$L_{\infty}$ Error - $\sigma_{11}$]
			\addplot+[green,mark=square*,mark options={fill=green}]   table [x=NodeNum, y=L-INF-S11-GFD
			, col sep=comma] {MethodCompCylinder.csv};
			\addplot+[red,mark=diamond*,mark options={fill=red}]   table [x=NodeNum, y=L-INF-S11-DCPSE1
			, col sep=comma] {MethodCompCylinder.csv};
			\addplot+[blue,mark=triangle*,mark options={fill=blue}]   table [x=NodeNum, y=L-INF-S11-MLS1
			, col sep=comma] {MethodCompCylinder.csv};
			\addplot+[yellow,mark=*,mark options={fill=yellow}]   table [x=NodeNum, y=L-INF-S11-IMLS
			, col sep=comma] {MethodCompCylinder.csv};
			\addplot+[purple,mark=+,mark options={fill=purple}]   table [x=NodeNum, y=L-INF-S11-RBF-FD
			, col sep=comma] {MethodCompCylinder.csv};
			\addplot+[orange,mark=pentagon*,mark options={fill=orange}]   table [x=NodeNum, y=L-INF-S11-FEA
			, col sep=comma] {MethodCompCylinder.csv};
			\logLogSlopeTriangle{0.9}{0.1}{0.13}{1.0}{black};
			\end{axis}
			\end{tikzpicture} \\
			
			\begin{tikzpicture}
			\begin{axis}[height=7cm,width=7.5cm,ymode=log, ymin=0.0000003,ymax=0.002, xmin=1000,xmode=log, legend entries={GFD,DCPSE1,MLS,IMLS,RBF-FD,FEA},legend style={ at={(0.5,-0.2)},anchor=south west,legend columns=3, cells={anchor=west},  font=\footnotesize, rounded corners=2pt,}, legend pos=north east,xlabel=Number of Nodes,ylabel=$L_2$ Error - $\sigma_{12}$]
			\addplot+[green,mark=square*,mark options={fill=green}]   table [x=NodeNum, y=L2-REL-S12-GFD
			, col sep=comma] {MethodCompCylinder.csv};
			\addplot+[red,mark=diamond*,mark options={fill=red}]   table [x=NodeNum, y=L2-REL-S12-DCPSE1
			, col sep=comma] {MethodCompCylinder.csv};
			\addplot+[blue,mark=triangle*,mark options={fill=blue}]   table [x=NodeNum, y=L2-REL-S12-MLS1
			, col sep=comma] {MethodCompCylinder.csv};
			\addplot+[yellow,mark=*,mark options={fill=yellow}]   table [x=NodeNum, y=L2-REL-S12-IMLS
			, col sep=comma] {MethodCompCylinder.csv};
			\addplot+[purple,mark=+,mark options={fill=purple}]   table [x=NodeNum, y=L2-REL-S12-RBF-FD
			, col sep=comma] {MethodCompCylinder.csv};
			\addplot+[orange,mark=pentagon*,mark options={fill=orange}]   table [x=NodeNum, y=L2-REL-S12-FEA
			, col sep=comma] {MethodCompCylinder.csv};
			\logLogSlopeTriangle{0.9}{0.1}{0.09}{1.6}{black};
			\end{axis}
			\end{tikzpicture}&
			
			\begin{tikzpicture}
			\begin{axis}[height=7cm,width=7.5cm,ymode=log, ymin=0.0001,ymax=0.6, xmin=1000,xmode=log, legend entries={GFD,DCPSE1,MLS,IMLS,RBF-FD,FEA},legend style={ at={(0.5,-0.2)},anchor=south west,legend columns=3, cells={anchor=west},  font=\footnotesize, rounded corners=2pt,}, legend pos=north east,xlabel=Number of Nodes,ylabel=$L_{\infty}$ Error - $\sigma_{12}$]
			\addplot+[green,mark=square*,mark options={fill=green}]   table [x=NodeNum, y=L-INF-S12-GFD
			, col sep=comma] {MethodCompCylinder.csv};
			\addplot+[red,mark=diamond*,mark options={fill=red}]   table [x=NodeNum, y=L-INF-S12-DCPSE1
			, col sep=comma] {MethodCompCylinder.csv};
			\addplot+[blue,mark=triangle*,mark options={fill=blue}]   table [x=NodeNum, y=L-INF-S12-MLS1
			, col sep=comma] {MethodCompCylinder.csv};
			\addplot+[yellow,mark=*,mark options={fill=yellow}]   table [x=NodeNum, y=L-INF-S12-IMLS
			, col sep=comma] {MethodCompCylinder.csv};
			\addplot+[purple,mark=+,mark options={fill=purple}]   table [x=NodeNum, y=L-INF-S12-RBF-FD
			, col sep=comma] {MethodCompCylinder.csv};
			\addplot+[orange,mark=pentagon*,mark options={fill=orange}]   table [x=NodeNum, y=L-INF-S12-FEA, col sep=comma] {MethodCompCylinder.csv};
			\logLogSlopeTriangle{0.9}{0.1}{0.11}{1.0}{black};
			\end{axis}
			\end{tikzpicture}\\
		\end{tabular}
		\caption{Methods comparison - 2D Cylinder. $L_2$ (Left) and $L_{\infty}$ (Right) errors as a function of the number of nodes for various collocation methods (i.e. GFD, DCPSE1, MLS, IMLS and RBF-FD) and for the FEM. The lowest error is obtained with the DCPSE1 method.}
		\label{MethodCompResultsCylinder1}
	\end{figure}

	From Figure \ref{MethodCompResultsCylinder1}, it can be observed that the DCPSE1 method leads to the lowest error in terms of $L_2$ and $L_{\infty}$ norms for the $\sigma_{11}$ and $\sigma_{12}$ stress components. The GFD method leads to a slightly higher error than the DC PSE method. The MLS and IMLS methods lead to very similar results. The IMLS method lead to an error constantly lower than the MLS method. Finally, the FEM and the RBF-FD method lead to the largest errors. The error obtained with the RBF-FD method does not monotonically decrease as the node density increases.
	
	\subsubsection{Results for the L-Shape Problem}
	
	The different methods are also compared for the 2D L-shape problem. Only the $L_2$ error is presented for this problem as the $L_{\infty}$ error diverges to infinity as the distance to the singular node tends to zero. The results are presented in Figure \ref{MethodCompResultsLShaped} below.
	\begin{figure}[H] 
		\centering
		\begin{tabular}{c:c}
			\begin{tikzpicture}
			\begin{axis}[height=7cm,width=7.5cm,ymode=log, ymin=0.00002,ymax=0.03, xmin=3000,xmax=1000000,xmode=log, legend entries={GFD,DCPSE1,MLS,IMLS,RBF-FD,FEA},legend style={ at={(0.5,-0.2)},anchor=south west,legend columns=3, cells={anchor=west},  font=\footnotesize, rounded corners=2pt,}, legend pos=north east,xlabel=Number of Nodes,ylabel=$L_2$ Error - $\sigma_{11}$]
			\addplot+[green,mark=square*,mark options={fill=green}]   table [x=NodeNum, y=L2-REL-S11-GFD
			, col sep=comma] {MethodCompLShaped.csv};
			\addplot+[red,mark=diamond*,mark options={fill=red}]   table [x=NodeNum, y=L2-REL-S11-DCPSE1
			, col sep=comma] {MethodCompLShaped.csv};
			\addplot+[blue,mark=triangle*,mark options={fill=blue}]   table [x=NodeNum, y=L2-REL-S11-MLS1
			, col sep=comma] {MethodCompLShaped.csv};
			\addplot+[yellow,mark=*,mark options={fill=yellow}]   table [x=NodeNum, y=L2-REL-S11-IMLS
			, col sep=comma] {MethodCompLShaped.csv};
			\addplot+[purple,mark=+,mark options={fill=purple}]   table [x=NodeNum, y=L2-REL-S11-RBF-FD
			, col sep=comma] {MethodCompLShaped.csv};
			\addplot+[orange,mark=pentagon*,mark options={fill=orange}]   table [x=NodeNum, y=L2-REL-S11-FEA
			, col sep=comma] {MethodCompLShaped.csv};
			\logLogSlopeTriangle{0.85}{0.1}{0.2}{0.6}{black};
			\end{axis}
			\end{tikzpicture} & 
			
			\begin{tikzpicture}
			\begin{axis}[height=7cm,width=7.5cm,ymode=log, ymin=0.00002,ymax=0.03, xmin=3000,xmax=1000000,xmode=log, legend entries={GFD,DCPSE1,MLS,IMLS,RBF-FD,FEA},legend style={ at={(0.5,-0.2)},anchor=south west,legend columns=3, cells={anchor=west},  font=\footnotesize, rounded corners=2pt,}, legend pos=north east,xlabel=Number of Nodes,ylabel=$L_2$ Error - $\sigma_{12}$]
			\addplot+[green,mark=square*,mark options={fill=green}]   table [x=NodeNum, y=L2-REL-S12-GFD
			, col sep=comma] {MethodCompLShaped.csv};
			\addplot+[red,mark=diamond*,mark options={fill=red}]   table [x=NodeNum, y=L2-REL-S12-DCPSE1
			, col sep=comma] {MethodCompLShaped.csv};
			\addplot+[blue,mark=triangle*,mark options={fill=blue}]   table [x=NodeNum, y=L2-REL-S12-MLS1
			, col sep=comma] {MethodCompLShaped.csv};
			\addplot+[yellow,mark=*,mark options={fill=yellow}]   table [x=NodeNum, y=L2-REL-S12-IMLS
			, col sep=comma] {MethodCompLShaped.csv};
			\addplot+[purple,mark=+,mark options={fill=purple}]   table [x=NodeNum, y=L2-REL-S12-RBF-FD
			, col sep=comma] {MethodCompLShaped.csv};
			\addplot+[orange,mark=pentagon*,mark options={fill=orange}]   table [x=NodeNum, y=L2-REL-S12-FEA
			, col sep=comma] {MethodCompLShaped.csv};
			\logLogSlopeTriangle{0.85}{0.1}{0.1}{0.6}{black};
			\end{axis}
			\end{tikzpicture} \\	
		\end{tabular}
		\caption{Methods comparison - 2D L-Shape. $L_2$ errors as a function of the number of nodes for various collocation methods (i.e. GFD, DCPSE1, MLS, IMLS and RBF-FD) and for the FEM. The lowest error is obtained with the FEM. The collocation method leading to the lowest error is the MLS method.}
		\label{MethodCompResultsLShaped}
	\end{figure}
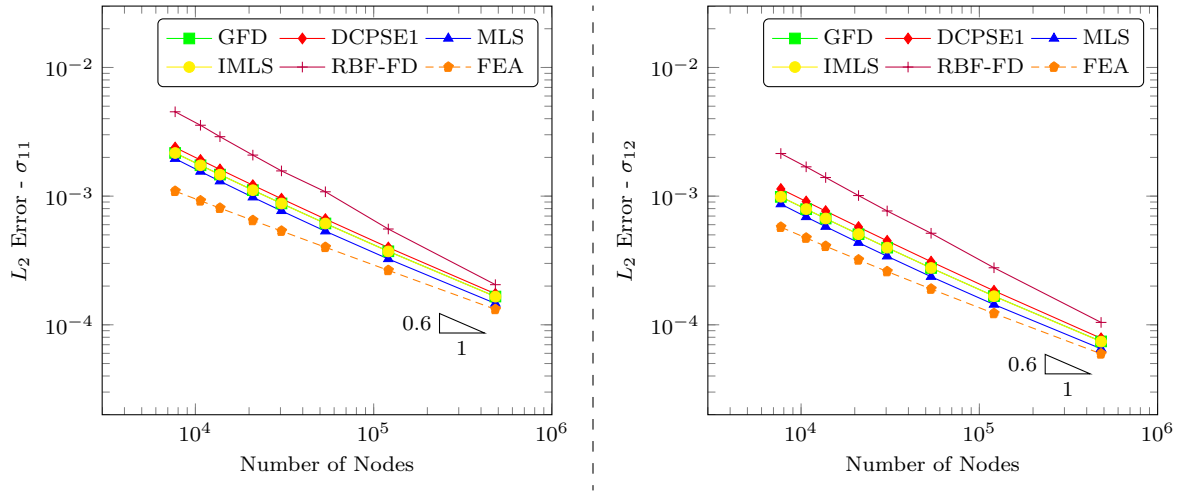
	We can see in Figure \ref{MethodCompResultsLShaped} that the results obtained with the FEM are the closest to the analytical solution. The MLS, IMLS, GFD and DC PSE methods lead to similar results. The trend is however opposite to the results presented in Figure \ref{MethodCompResultsCylinder1}. It shows hows the methods are affected by a rapid change in the field solution. The MLS method is the method leading to the lowest error among the collocation methods. Finally, the RBF-FD is the method leading to the highest error. 	
	
	\subsubsection{Convergence Rate and Computational Expense}
	
	The methods are also compared in terms of convergence rate and solution time. Results are summarized in Table \ref{MethodCompSummaryTable} below. It can be observed that the RBF-FD method has the lowest convergence rate for the 2D cylinder problem, and the largest convergence rate for the L-shape problem. The IMLS method shows the largest convergence rate for the 2D cylinder problem. The GFD method is the one having the lowest convergence rate for the L-shape problem.
	
	The right column of Table \ref{MethodCompSummaryTable} shows that the computation time for the MLS and IMLS method is significantly larger than the computation time for the other methods. This is due to the assembly step which is more time consuming for these methods, as the system presented in Equation (\ref{DerivativeSystem_MLS}) needs to be solved for each derivative.
	\begin{table}[h]
		\centering
		\caption{Method Comparison Summary}
		\label{MethodCompSummaryTable}
		\renewcommand{\arraystretch}{1.5}
		\begin{tabular}{|l|c|c|c|}
			\hline
			\multicolumn{1}{|c|}{\multirow{2}{*}{\centering \textbf{Method}}} & \multicolumn{2}{c|}{\begin{tabular}[x]{@{}c@{}}\textbf{Average $L_2$}\\ \textbf{Convergence Rate} \end{tabular} }&
			\multirow{2}{*}{\begin{tabular}[x]{@{}c@{}}\textbf{Computation Time $^\text{(1)}$} \end{tabular} }\\
			& \text{2D Cylinder} & \text{2D L-Shape} &  \\
			\hline
			\centering
			GFD & 1.5109 & 0.6244 & 9.3s \\
			DCPSE1  & 1.5592 & 0.6405 & 11.4s\\
			MLS & 1.5064 & 0.6249 & 19.7s \\
			IMLS & 1.5717 & 0.6279 & 20.5s \\
			RBF-FD & 1.0952 & 0.7401 & 10.1s \\
			\hline
			\multicolumn{4}{@{}l}{\footnotesize (1) Based on a 12,087 nodes 2D cylinder model solved with a direct solver.}
		\end{tabular}
	\end{table}
	
	\pagebreak
	
	The computation times presented in Table \ref{MethodCompSummaryTable} can be split into four main steps. These steps are:
	\begin{itemize}
		\item Problem initialization;
		\item Matrix assembly;
		\item Solution of the linear problem;
		\item Postprocessing and results output.
	\end{itemize}
	The initialization step consists in loading the problem from the input file and searching for the node neighbors (the nodes to be included in the support of each collocation nodes). During the assembly step, all the derivatives of the unknown field are approximated as a function of the field, and the linear problem is assembled. For the solution of the linear problem, while direct solvers can be used for 2D problems of a reasonable size, iterative solvers should be used for large 3D problems, as the matrix of the linear problem is significantly denser than for 2D problems. Finally, the postprocessing step consists of the computation of the quantities of interest (stress components based on the displacement field). In this work, the linear solver MUMPS \cite{MUMPS01,MUMPS02} and the iterative solver PETSc KSP \cite{petsc-user-ref,petsc-efficient} have been used. The analyses have been run using a C++ code developed in-house. The code was run on a machine equipped with an Intel Xeon E5-1650 processor at 3.2 GHz. A single process and a single thread has been used for the analyses.
	
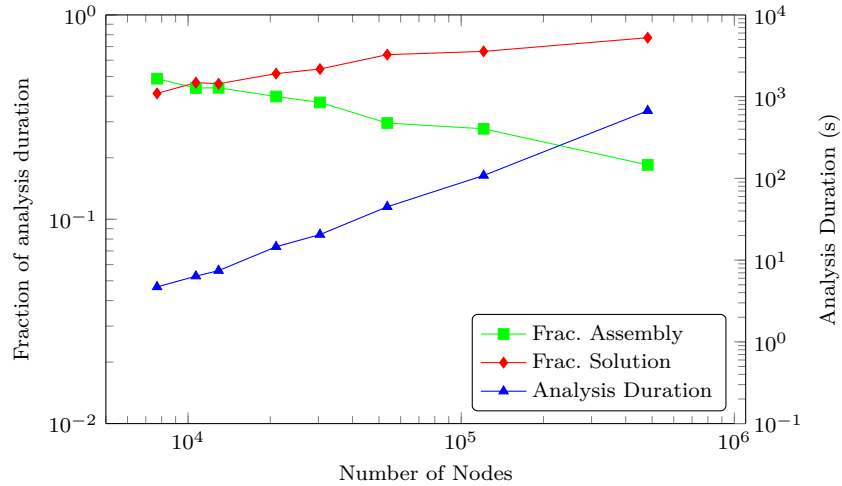
\begin{figure}[H] 
		\centering
		\begin{tikzpicture}
		\begin{axis}[height=7cm,width=10cm,ymode=log, ymin=0.01,ymax=1, xmin=5000,xmax=1100000,xmode=log, legend entries={Frac. Assembly,Frac. Solution,Analysis Duration},legend style={ at={(0.5,-0.2)},anchor=south east,legend columns=1, cells={anchor=west},  font=\footnotesize, rounded corners=2pt,}, axis y line*=left, legend pos=south east, xlabel=Number of Nodes,ylabel=Fraction of analysis duration]
		\addplot+[green,mark=square*,mark options={fill=green}]   table [x=X, y=GFD-Build, col sep=comma] {TimeSplitLShaped.csv};
		\addplot+[red,mark=diamond*,mark options={fill=red}]   table [x=X, y=GFD-Solve, col sep=comma] {TimeSplitLShaped.csv};
		\addplot+[blue,mark=triangle*,mark options={fill=blue}]   table [x=X, y=X, col sep=comma] {TimeSplitLShaped.csv};
		\end{axis}
		\begin{axis}[height=7cm,width=10cm,ymode=log, ymin=0.1,ymax=10000, xmin=5000,xmax=1100000,xmode=log, legend entries={},legend style={ at={(0.5,-0.2)},anchor=south east,legend columns=1, cells={anchor=west},  font=\footnotesize, rounded corners=2pt,}, axis x line=none, legend pos=south west,axis y line*=right,ylabel=Analysis Duration (s)]
		\addplot+[blue,mark=triangle*,mark options={fill=blue}]   table [x=X, y=GFD-Time, col sep=comma] {TimeSplitLShaped.csv};
		\end{axis}
		\end{tikzpicture}
		\caption{Computation time split and analysis duration - GFD method. Impact of the number of nodes on the fraction of the analysis spent in each step and total analysis time.}
		\label{TimeSplitLShapedGFD}
	\end{figure}

\begin{figure}[H] 
		\centering
		\begin{tikzpicture}
		\begin{axis}[height=7cm,width=10cm,ymode=log, ymin=0.01,ymax=1, xmin=5000,xmax=1100000,xmode=log, legend entries={Frac. Assembly,Frac. Solution,Analysis Duration},legend style={ at={(0.5,-0.2)},anchor=south east,legend columns=2, cells={anchor=west},  font=\footnotesize, rounded corners=2pt,}, axis y line*=left, legend pos=south east, xlabel=Number of Nodes,ylabel=Fraction of analysis duration]
		\addplot+[green,mark=square*,mark options={fill=green}]   table [x=X, y=IMLS-Build, col sep=comma] {TimeSplitLShaped.csv};
		\addplot+[red,mark=diamond*,mark options={fill=red}]   table [x=X, y=IMLS-Solve, col sep=comma] {TimeSplitLShaped.csv};
		\addplot+[blue,mark=triangle*,mark options={fill=blue}]   table [x=X, y=X, col sep=comma] {TimeSplitLShaped.csv};
		\end{axis}
		\begin{axis}[height=7cm,width=10cm,ymode=log, ymin=0.1,ymax=10000, xmin=5000,xmax=1100000,xmode=log, legend entries={},legend style={ at={(0.5,-0.2)},anchor=south east,legend columns=1, cells={anchor=west},  font=\footnotesize, rounded corners=2pt,}, axis x line=none, legend pos=south west,axis y line*=right,ylabel=Analysis Duration (s)]
		\addplot+[blue,mark=triangle*,mark options={fill=blue}]   table [x=X, y=IMLS-Time, col sep=comma] {TimeSplitLShaped.csv};
		\end{axis}
		\end{tikzpicture}
		\caption{Computation time split and analysis duration - IMLS method. Impact of the number of nodes on the fraction of the analysis spent in each step and total analysis time.}
		\label{TimeSplitLShapedIMLS}
	\end{figure}
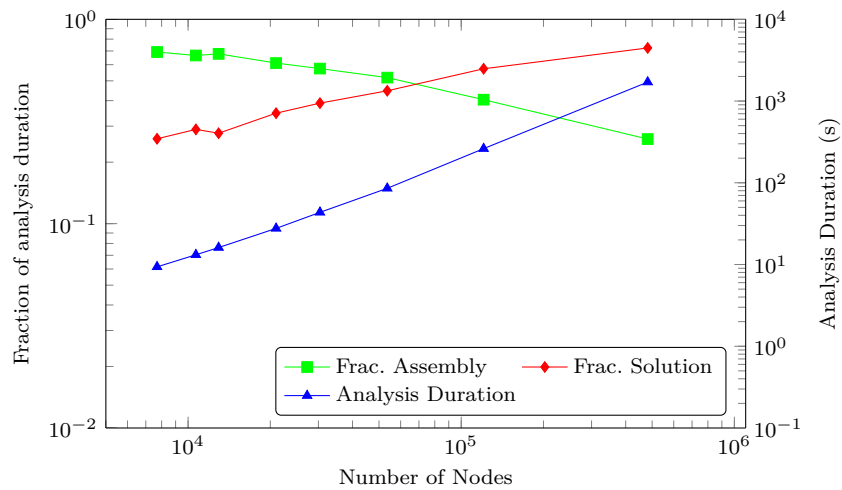

	For various node densities, the fraction of the total analysis time of the matrix assembly step and of the solution step is presented in Figure \ref{TimeSplitLShapedGFD} and Figure \ref{TimeSplitLShapedIMLS} for the L-shape problem. The duration of the problem initialization and of the postprocessing steps is negligible compare to the two other steps. Is has not been presented in those figures. The total analysis duration has also been presented in these figures on a secondary axis. The results are presented for the GFD and for the IMLS methods as these methods are, respectively, the fastest and the slowest methods for the node density presented in Table \ref{MethodCompSummaryTable}.
	It can be observed from Figure \ref{TimeSplitLShapedGFD} and Figure \ref{TimeSplitLShapedIMLS} that for both methods the trend of the results is similar. The assembly step represents the largest proportion of the total analysis time for 2D problems of small dimensions. As the number of degrees of freedom increases, the fraction of the solution step (here with a direct solver) in the analysis time increases and becomes larger than the fraction of the assembly step. This result is expected as the assembly time increases linearly with the number of nodes while the solution time increases exponentially with the number of degrees of freedom. The increased computation effort required by the IMLS method during the assembly step is observed on this graph as the fraction of the analysis spent assembling the matrix is larger than for the GFD method.
	
	\section{Three Dimensional Problems} \label{3DResults}
	
	The parametric study presented earlier allowed the selection of ``optimal" weight functions and support sizes for solving problems from different fields of application. In this section, we present the results from the stress analysis of various three dimensional problems. The results in terms of von Mises stress obtained with the GFD method are compared to the results obtained with ABAQUS. For a consistent comparison, the same discretization has been used for both methods. The results obtained from the FEA are extrapolated to the nodes. This tends to overestimate the error for the FEA method as the stresses are less accurate at the nodes than at the integration points. This allows, however, a comparison, at each node, of the results derived with the FE method to the results derived with the GFD method.
	In order to obtain convergence of the 3D problems with the collocation method, a relatively large number of nodes is used. Such a high density is required to capture the details of the geometry.
	
	\paragraph{Flange Model} \
	
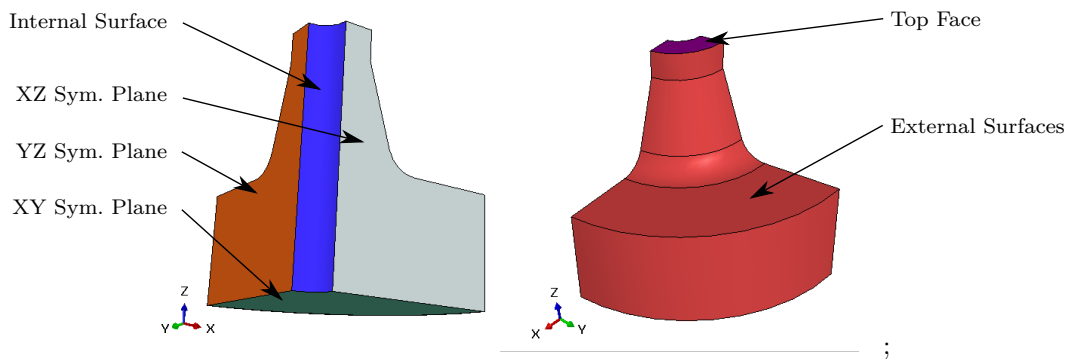
\begin{figure}[H] 
		\centering
		\begin{tikzpicture}
		\def\svgwidth{10cm}
		\node at (0,0) {\includegraphics{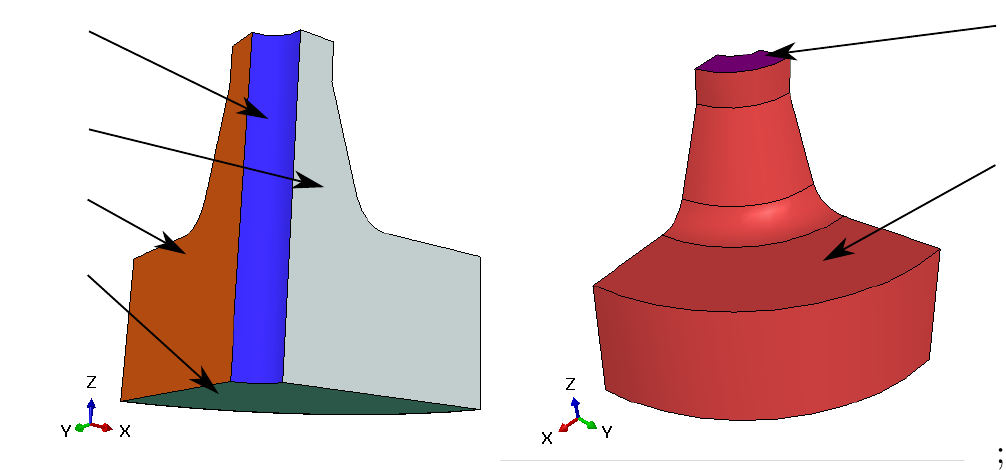}};
		\node[color=black] at (-4.3,2.1) [left] {Internal Surface};
		\node[color=black] at (-4.3,1.1) [left] {XZ Sym. Plane};
		\node[color=black] at (-4.3,0.35) [left] {YZ Sym. Plane};
		\node[color=black] at (-4.3,-0.4) [left] {XY Sym. Plane};
		\node[color=black] at (5,2.1) [right] {Top Face};
		\node[color=black] at (5,0.7) [right] {External Surfaces};
		\end{tikzpicture}
		\caption{Flange Model and Boundary Conditions}
		\label{FlangeBCs}
	\end{figure}

	The first problem considered in this section is an ISO flange. Only a quarter of the flange has been modeled due to the symmetries of the domain. The model and the various surfaces, on which boundary conditions are applied, are presented in Figure \ref{FlangeBCs}. Two load cases have been considered for this model: an internal pressure loading and an axial traction imposed by the connected pipe. The boundary conditions associated to each load case are presented in Table \ref{FlangeLCsBCs}.
	
	\begin{table}[h]
		\centering
		\caption{Boundary conditions applied to the flange for the pressure and displacement load cases. The surfaces are highlighted in Figure \ref{FlangeBCs}.}
		\label{FlangeLCsBCs}
		\renewcommand{\arraystretch}{1.5}
		\begin{tabular}{|l|l|l|}
			\hline
			\multicolumn{1}{|c|}{\multirow{2}{*}{\centering \textbf{Surface}}} & \multicolumn{2}{c|}{\textbf{Boundary Conditions}}\\
			& \multicolumn{1}{c|}{\text{Pressure Loading}} & \multicolumn{1}{c|}{\text{Diplacement loading}} \\
			\hline
			\centering
			XY Sym. Plane & Constrained in the Z direction & Constrained in the Z direction \\
			XZ Sym. Plane  & Constrained in the Y direction & Constrained in the Y direction\\
			YZ Sym. Plane & Constrained in the X direction & Constrained in the X direction\\
			Internal Surface & Constant pressure of 1.0 & Stress free\\
			External Surface & Stress free  & Stress free\\
			Top Face & Constrained in the Z direction & Applied displacement of 6.2e-04 in the Z direction\\
			\hline
		\end{tabular}
	\end{table}
	
\begin{figure}[H] 
		\centering
		\subfloat[]{\includegraphics[scale=0.35]{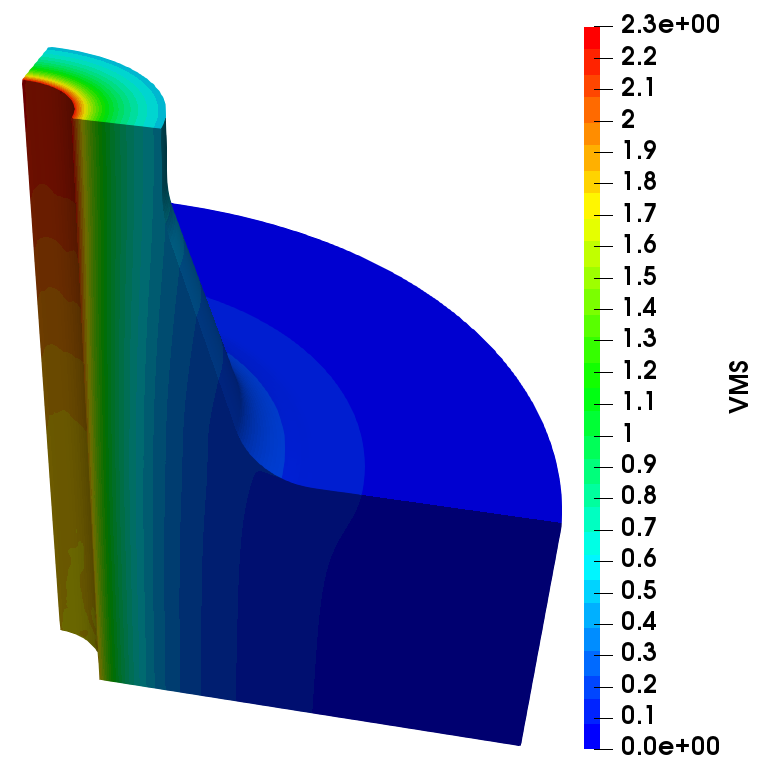}}
		\qquad
		\subfloat[]{\includegraphics[scale=0.35]{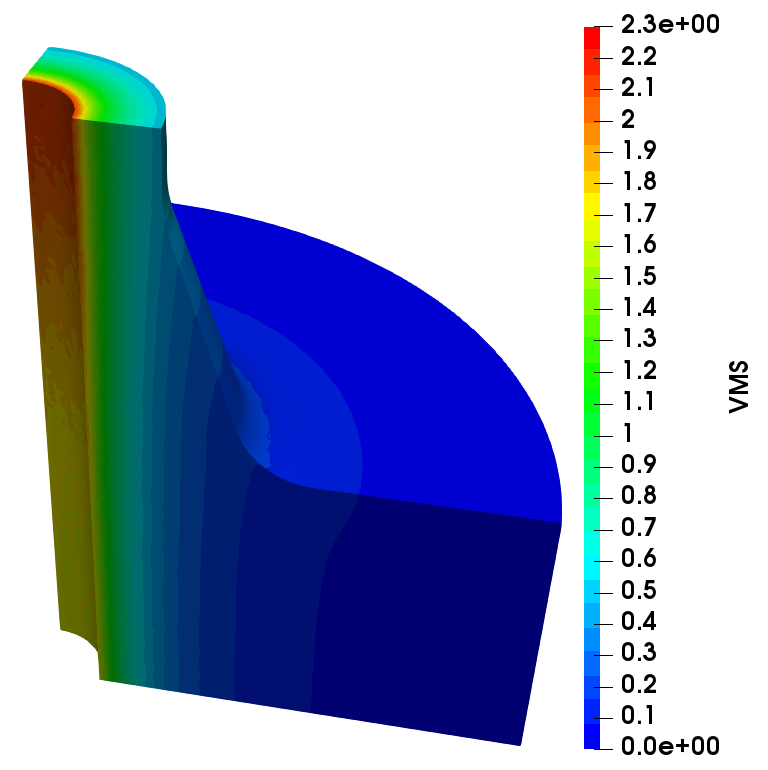}}
		\caption{Flange ISO PN50 DN25 subject to an internal pressure - von Mises stress results from the GFD method (a) and FEM (b) (548,648 nodes). The results from both models are very similar. The stress on the inner surface of the flange is larger for the GFD method.}
		\label{ISO_Flange_Pressure}
	\end{figure}

\begin{figure}[H] 
		\centering
		\subfloat[]{\includegraphics[scale=0.35]{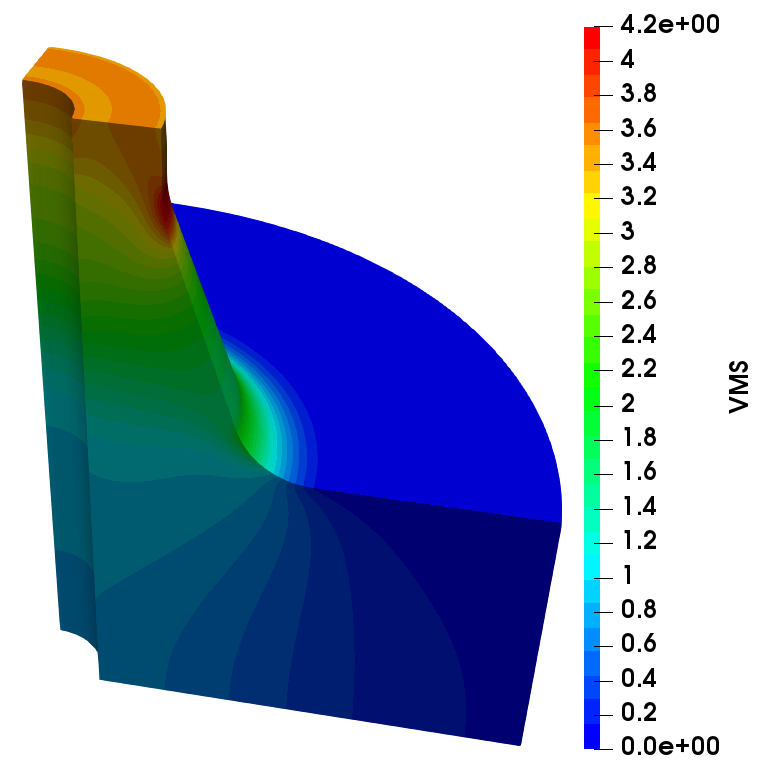}}
		\qquad
		\subfloat[]{\includegraphics[scale=0.35]{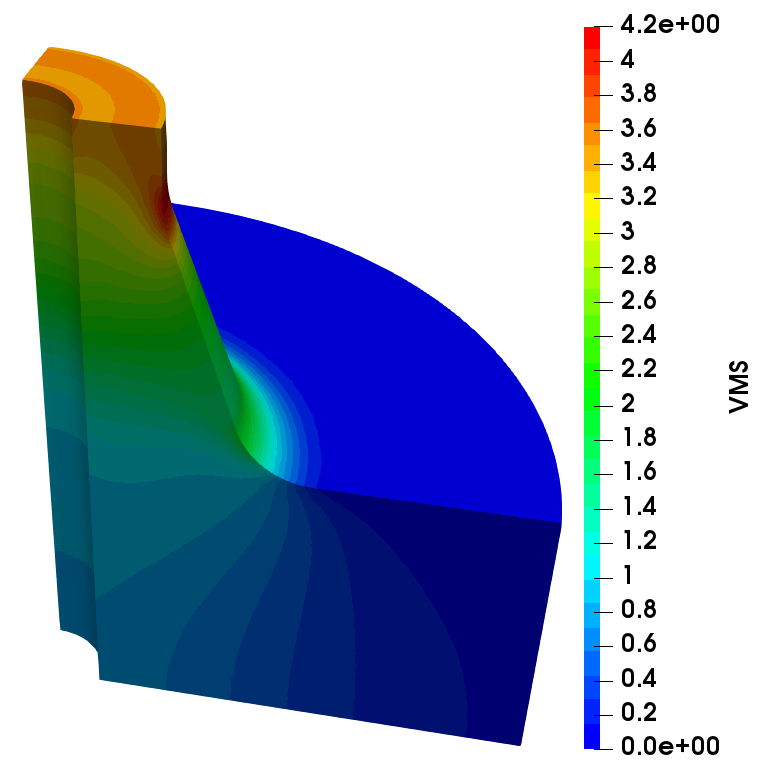}}
		\caption{Flange ISO PN50 DN25 under traction - von Mises stress results from the GFD method (a) and FEM (b) (548,648 nodes). The results from both models are very similar. The stress in the neck of the flange is slightly larger for the GFD method.}
		\label{ISO_Flange_Traction}
	\end{figure}

	Figure \ref{ISO_Flange_Pressure} and Figure \ref{ISO_Flange_Traction} show the von Mises stresses obtained with the GFD and the FE methods, respectively, for the internal pressure and traction load cases.
	
\begin{figure}[H] 
		\centering
		\subfloat[]{\includegraphics[scale=0.35]{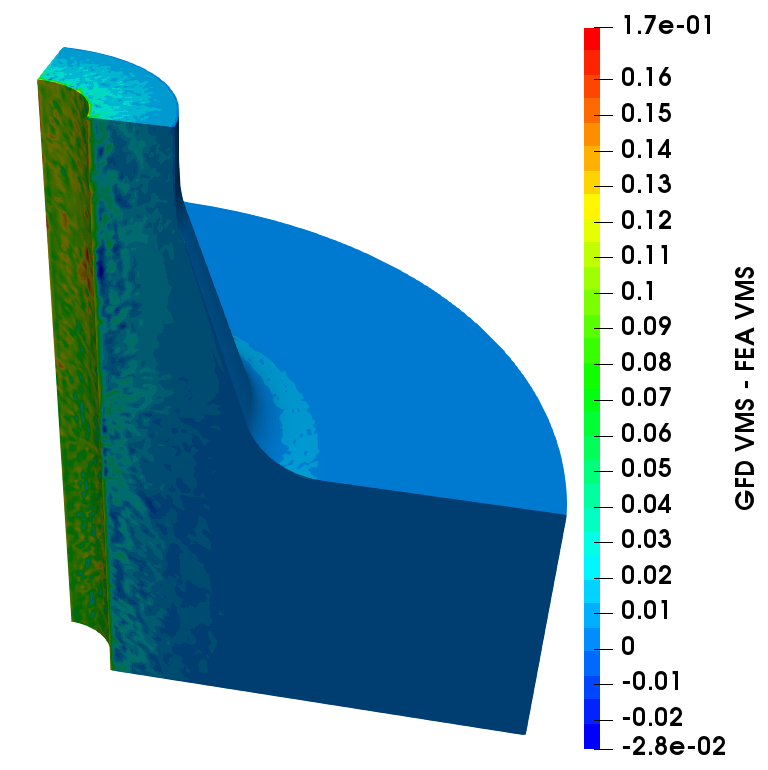}}
		\qquad
		\subfloat[]{\includegraphics[scale=0.35]{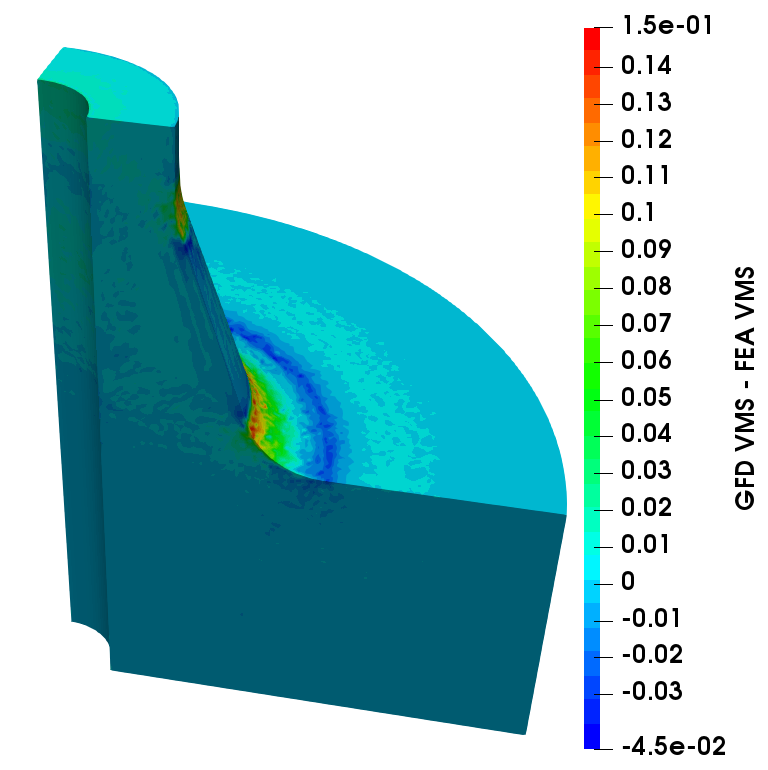}}
		\caption{Flange ISO PN50 DN25 under traction - Difference between von Mises stress results obtained from GFD method and FEM for the internal pressure load case (a) and the traction load case (b) (548,648 nodes). The von Mises stress results are larger for the GFD method on the the inner surface for the flange under internal pressure, and in the neck and in the cone bottom section for the flange under traction.}
		\label{ISO_Flange_CollocFEA}
	\end{figure}

	The results obtained from the GFD method and from the FEM are very close. In order to visually assess the difference between the two solutions, the difference between the von Mises stress results obtained from the GFD model and from the FE model are presented in Figure \ref{ISO_Flange_CollocFEA} for both load cases.
	
	It can be observed from Figure \ref{ISO_Flange_CollocFEA}(a) that the stresses on the inner surface of the flange are larger for the GFD method. From Figure \ref{ISO_Flange_CollocFEA}(b), it can be observed that the stresses obtained in the neck and in the bottom of the conical section are larger for the GFD method.
	
	\paragraph{Blade Model} \
	
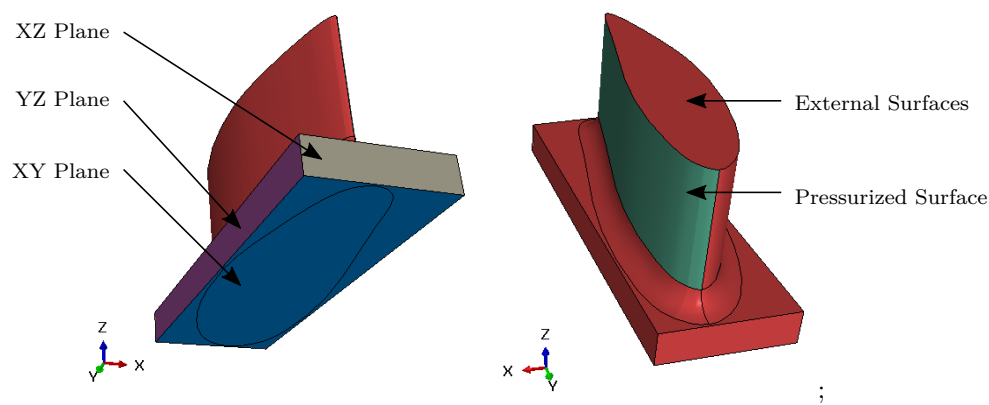
\begin{figure}[H] 
		\centering
		\begin{tikzpicture}
		\def\svgwidth{10cm}
		\node at (0,0) {\includegraphics{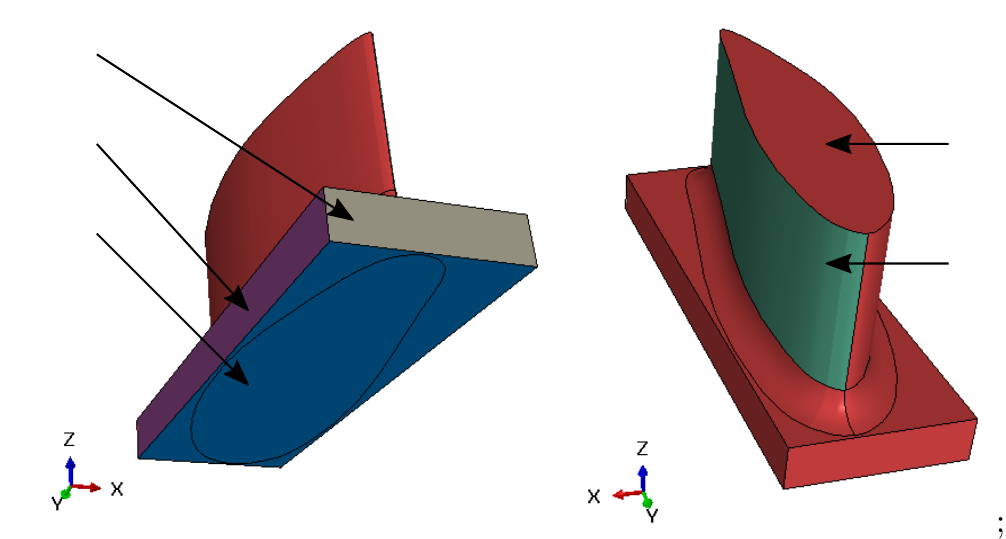}};
		\node[color=black] at (-4.2,2.2) [left] {XZ Plane};
		\node[color=black] at (-4.2,1.3) [left] {YZ Plane};
		\node[color=black] at (-4.2,0.35) [left] {XY Plane};
		\node[color=black] at (4.6,1.25) [right] {External Surfaces};
		\node[color=black] at (4.6,0.0) [right] {Pressurized Surface};
		\end{tikzpicture}
		\caption{Blade model and boundary conditions}
		\label{BladeBCs}
	\end{figure}

	Figure \ref{BladeBCs} presents a simplified model of a high pressure blade. The surfaces, on which the boundary conditions are applied, are presented in this figure. The nodes in the planes YZ, XZ and XY are, respectively, fixed in the X, Y and Z directions. A constant pressure resulting from a gas flow is applied on the pressurized surface. The remaining surfaces of the blade are considered stress free.
	
\begin{figure}[H] 
		\centering
		\subfloat[]{\includegraphics[scale=0.6]{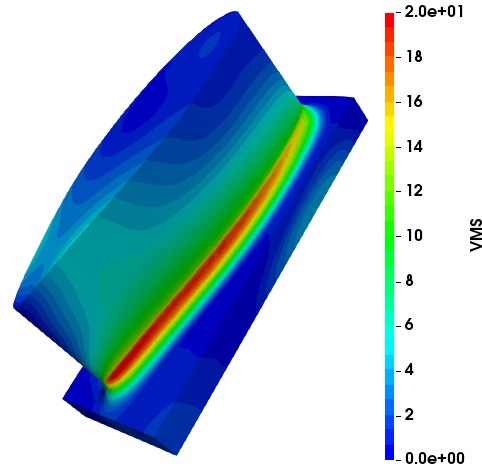}}
		\subfloat[]{\includegraphics[scale=0.6]{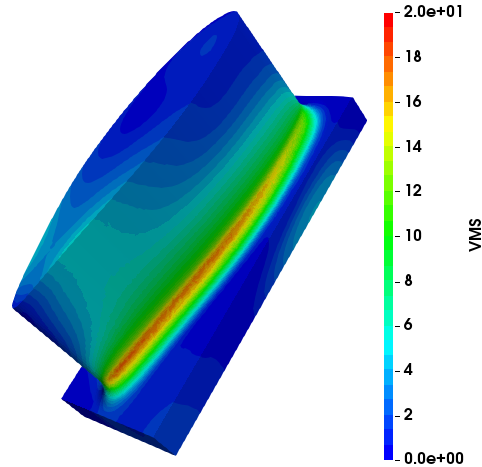}}
		\caption{Simplified high pressure blade subjected to a uniform pressure on one face - The von Mises stress results from the GFD method (a) and FEM (b) (484,238 nodes).}
		\label{SmallBlade}
	\end{figure}

\begin{figure}[H] 
		\centering
		\subfloat[]{\includegraphics[scale=0.4]{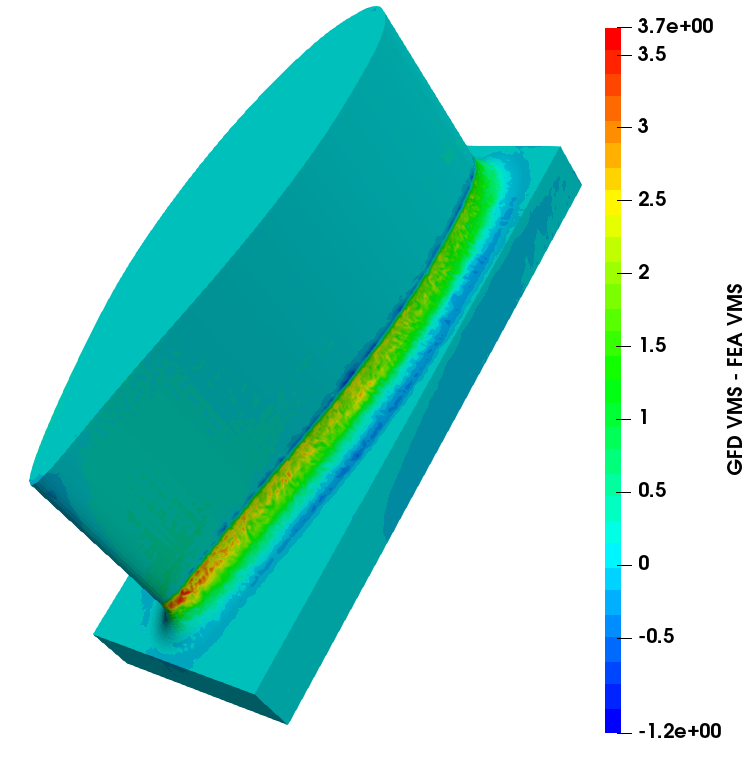}}
		\subfloat[]{\includegraphics[scale=0.4]{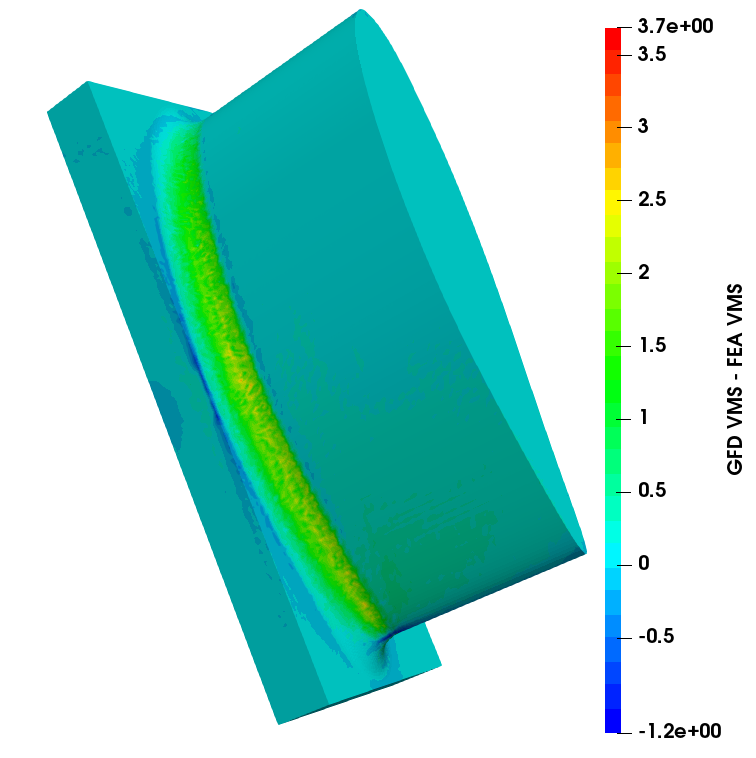}}
		\caption{Simplified high pressure blade subjected to a uniform pressure on one face - Difference between von Mises stress results obtained from the GFD method and FEM (484,238 nodes). The stress concentration at the interface between the blade and the support is larger for the GFD model than for the FE model.}
		\label{BladeResultsDifference}
	\end{figure}

	Figure \ref{SmallBlade} shows the von Mises stress results for the GFD and FE methods. As for the flange problem, the difference between the two von Mises stress solutions is presented in Figure \ref{BladeResultsDifference}.
	
	If can be observed from Figure \ref{SmallBlade} and Figure \ref{BladeResultsDifference} that the stress concentration in the zone between the blade and the support is larger for the GFD model.
	
	\paragraph{Horseshoe Model} \
	
\begin{figure}[H] 
		\centering
		\begin{tikzpicture}
		\def\svgwidth{8cm}
		\node at (0,0) {\includegraphics{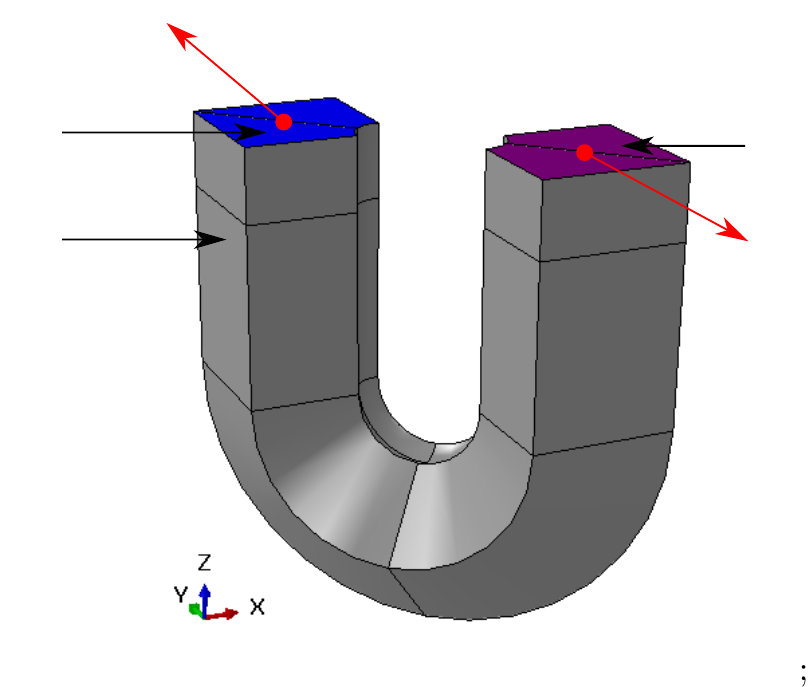}};
		\node[color=black] at (-3.7,2.1) [left] {Top Left Plane};
		\node[color=black] at (-3.7,1.05) [left] {External Surfaces};
		\node[color=black] at (3.7,2) [right] {Top Right Plane};
		
		\end{tikzpicture}
		\caption{Horseshoe model and boundary conditions}
		\label{HorseshoeBCs}
	\end{figure}
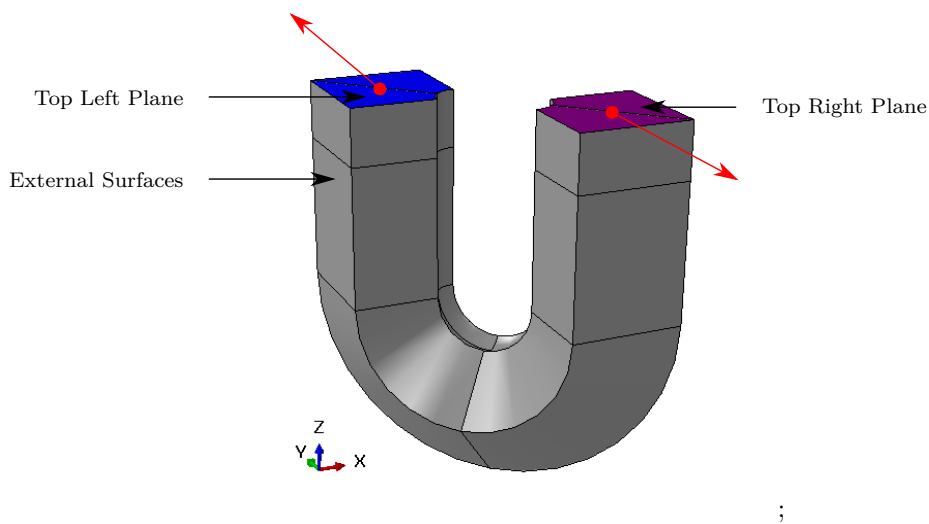

	In 2005, the horseshoe model was solved by Hughes et al \cite{Hughes2005} using the IGA method. This model has been reproduced and is presented in Figure \ref{HorseshoeBCs}. The nodes of the top left plane are fixed in the X and Z directions. The nodes of the top left plane and of the top right plane are, respectively, subjected to a positive and a negative displacement applied in the Y direction. These displacements are equal in absolute value. The external surfaces of the horseshoe are considered stress free. 
	
\begin{figure}[H] 
		\centering
		\subfloat[]{\includegraphics[scale=0.39]{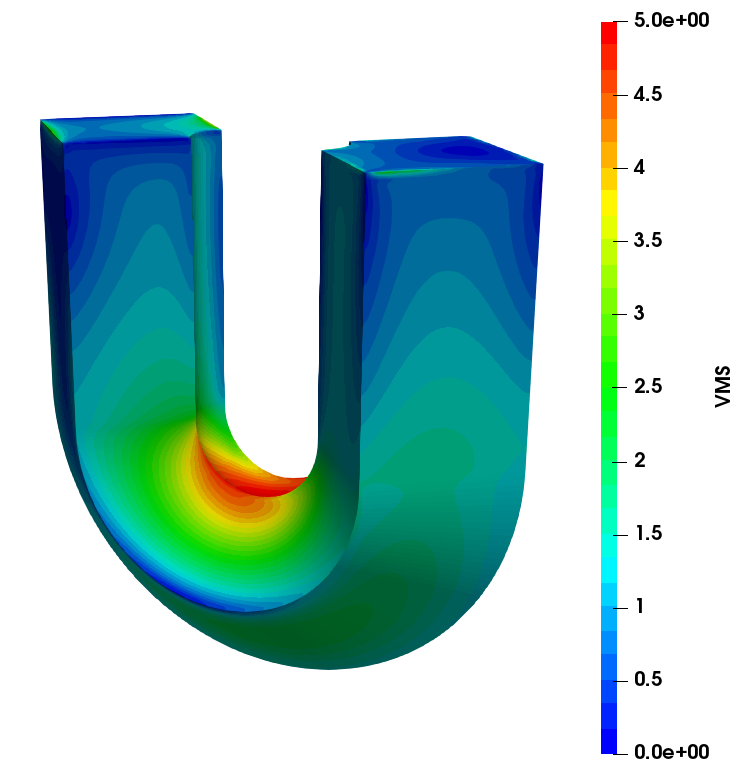}}
		\qquad
		\subfloat[]{\includegraphics[scale=0.39]{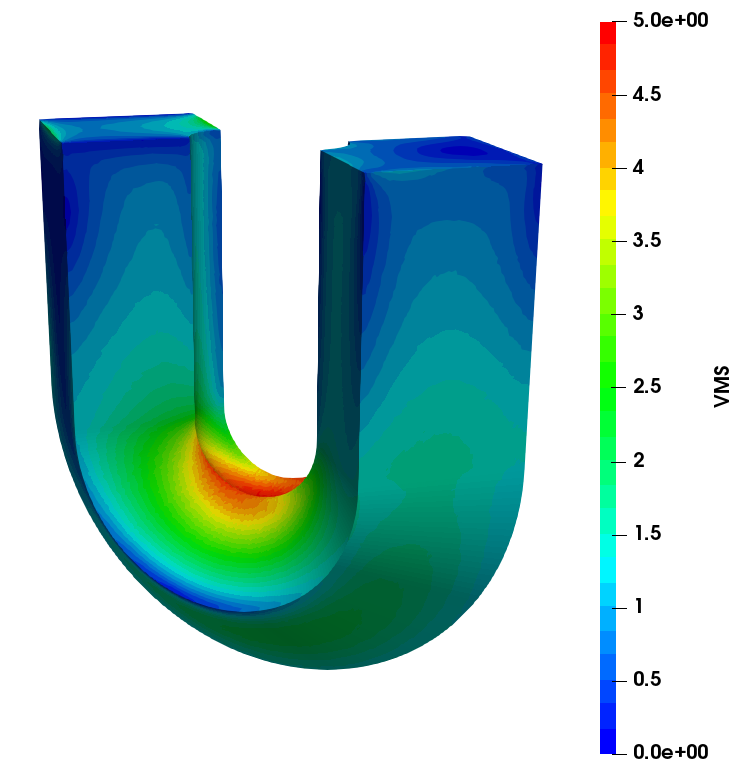}}\\
		\subfloat[]{\includegraphics[scale=0.39]{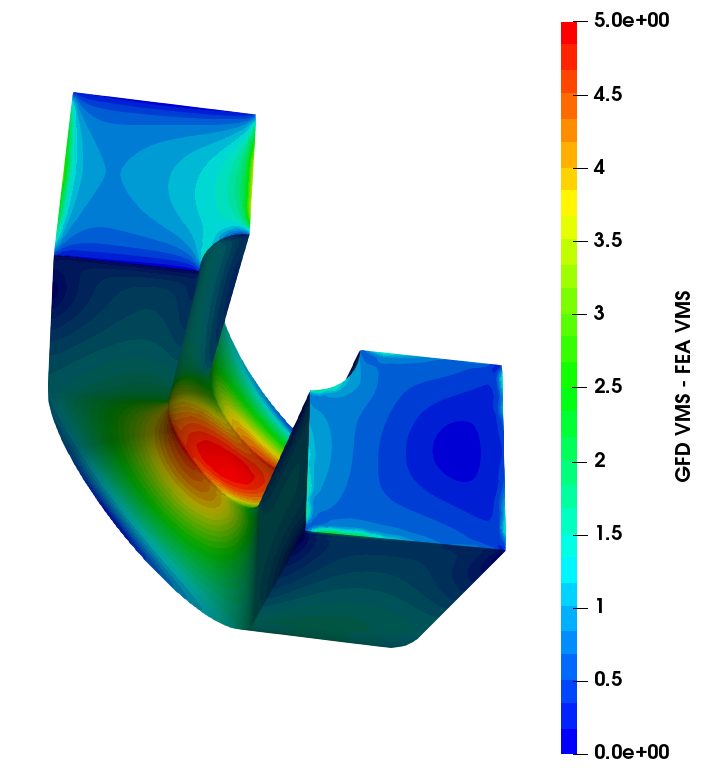}}
		\qquad
		\subfloat[]{\includegraphics[scale=0.39]{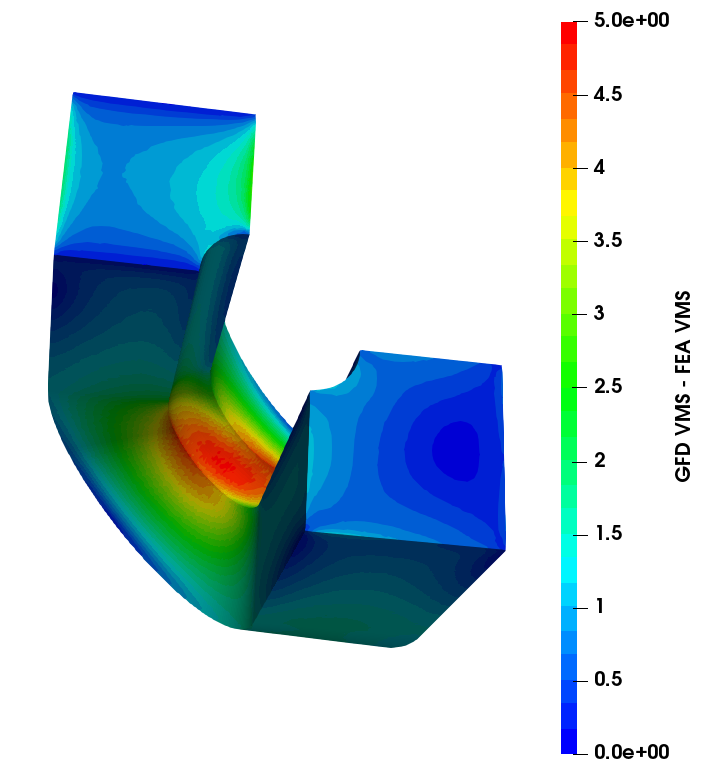}}
		\caption{Horseshoe under shear loading - von Mises stress results from GFD method (a) (c) and FEM (b) (d) (521,326 nodes). The stress concentration in the inner surface of the horeshoe is slightly larger for the GFD method.}
		\label{Horseshoe_CollocFEA}
	\end{figure}

	Figure \ref{Horseshoe_CollocFEA} shows the von Mises stress results for the GFD and FE models. The two figures on the left show two different views of the solution of the problem solved with the GFD method. The two figures on the right show the solution of the problem solved with the FEM. It can be observed form Figure \ref{Horseshoe_CollocFEA} that the stress concentration in the bottom of the horseshoe is larger in the GFD solution. A higher stress concentration is also observed at the edges of the top left plane in the GFD solution.
	
\begin{figure}[H]
		\centering
		\subfloat[]{\includegraphics[scale=0.39]{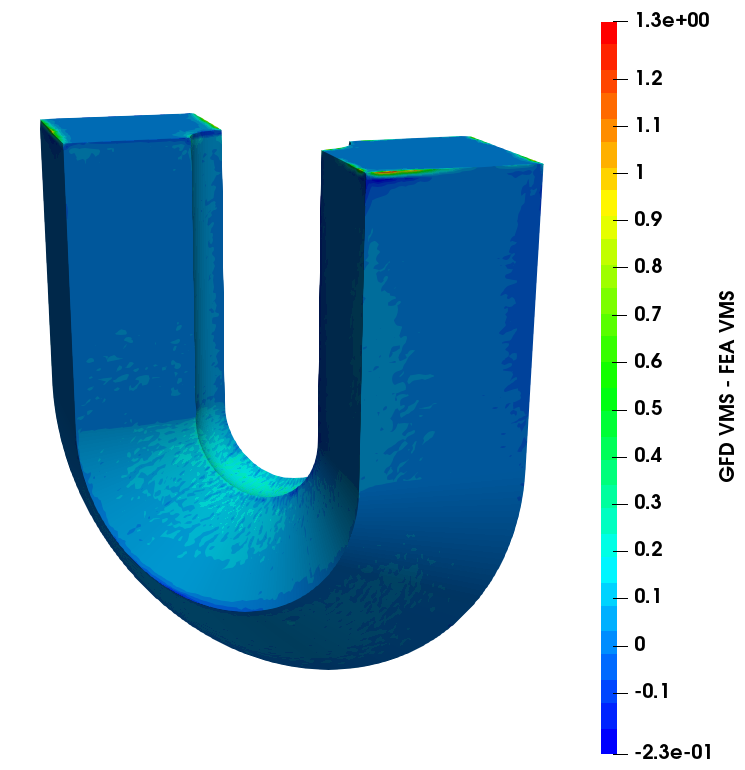}}
		\qquad
		\subfloat[]{\includegraphics[scale=0.39]{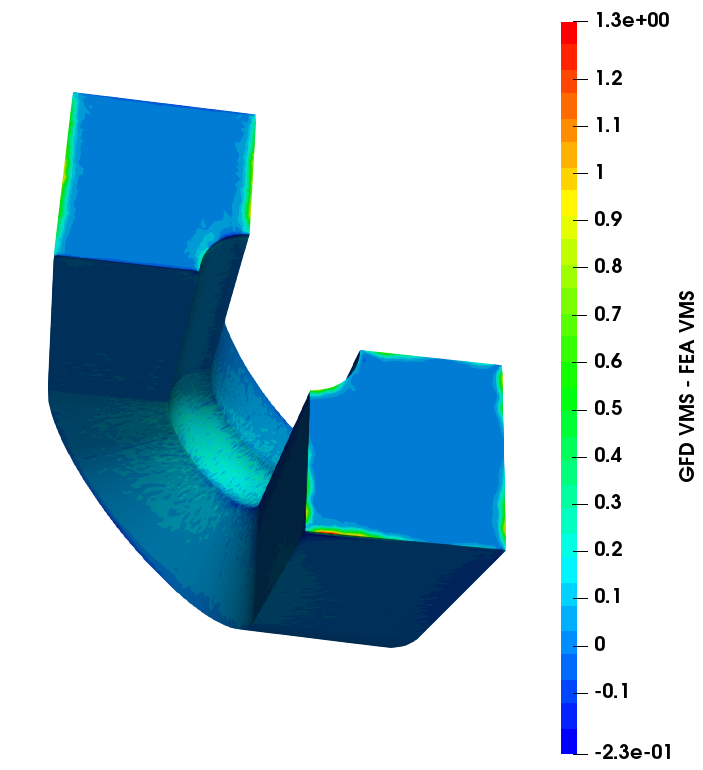}}
		\caption{Horseshoe under shear loading - von Mises stress results from GFD method (a) and FEM (b) (521,326 nodes). The stress concentration in the inner surface of the horseshoe is slightly larger for the GFD method.}
		\label{HorseshoeTop_CollocFEA}
	\end{figure}

	It should be noted that the computation time was lower for the FE method than for the GFD method. This is due to the loading, which creates a singularity at the edges of the top left plane and of the top right plane. The impact of this singularity affects more the GFD model as it is solved by collocation (strong form of the PDE). 
	The difference between the two von Mises stress solutions is presented in Figure \ref{HorseshoeTop_CollocFEA}. This figure confirms that higher stress concentrations are observed for the GFD model at the top planes edges and in the bottom section of the horseshoe.
	
	\paragraph{Fichera's Corner Model} \
	
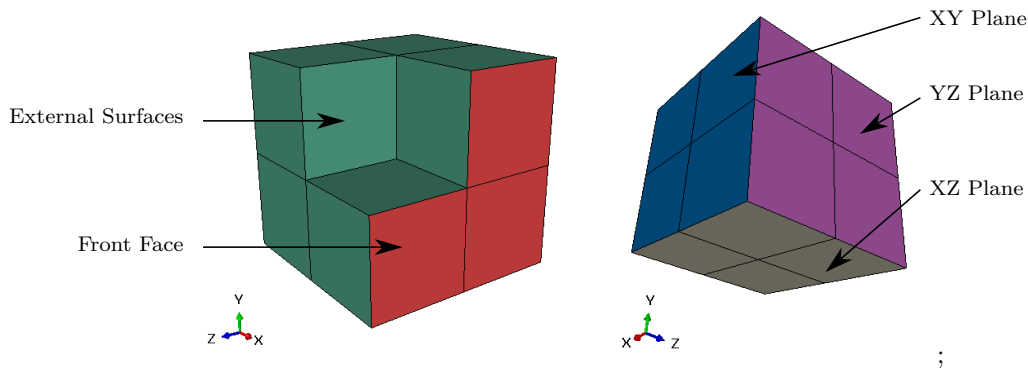
\begin{figure}[H] 
		\centering
		\begin{tikzpicture}
		\def\svgwidth{10cm}
		\node at (0,0) {\includegraphics{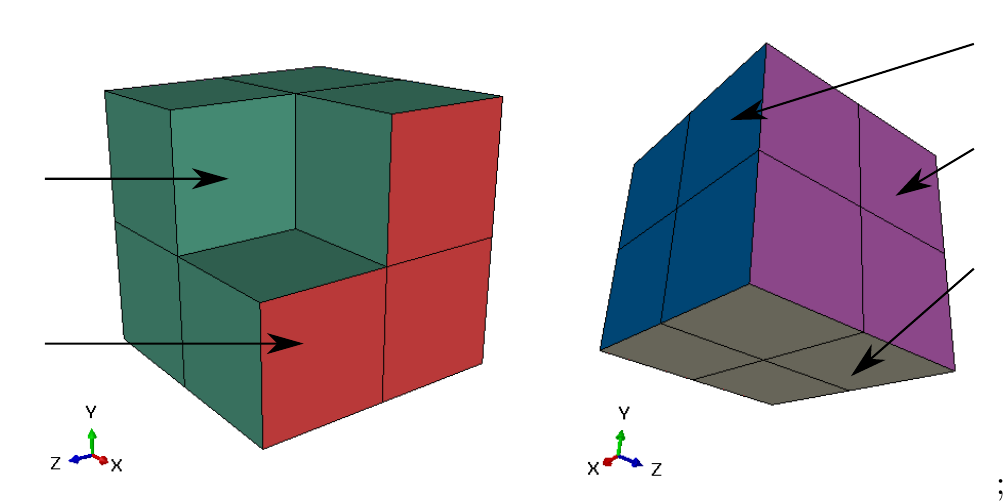}};
		\node[color=black] at (-4.8,0.8) [left] {External Surfaces};
		\node[color=black] at (-4.8,-0.9) [left] {Front Face};
		\node[color=black] at (4.8,2.1) [right] {XY Plane};
		\node[color=black] at (4.8,1.1) [right] {YZ Plane};
		\node[color=black] at (4.8,-0.15) [right] {XZ Plane};
		\end{tikzpicture}
		\caption{Fichera's corner model}
		\label{FicheraBCs}
	\end{figure}

\begin{figure}[H] 
		\centering
		\def\svgwidth{8.5cm}
		\subfloat[]{\includegraphics{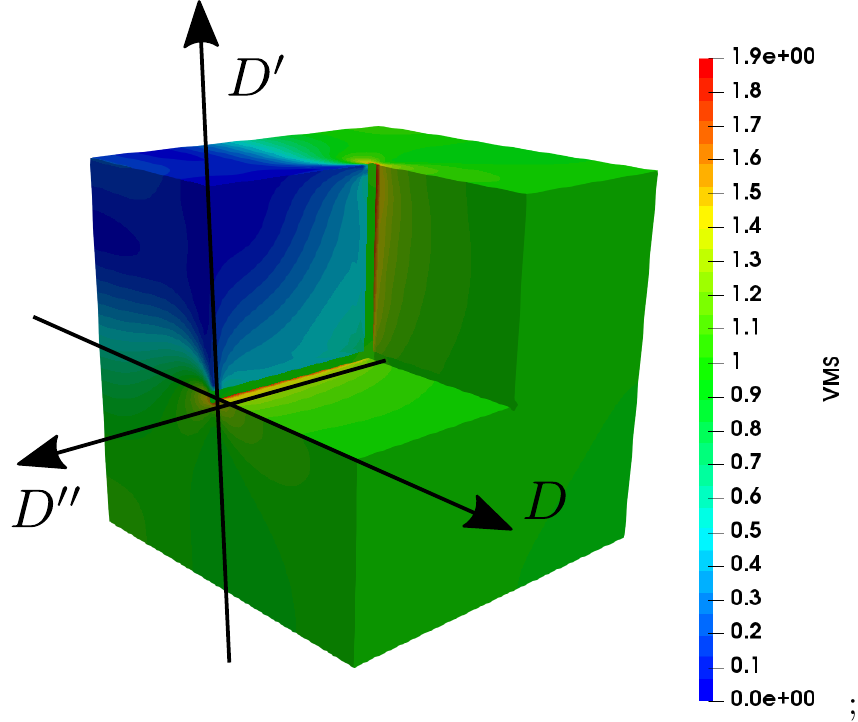}}
		\subfloat[]{\includegraphics[scale=0.327]{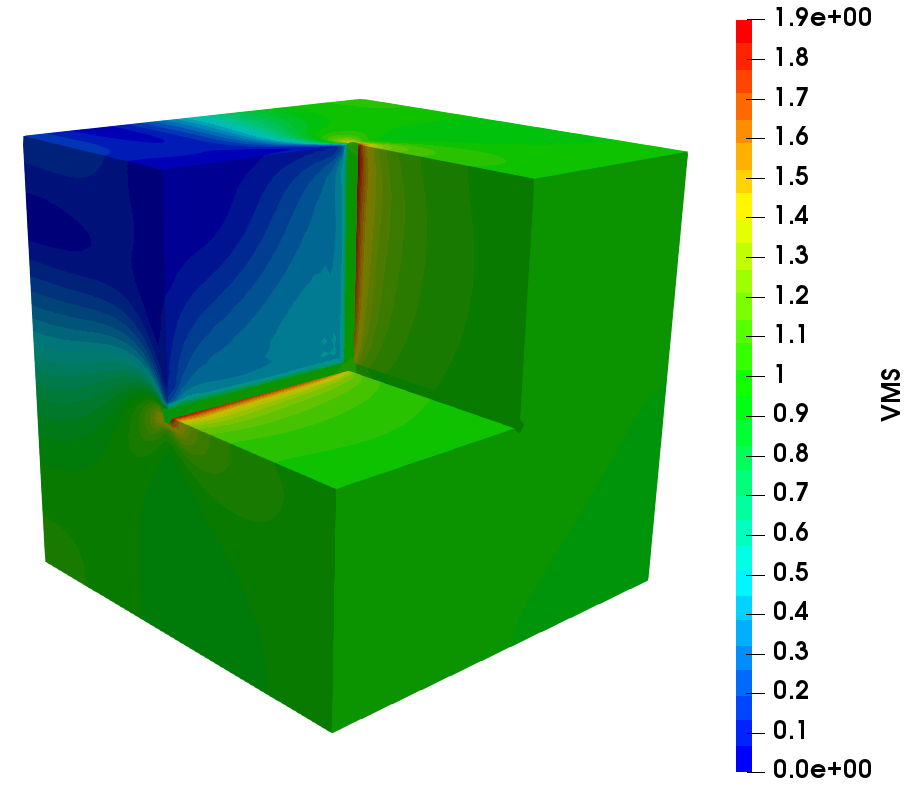}}
		\caption{Fichera corner subjected to a uniform traction on the front side - von Mises stress results from the GFD method (a) and FEM (b) (264,726 nodes).}
		\label{FicheraCorner}
	\end{figure}

\begin{figure}[H] 
		\centering
		\subfloat[]{\includegraphics[scale=0.33]{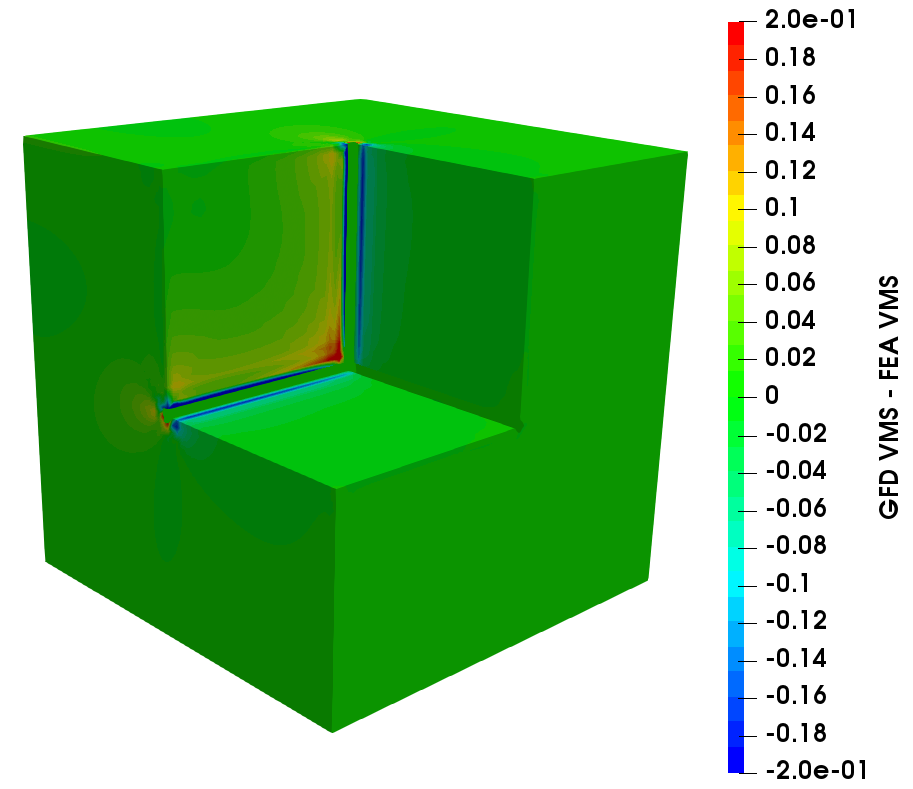}}
		\subfloat[]{\includegraphics[scale=0.33]{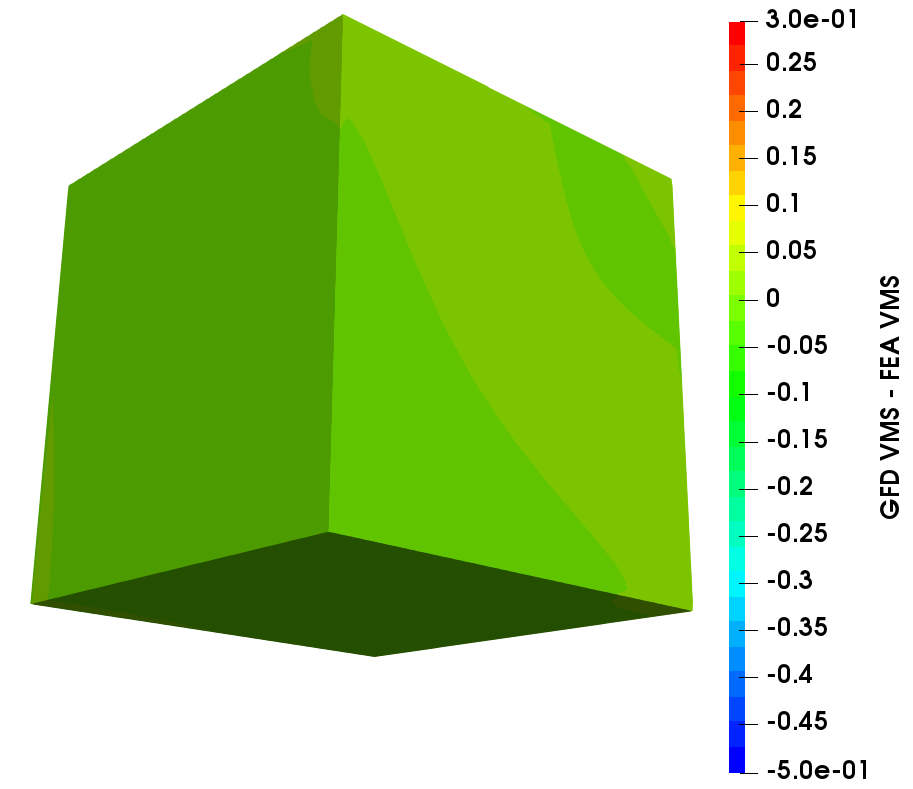}}
		\caption{Fichera corner subjected to a uniform traction on the front side - Difference between von Mises stress results obtained from the GFD method and FEM (264,726 nodes). The stress concentration in the internal corners is larger for the FE model.}
		\label{LShapeResultsDifference}
	\end{figure}

\begin{figure}[H] 
		\centering
		\begin{tabular}{c:c}
			\begin{tikzpicture}
			\begin{axis}[height=7cm,width=7.5cm, ymin=0.4,ymax=2, xmin=-2,xmax=2, legend entries={GFD,FEA},legend style={ at={(0.5,-0.2)},anchor=south west,legend columns=1, cells={anchor=west},  font=\footnotesize, rounded corners=2pt,}, legend pos=north east,xlabel=x coordinate along axis D,ylabel=von Mises Stress]
			\addplot+[green,mark=square*,mark options={fill=green}]   table [x=X1-1, y=X1-1-Colloc
			, col sep=comma] {LShapeCollocFEAComp.csv};
			\addplot+[red,mark=diamond*,mark options={fill=red}]   table [x=X1-1, y=X1-1-FEA
			, col sep=comma] {LShapeCollocFEAComp.csv};
			\addplot+[green,mark=square*,mark options={fill=green}]   table [x=X1-2, y=X1-2-Colloc
			, col sep=comma] {LShapeCollocFEAComp.csv};
			\addplot+[red,mark=diamond*,mark options={fill=red}]   table [x=X1-2, y=X1-2-FEA
			, col sep=comma] {LShapeCollocFEAComp.csv};
			\end{axis}
			\end{tikzpicture} & 
			
			\begin{tikzpicture}
			\begin{axis}[height=7cm,width=7.5cm, ymin=-0.10,ymax=2, xmin=-2,xmax=2, legend entries={GFD,FEA},legend style={ at={(0.5,-0.2)},anchor=south west,legend columns=1, cells={anchor=west},  font=\footnotesize, rounded corners=2pt,}, legend pos=north east,xlabel=y coordinate along axis D',ylabel=von Mises Stress]
			\addplot+[green,mark=square*,mark options={fill=green}]   table [x=X2-1, y=X2-1-Colloc
			, col sep=comma] {LShapeCollocFEAComp.csv};
			\addplot+[red,mark=diamond*,mark options={fill=red}]   table [x=X2-1, y=X2-1-FEA
			, col sep=comma] {LShapeCollocFEAComp.csv};
			\addplot+[green,mark=square*,mark options={fill=green}]   table [x=X2-2, y=X2-2-Colloc
			, col sep=comma] {LShapeCollocFEAComp.csv};
			\addplot+[red,mark=diamond*,mark options={fill=red}]   table [x=X2-2, y=X2-2-FEA
			, col sep=comma] {LShapeCollocFEAComp.csv};
			\end{axis}
			\end{tikzpicture} \\
			(a) & (b) \\
		\end{tabular}
		\begin{tabular}{c}
			\begin{tikzpicture}
			\begin{axis}[height=7cm,width=7.5cm, ymin=0.7,ymax=2, xmin=-2,xmax=2, legend entries={GFD,FEA},legend style={ at={(0.5,-0.2)},anchor=south west,legend columns=1, cells={anchor=west},  font=\footnotesize, rounded corners=2pt,}, legend pos=north west,xlabel=z coordinate along axis D",ylabel=von Mises Stress]
			\addplot+[green,mark=square*,mark options={fill=green}]   table [x=X3-1, y=X3-1-Colloc
			, col sep=comma] {LShapeCollocFEAComp.csv};
			\addplot+[red,mark=diamond*,mark options={fill=red}]   table [x=X3-1, y=X3-1-FEA
			, col sep=comma] {LShapeCollocFEAComp.csv};
			\addplot+[green,mark=square*,mark options={fill=green}]   table [x=X3-2, y=X3-2-Colloc
			, col sep=comma] {LShapeCollocFEAComp.csv};
			\addplot+[red,mark=diamond*,mark options={fill=red}]   table [x=X3-2, y=X3-2-FEA
			, col sep=comma] {LShapeCollocFEAComp.csv};
			\end{axis}
			\end{tikzpicture}\\	
			(c) \\
		\end{tabular}
		\caption{The von Mises stress results comparison along the axes D, D' and D", as presented in Figure \ref{FicheraCorner}(a). Depending on the considered axis, either the GFD method or the FE method leads to the maximum observed stress. The largest stress is observed in the subfigure (c) for the FE method.}
		\label{ResultsCompFichera}
	\end{figure}
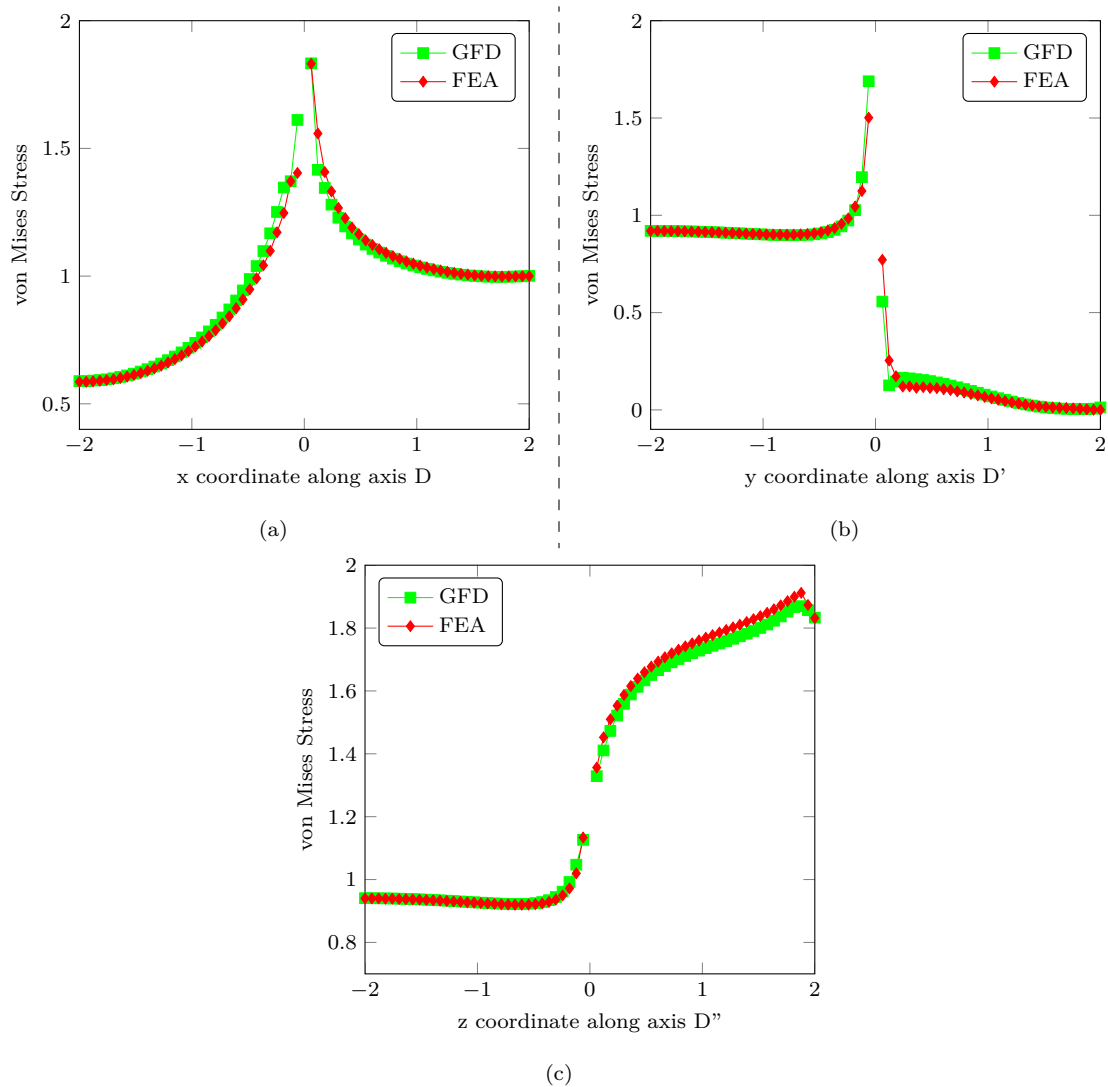

	The Fichera's corner model analyzed in \cite{Dimitrov2001,Rachowicz2006,Zander2016} is presented in Figure \ref{FicheraBCs}. The characteristic planes, on which boundary conditions have been applied, are highlighted and labeled in this figure. The nodes in the planes YZ, XZ and XY are, respectively, fixed in the X, Y and Z directions. A uniform traction is applied on the front face of the truncated cube. The rest of the surfaces are considered stress free. When solved using the GFD method, the internal corner nodes have not been included in the Fichera's corner model. These nodes lead to the divergence of the solution as the stress is infinite at these locations. The visibility criterion presented in Section \ref{SupNodeSelectionGle} is applied to this problem.
	Figure \ref{FicheraCorner} shows the von Mises stress results for the GFD and FE methods. As for the three other problems, the difference between the two von Mises stress solutions is presented in Figure \ref{LShapeResultsDifference}.
	It can be observed from Figure \ref{FicheraCorner} and Figure \ref{LShapeResultsDifference} that the stress concentration near the internal edges is larger for the FE model than for the GFD model. The GFD model leads to larger results only in the center of the corner. In order to visualize more precisely the results of this analysis, von Mises stress results are plotted in Figure \ref{ResultsCompFichera} along the axes D, D' and D" presented in Figure \ref{FicheraCorner}(a). The axes D and D' follow the edge of the model, while the axis D" is slightly offset from the internal corner as results are not available at the corner nodes for the GFD method. The truncated cube has an edge length of 4. The coordinate of the internal corner is (0,0,0).
	Figure \ref{ResultsCompFichera} shows that, depending on the considered axis, either the GFD method or the FE method leads to the largest results. For this model, the FE method is expected to give a higher stress concentration as the corner nodes have not been included in the GFD method but have been included in the  FE model.
	
	\paragraph{Discussion} \
	
	From the figures presented in this section for the various problems considered, it can be observed that the results obtained with the GFD method are very close to the results obtained with the FEM. The stress concentrations are slightly larger for the flange, the blade and the horseshoe problems when the GFD method is used. This might be due to the use of the strong form of the equations, which allows solving the loading equations on the boundary of the domain. In the FEM, the problem is solved in a weak form using an integration over the domain. For the Fichera's corner problem, the largest von Mises stress concentration is observed for the FE method. This is due to the internal corner nodes which have not been included in the GFD method. The FEM is thus deemed more accurate.
	
	\section{Conclusions} \label{Conclusions}
	
	The aim of the paper was three-fold:
	
	\begin{description}
		\item[Brief review and primer] We briefly reviewed Taylor-series expansion based collocation/meshfree methods. Our aim here was not to be exhaustive, but to cover the main material available, to our knowledge. We also presented detailed derivation of the system matrices, in order to facilitate the entry of  newcomers into the field and attempted to unify the generalised finite difference method and the discrete correction particle strength exchange method under one umbrella. 
		\item[Performance benchmarking] We provided a detailed benchmarking strategy for Taylor-series expansion based collocation methods as well as all data files including all input files, boundary conditions, point distribution and solution fields, so as to facilitate future benchmarking of new methods.
		\item[New methods for non-smooth solutions] We proposed a few improvements to the original methods, both DCPSE and GFD, in order to treat problem with non-smooth solutions, including discontinuities, singularities or sharp gradients.
	\end{description}
	
	We noted that the various parameters involved in the methods have a significant impact on the solution, and that they should therefore be carefully chosen. In itself, this is a drawback compared to more parameter-robust methods, in particular the finite element method. Another main conclusion of this work is that common approaches used in practice to improve collocation methods must be used with caution as they do not always lead to the reduction of the overall error. We observed the following:
	
	\begin{enumerate}[(1)]
		\item For the GFD method, the weight function based on the 4$^{\text{th}}$ order spline leads to the minimum error for problems with a polynomial solution such as the pressurized cylinder. For singular problems, such as the L-shape in mode I loading, both linear and 4$^{\text{th}}$ order spline weight functions lead to a minimum error.
		\item For the DC PSE method, the weight function based on the exponential functions leads to the minimum error for both polynomial and singular problems.
		\item For the problems with a polynomial solution, a polynomial correction function basis leads to an error approximately fifteen time lower than with an exponential basis function. For the singular problems, an exponential correction function basis lead to an error approximately 5\% lower than with an exponential basis function. A polynomial correction function is recommended for most problems as the solution type in not known a priori.
		\item Increased size of the node supports on the boundary helps decreasing the overall error, while increasing only slightly the number of non zero elements in the system matrix. For the polynomial problem considered, a reduction of a factor one hundred is observed between boundary support nodes composed of thirteen nodes and eighteen nodes. For the singular problem, no significant error reduction is observed.
		\item Voronoi diagrams can be used to give additional information to the collocation methods on the spatial arrangement of the nodes over the domain.
		\begin{enumerate}[(a)]
			\item For the GFD method, Voronoi diagrams allow the selection of weights which depend on the node placement over the collocation node support.
			\item For the DC PSE method, Voronoi diagrams are expected to improve the accuracy of the convolution, but they may also lead to an increased error for some node distributions.
		\end{enumerate}	
		The use of Voronoi diagrams helps in reducing the error for the 2D cylinder problem with a free node distribution (based on Delaunay triangulation). A reduction of up to 17\% is observed for the GFD method and of up to 10\% for the DC PSE methods. For a regular node distribution an error increase of 3\% is observed for the considered problems when Voronoi diagrams are used. For the L-shape problem, the use of Voronoi diagrams has no significant impact on the error. It can be concluded that the use of Voronoi diagrams does not allow a significant error reduction for the considered node arrangements, and their use is not recommended in the general case.
		\item The stabilization method reduces the error for the L-shape problem by respectively 25\% and 35\% for the GFD and DC PSE methods. A large error increase (up to a factor 30) is observed for both methods for the pressurized cylinder problem. This difference is due to the type of boundary conditions imposed. The stabilization method is more suitable to Dirichlet loaded problems than to Neumann loaded problems.
		\item For problems with singularities and concave geometries, the visibility criterion improves the convergence of the solution when solved with iterative solvers. Moreover, it allows to significantly reduce the observed error. A reduction of up to 50\% and 55\% are respectively observed for the GFD and DC PSE methods when the visibility criterion is used. The use of this criterion for support node selection is recommended for all singular and concave problems.
		\item Compared to other typical collocation methods (e.g., MLS, IMLS, RBF-FD), the GFD and DC PSE methods have shown good performance both in terms of observed error and computation time. The results obtained with the collocation methods are close to the results obtained using the FEM and more accurate for some problems. For large 3D problems, the GFD method leads to very similar results as those obtained using FEM. 
		\item A slightly larger stress concentration has been observed for the flange, the blade and the horseshoe problems when solved with the GFD method compared to the results from FEA. For the Fichera's corner problem, the stress concentration obtained with the FEM is slightly more pronounced than that obtained using GFD method. This is due to the fact that the corner nodes have not been included in the GFD method, and thereby, the FEM results represent more accurately the actual solution.
	\end{enumerate}
	
	To summarize, we have proposed for the GFD and DC PSE methods a set of optimal parameters that can be used to readily solve any linear elastic problem. We also showed that point collocation methods may be used effectively for problems with singularities and 3D problems of industrial size. Using the visibility criterion for concave and singular problems improves the convergence of the methods and leads to a significant error reduction. A logical next step in collocation methods is to investigate the use of enriched weight functions and enriched stencils near singularities in order to improve the results obtained in regions of rapid field change. Similarly, a posteriori error estimation driven local refinement, vastly simplified in collocation methods, should be investigated, which is the topic of ongoing work in our teams. Finally, a massively parallel approach, if possible based on graphical processing units, should be investigated to accelerate the solution scheme. 
	
	\begin{acknowledgements}
		St\'ephane P.A. Bordas and Satyendra Tomar thank partial funding for their time provided by the European Research Council Starting Independent Research Grant (ERC Stg grant agreement No. 279578) RealTCut ``Towards real time multiscale simulation of cutting in non-linear materials with applications to surgical simulation and computer guided surgery". The authors are also grateful for the funding from the Luxembourg National Research Fund (INTER/FWO/15/10318764). This is a pre-print of an article published in Archives of Computational Methods in Engineering. The final authenticated version is available online at: https://doi.org/10.1007/s11831-019-09357-5.
	\end{acknowledgements}
	
	\nocite{*}
	\bibliographystyle{unsrt}
	\bibliography{Bibliography.bib}

\begin{thebibliography}{100}

\bibitem{Runge1908}
C.~Runge.
\newblock {\em Z. Math. u. Physik}, 50:255, 1908.

\bibitem{Ritz1908}
W.~Ritz.
\newblock {\"{U}}ber eine neue methode zur l{\"{o}}sung gewisser
  variationsprobleme der mathematischen physik.
\newblock {\em Journal f{\"{u}}r die Reine und Angewandte Mathematik},
  135:1--61, 1908.

\bibitem{Galerkin1915}
B.~G. Galerkin.
\newblock Rods and plates. series occurring in various questions concerning the
  elastic equilibrium of rods and plates.
\newblock {\em Vestnik Inzh.}, 19:897--908, 1915.

\bibitem{Monaghan1992}
J.J. Monaghan.
\newblock Smoothed particle hydrodynamics.
\newblock {\em Annual Review of Astronomy and Astrophysics}, 30(1):543--574,
  sep 1992.

\bibitem{Belytschko1994}
T.~Belytschko, Y.Y. Lu, and L.~Gu.
\newblock Element-free galerkin methods.
\newblock {\em International Journal for Numerical Methods in Engineering},
  37(2):229--256, jan 1994.

\bibitem{Liu1995}
W.~K. Liu, S.~Jun, and Y.~F. Zhang.
\newblock Reproducing kernel particle methods.
\newblock {\em International Journal for Numerical Methods in Fluids},
  20(8-9):1081--1106, apr 1995.

\bibitem{Duarte1996}
C.~Armando Duarte and J.~Tinsley Oden.
\newblock H-p clouds{\textemdash}anh-p meshless method.
\newblock {\em Numerical Methods for Partial Differential Equations},
  12(6):673--705, nov 1996.

\bibitem{Atluri1998}
S.~N. Atluri and T.~Zhu.
\newblock A new meshless local petrov-galerkin ({MLPG}) approach in
  computational mechanics.
\newblock {\em Computational Mechanics}, 22(2):117--127, aug 1998.

\bibitem{De2000}
S.~De and K.~J. Bathe.
\newblock The method of finite spheres.
\newblock {\em Computational Mechanics}, 25(4):329--345, apr 2000.

\bibitem{Chen2000}
J.S. Chen, C.T. Wu, S.~Yoon, and Y.~You.
\newblock A stabilized conforming nodal integration for galerkin mesh-free
  methods.
\newblock {\em International Journal for Numerical Methods in Engineering},
  50(2):435--466, 2000.

\bibitem{Nguyen2008}
V.P. Nguyen, T.~Rabczuk, S.~Bordas, and M.~Duflot.
\newblock Meshless methods: A review and computer implementation aspects.
\newblock {\em Mathematics and Computers in Simulation}, 79(3):763--813, dec
  2008.

\bibitem{Babuska1995}
I.~Babu{\v{s}}ka and J.~M. Melenk.
\newblock The partition of unity finite element method.
\newblock Technical report, 1995.

\bibitem{Babuska1997}
I.~Babu{\v{s}}ka and J.~M. Melenk.
\newblock The partition of unity mehtod.
\newblock {\em International Journal for Numerical Methods in Engineering},
  40(4):727--758, feb 1997.

\bibitem{Strouboulis2001}
T.~Strouboulis, K.~Copps, and I.~Babu{\v{s}}ka.
\newblock The generalized finite element method.
\newblock {\em Computer Methods in Applied Mechanics and Engineering},
  190(32-33):4081--4193, may 2001.

\bibitem{Mos1999}
N.~Mo{\"{e}}s, J.~Dolbow, and T.~Belytschko.
\newblock A finite element method for crack growth without remeshing.
\newblock {\em International Journal for Numerical Methods in Engineering},
  46(1):131--150, sep 1999.

\bibitem{Sukumar2000}
N.~Sukumar, N.~Mo{\"{e}}s, B.~Moran, and T.~Belytschko.
\newblock Extended finite element method for three-dimensional crack modelling.
\newblock {\em International Journal for Numerical Methods in Engineering},
  48(11):1549--1570, 2000.

\bibitem{Dolbow2000}
J.~Dolbow, N.~Mo{\"{e}}s, and T.~Belytschko.
\newblock Modeling fracture in mindlin{\textendash}reissner plates with the
  extended finite element method.
\newblock {\em International Journal of Solids and Structures},
  37(48-50):7161--7183, nov 2000.

\bibitem{Dolbow2001}
J.~Dolbow, N.~Mo{\"{e}}s, and T.~Belytschko.
\newblock An extended finite element method for modeling crack growth with
  frictional contact.
\newblock {\em Computer Methods in Applied Mechanics and Engineering},
  190(51-52):6825--6846, oct 2001.

\bibitem{Sukumar2001}
N.~Sukumar, D.L. Chopp, N.~Mo{\"{e}}s, and T.~Belytschko.
\newblock Modeling holes and inclusions by level sets in the extended
  finite-element method.
\newblock {\em Computer Methods in Applied Mechanics and Engineering},
  190(46-47):6183--6200, sep 2001.

\bibitem{Mos2002}
N.~Mo{\"{e}}s and T.~Belytschko.
\newblock Extended finite element method for cohesive crack growth.
\newblock {\em Engineering Fracture Mechanics}, 69(7):813--833, may 2002.

\bibitem{Ji2004}
H.~Ji and J.~Dolbow.
\newblock On strategies for enforcing interfacial constraints and evaluating
  jump conditions with the extended finite element method.
\newblock {\em International Journal for Numerical Methods in Engineering},
  61(14):2508--2535, 2004.

\bibitem{Duflot2008}
M.~Duflot and S.~Bordas.
\newblock A posteriorierror estimation for extended finite elements by an
  extended global recovery.
\newblock {\em International Journal for Numerical Methods in Engineering},
  76(8):1123--1138, nov 2008.

\bibitem{Rabczuk2007}
T.~Rabczuk, T.~Belytschko, S.~Bordas, and G.~Zi.
\newblock Enriched meshfree methods for crack problems.
\newblock {\em 9th National Congress on Computational Mechanics}, jun 2007.

\bibitem{Rabczuk2007Sec}
T.~Rabczuk, S.~Bordas, and G.~Zi.
\newblock Initiation, nucleation and propagation of cracks in a cohesive way
  without mesh enriched in the quasi-sensitive materials: large strains,
  quasi-static and dynamic, 2007.

\bibitem{Rabczuk2007Thi}
T.~Rabczuk, G.~Zi, and S.~Bordas.
\newblock Enriched finite element and meshfree methods for dynamic crack
  propagation problems.
\newblock {\em 5th Australasian Congress on Applied Mechanics, ACAM 2007},
  2007.

\bibitem{Bordas2007}
S.~Bordas, G.~Zi, and T.~Rabczuk.
\newblock Three-dimensional non-linear fracture mechanics by enriched meshfree
  methods without asymptotic enrichment.
\newblock In {\em {IUTAM} Symposium on Discretization Methods for Evolving
  Discontinuities}, pages 21--36. Springer Netherlands, 2007.

\bibitem{Bordas2008}
S.~Bordas, T.~Rabczuk, and G.~Zi.
\newblock Three-dimensional crack initiation, propagation, branching and
  junction in non-linear materials by an extended meshfree method without
  asymptotic enrichment.
\newblock {\em Engineering Fracture Mechanics}, 75(5):943--960, mar 2008.

\bibitem{Talebi2011}
H.~Talebi, C.~Samaniego, E.~Samaniego, and T.~Rabczuk.
\newblock On the numerical stability and mass-lumping schemes for explicit
  enriched meshfree methods.
\newblock {\em International Journal for Numerical Methods in Engineering},
  89(8):1009--1027, nov 2011.

\bibitem{Natarajan2011}
S.~Natarajan, P.~Kerfriden, S.~Bordas, D.R. Mahapatra, and T.~Rabczuk.
\newblock Enriched element free galerkin method for gradient elasticity.
\newblock {\em XFEM 2011}, jun 2011.

\bibitem{Belytschko2002}
T.~Belytschko, C.~Parimi, N.~Mo{\"{e}}s, N.~Sukumar, and S.~Usui.
\newblock Structured extended finite element methods for solids defined by
  implicit surfaces.
\newblock {\em International Journal for Numerical Methods in Engineering},
  56(4):609--635, nov 2002.

\bibitem{Moumnassi2014}
M.~Moumnassi, S.~Bordas, R.~Figueredo, and P.~Sansen.
\newblock Analysis using higher-order {XFEM}: implicit representation of
  geometrical features from a given parametric representation.
\newblock {\em Mechanics {\&} Industry}, 15(5):443--448, 2014.

\bibitem{Rabczuk2010}
T.~Rabczuk, S.~Bordas, and Goangseup Zi.
\newblock On three-dimensional modeling of crack growth using partition of
  unity methods.
\newblock {\em Computers {\&} Structures}, 88(23-24):1391--1411, dec 2010.

\bibitem{Hughes2005}
T.J.R. Hughes, J.A. Cottrell, and Y.~Bazilevs.
\newblock Isogeometric analysis: {CAD}, finite elements, {NURBS}, exact
  geometry and mesh refinement.
\newblock {\em Computer Methods in Applied Mechanics and Engineering},
  194(39-41):4135--4195, oct 2005.

\bibitem{Simpson2012}
R.N. Simpson, S.~Bordas, J.~Trevelyan, and T.~Rabczuk.
\newblock A two-dimensional isogeometric boundary element method for
  elastostatic analysis.
\newblock {\em Computer Methods in Applied Mechanics and Engineering},
  209-212:87--100, feb 2012.

\bibitem{Simpson2013}
R.N. Simpson, S.~Bordas, H.~Lian, and J.~Trevelyan.
\newblock An isogeometric boundary element method for elastostatic analysis: 2d
  implementation aspects.
\newblock {\em Computers {\&} Structures}, 118:2--12, mar 2013.

\bibitem{Scott2013}
M.A. Scott, R.N. Simpson, J.A. Evans, S.~Lipton, S.~Bordas, T.J.R. Hughes, and
  T.W. Sederberg.
\newblock Isogeometric boundary element analysis using unstructured
  {T}-splines.
\newblock {\em Computer Methods in Applied Mechanics and Engineering},
  254:197--221, feb 2013.

\bibitem{Lian2013}
H.~Lian, R.N. Simpson, and S.~Bordas.
\newblock Stress analysis without meshing: isogeometric boundary-element
  method.
\newblock {\em Proceedings of the Institution of Civil Engineers - Engineering
  and Computational Mechanics}, 166(2):88--99, jun 2013.

\bibitem{Peng2014}
X.~Peng, E.~Atroshchenko, and S.~Bordas.
\newblock Damage tolerance assessment directly from cad: (extended)
  isogeometric boundary element methods (xigabem).
\newblock {\em 6th International Conference on Advanced Computational Methods
  in Engineering}, 2014.

\bibitem{Atroshchenko2015}
E.~Atroshchenko and S.~Bordas.
\newblock Fundamental solutions and dual boundary element methods for fracture
  in plane cosserat elasticity.
\newblock {\em Proceedings of the Royal Society A: Mathematical, Physical and
  Engineering Science}, 471(2179):20150216, jul 2015.

\bibitem{Lian2016}
H.~Lian, P.~Kerfriden, and S.~Bordas.
\newblock Implementation of regularized isogeometric boundary element methods
  for gradient-based shape optimization in two-dimensional linear elasticity.
\newblock {\em International Journal for Numerical Methods in Engineering},
  106(12):972--1017, apr 2016.

\bibitem{Peng2017}
X.~Peng, E.~Atroshchenko, P.~Kerfriden, and S.~Bordas.
\newblock Isogeometric boundary element methods for three dimensional static
  fracture and fatigue crack growth.
\newblock {\em Computer Methods in Applied Mechanics and Engineering},
  316:151--185, apr 2017.

\bibitem{Lian2017}
H.~Lian, P.~Kerfriden, and S.~Bordas.
\newblock Shape optimization directly from {CAD}: An isogeometric boundary
  element approach using t-splines.
\newblock {\em Computer Methods in Applied Mechanics and Engineering},
  317:1--41, apr 2017.

\bibitem{Atroshchenko2017}
E.~Atroshchenko, J.S. Hale, Javier~A. Videla, S.~Potapenko, and S.~Bordas.
\newblock Micro-structured materials: Inhomogeneities and imperfect interfaces
  in plane micropolar elasticity, a boundary element approach.
\newblock {\em Engineering Analysis with Boundary Elements}, 83:195--203, oct
  2017.

\bibitem{Xu2011}
G.~Xu, B.~Mourrain, R.~Duvigneau, and A.~Galligo.
\newblock Parameterization of computational domain in isogeometric analysis:
  Methods and comparison.
\newblock {\em Computer Methods in Applied Mechanics and Engineering},
  200(23-24):2021--2031, jun 2011.

\bibitem{Xu2013}
G.~Xu, B.~Mourrain, R.~Duvigneau, and A.~Galligo.
\newblock Analysis-suitable volume parameterization of multi-block
  computational domain in isogeometric applications.
\newblock {\em Computer-Aided Design}, 45(2):395--404, feb 2013.

\bibitem{Xu2018}
G.~Xu, M.~Li, B.~Mourrain, T.~Rabczuk, J.~Xu, and S.~Bordas.
\newblock Constructing {IGA}-suitable planar parameterization from complex
  {CAD} boundary by domain partition and global/local optimization.
\newblock {\em Computer Methods in Applied Mechanics and Engineering},
  328:175--200, jan 2018.

\bibitem{Nguyen2015}
V.P. Nguyen and S.~Bordas.
\newblock Extended isogeometric analysis for strong and weak discontinuities.
\newblock In {\em Isogeometric Methods for Numerical Simulation}, pages
  21--120. Springer Vienna, 2015.

\bibitem{NguyenImplementation2015}
V.P. Nguyen, C.~Anitescu, S.~Bordas, and T.~Rabczuk.
\newblock Isogeometric analysis: An overview and computer implementation
  aspects.
\newblock {\em Mathematics and Computers in Simulation}, 117:89--116, nov 2015.

\bibitem{Atroshchenko2018}
E.~Atroshchenko, S.~Tomar, G.~Xu, and S.~Bordas.
\newblock Weakening the tight coupling between geometry and simulation in
  isogeometric analysis: From sub- and super-geometric analysis to
  geometry-independent field {approximaTion} ({GIFT}).
\newblock {\em International Journal for Numerical Methods in Engineering},
  114(10):1131--1159, mar 2018.

\bibitem{Burman2010}
E.~Burman and P.~Hansbo.
\newblock Fictitious domain finite element methods using cut elements: {I}. a
  stabilized lagrange multiplier method.
\newblock {\em Computer Methods in Applied Mechanics and Engineering},
  199(41-44):2680--2686, oct 2010.

\bibitem{Burman2012}
E.~Burman and P.~Hansbo.
\newblock Fictitious domain finite element methods using cut elements: {II}. a
  stabilized nitsche method.
\newblock {\em Applied Numerical Mathematics}, 62(4):328--341, apr 2012.

\bibitem{Burman2014}
E.~Burman and P.~Hansbo.
\newblock Fictitious domain methods using cut elements: {III}. a stabilized
  nitsche method for stokes' problem.
\newblock {\em {ESAIM}: Mathematical Modelling and Numerical Analysis},
  48(3):859--874, apr 2014.

\bibitem{Burman2014_2}
E.~Burman, S.~Claus, P.~Hansbo, M.G. Larson, and A.~Massing.
\newblock {CutFEM}: Discretizing geometry and partial differential equations.
\newblock {\em International Journal for Numerical Methods in Engineering},
  104(7):472--501, dec 2014.

\bibitem{Hansbo2014}
P.~Hansbo, M.G. Larson, and S.~Zahedi.
\newblock A cut finite element method for a stokes interface problem.
\newblock {\em Applied Numerical Mathematics}, 85:90--114, nov 2014.

\bibitem{Burman2015}
E.~Burman, P.~Hansbo, and M.G. Larson.
\newblock A stabilized cut finite element method for partial differential
  equations on surfaces: The laplace{\textendash}beltrami operator.
\newblock {\em Computer Methods in Applied Mechanics and Engineering},
  285:188--207, mar 2015.

\bibitem{Claus2015}
S.~Claus, E.~Burman, and A.~Massing.
\newblock {CutFEM}: a stabilised nitsche {XFEM} method for multi-physics
  problems, 2015.

\bibitem{Claus2017}
S.~Claus and P.~Kerfriden.
\newblock A stable and optimally convergent {LaTIn}-{CutFEM} algorithm for
  multiple unilateral contact problems.
\newblock {\em International Journal for Numerical Methods in Engineering},
  113(6):938--966, oct 2017.

\bibitem{Claus2018}
S.~Claus, S.~Bigot, and P.~Kerfriden.
\newblock {CutFEM} method for stefan--signorini problems with application in
  pulsed laser ablation.
\newblock {\em {SIAM} Journal on Scientific Computing}, 40(5):B1444--B1469, jan
  2018.

\bibitem{BordasUnfitted2017}
S.~Bordas, E.~Burman, M.G. Larson, and M.A. Olshanskii, editors.
\newblock {\em Geometrically Unfitted Finite Element Methods and Applications}.
\newblock Springer International Publishing, 2017.

\bibitem{Lancaster1981}
P.~Lancaster and K.~Salkauskas.
\newblock Surfaces generated by moving least squares methods.
\newblock {\em Mathematics of Computation}, 37(155):141--141, sep 1981.

\bibitem{Shepard1968}
D.~Shepard.
\newblock A two-dimensional interpolation function for irregularly-spaced data.
\newblock In {\em Proceedings of the 1968 23rd {ACM} national conference}.
  {ACM} Press, 1968.

\bibitem{Macneal953}
R.H. Macneal.
\newblock An asymmetrical finite difference network.
\newblock {\em Quarterly of Applied Mathematics}, 11(3):295--310, 1953.

\bibitem{Forsythe1960}
G.E. Forsythe and W.R. Wasow.
\newblock {\em Finite Difference Methods for Partial Differential Equations}.
\newblock Wiley, 1960.

\bibitem{Jensen1972}
P.S. Jensen.
\newblock Finite difference techniques for variable grids.
\newblock {\em Computers {\&} Structures}, 2(1-2):17--29, feb 1972.

\bibitem{Liszka1980}
T.~Liszka and J.~Orkisz.
\newblock The finite difference method at arbitrary irregular grids and its
  application in applied mechanics.
\newblock {\em Computers {\&} Structures}, 11(1-2):83--95, feb 1980.

\bibitem{Orkisz1998}
J.~Orkisz.
\newblock Finite difference method ({Part III}).
\newblock In {\em Handbook of Computational Solid Mechanics}, pages 335--432.
  Springer-Verlag, 1998.

\bibitem{Degond1989}
P.~Degond and S.~Mas-Gallic.
\newblock The weighted particle method for convection-diffusion equations. part
  1: The case of an isotropic viscosity.
\newblock {\em Mathematics of Computation}, 53(188):485, oct 1989.

\bibitem{Eldredge2002}
J.D. Eldredge, A.~Leonard, and T.~Colonius.
\newblock A general deterministic treatment of derivatives in particle methods.
\newblock {\em Journal of Computational Physics}, 180(2):686--709, aug 2002.

\bibitem{Schrader2010}
B.~Schrader, S.~Reboux, and I.F. Sbalzarini.
\newblock Discretization correction of general integral {PSE} operators for
  particle methods.
\newblock {\em Journal of Computational Physics}, 229(11):4159--4182, jun 2010.

\bibitem{Abaqus2017}
Dassault Systemes.
\newblock Abaqus 2017.
\newblock 2017.

\bibitem{Schrader2012}
B.~Schrader, S.~Reboux, and I.F. Sbalzarini.
\newblock Choosing the best kernel: Performance models for diffusion operators
  in particle methods.
\newblock {\em {SIAM} Journal on Scientific Computing}, 34(3):A1607--A1634, jan
  2012.

\bibitem{Kansa1990a}
E.J. Kansa.
\newblock Multiquadrics{\textemdash}a scattered data approximation scheme with
  applications to computational fluid-dynamics{\textemdash}i surface
  approximations and partial derivative estimates.
\newblock {\em Computers {\&} Mathematics with Applications}, 19(8-9):127--145,
  1990.

\bibitem{Kansa1990b}
E.J. Kansa.
\newblock Multiquadrics{\textemdash}a scattered data approximation scheme with
  applications to computational fluid-dynamics{\textemdash}{II} solutions to
  parabolic, hyperbolic and elliptic partial differential equations.
\newblock {\em Computers {\&} Mathematics with Applications}, 19(8-9):147--161,
  1990.

\bibitem{Driscoll2002}
T.A. Driscoll and B.~Fornberg.
\newblock Interpolation in the limit of increasingly flat radial basis
  functions.
\newblock {\em Computers {\&} Mathematics with Applications}, 43(3-5):413--422,
  feb 2002.

\bibitem{Shu2003}
C.~Shu, H.~Ding, and K.S Yeo.
\newblock Local radial basis function-based differential quadrature method and
  its application to solve two-dimensional incompressible
  navier{\textendash}stokes equations.
\newblock {\em Computer Methods in Applied Mechanics and Engineering},
  192(7-8):941--954, feb 2003.

\bibitem{Fornberg2011}
B.~Fornberg and E.~Lehto.
\newblock Stabilization of {RBF}-generated finite difference methods for
  convective {PDEs}.
\newblock {\em Journal of Computational Physics}, 230(6):2270--2285, mar 2011.

\bibitem{Fornberg2013}
B.~Fornberg, E.~Lehto, and C.~Powell.
\newblock Stable calculation of gaussian-based {RBF}-{FD} stencils.
\newblock {\em Computers {\&} Mathematics with Applications}, 65(4):627--637,
  feb 2013.

\bibitem{Davydov2011}
O.~Davydov and D.~Thi Oanh.
\newblock On the optimal shape parameter for gaussian radial basis function
  finite difference approximation of the poisson equation.
\newblock {\em Computers {\&} Mathematics with Applications}, 62(5):2143--2161,
  sep 2011.

\bibitem{Davydov2011a}
O.~Davydov and D.~Thi Oanh.
\newblock Adaptive meshless centres and {RBF} stencils for poisson equation.
\newblock {\em Journal of Computational Physics}, 230(2):287--304, jan 2011.

\bibitem{Kee2008}
B.B.T. Kee, G.R. Liu, and C.~Lu.
\newblock A least-square radial point collocation method for adaptive analysis
  in linear elasticity.
\newblock {\em Engineering Analysis with Boundary Elements}, 32(6):440--460,
  jun 2008.

\bibitem{Onate1996}
E.~O{\~{n}}ate, S.~Idelsohn, O.C. Zienkiewick, and R.L. Taylor.
\newblock A finite point method in computational mechanics. application to
  convective transport. and fluid flow.
\newblock {\em International Journal for Numerical Methods in Engineering},
  39(22):3839--3866, November 1996.

\bibitem{Lancaster1986}
P.~Lancaster and K.~Salkauskas.
\newblock {\em Curve \& Surface Fitting}.
\newblock Academic Press, 1986.

\bibitem{Ishida1999}
T.~Ishida and G.~C. Schatz.
\newblock A local interpolation scheme using no derivatives in quantum-chemical
  calculations.
\newblock {\em Chemical Physics Letters}, 314(3-4):369--375, dec 1999.

\bibitem{Maisuradze2003}
G.G. Maisuradze, D.L. Thompson, A.F. Wagner, and M.~Minkoff.
\newblock Interpolating moving least-squares methods for fitting potential
  energy surfaces: Detailed analysis of one-dimensional applications.
\newblock {\em The Journal of Chemical Physics}, 119(19):10002--10014, nov
  2003.

\bibitem{Sukumar2003}
N.~Sukumar.
\newblock Voronoi cell finite difference method for the diffusion operator on
  arbitrary unstructured grids.
\newblock {\em International Journal for Numerical Methods in Engineering},
  57(1):1--34, 2003.

\bibitem{Zhou2007}
J.X. Zhou, M.E. Li, Z.Q. Zhang, W.~Zou, and L.~Zhang.
\newblock A subdomain collocation method based on voronoi domain partition and
  reproducing kernel approximation.
\newblock {\em Computer Methods in Applied Mechanics and Engineering},
  196(13-16):1958--1967, mar 2007.

\bibitem{Oate1998}
E.~O{\~{n}}ate.
\newblock Derivation of stabilized equations for numerical solution of
  advective-diffusive transport and fluid flow problems.
\newblock {\em Computer Methods in Applied Mechanics and Engineering},
  151(1-2):233--265, jan 1998.

\bibitem{Oate2001}
E.~O{\~{n}}ate, F.~Perazzo, and J.~Miquel.
\newblock A finite point method for elasticity problems.
\newblock {\em Computers {\&} Structures}, 79(22-25):2151--2163, sep 2001.

\bibitem{Duflot2004}
M.~Duflot.
\newblock {\em Application des m\'ethodes sans maillage en m\'ecanique de la
  rupture}.
\newblock PhD thesis, 2004.

\bibitem{Organ1996}
D.~Organ, M.~Fleming, T.~Terry, and T.~Belytschko.
\newblock Continuous meshless approximations for nonconvex bodies by
  diffraction and transparency.
\newblock {\em Computational Mechanics}, 18(3):225--235, jul 1996.

\bibitem{MUMPS01}
P.R. Amestoy, I.S. Duff, J.Y. L'Excellent, and J.~Koster.
\newblock A fully asynchronous multifrontal solver using distributed dynamic
  scheduling.
\newblock {\em SIAM Journal on Matrix Analysis and Applications}, 23(1):15--41,
  2001.

\bibitem{MUMPS02}
P.R. Amestoy, A.~Guermouche, J.Y. L'Excellent, and S.~Pralet.
\newblock Hybrid scheduling for the parallel solution of linear systems.
\newblock {\em Parallel Computing}, 32(2):136--156, 2006.

\bibitem{petsc-user-ref}
S.~Balay, S.~Abhyankar, M.F. Adams, J.~Brown, P.~Brune, K.~Buschelman,
  L.~Dalcin, V.~Eijkhout, W.D. Gropp, D.~Kaushik, M.G. Knepley, D.A. May,
  L.~Curfman McInnes, R.~Tran Mills, T.~Munson, K.~Rupp, P.~Sanan, B.F. Smith,
  S.~Zampini, H.~Zhang, and H.~Zhang.
\newblock {PETS}c users manual.
\newblock Technical Report ANL-95/11 - Revision 3.9, Argonne National
  Laboratory, 2018.

\bibitem{petsc-efficient}
S.~Balay, W.D. Gropp, L.~Curfman McInnes, and B.F. Smith.
\newblock Efficient management of parallelism in object oriented numerical
  software libraries.
\newblock In E.~Arge, A.M. Bruaset, and H.P. Langtangen, editors, {\em Modern
  Software Tools in Scientific Computing}, pages 163--202. Birkh{\"{a}}user
  Press, 1997.

\bibitem{Dimitrov2001}
A.~Dimitrov, H.~Andr\"{a}, and E.~Schnack.
\newblock Efficient computation of order and mode of corner singularities in
  3d-elasticity.
\newblock {\em International Journal for Numerical Methods in Engineering},
  52(8):805--827, nov 2001.

\bibitem{Rachowicz2006}
W.~Rachowicz, D.~Pardo, and L.~Demkowicz.
\newblock Fully automatic hp-adaptivity in three dimensions.
\newblock {\em Computer Methods in Applied Mechanics and Engineering},
  195(37-40):4816--4842, jul 2006.

\bibitem{Zander2016}
N.~Zander, T.~Bog, M.~Elhaddad, F.~Frischmann, S.~Kollmannsberger, and E.~Rank.
\newblock The multi-level hp-method for three-dimensional problems: Dynamically
  changing high-order mesh refinement with arbitrary hanging nodes.
\newblock {\em Computer Methods in Applied Mechanics and Engineering},
  310:252--277, oct 2016.

\bibitem{Hrennikoff1941}
A.~Hrennikoff.
\newblock Solution of problems of elasticity by the framework method.
\newblock {\em Journal of applied mechanics}, 8.4:169--175, 1941.

\bibitem{Courant1943}
R.~Courant.
\newblock Variational methods for the solution of problems of equilibrium and
  vibrations.
\newblock {\em Bulletin of the American Mathematical Society}, 49:1--23, 1943.

\bibitem{Strang1973}
G.~Strang and G.~Fix.
\newblock {\em An Analysis of the Finite Element Method}.
\newblock Wellesley-Cambridge Press, 2008.

\bibitem{Aluru2000}
N.R. Aluru.
\newblock A point collocation method based on reproducing kernel
  approximations.
\newblock {\em International Journal for Numerical Methods in Engineering},
  47(6):1083--1121, feb 2000.

\bibitem{Perrone1975}
N.~Perrone and R.~Kao.
\newblock A general finite difference method for arbitrary meshes.
\newblock {\em Computers {\&} Structures}, 5(1):45--57, apr 1975.

\bibitem{Cottet1990}
G.H. Cottet.
\newblock A particle-grid superposition method for the navier-stokes equations.
\newblock {\em Journal of Computational Physics}, 89(2):301--318, aug 1990.

\bibitem{Sbalzarini2007}
I.F. Sbalzarini.
\newblock Particle methods for the simulation of diffusion processes in space.
\newblock 2007.

\bibitem{Schrader2011}
B.~Schrader.
\newblock {\em Discretization-Corrected PSE Operators for Adaptive
  Multiresolution Particle Methods}.
\newblock PhD thesis, 2011.

\bibitem{Bourantas2016}
G.C. Bourantas, B.L. Cheeseman, R.~Ramaswamy, and I.F. Sbalzarini.
\newblock Using {DC} {PSE} operator discretization in eulerian meshless
  collocation methods improves their robustness in complex geometries.
\newblock {\em Computers {\&} Fluids}, 136:285--300, sep 2016.

\bibitem{Flyer2016}
N.~Flyer, G.A. Barnett, and L.J. Wicker.
\newblock Enhancing finite differences with radial basis functions: Experiments
  on the navier{\textendash}stokes equations.
\newblock {\em Journal of Computational Physics}, 316:39--62, jul 2016.

\bibitem{Yensiri2017}
S.~Yensiri and R.~Skulkhu.
\newblock An investigation of radial basis function-finite difference
  ({RBF}-{FD}) method for numerical solution of elliptic partial differential
  equations.
\newblock {\em Mathematics}, 5(4):54, oct 2017.

\bibitem{Belytschko1996}
T.~Belytschko, Y.~Krongauz, D.~Organ, M.~Fleming, and P.~Krysl.
\newblock Meshless methods: An overview and recent developments.
\newblock {\em Computer Methods in Applied Mechanics and Engineering},
  139(1-4):3--47, dec 1996.

\bibitem{Szab1991}
B.~Szab\'o and I.~Babu{\v{s}}ka.
\newblock {\em Finite Element Analysis}.
\newblock Wiley-Interscience, 1991.

\bibitem{eigenweb}
G.~Guennebaud, B.~Jacob, et~al.
\newblock Eigen v3.
\newblock http://eigen.tuxfamily.org, 2010.

\bibitem{cgal:s-gkd-18b}
M.~Seel.
\newblock {dD} geometry kernel.
\newblock In {\em {CGAL} User and Reference Manual}. {CGAL Editorial Board},
  {4.13} edition, 2018.

\bibitem{cgal:hdj-t-18b}
O.~Devillers, S.~Hornus, and C.~Jamin.
\newblock {dD} triangulations.
\newblock In {\em {CGAL} User and Reference Manual}. {CGAL Editorial Board},
  {4.13} edition, 2018.

\bibitem{Rycroft2009}
C.H. Rycroft.
\newblock Voro++: A three-dimensional \textsc{V}oronoi cell library in
  \textsc{C}++.
\newblock {\em Chaos: An Interdisciplinary Journal of Nonlinear Science},
  19(4):041111, dec 2009.

\bibitem{Ahrens2005}
J.~Ahrens, B.~Geveci, and C.~Law.
\newblock Paraview: An end-user tool for large data visualization.
\newblock {\em Visualization Handbook}, 2005.

\bibitem{Ayachit2015}
U.~Ayachit.
\newblock {\em The ParaView Guide: A Parallel Visualization Application}.
\newblock Kitware, Incorporated, 2015.

\end{thebibliography}
	\pagebreak
	
	\appendix
	\setcounter{equation}{0}
	\renewcommand{\theequation}{A\arabic{equation}}
	\setcounter{figure}{0}
	\renewcommand{\theequation}{A\arabic{figure}}
	
	\section{GFD and DC PSE Methods Comparison for 1D Problems} \label{app:a}
	
	The purpose of this appendix is to present a detailed exposition of the GFD and the DC PSE methods for a simple 1D problem. We have selected the case of a second order PDE for illustration purposes.
	
	Considering a differential operator $\mathcal{A}$ and a field $f: \rm I\!R \rightarrow \rm  I\!R$, the following PDE can be written over the domain $\Omega$:
	\begin{equation} \label{1DProblemEquations}
	\mathcal{A}(f)=0 \text{ \quad in \quad} \Omega .
	\end{equation}
	The field $f$ shall also verify the conditions imposed on the Dirichlet and Neumann boundaries, which are denoted by $\Gamma_u$ and $\Gamma_t$, respectively. The field values are set to $\overline{f}$ on $\Gamma_u$. On $\Gamma_t$, the field shall verify a lower order PDE defined by an operator $\mathcal{B}$. 
	\begin{equation} \label{1DProblemBCs}
	\begin{aligned}
	f-\overline{f}&=0 \text{ \quad on \quad} \Gamma_u, \\
	\mathcal{B}(f) &=0 \text{ \quad on \quad} \Gamma_t.
	\end{aligned}
	\end{equation}
	\begin{figure}[H] 
		\centering
		\begin{tikzpicture}
		\def\svgwidth{6cm}
		\node at (0,0) {\includegraphics{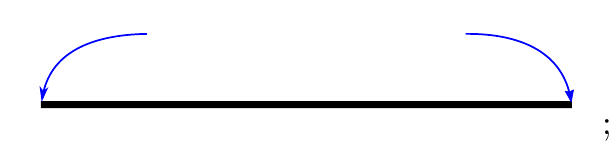}};
		\node[color=black] at (0,-0.1) [left] {$\Omega$};
		\node[color=blue] at (-0.8,0.4) [left] {$\Gamma_u$};
		\node[color=blue] at (1.5,0.4) [left] {$\Gamma_t$};
		\end{tikzpicture}
		\caption{1D Domain $\Omega$ with $\Gamma_u$ and $\Gamma_t$ boundaries.}
		\label{1D_DomainDrawing}
	\end{figure}
	In order to solve this problem by collocation, we need to transform it into a linear system of the form $\mathbf{AF=B}$, where $\mathbf{A}$ is the problem matrix, $\mathbf{F}$ is a vector containing the field values at each node of the domain, and $\mathbf{B}$ is a vector containing various constraints of the problem. The field derivatives at the collocation centers need to be approximated as a function of the values at the nodes. Collocation is typically performed at the nodes but can also be performed in other locations.
	
	\paragraph{Principle}\
	
	Both the GFD and the DC PSE methods are based on a Taylor's series expansion of the unknown field around the collocation node.
	
	\begin{minipage}[t]{.5\textwidth}
		\begin{center}
			\textbf{GFD}
		\end{center}
		\begin{adjustwidth}{0.0in}{0.075in}
			The field derivatives at the collocation node are determined so that the field values at the support nodes can be reproduced using a Taylor's series expansion. \\
			The field derivatives are approximated simultaneously at each node of the domain.
		\end{adjustwidth}		
	\end{minipage}
	\vline
	\begin{minipage}[t]{.5\textwidth}
		\begin{center}
			\textbf{DC PSE}
		\end{center}
		\begin{adjustwidth}{0.075in}{0.0in}
			A convolution function is selected so that the approximated derivative in the Taylor's series expansion only depends on the field values at the support nodes. All the other unknown derivatives are canceled by the selected function.\\
			Different convolution functions are used to approximate the derivatives of various orders at a collocation node.
		\end{adjustwidth}	
	\end{minipage}
	
	\quad \\
	
	The steps associated to each method are presented below for the case of a second order PDE. In the sections below, the nodes are labeled $X_i$ and the 1D coordinate associated to the node is written $x_i$.
	
	\paragraph{Step 1: Taylor's Series Approximation}\	
	
	For both methods, the first step consists in writing an approximation of the Taylor's series expansion up to the desired order. The approximation order shall be of at the least the highest derivative order of the differential operator $\mathcal{A}$. In 1D, the Taylors's series expansion at a point $X_{pi}$ in the vicinity of $X_c$ is:
	\begin{equation} \label{Taylor1D_AllTerms}
	f(X_{pi})=f(X_c) + \sum_{i=1}^{+\infty} \frac{(x_{pi} - x_c)^i}{i!} \frac{d^{i}f(X_c)}{d x^i}.
	\end{equation}
	Denoting the second order approximation of this expansion by $f_h(X_{pi})$, Equation (\ref{Taylor1D_AllTerms}) becomes:
	\begin{equation} \label{Taylor1D_SecondOrderApprox}
	f_h(X_{pi})= f(X_c) +(x_{pi} - x_c)\frac{d f (X_c)}{d x} + \frac{(x_{pi} - x_c)^{2}}{2!}\frac{d^2 f (X_c)}{d x^2}.
	\end{equation}
	
	\paragraph{Step 2: Support Node Selection}\
	
	The nodes in the vicinity of the collocation nodes are used to approximate the desired derivatives. These nodes are selected within a distance or radius $R_c$ to the collocation nodes $X_c$.
	\begin{figure}[H] 
		\centering
		\begin{tikzpicture}
		\def\svgwidth{10cm}
		\node at (0,0) {\includegraphics{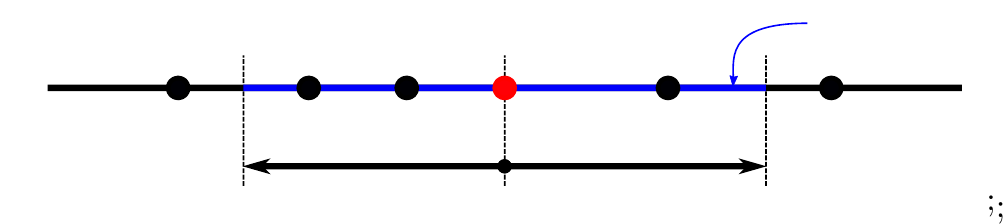}};
		\node[color=black] at (0.7,0.65) [left] {$X_c$};
		\node[color=black] at (-1.7,0.65) [left] {$X_{p1}$};
		\node[color=black] at (-0.7,0.65) [left] {$X_{p2}$};
		\node[color=black] at (2,0.65) [left] {$X_{p3}$};
		\node[color=black] at (-0.9,-0.4) [left] {$R_c$};
		\node[color=black] at (1.75,-0.4) [left] {$R_c$};
		\node[color=blue] at (3.8,0.9) [left] {$\Omega_c$};
		\node[color=black] at (4.7,0.4) [left] {$\Omega$};
		\end{tikzpicture}
		\caption{1D support $\Omega_c$ of the collocation node $X_c$. The radius of the support is $R_c$. The nodes $X_{p1}$, $X_{p2}$ and $X_{p3}$ are in the support of $X_c$.}
		\label{1D_SupportDrawing}
	\end{figure}
	
	\paragraph{Step 3: Derivatives Approximation}\
	
	\underline{\textbf{GFD}}
	
	The Taylor's series expansion presented in Equation (\ref{Taylor1D_SecondOrderApprox}) can be written for each node of the collocation node support. For the example presented in Figure \ref{1D_SupportDrawing}, three nodes are present in the support of the collocation node. The following system is obtained:
	\begin{equation}
	\left \{
	\begin{array}{ll}
	f_h(X_{p1})= f(X_c) +(x_{p1} - x_c)\frac{d f (X_c)}{d x} + \frac{(x_{p1} - x_c)^{2}}{2!}\frac{d^2 f (X_c)}{d x^2}\\
	f_h(X_{p2})= f(X_c) +(x_{p2} - x_c)\frac{d f (X_c)}{d x} + \frac{(x_{p2} - x_c)^{2}}{2!}\frac{d^2 f (X_c)}{d x^2}\\
	f_h(X_{p3})= f(X_c) +(x_{p3} - x_c)\frac{d f (X_c)}{d x} + \frac{(x_{p3} - x_c)^{2}}{2!}\frac{d^2 f (X_c)}{d x^2}.\\
	\end{array}
	\right.
	\end{equation}
	If the number of support nodes is larger than the approximated derivative order (two for the selected example), the system is overdetermined. In that case, the derivatives that best reproduce the field values are determined using a mean least square method. The associated minimization problem is written in the form of a functional $B$. A weight function $w$ is used to balance the contribution of the support nodes as a function of their distance to the collocation node.
	\begin{equation}
	B(X_c)=\sum_{i=1}^3 w(X_{pi} - X_c) \Big[ f(X_c) - f(X_{pi}) + (x_{pi} - x_c)\frac{d f (X_c)}{d x} + \frac{(x_{pi} - x_c)^{2}}{2!}\frac{d^2 f (X_c)}{d x^2} \Big]^2.
	\end{equation}
	The derivatives $\mathbf{Df}(X)= \Big[ \frac{d f (X)}{d x}, \frac{d^2 f (X)}{d x^2} \Big]^T $, that best approximate the known field values based on the Taylor's series expansion, minimize $B$ when:
	\begin{equation} \label{DerivativeFunctionalB_1D_GFD}
	\frac{d B(X)}{d \mathbf{Df}(X)}\biggr\rvert_{X=X_c}=0,
	\end{equation}
	\begin{equation}
	\left \{
	\begin{array}{ll}
	\begin{aligned}
	&\sum_{i=1}^3 w(X_{pi} - X_c) (x_{pi} - x_c) \Big[ f(X_c) - f(X_{pi}) + (x_{pi} - x_c)\frac{d f (X_c)}{d x} + \frac{(x_{pi} - x_c)^{2}}{2!}\frac{d^2 f (X_c)}{d x^2} \Big]&=0\\
	&\sum_{i=1}^3 w(X_{pi} - X_c) \frac{(x_{pi} - x_c)^{2}}{2!} \Big[ f(X_c) - f(X_{pi}) + (x_{pi} - x_c)\frac{d f (X_c)}{d x} + \frac{(x_{pi} - x_c)^{2}}{2!}\frac{d^2 f (X_c)}{d x^2} \Big]&=0.\\
	\end{aligned}
	\end{array}
	\right.
	\end{equation}
	This system can be rearranged in a matrix form as follows:
	\begin{equation}\label{GFD_1D_DerivApprox}
	\begin{bmatrix}
	m_{11} & m_{12} \\
	m_{21} & m_{22} \\
	\end{bmatrix}
	\begin{bmatrix}
	\frac{d f (X_c)}{d x} \\
	\frac{d^2 f (X_c)}{d x^2} \\
	\end{bmatrix}=
	\begin{bmatrix}
	-m_{01} & m_{01,1} & m_{01,2} & m_{01,3} \\
	-m_{02} & m_{02,1} & m_{02,2} & m_{02,3} \\
	\end{bmatrix}
	\begin{bmatrix}
	f(X_{c}) \\
	f(X_{p1}) \\
	f(X_{p2}) \\
	f(X_{p3}) \\
	\end{bmatrix},
	\end{equation}
	where the moments $m_{ij,k}$ and $m_{ij}$ and the matrix $\mathbf{P}(X_c) \in \rm I\!R^{3 \times 3}$ correspond to:
	\begin{equation}\label{Moments_GFD1D}
	\begin{aligned}
	m_{ij,k} &= w(X_{pk} - X_c) P_{(i+1)k}(X_c) P_{(j+1)k}(X_c), \\
	m_{ij} &= \sum_{k=1}^3 {m_{ij,k}},\\
	\mathbf{P}(X_c)&=\begin{bmatrix}
	1 &1 & 1 \\
	(x_{p1} - x_c) & (x_{p2} - x_c) & (x_{p3} - x_c)  \\
	\frac{(x_{p1} - x_c)^2}{2!} & \frac{(x_{p2} - x_c)^2}{2!} & \frac{(x_{p3} - x_c)^2}{2!}  \\
	\end{bmatrix}.\\
	\end{aligned}
	\end{equation}
	The Equation (\ref{GFD_1D_DerivApprox}) can be represented in the form $\mathbf{A}(X_c) \mathbf{Df}(X_c) = \mathbf{E}(X_c) \mathbf{F}(X_c)$.
	
	\underline{\textbf{DC PSE}}
	
	The Taylor's series expansion presented in Equation (\ref{Taylor1D_SecondOrderApprox}) can be convoluted by a function $\eta$ over the support $\Omega_c$ of the collocation node $X_c$:
	\begin{equation} \label{Convolution_SecondOrderApprox1D}
	\begin{aligned}
	\int_{\Omega_c} {f_h(X_p)} \eta(X_p-X_c) dX_p = &\int_{\Omega_c} {f(X_c)} \eta(X_p-X_c) dX_p + \int_{\Omega_c} {\frac{d f(X_c)}{d x}} (x_p - x_c) \eta(X_p-X_c) dX_p \\
	&+ \int_{\Omega_c} {\frac{d^2 f(X_c)}{d x^2}} \frac{(x_p - x_c)^{2}}{2!} \eta(X_p-X_c) dX_p. \\
	\end{aligned}
	\end{equation}
	The integral can be approximated by a sum, assuming that the nodes are regularly distributed over the support and that the field $f$ is sufficiently smooth.
	\begin{equation} \label{ConvWithMoment_SecondOrderApprox1D}
	\sum_{i=1}^3 {f_h(X_{pi})} \eta(X_{pi}-X_c) = f(X_c) M_{0}(X_c) + {\frac{d f(X_c)}{d x}} M_{1}(X_c)+ {\frac{d f^2(X_c)}{d x^2}} M_{2}(X_c),
	\end{equation}
	where the moments $M_j$ are:
	\begin{equation} \label{DiscreteMomentEquation_UnitVol1D}
	M_{j}(X_c)= \sum_{i=1}^3 \frac{(x_{pi} - x_c)^j}{j!} \eta(X_{pi}-X_c).
	\end{equation}
	The convolution function is chosen so that all the moments in Equation (\ref{ConvWithMoment_SecondOrderApprox1D}) are null except the one multiplying the approximated derivative of order $n_x$, which is denoted by $D^{n_x}f(X_c)$. This moment is set to unity.
	\begin{equation}\label{DC PSE Operator1D}
	\left \{
	\begin{aligned}
	&D^{n_x}f(X_c) = \sum_{i=1}^3 {f_h(X_{pi})} \eta(X_{pi}-X_c)\\
	&\begin{array}{ll}
	\text{with} &M_{n_x}(X_c)=1\\
	&M_{i}(X_c)=0 \text{ \quad } \text{if} \ i \ne n_x.\\
	\end{array}\\
	\end{aligned}
	\right.	
	\end{equation}
	In order to satisfy this moment condition, the convolution function is chosen as the product of two functions: a correction function $K$ and a weight function $w$. The correction function is typically the product of a coefficient vector $\mathbf{a}$ and a polynomial basis $\mathbf{P}$. For the 1D second order approximation, the polynomial basis $\mathbf{P}=[1, x, x^2]^T$ can be selected. The correction function can then be written as follows:
	\begin{equation} \label{SimpleKernel1D}
	\eta(X_{p}-X_c) ={\mathbf{P}(X_{p}-X_c)}^T \mathbf{a} \ w(X_{p}-X_c).
	\end{equation}
	The coefficients of the vector $\mathbf{a}$ are determined in order to satisfy the moment condition set in Equation (\ref{DC PSE Operator1D}). For instance, the moment condition associated to the second order derivative approximation is:
	\begin{equation} \label{DCPSE_LinerSyst1D}
	\left \{
	\begin{array}{ll}
	\begin{aligned}
	&M_{0}(X_c)=0 &\Leftrightarrow \quad &\sum_{i=1}^3 {\mathbf{P}(X_{pi}-X_c)}^T \mathbf{a} w(X_{pi}-X_c) = 0 \\
	&M_{1}(X_c)=0 &\Leftrightarrow \quad &\sum_{i=1}^3 (x_{pi} - x_c) {\mathbf{P}(X_{pi}-X_c)}^T \mathbf{a} w(X_{pi}-X_c) = 0 \\
	&M_{2}(X_c)=1 &\Leftrightarrow \quad &\sum_{i=1}^3 \frac{(x_{pi} - x_c)^2}{2!} {\mathbf{P}(X_{pi}-X_c)}^T \mathbf{a} w(X_{pi}-X_c) = 1. \\
	\end{aligned}
	\end{array}
	\right.
	\end{equation}
	The system of equations can be put in a matrix form as follows:
	\begin{equation}\label{DCPSE_1D_DerivApprox}
	\begin{bmatrix}
	A_{11} & A_{12} & A_{13} \\
	A_{21} & A_{22} & A_{23} \\
	A_{31} & A_{32} & A_{33} \\
	\end{bmatrix}
	\begin{bmatrix}
	a_1 \\
	a_2 \\
	a_3 \\
	\end{bmatrix}=
	\begin{bmatrix}
	0 \\
	0 \\
	1 \\
	\end{bmatrix}.
	\end{equation}
	Considering the vector $\mathbf{Q}(X_c,X_p)=[1,(x_p - x_c),\frac{(x_p - x_c)^2}{2!}]^T$, the correction function basis $\mathbf{P}$ and the weight function $w$, the coefficients of the matrix $\mathbf{A} \in \rm I\!R^{3 \times 3}$ can be written as:
	\begin{equation}\label{DCPSE_ACoefficients1D}
	A_{i,j}(X_c)=\sum_{i=1}^3 Q_i(X_c,X_{pi})P_j(X_{pi}-X_c)w(X_{pi}-X_c).
	\end{equation}
	
	\paragraph{Step 4: Solution of the Collocation Linear Systems}\
	
	\begin{minipage}[t]{.5\textwidth}
		\begin{center}
			\textbf{GFD}
		\end{center}
		\begin{adjustwidth}{0.0in}{0.075in}
			The system presented in Equation (\ref{GFD_1D_DerivApprox}) can be solved in order to obtain the derivatives $\mathbf{Df}(X_c)$ as a function of the field values $\mathbf{F}(X_c)$.
		\end{adjustwidth}		
	\end{minipage}
	\vline
	\begin{minipage}[t]{.5\textwidth}
		\begin{center}
			\textbf{DC PSE}
		\end{center}
		\begin{adjustwidth}{0.075in}{0.0in}
			The system presented in Equation (\ref{DCPSE_1D_DerivApprox}) is solved in order to obtain the coefficients of the correction function $\eta$. Once these coefficients are obtained, the convolution function presented in Equation (\ref{SimpleKernel1D}) can be calculated and the derivative $D^{n_x}f(X_c)$ presented in Equation (\ref{DC PSE Operator1D}) can be approximated. \\
			The solution of Equation (\ref{DCPSE_1D_DerivApprox}) needs to be performed for all the moment conditions associated to the approximated derivatives in the differential operators $\mathcal{A}$ and $\mathcal{B}$.
		\end{adjustwidth}	
	\end{minipage}
	
	\paragraph{Step 5: Assembly of the Linear Problem}\
	
	The steps 1 to 4 allowed the approximation of the derivatives at the collocation nodes as a function of the field values at theses nodes and at their support nodes. Sets of coefficients $C_{n_x}(X_c)$ are obtained for each derivatives so that $D^{n_x}f(X_c)=\mathbf{C_{n_x}}(X_c) \ \mathbf{F}(X_c)$. Finally, based on the differential operator and on the boundary conditions, the problem matrices $\mathbf{A}$ and $\mathbf{B}$ are assembled.
	
\end{document}